\documentclass[a4paper,leqno,12pt]{amsart}

\usepackage{fancyheadings}
\usepackage{times}
\usepackage{amssymb,latexsym}

\providecommand{\prenom}{} 
\providecommand{\nomfam}{F. Gavarini}
\providecommand{\adresse}{Universit\`a degli Studi di Roma
                          ``Tor Vergata'' \\
                          Dipartimento di Matematica  \\
                          Via della Ricerca Scientifica 1\\
                          I-00133 Roma \\
                          Italy}
\providecommand{\celec}{gavarini@mat.uniroma2.it} 
\providecommand{\titreexp}{The global quantum duality principle:  \\  
                           a survey through examples}
\providecommand{\titreexpo}{The global quantum duality principle:  
                            a survey through examples}
\providecommand{\annee}{2002}
\providecommand{\noexpose}{II}

\renewcommand{\thepage}{\noexpose-\arabic{page}}
\providecommand{\HRule}{\rule{\linewidth}{1mm}}

\newtheorem{thm}{Th\'eor\`eme}[section]

\newtheorem{hyp}[thm]{Hypoth\`ese}

\def \loongrightarrow {\relbar\joinrel\relbar\joinrel\rightarrow}
\def \llongrightarrow
{\relbar\joinrel\relbar\joinrel\relbar\joinrel\rightarrow}
\def \longtwoheadrightarrow
{\relbar\joinrel\relbar\joinrel\twoheadrightarrow}
\def \llongtwoheadrightarrow
{\relbar\joinrel\relbar\joinrel\relbar\joinrel\twoheadrightarrow}
\def \gerg {\mathfrak g}
\def \gere {\mathfrak e}
\def \gerh {\mathfrak h}
\def \gerk {\mathfrak k}
\def \germ {\mathfrak m}
\def \gersl {{\mathfrak{sl}}}
\def \und1 {\underline{1}}

\def \u {{\text{\bf u}}}
\def \hyp {{\text{\sl Hyp}\,}}
\def \h {\hbar}

\def \id {\text{\rm id}}

\def \L {\mathcal{L}}
\def \U {\mathcal{U}}
\def \calM {\mathcal{M}}

\def \HA {\mathcal{HA}}
\def \H {\mathbb{H}}
\def \N {\mathbb{N}}
\def \Z {\mathbb{Z}}

\def \F {\mathbb{F}}
\def \A {\mathbb{A}}

\def \Char {\hbox{\it Char}\,}
\def \QrUEA {\mathcal{Q}{\hskip1pt}r\mathcal{UE{\hskip-1,1pt}A}}
\def \QFA {\mathcal{QF{\hskip-1,7pt}A}}
\def \PrUEA {\mathcal{P}{\hskip0,8pt}r\mathcal{UE{\hskip-1,1pt}A}}
\def \PFA {\mathcal{PF{\hskip-1,7pt}A}}

\begin{document}  

\thispagestyle{empty}

\vspace*{\stretch{1}}
\noindent\HRule
\begin{flushright}
\Huge
\prenom\ \textsc{\nomfam}\\[8mm]
\emph{\titreexp}
\end{flushright}
\noindent\HRule
\vspace*{\stretch{2}}
\begin{center}
\Large\textsc{Rencontres Math\'ematiques de Glanon}\\ 
\Large\textsc{--- $\mathfrak{\annee}$ ---}
\end{center}

\newpage
\pagestyle{fancy}
\thispagestyle{empty}

\vspace*{\stretch{1}}
\noindent\Large\textsc{Rencontres Math\'ematiques de Glanon \annee}\\
\hskip2cm\Large\textsc{Contribution \noexpose}\\
\\ 
\\ 
 
\noindent\Large \prenom\ \textsc{\nomfam}\large\\[8mm]
\ \hskip6cm\ \adresse
\\ \ \hskip6cm\ {\tt\celec}\\
\\ 
\\ 
\ \hskip6cm\ \Large\emph{\titreexpo}\\ 
\ 
\vspace*{\stretch{2}}

\normalsize
\setcounter{page}{0}
\newpage

\thispagestyle{empty}

\lhead[\textsc{Rencontres Math\'ematiques de Glanon 2002}]{\prenom\
\textsc{\nomfam}}
\rhead[\textsc{Contribution \noexpose}]{\it The global quantum duality 
                                            principle...}

\ 

\vskip3cm

\begin{center}
{\Large\textbf{\titreexp}}\\ 
\ \\ 
{\large\prenom\ \textsc{\nomfam}}\\ 
\adresse\\ 
{\tt\celec}
\end{center} 
\normalsize

\vskip1.5cm
\setcounter{page}{1}

\section*{Abstract}
\begin{quotation}
\begin{sloppy}
   Let  $ R $  be a 1-dimensional integral domain, let  $ \, \h \in
R \setminus \{0\} \, $  be prime, and let  $ \HA $  be the category
of torsionless Hopf algebras over  $ R $.  We call  $ \, H \in \HA
\, $  a  {\it ``quantized function algebra'' (=QFA)}, resp.~{\it
``quantized restricted universal enveloping algebras'' (=QrUEA)},
at  $ \h \, $  if  $ \, H \big/ \h \, H \, $  is the function algebra
of a connected Poisson group, resp.~the (restricted, if  $ \, R \big/
\h \, R \, $  has positive characteristic) universal enveloping algebra
of a (restricted) Lie bialgebra.
%
%
                                               \par
   An ``inner'' Galois correspondence on  $ \HA $  is established
via the definition of two endofunctors,  $ (\ )^\vee $  and
$ (\ )' $,  of  $ \HA $  such that:  {\it (a)} \, the image of
$ (\ )^\vee $,  resp.~of  $ (\ )' $,  is the full subcategory
of all QrUEAs, resp.~QFAs,  at  $ \h \, $;  {\it (b)} \, if
$ \, p := \text{\it Char}\big(R \big/ \h \, R \big) = 0 \, $,
the restrictions  $ (\ )^\vee{\big|}_{\text{QFAs}} $  and
$ (\ )'{\big|}_{\text{QrUEAs}} $  yield equivalences inverse
to each  other;  {\it (c)} \, if  $ \, p = 0 \, $,  starting from
a QFA over a Poisson group  $ G $,  resp.~from a QrUEA over a Lie
bialgebra  $ \gerg $,  the functor  $ (\ )^\vee $,  resp.~$ (\ )' $,
gives a QrUEA, resp.~a QFA, over the dual Lie bialgebra, resp.~a dual
Poisson group.  In particular,  {\it (a)\/}  provides a machine to
produce quantum groups of both types (either QFAs or QrUEAs),  {\it
(b)\/}  gives a characterization of them among objects of  $ \HA $,
and  {\it (c)\/}  gives a ``global'' version of the so-called
``quantum duality principle'' (after Drinfeld's, cf.~[Dr]).
                                               \par
   This result applies in particular to Hopf algebras of the form  $ \,
\Bbbk[\h] \otimes_\Bbbk H \, $  where  $ H $  is a Hopf algebra over
the field  $ \Bbbk $:  this yields quantum groups, hence ``classical''
geometrical symmetries of Poisson type (Poisson groups or Lie bialgebras,
via specialization) associated to the ``generalized'' symmetry encoded
by  $ H $.  Both our main result and the above mentioned application
are illustrated by means of several examples, which are studied in
some detail.
                                               \par
   These notes draw a sketch of the theoretical construction leading
to the ``global quantum duality principle''.  Besides, the principle
itself, and in particular the above mentioned application, is
illustrated by means of several examples: group algebras, the
standard quantization of the Kostant-Kirillov structure on any
Lie algebra, the quantum semisimple groups, the quantum Euclidean
group and the quantum Heisenberg group.
\end{sloppy}
\end{quotation}

\vskip1.5cm


\centerline {\bf Introduction }

\vskip10pt

   The most general notion of ``symmetry'' in mathematics is encoded
in the definition of Hopf algebra.  Among Hopf algebras  $ H $  over a
field, the commutative and the cocommutative ones encode ``geometrical''
symmetries, in that they correspond, under some technical conditions, to
algebraic groups and to (restricted, if the ground field has positive
characteristic) Lie algebras respectively: in the first case  $ H $
is the algebra  $ F[G] $  of regular functions over an algebraic group
$ G $,  whereas in the second case it is the (restricted) universal
enveloping algebra  $ U(\gerg) $  ($ \u(\gerg) $  in the restricted case)
of a (restricted) Lie algebra  $ \gerg \, $.  A popular generalization
of these two types of ``geometrical symmetry'' is given by quantum groups:
roughly, these are Hopf algebras  $ H $  depending on a parameter  $ \h $
such that setting  $ \h = 0 \, $  the Hopf algebra one gets is either of
the type  $ F[G] $   --- hence  $ H $  is a  {\sl quantized function
algebra}, in short QFA ---   or of the type  $ U(\gerg) $  or
$ \u(\gerg) $  (according to the characteristic of the ground field)
--- hence  $ H $  is a  {\sl quantized (restricted) universal enveloping
algebra}, in short QrUEA.  When a QFA exists whose specialization
(i.e.~its ``value'' at  $ \, \h = 0 \, $)  is  $ F[G] $,  the group
$ G $  inherits from this ``quantization'' a Poisson bracket, which
makes it a Poisson (algebraic) group; similarly, if a QrUEA exists
whose specialization is  $ U(\gerg) $  or  $ \u(\gerg) $,  the
(restricted) Lie algebra  $ \gerg $  inherits a Lie cobracket
which makes it a Lie bialgebra.  Then by Poisson group theory
one has Poisson groups  $ G^* $  dual to  $ G $  and a Lie bialgebra
$ \gerg^* $  dual to  $ \gerg \, $,  so other geometrical symmetries
are related to the initial ones.
                                            \par
   The dependence of a Hopf algebra on  $ \h $  can be described
as saying that it is defined over a ring  $ R $  and  $ \, \h \in
R \, $:  \, so one is lead to dwell upon the category  $ \HA $  of
Hopf  $ R $--algebras  (maybe with some further conditions), and
then raises three basic questions:

\vskip3pt

   {\bf --- (1)}  {\sl How can we produce quantum groups?}
                                            \par
   {\bf --- (2)}  {\sl How can we characterize quantum groups
(of either kind) within  $ \HA $?}
                                            \par
   {\bf --- (3)}  {\sl What kind of relationship, if any, does
exist between quantum groups over mutually dual Poisson groups,
or mutually dual Lie bialgebras?}

\vskip3pt

   A first answer to question  {\bf (1)}  and  {\bf (3)}  together
is given, in characteristic zero, by the so-called ``quantum duality
principle'', known in literature with at least two formulations.  One
claims that quantum function algebras associated to dual Poisson groups
can be taken to be dual   --- in the Hopf sense ---   to each other;
and similarly for quantum enveloping algebras (cf.~[FRT1] and [Se]).
The second one, formulated by Drinfeld in local terms (i.e., using
formal groups, rather than algebraic groups, and complete topological
Hopf algebras; cf.~[Dr], \S 7, and see [Ga4] for a proof), provides a
recipe to get, out of a QFA over  $ G $,  a QrUEA over  $ \gerg^* $,
and, conversely, to get a QFA over  $ G^* $  out of a QrUEA over
$ \gerg \, $.  More precisely, Drinfeld defines two functors, inverse
to each other, from the category of quantized universal enveloping
algebras (in his sense) to the category of quantum formal series
Hopf algebras (his analogue of QFAs) and viceversa, such that  $ \,
U_\h(\gerg) \mapsto F_\h[[G^*]] \, $  and  $ \, F_\h[[G]] \mapsto
U_\h(\gerg^*) \, $  (in his notation, where the subscript  $ \h $
stands as a reminder for ``quantized'' and the double square
brackets stand for ``formal series Hopf algebra'').
                                            \par
   In this paper we present a  {\sl global\/}  version of the quantum
duality principle which gives a complete answer to questions  {\bf (1)}
through  {\bf (3)}.  The idea is to push as far as possible Drinfeld's
original method so to apply it to the category  $ \HA $  of all Hopf
algebras which are torsion-free modules over some 1-dimensional domain
(in short, 1dD), say  $ R $,  and to do it for each non-zero prime
element  $ \h $  in  $ R \, $.
                                            \par
   To be precise, we extend Drinfeld's recipe so to define functors
from  $ \HA $  to itself.  We show that the image of these
``generalized'' Drinfeld's functors is contained in a category
of quantum groups   --- one gives QFAs, the other QrUEAs ---   so
we answer question  {\bf (1)}.  Then, in the zero characteristic case,
we prove that when restricted to quantum groups these functors yield
equivalences inverse to each other.  Moreover, we show that these
equivalences exchange the types of quantum group (switching QFA with
QrUEA) and the underlying Poisson symmetries (interchanging  $ G $
or  $ \gerg $  with  $ G^* $  or  $ \gerg^* $),  thus solving  {\bf
(3)}.  Other details enter the picture to show that these functors
endow  $ \HA $  with sort of a (inner) ``Galois correspondence'',
in which QFAs on one side and QrUEAs on the other side are the
subcategories (in  $ \HA $)  of ``fixed points'' for the composition
of both Drinfeld's functors (in the suitable order): in particular,
this answers question  {\bf (2)}.  It is worth stressing that, since
our ``Drinfeld's functors'' are defined for each non-trivial point
$ (\h) $  of  $ \text{\it Spec}\,(R) $,  \, for any such  $ (\h) $
and for any  $ H $  in  $ \HA $  they yield two quantum groups,
namely a QFA and a QrUEA,  w.r.t.~$ \h $  itself.
   \hbox{Thus we have a method to get, out of  any single
$ H \hskip-0pt \in \hskip-0pt \HA \, $, several quantum groups.}
                                             \par
   Therefore the ``global'' in the title is meant in several respects:
geometrically, we consider global objects (namely Poisson groups rather
than Poisson  {\sl formal\/}  groups, as in Drinfeld's approach);
algebraically we consider quantum groups over any 1dD  $ R $,  so there
may be several different ``semiclassical limits'' (=specialization)
to consider, one for each non-trivial point in the spectrum of  $ R $
(whereas in Drinfeld's context  $ \, R = \Bbb\Bbbk[[\h]] \, $  so one
can specialise only at  $ \, \h = 0 \, $);  more generally, our recipe
applies to  {\sl any\/}  Hopf algebra, i.e.~not only to quantum groups;
finally, most of our results are characteristic-free, i.e.~they hold
not only in zero characteristic (as in Drinfeld's case) but also in
positive characteristic.  As a further outcome, this ``global version''
of the quantum duality principle leads to formulate a ``quantum duality
principle for subgroups and homogeneous spaces'', see [CG].
                                             \par
   A key, long-ranging application of our  {\sl global quantum duality
principle\/}  (GQDP) is the following.  Take as  $ R $  the polynomial
ring  $ \, R = \Bbbk[\h\,] \, $,  where  $ \Bbbk $  is a field: then
for any Hopf algebra over  $ \Bbbk $  we have that  $ \, H[\h\,] := R
\otimes_\Bbbk H \, $  is a torsion-free Hopf  $ R $--algebra,  hence
we can apply Drinfeld's functors to it.  The outcome of this procedure
is the  {\sl crystal duality principle\/}  (CDP), whose statement strictly
resembles that of the GQDP: now Hopf  $ \Bbbk $--algebras  are looked at
instead of torsionless Hopf  $ R $--algebras,  and quantum groups are
replaced by Hopf algebras with canonical filtrations such that the
associated graded Hopf algebra is either commutative or cocommutative.
Correspondingly, we have a method to associate to  $ H $  a Poisson
group  $ G $  and a Lie bialgebra  $ \gerk $  such that  $ G $  is an
affine space (as an algebraic variety) and  $ \gerk $  is graded (as
a Lie algebra); in both cases, the ``geometrical'' Hopf algebra can be
attained   --- roughly ---   through a continuous 1-parameter deformation
process.  This result can also be formulated in purely classical   ---
i.e.~``non-quantum'' ---   terms and proved by purely classical means.
However, the approach via the GQDP also yields further possibilities to
deform  $ H $  into other Hopf algebras of geometrical type, which is
out of reach of any classical approach.
                                             \par
   The purpose of these notes is to illustrate the global quantum
duality principle in some detail through some relevant examples,
namely the application to the ``Crystal Duality Principle'' (\S 3) and
some quantum groups: the standard quantization of the Kostant-Kirillov
structure on a Lie algebra (\S 4), the quantum semisimple groups (\S 5),
the three dimensional quantum Euclidean group (\S 6), the quantum
Heisenberg group.  All details and technicalities which are skipped
in the present paper can be found in [Ga5], together with another
relevant example (see also [Ga6] and [Ga7]).

\vskip5pt

   These notes are the written version of the author's contribution
to the conference  {\it ``Rencontres Math\'ematiques de Glanon'',
6th edition\/} (1--5 july 2002) held in Glanon (France).  The
author's heartily thanks the organizers   --- especially Gilles
Halbout ---   for kindly inviting him.  Il remercie aussi tous
les Glanonnets pour leur charmante hospitalit\'e.

\vskip11pt

\centerline{ \sc acknowledgements }

  The author thanks C.~Gasbarri, A.~Frabetti, L.~Foissy, B.~Di
Blasio, A.~D'An\-drea, I.~Damiani, N.~Ciccoli, G.~Carnovale,
D.~Fiorenza, E.~Taft and P.~Baumann for many helpful discussions.

\vskip1,3truecm

\centerline {\bf \S \; 1 \ Notation and terminology }

\vskip10pt

  {\bf 1.1 The classical setting.} \, Let  $ \Bbbk $  be a fixed
field of any characteristic.  We call ``algebraic group'' the
maximal spectrum  $ G $  associated to any commutative Hopf
$ \Bbbk $--algebra  $ H $  (in particular, we deal with  {\sl
proaffine\/}  as well as  {\sl affine\/}  algebraic groups); then
$ H $  is called the algebra of regular function on  $ G $,  denoted
with  $ F[G] $.  We say that  $ G $  is connected if  $ F[G] $  has
no non-trivial idempotents; this is equivalent to the classical
topological notion when
%
%
$ \dim(G) $
is finite.  If  $ G $  is an algebraic group, we denote by $ \germ_e $
the defining ideal of the unit element  $ \, e \in G \, $  (in fact
$ \germ_e $  is the augmentation ideal of  $ F[G] \, $).  The cotangent
space of  $ G $  at  $ e $  is  $ \, \gerg^\times := \germ_e  \Big/
{\germ_e}^{\!2} \, $,  \, endowed with its weak topology, which is
naturally a Lie coalgebra.  By  $ \gerg $  we mean the tangent space
of  $ G $  at  $ e $,  realized as the topological dual  $ \, \gerg
:= {\big( \gerg^\times \big)}^\star \, $  of  $ \gerg^\times \, $: 
\, this is the tangent Lie algebra of  $ G $.  By  $ U(\gerg) $  we
mean the universal enveloping algebra of  $ \gerg $:  this is a
connected cocommutative Hopf algebra, and there is a natural Hopf
pairing (see  \S 1.2{\it (a)\/})  between  $ F[G] $  and  $ U(\gerg) $. 
If  $ \, \hbox{\it Char} \,(\Bbbk) = p > 0 \, $,  \, then  $ \gerg $ 
is a restricted Lie algebra, and  $ \, \u(\gerg) := U(\gerg) \Big/ 
\big( \big\{\, x^p - x^{[p\hskip0,7pt]} \,\big|\, x \in \gerg \,\big\}
\big) \, $  is the restricted universal enveloping algebra of  $ \gerg
\, $.  In the sequel, in order to unify notation and terminology,
when  $ \, \Char(\Bbbk) = 0 \, $  we call any Lie algebra  $ \gerg $ 
``restricted'', and its ``restricted universal enveloping algebra'' will
be  $ U(\gerg) $,  and we write  $ \, \U(\gerg) := U(\gerg) \, $  if 
$ \, \Char(\Bbbk) = 0 \, $  and  $ \, \U(\gerg) := \u(\gerg) \, $  if 
$ \, \Char(\Bbbk) > 0 \, $.   
                                            \par
   We shall also consider  $ \, \hyp(G) := {\big( {F[G]}^\circ
\big)}_\varepsilon = \big\{\, f \in {F[G]}^\circ \,\big|\,
f({\germ_e}^{\!n}) = 0 \; \forall \, n \geqq 0 \,\big\} \, $,  \,
i.e.~the connected component of the Hopf algebra  $ {F[G]}^\circ $
dual to  $ F[G] $;  this is called the  {\sl hyperalgebra\/}  of
$ G $.  By construction  $ \hyp(G) $  is a connected Hopf algebra,
containing  $ \, \gerg = \text{\sl Lie\/}(G) \, $;  if  $ \, \text{\it
Char}\,(\Bbbk) = 0 \, $  one has  $ \, \hyp(G) = U(\gerg) \, $,  \,
whereas if  $ \, \text{\it Char}\,(\Bbbk) > 0 \, $  one has a sequence
of Hopf algebra morphisms  $ \; U(\gerg) \longtwoheadrightarrow \u(\gerg)
\; {\lhook\joinrel\relbar\joinrel\relbar\joinrel\rightarrow} \, \hyp(G)
\; $.  In any case, there exists a natural perfect (= non-degenerate)
Hopf pairing between  $ F[G] $  and  $ \hyp(G) $.
                                            \par
   Now assume  $ G $  is a Poisson group (for this and other notions
hereafter see, e.g., [CP], but within an  {\sl algebraic geometry\/}
setting): then  $ F[G] $  is a Poisson Hopf algebra, and its Poisson
bracket induces on  $ \gerg^\times $  a Lie bracket which makes it
into a Lie bialgebra; so  $ U(\gerg^\times) $  and  $ \U(\gerg^\times) $ 
are co-Poisson Hopf algebras too.  On the other hand,  $ \gerg $  turns
into a Lie bialgebra   --- maybe in topological sense, if  $ G $  is
infinite dimensional ---   and  $ U(\gerg) $  and  $ \U(\gerg) $  are
(maybe topological) co-Poisson Hopf algebras.  The Hopf pairing above
between  $ F[G] $  and  $ \U(\gerg) $  then is compatible with these
additional co-Poisson and Poisson structures.  Similarly,  $ \hyp(G) $ 
is a co-Poisson Hopf algebra as well and the Hopf pairing between 
$ F[G] $  and  $ \hyp(G) $  is compatible with the additional structures. 
Moreover, the perfect pairing  $ \, \gerg \times \gerg^\times \!
\longrightarrow \Bbbk \, $  given by evaluation is compatible with
the Lie bialgebra structure on either side (see  \S 1.2{\it (b)\/}): 
so  $ \gerg $  and  $ \gerg^\times $  are Lie bialgebras  {\sl dual to
each other}.  In the sequel, we denote by  $ G^\star $  any connected
algebraic Poisson group with  $ \gerg $  as cotangent Lie bialgebra,
   \hbox{and say it is  {\sl (Poisson) dual\/}  to  $ G \, $.}
                                            \par
   For the Hopf operations in any Hopf algebra we shall use standard
notation, as in [Ab].

\vskip7pt

\noindent
{\bf Definition 1.2.}{\em
                                   \hfill\break
   \indent   (a) \, Let  $ H $,  $ K $  be Hopf algebras (in any
category).  A pairing  $ \; \langle \,\ , \,\ \rangle \, \colon \, H
\times K \loongrightarrow R \; $  (where  $ R $  is the ground ring)
is a  {\sl Hopf (algebra) pairing\/}  if  $ \;\; \big\langle x, y_1
\cdot y_2 \big\rangle = \big\langle \Delta(x), y_1 \otimes y_2
\big\rangle := \sum_{(x)} \big\langle x_{(1)}, y_1 \big\rangle \cdot
\big\langle x_{(2)}, y_2 \big\rangle \, $,  $ \; \big\langle x_1
\cdot x_2, y \big\rangle = \big\langle x_1 \otimes x_2, \Delta(y)
\big\rangle := \sum_{(y)} \big\langle x_1, y_{(1)} \big\rangle
\cdot \big\langle x_2, y_{(2)} \big\rangle  \, $,  $ \, \langle
x, 1 \rangle = \epsilon(x) \, $,  $ \; \langle 1, y \rangle =
\epsilon(y) \, $,  $ \; \big\langle S(x), y \big\rangle =
\big\langle x, S(y) \big\rangle \, $,  for all  $ \, x, x_1,
x_2 \in H $,  $ \, y, y_1, y_2 \in K $.
                                   \hfill\break
   \indent   (b) \, Let  $ \gerg $,  $ \gerh $  be Lie bialgebras (in
any category).  A pairing  $ \; \langle \,\ , \,\ \rangle \, \colon \,
\gerg \times \gerh \loongrightarrow \Bbbk \; $  (where  $ \Bbbk $  is the
ground ring) is  called a  {\sl Lie bialgebra pairing\/}  if  $ \;\;
\big\langle x, [y_1,y_2] \big\rangle = \big\langle \delta(x), y_1
\otimes y_2 \big\rangle := \sum_{[x]} \big\langle x_{[1]}, y_1
\big\rangle \cdot \big\langle x_{[2]}, y_2 \big\rangle  \, $,
$ \; \big\langle [x_1,x_2], y \big\rangle = \big\langle x_1 \otimes
x_2, \delta(y) \big\rangle := \sum_{[y]} \big\langle x_1, y_{[1]}
\big\rangle \cdot \big\langle x_2, y_{[2]} \big\rangle \, $,  \,
for all  $ \, x, x_1, x_2 \in \gerg \, $  and  $ \, y, y_1, y_2
\in \gerh $,  \, with  $ \, \delta(x) = \sum_{[x]} x_{[1]} \otimes
x_{[2]} \, $  and  $ \, \delta(y) = \sum_{[x]} y_{[1]} \otimes
y_{[2]} \, $.
}

\vskip7pt

  {\bf 1.3 The quantum setting.} \, Let  $ R $  be a 1-dimensional
(integral) domain (=1dD), and let  $ \, F = F(R) \, $  be its
quotient field.  Denote by  $ \calM $  the category of torsion-free
$ R $--modules,  and by  $ \HA $  the category of all Hopf algebras
in  $ \calM $.  Let  $ \calM_F $  be the category of  $ F $--vector
spaces, and  $ \HA_F $  be the category of all Hopf algebras in
$ \calM_F \, $.  For any  $ \, M \in \calM \, $,  \, set  $ \, M_F
:= F(R) \otimes_R M \, $.  Scalar extension gives a functor  $ \;
\calM \longrightarrow \calM_F \, $,  $ \, M \mapsto M_F \, $,  \,
which restricts to a functor  $ \; \HA \longrightarrow \HA_F \, $.
                                            \par
   Let  $ \, \h \in R \, $  be a non-zero prime element (which will be
fixed throughout), and  $ \, \Bbbk := R \big/ (\h) = R \big/ \h \, R
\, $  the corresponding quotient field.  For any  $ R $--module  $ M $,
we set  $ \, M_\h{\Big|}_{\h=0} \!\! := M \big/ \! \h \, M = \Bbbk
\otimes_R M \, $:  this is a  $ \Bbbk $--module  (via scalar
restriction  $ \, R \rightarrow R \big/ \h \, R =: \Bbbk \, $),
which we call the  {\sl specialization\/}  of  $ M $  at  $ \,
\h = 0 \, $;  we use also notation  $ \, M \,{\buildrel \, \h
\rightarrow 0 \, \over \llongrightarrow}\, \overline{N} \, $  to
mean that  $ \, M_\h{\Big|}_{\h=0} \hskip-3pt \cong \overline{N} \, $.
Moreover, set  $ \, M_\infty := \bigcap_{n=0}^{+\infty} \h^n M \, $ 
(this is the closure of  $ \{0\} $  in the  $ \h $--adic  topology
of  $ M $).  In addition, for any  $ H \! \in \! \HA \, $,  let 
$ \, I_{\! \scriptscriptstyle H} \! := \! \text{\sl Ker} \Big(
H \,{\buildrel \epsilon \over \twoheadrightarrow}\,
R \,{\buildrel {\h \mapsto 0} \over
{\relbar\joinrel\relbar\joinrel\twoheadrightarrow}}\,
\Bbbk \Big) \, $  and set  $ \, {I_{\!
\scriptscriptstyle H}}^{\!\!\infty} \!\! :=
\bigcap_{n=0}^{+\infty} {I_{\! \scriptscriptstyle H}}^{\!\!n} $.   
                                            \par
   Finally, given  $ \Bbb{H} $  in  $ \HA_F $,  a subset
$ \overline{H} $  of  $ \Bbb{H} $  is called  {\it an
$ R $--integer  form}  (or simply  {\it an  $ R $--form})  {\it of}
$ \, \Bbb{H} $  if  $ \; \overline{H} \, $  is a Hopf  $ R $--subalgebra
of  $ \, \Bbb{H} \, $  (hence  $ \, \overline{H} \in \HA \, $)  and 
$ \; H_F := F(R) \otimes_R \overline{H} = \Bbb{H} \, $.

\vskip7pt

\noindent
{\bf Definition 1.4.} {\em ({\sl ``Global quantum groups'' [or
``algebras'']})  Let  $ \, \h \in R \setminus \{0\} \, $  be a prime.
                                          \hfill\break
  \indent  (a) \, We call  {\sl quantized restricted universal
enveloping algebra (at  $ \h $)}  (in short,  {\sl QrUEA\/})  any
$ \U_\h \in \HA \, $  such that  $ \, \U_\h{\big|}_{\h=0} \! := \U_\h
\big/ \h \, \U_\h \, $  is (isomorphic to) the restricted universal
enveloping algebra  $ \U(\gerg) $  of some restricted Lie algebra
$ \gerg \, $.
                                         \hfill\break
   \indent   We call  $ \, \QrUEA \, $  the full subcategory
of  $ \HA $  whose objects are all the QrUEAs (at  $ \h $).
                                         \hfill\break
  \indent  (b) \, We call  {\sl quantized function algebra (at
$ \h $)}  (in short,  {\sl QFA\/})  any  $ \, F_\h \in \HA $
such that  $ \, {(F_\h)}_\infty = {I_{\!\scriptscriptstyle
F_\h}}^{\!\!\infty} \, $  (notation of \S 1.3)  and  $ \,
F_\h{\big|}_{\h=0} \! := F_\h \big/ \h \, F_\h \, $  is
(isomorphic to) the algebra of regular functions  $ F[G] $
of some connected algebraic group  $ G $.
                                         \hfill\break
   \indent   We call  $ \, \QFA \, $  the full subcategory
of  $ \HA $  whose objects are all the QFAs (at  $ \h $).
}

\vskip7pt

  {\bf Remark 1.5.} \, If  $ \, \U_\h \, $  is a QrUEA (at  $ \h \, $,
that is w.r.t.~to  $ \h \, $)  then  $ \, \U_\h{\big|}_{\h=0} \, $  is
a co-Poisson Hopf algebra, w.r.t.~the Poisson cobracket  $ \delta $
defined as follows: if  $ \, x \in \U_\h{\big|}_{\h=0} \, $  and  $ \,
x' \in \U_\h \, $  gives  $ \, x = x' \mod \h \, \U_\h \, $,  \,
then  $ \, \delta(x) := \big( \h^{-1} \, \big( \Delta(x') -
\Delta^{\text{op}}(x') \big) \big) \mod \h \, \big( \U_\h \otimes
\U_\h \big) \, $.  So  $ \, \U_\h{\big|}_{\h=0} \cong \U(\gerg) \, $
for some Lie algebra  $ \gerg $,  and by [Dr], \S 3, the restriction
of  $ \delta $  makes  $ \gerg $  into a  {\sl Lie bialgebra\/}  (the
isomorphism  $ \, \U_\h{\big|}_{\h=0} \cong \U(\gerg) \, $  being one
of  {\sl co-Poisson\/}  Hopf algebras); in this case we write  $ \,
\U_\h = \U_\h(\gerg) \, $.  Similarly, if  $ F_\h $  is a QFA at
$ \h $,  then  $ \, F_\h{\big|}_{\h=0} \, $  is a  {\sl Poisson\/}
Hopf algebra, w.r.t.~the Poisson bracket  $ \{\,\ ,\ \} $  defined
as follows: if  $ \, x $,  $ y \in F_\h{\big|}_{\h=0} \, $  and
$ \, x' $,  $ y' \in F_\h \, $  give  $ \, x = x' \mod \h \, F_\h
\, $,  $ \, y = y' \mod \h \, F_\h \, $,  \, then  $ \, \{x,y\} :=
\big( \h^{-1} (x' \, y' - y' \, x') \big) \mod \h \, F_\h \, $. 
So  $ \, F_\h{\big|}_{\h=0} \cong F[G] \, $  for some connected 
{\sl Poisson\/}  algebraic group  $ G $  (the isomorphism being
one of  {\sl Poisson\/}  Hopf algebras): in this case we write 
$ \, F_\h = F_\h[G] \, $.

\vskip7pt

\noindent {\bf Definition 1.6.}{\em 
                                           \hfill\break
   \indent   (a) \, Let  $ R $  be any (integral) domain, and let
$ F $  be its field of fractions.  Given two  $ F $--modules
$ \Bbb{A} $,  $ \Bbb{B} $,  and an  $ F $--bilinear  pairing
$ \, \Bbb{A} \times \Bbb{B} \longrightarrow F \, $,  for any
$ R $--submodule  $ \, A \subseteq \Bbb{A} \, $  and  $ \, B
\subseteq \Bbb{B} \, $  we set  $ \; \displaystyle{ A^\bullet
\, := \Big\{\, b \in \Bbb{B} \;\Big\vert\; \big\langle A, \, b
\big\rangle \subseteq R \Big\} } \; $  and  $ \; \displaystyle{
B^\bullet \, := \Big\{\, a \in \Bbb{A} \;\Big\vert\; \big\langle a,
B \big\rangle \subseteq R \Big\} } \, $.
                                           \hfill\break
   \indent   (b) \, Let  $ R $  be a 1dD.  Given  $ \, H $,  $ K \in
\HA \, $,  \, we say that  {\sl  $ H $  and  $ K $  are dual to each
other}  if there exists a perfect Hopf pairing between them for which
$ \, H = K^\bullet \, $  and  $ \, K = H^\bullet \, $.
}

\vskip1,1truecm

\centerline {\bf \S \; 2 \  The global quantum duality principle }

\vskip10pt

  {\bf 2.1 Drinfeld's functors.} \,  (Cf.~[Dr], \S 7) Let  $ R $,
$ \HA $  and  $ \, \h \in R \, $  be as in \S 1.3.  For any
$ \, H \in \HA \, $,  \, let  $ \, I = I_{\scriptscriptstyle
\! H} := \hbox{\sl Ker} \Big( H \,{\buildrel \epsilon \over
{\relbar\joinrel\twoheadrightarrow}}\, R \,{\buildrel {\h \mapsto 0}
\over \llongtwoheadrightarrow}\, R \big/ \h \, R = \Bbbk \Big) =
\hbox{\sl Ker} \Big( H \,{\buildrel {\h \mapsto 0} \over
\llongtwoheadrightarrow}\, H \big/ \h \, H \,{\buildrel
\bar{\epsilon} \over {\relbar\joinrel\twoheadrightarrow}}\,
\Bbbk \Big) \, $  (as in \S 1.3), a maximal Hopf ideal of  $ H $
(where  $ \bar{\epsilon} $  is the counit of  $ \, H{\big|}_{\h=0}
\, $,  \, and the two composed maps clearly coincide): we define
  $$  H^\vee \; := \; {\textstyle \sum\limits_{n \geq 0}} \, \h^{-n}
I^n \; = \; {\textstyle \sum\limits_{n \geq 0}} \, {\big( \h^{-1} I
\, \big)}^n \; = \; {\textstyle \bigcup\limits_{n \geq 0}} \, {\big(
\h^{-1} I \, \big)}^n  \quad  \big( \! \subseteq H_F \, \big) \; .  $$  
If  $ \, J = J_{\scriptscriptstyle \! H} := \hbox{\sl Ker}\,
(\epsilon_{\scriptscriptstyle \! H}) \, $  then  $ \, I = J +
\h \cdot 1_{\scriptscriptstyle H} \, $,  \, thus  $ \; H^\vee =
\sum_{n \geq 0} \h^{-n} J^n = \sum_{n \geq 0} {\big( \h^{-1} J \,
\big)}^n \; $  too.
                                              \par
  Given any Hopf algebra  $ H $,  for every  $ \, n \in \N \, $  define
$ \; \Delta^n \colon H \longrightarrow H^{\otimes n} \; $  by  $ \,
\Delta^0 := \epsilon \, $,  $ \, \Delta^1 := \id_{\scriptscriptstyle
H} $,  \, and  $ \, \Delta^n := \big( \Delta \otimes
\id_{\scriptscriptstyle H}^{\otimes (n-2)} \big)
\circ \Delta^{n-1} \, $  if  $ \, n > 2 $.  For any
ordered subset  $ \, \Sigma = \{i_1, \dots, i_k\} \subseteq
\{1, \dots, n\} \, $  with  $ \, i_1 < \dots < i_k \, $,  \, define
the morphism  $ \; j_{\scriptscriptstyle \Sigma} : H^{\otimes k}
\longrightarrow H^{\otimes n} \; $  by  $ \; j_{\scriptscriptstyle
\Sigma} (a_1 \otimes \cdots \otimes a_k) := b_1 \otimes \cdots
\otimes b_n \; $  with  $ \, b_i := 1 \, $  if  $ \, i \notin
\Sigma \, $  and  $ \, b_{i_m} := a_m \, $  for  $ \, 1 \leq m
\leq k $~;  then set  $ \; \Delta_\Sigma := j_{\scriptscriptstyle
\Sigma} \circ \Delta^k \, $,  $ \, \Delta_\emptyset := \Delta^0
\, $,  and  $ \; \delta_\Sigma := \sum_{\Sigma' \subset \Sigma}
{(-1)}^{n- \left| \Sigma' \right|} \Delta_{\Sigma'} \, $,
$ \; \delta_\emptyset := \epsilon \, $.  By the inclusion-exclusion
principle, this definition admits the inverse formula  $ \; \Delta_\Sigma
= \sum_{\Psi \subseteq \Sigma} \delta_\Psi \, $.  We shall also use the
notation  $ \, \delta_0 := \delta_\emptyset \, $,  $ \, \delta_n :=
\delta_{\{1, 2, \dots, n\}} \, $,  and the useful formula  $ \; \delta_n
= {(\id_{\scriptscriptstyle H} - \epsilon)}^{\otimes n} \circ \Delta^n
\, $,  \, for all  $ \, n \in \N_+ \, $.
                                              \par
   Now consider any  $ \, H \in \HA \, $  and  $ \, \h \in R \, $  as
in \S 1.3: we define
  $$  H' \; := \; \big\{\, a \in H \,\big\vert\, \delta_n(a) \in
\h^n H^{\otimes n} ,  \; \forall \,\, n \in \N \, \big\}  \quad
\big( \! \subseteq H \, \big) \, .  $$

\vskip7pt

\noindent
 {\bf Theorem 2.2} \, 
{\em ({\sl ``The Global Quantum Duality Principle''})
                                        \hfill\break
  \indent   (a) \, The assignment  $ \, H \mapsto H^\vee \, $,
resp.~$ \, H \mapsto H' \, $,  defines a functor  $ \; {(\ )}^\vee
\colon \, \HA \, \longrightarrow \, \HA \, $,  \, resp.~$ \; {(\ )}'
\colon \, \HA \, \longrightarrow \, \HA \, $,  \, whose image lies in
$ \QrUEA $,  resp.~in  $ \QFA $.  In particular, when  $ \, \hbox{\it
Char}\,(\Bbbk) > 0 \, $  the algebraic Poisson group  $ G $  such that
$ \; H'{\big|}_{\h=0} = F[G] \; $  is zero-dimensional of height 1.
Moreover, for all  $ \, H \in \HA \, $  we have  $ \, H \subseteq
{\big( H^\vee \big)}' \, $  and  $ \, H \supseteq {\big( H'
\big)}^{\!\vee} \! $,  hence also  $ \, H^\vee = \big(
\big(H^\vee\big)' \,\big)^{\!\vee} $  and  $ \, H' =
\big( \big(H'\big)^{\!\vee} \big)' $.
                                        \hfill\break
   \indent   (b) \, Let  $ \, \hbox{\it Char}\,(\Bbbk) = 0 \, $.
Then for any  $ \, H \in \HA \, $  one has
 \vskip1pt
   \centerline{ $ \displaystyle{ H = {\big(H^\vee\big)}'
\,\Longleftrightarrow\, H \in \QFA  \qquad  \hbox{and}
\qquad  H = {\big( H' \big)}^{\!\vee} \,\Longleftrightarrow\,
H \in \QrUEA \, , } $ }
 \vskip1pt
\noindent   thus we have two induced equivalences, namely
$ \; {(\ )}^\vee \colon \, \QFA \, \llongrightarrow
\, \QrUEA \, $,  $ \, H \mapsto H^\vee \, $,  \; and
$ \; {(\ )}' \colon \, \QrUEA \, \llongrightarrow \, \QFA \, $,
$ \, H \mapsto H' \, $,  \; which are inverse to each other.
                                        \hfill\break
  \indent   (c) \, (``Quantum Duality Principle'') Let
$ \, \hbox{\it Char}\,(\Bbbk) = 0 \, $.  Then
 \vskip5pt
   \centerline{ $ {F_\h[G]}^\vee{\Big|}_{\h=0} := {F_\h[G]}^\vee \Big/
\h \, {F_\h[G]}^\vee \, = \, U(\gerg^\times) \; ,  \quad
\displaystyle{ {U_\h(\gerg)}'{\Big|}_{\h=0} := {U_\h(\gerg)}' \Big/
\h \, {U_\h(\gerg)}' \, = \, F\big[G^\star\big] } $ }
 \vskip5pt
\noindent   (with  $ G $,  $ \gerg $,  $ \gerg^\times $,
$ \gerg^\star $  and  $ G^\star $  as in \S 1.1, and  $ U_\h(\gerg) $
has the obvious meaning, cf.~\S 1.5) where the choice of the group
$ G^\star $  (among all the connected Poisson algebraic groups with
tangent Lie bialgebra  $ \gerg^\star $)  depends on the choice of the
QrUEA  $ U_\h(\gerg) $.  In other words,  $ \, {F_\h[G]}^\vee \, $
is a QrUEA for the Lie bialgebra  $ \gerg^\times $,  and  $ \,
{U_\h(\gerg)}' \, $  is a QFA for the Poisson group  $ G^\star $.
                                        \hfill\break
  \indent   (d) \, Let  $ \, \text{\it Char}\,(\Bbbk) = 0 \, $.  Let
$ \, F_\h \in \QFA \, $,  $ \, U_\h \in \QrUEA \, $  be dual to each
other (with respect to some pairing).  Then  $ {F_\h}^{\!\vee} $  and
$ {U_\h}' $  are dual to each other (w.r.t.~the  {\sl same}  pairing).
                                        \hfill\break
  \indent   (e) \, Let  $ \, \hbox{\it Char}\,(\Bbbk) = 0 \, $.
Then for all  $ \, \H \in \HA_F \, $  the following are equivalent:
                                        \hfill\break
   \indent \indent  \;  $ \H $  has an  $ R $--integer  form
$ H_{(f)} $  which is a QFA at  $ \h \, $;
                                        \hfill\break
   \indent \indent  \;  $ \H $  has an  $ R $--integer  form
$ H_{(u)} $  which is a QrUEA at  $ \h \, $.
}

\vskip7pt

{\bf Remarks 2.3.}  After stating our main theorem, some
comments are in order.
                                                 \par
   {\it (a)} \, {\sl The Global Quantum Duality Principle as
a ``Galois correspondence'' type theorem.}
                                                 \par
\noindent   Let  $ \, L \subseteq E \, $  be a Galois (not
necessarily finite) field extension, and let  $ \, G := \hbox{\it
Gal}\,\big(E/L\big) \, $  be its Galois group.  Let  $ \, \mathcal{F}
\, $  be the set of intermediate extensions (i.e.~all fields
$ F $  such that  $ \, L \subseteq F \subseteq E \, $),  \, let
$ \, \mathcal{S} \, $  be the set of all subgroups of  $ G $  and let
$ \, \mathcal{S}^c \, $  be the set of all subgroups of  $ G $  which
are  {\sl closed\/}  w.r.t.~the Krull topology of  $ G $.  Note that 
$ \mathcal{F} $,  $ \mathcal{S} $  and  $ \mathcal{S}^c $  can
all be seen as lattices w.r.t.~set-theoretical inclusion   ---
$ \mathcal{S}^c $  being a sublattice of  $ \mathcal{S} $  ---   hence
as categories too.  The celebrated Galois Theorem provides two maps,
namely  $ \; \varPhi \, \colon \, \mathcal{F} \loongrightarrow \mathcal{S}
\, $,  $ \, F \mapsto \hbox{\it Gal}\,\big(E/F\big) := \big\{\, \gamma
\in G \;\big\vert\;\, \gamma{\big\vert}_F = \hbox{id}_F \,\big\} \, $,
and  $ \; \varPsi \, \colon \, \mathcal{S} \loongrightarrow \mathcal{F}
\, $,  $ \, H \mapsto E^H := \big\{\, e \in E \;\big\vert\; \eta(e) = e
\;\; \forall \; \eta \in H \,\big\} \, $,
%
%
%
%
such that:
                                      \par
   {\it --- 1)} \;  $ \varPhi $  and  $ \varPsi $  are  {\sl
contravariant\/}  functors (that is, they are order-reversing maps
of lattices, i.e.{} lattice antimorphisms); moreover, the image of 
$ \varPhi $  lies in the subcategory  $ \, \mathcal{S}^c \, $;
                                      \par
   {\it --- 2)} \; for  $ \, H \in \mathcal{S} \, $  one has  $ \,
\varPhi\big(\varPsi(H)\big) = \overline{H} \, $,  \, the  {\sl
closure\/}  of  $ H $  w.r.t.~the Krull topology: thus  $ \, H
\subseteq \varPhi \big( \varPsi(H) \big) \, $,  \, and  $ \,
\varPhi \circ \varPsi \, $  is a  {\sl closure operator},  so
that   $ \, H \in \mathcal{S}^c \, $  iff  $ \, H = \varPhi
\big( \varPsi(H) \big) \, $;   
                                      \par
   {\it --- 3)} \; for  $ \, F \in \mathcal{F} \, $  one has 
$ \, \varPsi\big(\varPhi(F)\big) = F \, $;   
                                      \par
   {\it --- 4)} \;  $ \varPhi $  and  $ \varPsi $  restrict to
antiequivalences  $ \, \varPhi : \mathcal{F} \rightarrow \mathcal{S}^c
\, $  and  $ \, \varPsi : \mathcal{S}^c \rightarrow \mathcal{F} \, $ 
which are inverse to each other.
                                      \par
   Then one can see that Theorem 2.2 establishes a strikingly
similar result, which in addition is much more symmetric:  $ \HA $ 
plays the role of both  $ \mathcal{F} $  and  $ \mathcal{S} $,
whereas  $ {(\ )}' $  stands for  $ \varPsi $  and  $ {(\ )}^\vee $
stands for  $ \varPhi $.  $ \QFA $  plays the role of the
distinguished subcategory  $ \mathcal{S}^c $,  and symmetrically
we have the distinguished subcategory  $ \QrUEA $.  The composed
operator  $ \, {\big({(\ )}^\vee\big)}' = {(\ )}' \circ {(\ )}^\vee
\, $  plays the role of a ``closure operator'', and symmetrically
$ \, {\big( {(\ )}' \big)}^\vee = {(\ )}^\vee \circ {(\ )}' \, $
plays the role of a ``taking-the-interior operator'': in other
words, QFAs may be thought of as ``closed sets''  and QrUEAs
as ``open sets'' in  $ \HA \, $.  Yet note also that now all
involved functors are  {\sl covariant}.   
                                                 \par
   {\it (b)} \, {\sl Duality between Drinfeld's functors}.
For any  $ \, n \in \N \, $  let  $ \; \mu_n \, \colon
\, {J_{\scriptscriptstyle H}}^{\!\otimes n}
\lhook\joinrel\longrightarrow H^{\otimes n} \,{\buildrel
{\; m^n} \over \loongrightarrow}\, H \, $  be the composition of the
natural embedding of  $ {J_{\scriptscriptstyle H}}^{\!\otimes n} $
into  $ H^{\otimes n} $  with the  $ n $--fold  multiplication (in
$ H \, $):  then  $ \mu_n $  is the ``Hopf dual'' to  $ \delta_n \, $.
By construction we have  $ \, H^\vee = \sum_{n \in \N} \mu_n\big(
\h^{-n} {J_{\scriptscriptstyle H}}^{\!\otimes n}\big) \, $  and
$ \, H' = \bigcap_{n \in \N} {\delta_n}^{\!-1}\big(\h^{+n}
{J_{\scriptscriptstyle H}}^{\!\otimes n} \big) \, $:  \, this
shows that the two functors are built up as ``dual'' to each
other (see also part  {\it (d)\/}  of Theorem 2.2).
                                                 \par
   {\it (c)} \, {\sl Ambivalence \;
       \hbox{QrUEA  $ \leftrightarrow $
QFA  \; in  $ \HA_F \, $.}} \; Part  {\it (e)} of
Theorem 2.2 means that some Hopf algebras  {\sl over\/}  $ F(R) $ 
might be thought of  {\sl both\/}  as ``quantum function algebras'' 
{\sl and\/}  as ``quantum enveloping algebras'': examples are  $ U_F $ 
and  $ F_F $  for  $ \, U \in \QrUEA \, $  and  $ \, F \in \QFA \, $.
                                                 \par
   {\it (d)} \, {\sl Drinfeld's functors for algebras, coalgebras and
bialgebras}.  The definition of either of Drinfeld functors requires
only ``half of'' the notion of Hopf algebra.  In fact, one can define
$ (\ )^\vee $  for all ``augmented algebras'' (that is, roughly
speaking, ``algebras with a counit'') and  $ (\ )' $  for all
``coaugmented coalgebras'' (roughly, ``coalgebras with a unit''),
and in particular for bialgebras: this yields again nice functors,
and neat results extending the global quantum duality principle,
cf.~[Ga5], \S\S 3--4.   
                                                 \par
   {\it (e)} \, {\sl Relaxing the assumptions}.  We chose to work over
$ \HA $  for simplicity: in fact, this ensures that the specialization
functor  $ \, H \mapsto H{\big|}_{\h=0} \, $  yields Hopf algebras over
a field, so that we can use the more elementary geometric language of
algebraic groups and Lie algebras in the easiest sense.  Nevertheless,
what is really necessary to let the machine work is to consider any
(commutative, unital) ring  $ R $,  any  $ \, \h \in R \, $  and then
define Drinfeld's functors over Hopf  $ R $--algebras  {\sl which are
$ \h $--torsionless}.  For instance, this is   --- essentially ---  
what is done in [KT], where the ground ring is  $ \, R = \Bbbk[[u,v]]
\, $,  \, and the role of  $ \h $  is played by either  $ u $  or 
$ v \, $.  In general, working in such a more general setting amounts
to consider, at the semiclassical level (i.e.~after specialization), 
{\sl Poisson group schemes over  $ R \big/ \h \, R $}  (i.e.~over 
$ {\text{\it Spec}}\big( R / \h \, R \big) \, $)  and  {\sl Lie 
$ R \big/ \h \, R $--bialgebras},  where $ \, R / \h \, R \, $ 
might not be a field.
                                                 \par
   Similar considerations   --- about  $ R $  and  $ \h $  ---
hold w.r.t.~remark  {\it (d)\/}  above.

\vskip1,1truecm

\centerline {\bf \S \; 3 \  Application to trivial deformations:
the Crystal Duality Principle }

\vskip10pt

  {\bf 3.1 Drinfeld's functors on trivial deformations.} \, Let
$ \HA_\Bbbk $  be the category of all Hopf algebras over the field
$ \Bbbk $.  For all  $ \, n \in \N \, $,  \, let  $ \; J^n := {\big(
\text{\sl Ker}\,(\epsilon \colon \, H \longrightarrow \Bbbk) \big)}^n
\, $  and  $ \; D_n := \text{\sl Ker}\, \big( \delta_{n+1} \colon \, H
\longrightarrow H^{\otimes n} \big) \, $,  \; and set  $ \; \underline{J}
:= {\big\{ J^n \big\}}_{n \in \N} \, $,  $ \; \underline{D} := {\big\{
D_n \big\}}_{n \in \N} \, $.  Of course  $ \underline{J} $  is a
decreasing filtration of  $ H \, $  (maybe with  $ \, \bigcap_{n
\geq 0} J^n \supsetneqq \{0\} \, $),  and  $ \underline{D} $  is an
increasing filtration of  $ H \, $  (maybe with  $ \, \bigcup_{n \geq
0} D_n \subsetneqq H \, $), by coassociativity of the  $ \delta_n $'s.
                                      \par
   Let  $ \, R := \Bbbk[\h\,] \, $  be the polynomial ring in the
indeterminate  $ \h \, $:  then  $ R $  is a PID (= \, principal
ideal domain), hence a 1dD, and  $ \h $  is a non-zero prime in
$ R \, $.  Let  $ \, H_\h := H[\h] = R \otimes_{\Bbbk} H \, $,  the
scalar extension of  $ H \, $:  \, this is a torsion free Hopf algebra
over  $ R $,  hence one can apply Drinfeld's functors to  $ \, H_\h
\, $;  in this section we do that with respect to the prime  $ \, \h
\, $  itself.  We shall see that the outcome is quite neat, and can
be expressed purely in terms of Hopf algebras in  $ \HA_\Bbbk \, $:
because of the special relation between some features of  $ H $
--- namely, the filtrations  $ \underline{J} $  and  $ \underline{D}
\, $  ---   and some properties of Drinfeld's functors, we call this
result ``Crystal Duality Principle'', in that it is obtained through
sort of a ``crystallization'' process (bearing in mind, in a sense,
Kashiwara's motivation for the terminology ``crystal bases'' in the
context of quantum groups: see [CP], \S 14.1, and references therein).
Indeed, this theorem can also be proved almost entirely by using only
classical Hopf algebraic methods within  $ \HA_\Bbbk $,  i.e.~without
resorting to deformations: this is accomplished in [Ga6].  We first
discuss the general situation (\S\S 3.2--5), second we look at the
case of function algebras and enveloping algebras (\S\S 3.6--7), then
we state and prove the theorem of Crystal Duality Principle (\S 3.9). 
Eventually (\S\S 3.11--12) we dwell upon two other interesting
applications: hyperalgebras, and group algebras and their dual.
                                            \par
  Note that the same analysis and results (with only a few more
details to take care of) still hold if we take as  $ R $  any 1dD
and as  $ \, \h \, $  any prime element in  $ R $  such that  $ \,
R \big/ \h \, R = \Bbbk \, $  and  $ R $  carries a structure of 
$ \Bbbk $--algebra;  \, for instance, one can take  $ \, R = \Bbbk[[h]]
\, $  and  $ \, \h = h \, $,  or  $ \, R = \Bbbk \big[ q, q^{-1} \big]
\, $  and  $ \h = q-1 \, $.  Finally, in the sequel to be short we
perform our analysis for Hopf algebras only: however, as Drinfeld's
functors are defined not only for Hopf algebras but for augmented
algebras and coaugmented coalgebras too, we might do the same study
for them as well.  In particular, the Crystal Duality Principle has
a stronger version which concerns these more general objects too
(cf.~[Ga6]).

\vskip7pt

\noindent
{\bf Lemma 3.2}{\em
$$\begin{array}{llr}
   {H_\h}^{\hskip-3pt\vee} \;  &  = \; {\textstyle \sum\limits_{n
\geq 0}} \, R \cdot \h^{-n} J^n \, = \; R \cdot J^0 + R \cdot \h^{-1}
J^1 + \cdots + R \cdot \h^{-n} J^n + \cdots  \;  &   (3.1)  \\
   {H_\h}' \;  &  = \; {\textstyle \sum\limits_{n \geq 0}} \, R \cdot
\h^{+n} D_n \, = \; R \cdot D_0 + R \cdot \h^{+1} \, D_1 + \cdots +
R \cdot \h^{+n} D_n + \cdots  &   (3.2) 
\end{array}$$
}

\noindent
{\it Sketch of proof.} \, (3.1) follows directly from definitions, while
(3.2) is an easy exercise.   \qed

\vskip7pt

  {\bf 3.3 Rees Hopf algebras and their specializations.} \, Let
$ M $  be a module over a commutative unitary ring  $ R $,  and let
$ \; \underline{M} := {\{M_z\}}_{z \in \Z} = \Big( \cdots \subseteq
M_{-m} \subseteq \cdots \subseteq M_{-1} \subseteq M_0 \subseteq M_1
\subseteq \cdots \subseteq M_n \subseteq \cdots \Big) \; $  be a
bi-infinite filtration of  $ M $  by submodules  $ M_z $  ($ z \in
\Z $).  In particular, we consider increasing filtrations (i.e.,
those with  $ \, M_z = \{0\} \, $  for  $ \, z < 0 \, $)  and
decreasing filtrations (those with  $ \, M_z = \{0\} \, $  for
all  $ \, z > 0 \, $)  as special cases of bi-infinite filtrations.
First we define the associated  {\sl blowing module\/}  to be the
$ R $--submodule  $ \mathcal{B}_{\underline{M}}(M) $  of  $ M \big[ t,
t^{-1} \big] $  (where  $ t $  is any indeterminate) given by  $ \,
\mathcal{B}_{\underline{M}}(M) := \sum_{z \in \Z} t^z M_z \, $;  \,
this is isomorphic to the
         {\sl first graded module\,\footnote{
          I pick
this terminology from Serge Lang's textbook  {\it ``Algebra''}.}
associated to  $ M $},  namely  $ \, \bigoplus_{z \in \Z} M_z \, $.
Second, we define the associated  {\sl Rees module\/}  to be the
$ R[t] $--submodule  $ \mathcal{R}^t_{\underline{M}}(M) $  of  $ M \big[
t, t^{-1} \big] $  generated by  $ \mathcal{B}_{\underline{M}}(M) $;  \,
straight\-forward computations then give  $ R $--module  isomorphisms
  $$  \mathcal{R}^t_{\underline{M}}(M) \Big/ (t-1)
\, \mathcal{R}^t_{\underline{M}}(M) \; \cong \;
{\textstyle \bigcup\limits_{z \in \Z}} M_z \; , 
\qquad  \mathcal{R}^t_{\underline{M}}(M) \Big/ t \,
\mathcal{R}^t_{\underline{M}}(M) \; \cong \; G_{\underline{M}}(M)  $$
where  $ \, G_{\underline{M}}(M) := \bigoplus_{z \in Z} M_z \big/
M_{z-1} \, $  is the  {\sl second
                graded module$\,{}^1$  associated
to  $ M $}.  In other words,  $ \mathcal{R}^t_{\underline{M}}(M) $
is an  $ R[t] $--module  which specializes to  $ \, \bigcup_{z
\in \Z} M_z \, $  for  $ \, t = 1 \, $  and specializes to  $ \,
G_{\underline{M}}(M) \, $  for  $ \, t = 0 \, $;  \, therefore
the  $ R $--modules  $ \, \bigcup_{z \in \Z} M_z \, $  and  $ \,
G_{\underline{M}}(M) \, $  can be seen as 1-parameter (polynomial)
deformations of each other via the 1-parameter family of
$ R $--modules  given by  $ \mathcal{R}^t_{\underline{M}}(M) $.
   We can repeat this construction within the category of algebras,
coalgebras, bialgebras or Hopf algebras over  $ R $  with a filtration
in the proper sense:
%
%
then we'll end up with corresponding objects
$ \mathcal{B}_{\underline{M}}(M) $,  $ \mathcal{R}^t_{\underline{M}}(M) $,
etc.{} of the like type (algebras, coalgebras, etc.).  In particular
we'll deal with Rees Hopf algebras.

\vskip7pt

  {\bf 3.4 Drinfeld's functors on  $ H_\h \, $  and filtrations on
$ H \, $.} \, Lemma 3.2 sets a link between properties of  $ {H_\h}' $,
resp.~of  $ {H_\h}^{\!\vee} $,  and properties of the filtration
$ \underline{D} \, $,  resp.~$ \underline{J} \, $,  \, of  $ H \, $.
                                              \par
   First, (3.1) together with  $ \, {H_\h}^{\!\vee} \in \HA \, $
implies that  $ \, \underline{J} \, $  is a Hopf algebra filtration
of  $ H \, $;  \, conversely, if one proves that  $ \, \underline{J}
\, $  is a Hopf algebra filtration of  $ H \, $  (which is
straightforward) then from (3.1) we get a one-line proof that
$ \, {H_\h}^{\!\vee} \in \HA \, $.  Second, we can look at
$ \underline{J} $  as a bi-infinite filtration, reversing index
notation and extending trivially on positive indices,  $ \;
\underline{J} \, = \, \Big( \cdots \subseteq J^n \subseteq \cdots
J^2 \subseteq J \subseteq J^0 \big( = H \big) \subseteq H \subseteq
\cdots \subseteq H \subseteq \cdots \Big) \, $;  \; then the Rees Hopf
algebra  $ \mathcal{R}^\h_{\underline{J}}(H) $  is defined (see \S 3.3).
Now (3.1) give  $ \, {H_\h}^{\!\vee} = \mathcal{R}^\h_{\underline{J}}(H)
\, $,  \, so  $ \; {H_\h}^{\!\vee} \Big/ \h \, {H_\h}^{\!\vee} \cong
\mathcal{R}^\h_{\underline{J}}(H) \Big/ \h \, \mathcal{R}^\h_{\underline{J}}(H)
\cong G_{\underline{J}}(H) \, $.  Thus  $ \, G_{\underline{J}}(H) \, $
is cocommutative because  $ \, {H_\h}^{\!\vee} \Big/ \h \, {H_\h}^{\!
\vee} \, $  is; conversely, we get an easy proof of the cocommutativity
of  $ \, {H_\h}^{\!\vee} \Big/ \h \, {H_\h}^{\!\vee} \, $  once we
prove that  $ \, G_{\underline{J}}(H) \, $  is cocommutative, which
is straightforward.  Finally,  $ \, G_{\underline{J}}(H) \, $  is
generated by  $ \, Q(H) = J \big/ J^{\,2} \, $  whose elements are
primitive, so  {\sl a fortiori\/}  $ \, G_{\underline{J}}(H) \, $  is
generated by its primitive elements; then the latter holds for  $ \,
{H_\h}^{\!\vee} \Big/ \h \, {H_\h}^{\!\vee} \, $  as well.  To sum
up, as  $ {H_\h}^{\!\vee} \in \QrUEA \, $  we argue that  $ \,
G_{\underline{J}}(H) = \U(\gerg) \, $  for some restricted Lie
bialgebra  $ \gerg \, $;  conversely, we can get  $ {H_\h}^{\!\vee}
\in \QrUEA \, $  directly from the properties of the filtration
$ \underline{J} \, $  of  $ H $.  Moreover, since  $ \,
G_{\underline{J}}(H) = \U(\gerg) \, $  is graded,
$ \gerg $  as a restricted Lie algebra is 
{\sl graded\/}  too.
                                             \par
   On the other hand, it is easy to see that (3.2) and  $ \, {H_\h}'
\in \HA \, $  imply that  $ \, \underline{D} \, $  is a Hopf algebra
filtration of  $ H \, $;  \, conversely, if one shows that  $ \,
\underline{D} \, $  is a Hopf algebra filtration of  $ H \, $  (which
can be done) then (3.2) yields a direct proof that  $ \, {H_\h}' \in
\HA \, $.  Second, we can look at  $ \underline{D} $  as a bi-infinite
filtration, extending it trivially on negative indices, namely  $ \;
\underline{D} \, = \, \Big( \cdots \subseteq \{0\} \subseteq \cdots
\{0\} \subseteq \big( \{0\} = \big) D_0 \subseteq D_1 \subseteq \cdots
\subseteq D_n \subseteq \cdots \Big) \, $;  \; then the Rees Hopf
algebra  $ \mathcal{R}^\h_{\underline{D}}(H) $  is defined
(see \S 3.3).  Now (3.2) gives  $ \, {H_\h}' =
\mathcal{R}^\h_{\underline{D}}(H) \, $;  \, but
then  $ \; {H_\h}' \Big/ \h \, {H_\h}' \cong
\mathcal{R}^\h_{\underline{D}}(H) \Big/ \h \,
\mathcal{R}^\h_{\underline{D}}(H) \cong G_{\underline{D}}(H)
\, $.  Thus  $ \, G_{\underline{D}}(H) \, $  is commutative
because  $ \, {H_\h}' \Big/ \h \, {H_\h}' \, $  is; viceversa, we
get an easy proof of the commutativity of  $ \, {H_\h}' \Big/ \h
\, {H_\h}' \, $  once we prove that  $ \, G_{\underline{D}}(H)
\, $  is commutative (which can be done too).  Finally,  $ \,
G_{\underline{D}}(H) \, $  is graded with  $ 1 $-dimensional 
degree 0 component (by construction) hence it has no non-trivial
idempotents; so the latter is true for  $ \, {H_\h}' \Big/ \h
\, {H_\h}' \, $  too.  Note also that  $ \, {I_{\!{H_\h}'}
\phantom{\big|}}^{\hskip-11pt\infty} = \{0\} \, $  by construction
(because  $ H_\h $  is free over  $ R \, $).  To sum up, since 
$ {H_\h}' \in \QFA \, $  we get  $ \, G_{\underline{D}}(H) = F[G] \, $  for some connected algebraic Poisson group  $ G \, $;
conversely, we can argue that  $ \, {H_\h}' \in \QFA \, $  directly
from the properties of the filtration  $ \underline{D} \, $.
                                             \par
   In addition, since  $ \, G_{\underline{D}}(H) = F[G] \, $  is
graded, when  $ \, \Char(\Bbbk) = 0 \, $  the (pro)affine variety
$ \, G_{(\text{\it cl\/})} \, $  of closed points of  $ G $  is a
          (pro)affine space\,\footnote
          {\,For it is a  {\sl
cone}   --- since  $ H $  is graded ---   without vertex   ---
since  $ G_{(\text{\it cl\/})} $,  being a group, is smooth.},
that is  $ \, G_{(\text{\it cl\/})} \cong \Bbb{A}_{\,\Bbbk}^{\times
\mathcal{I}} = \Bbbk^{\mathcal{I}} \, $  for some index set 
$ \mathcal{I} \, $,  \, and so  $ \, F[G] = \Bbbk \big[
{\{x_i\}}_{i \in \mathcal{I}} \big] \, $  is a
polynomial algebra.   
                                             \par
   Finally, when  $ \, p := \Char(\Bbbk) > 0 \, $  the group  $ G $
has dimension 0 and height~1: indeed, we can see this as a consequence
of part of  Theorem 2.2{\it (a)\/}  via the identity  $ \; H'_\h \Big/
\h \, H'_\h = G_{\underline{D}}(H) \, $,  \, or conversely we can prove
the relevant part of  Theorem 2.2{\it (a)\/}   via this identity by
observing that  $ G $  has those properties (cf.~[Ga5], \S 5.4).  At
last, by general theory since  $ G $  has dimension 0 and height~1 the
function algebra  $ \, F[G] = G_{\underline{D}}(H) = H'_\h \Big/ \h \,
H'_\h \, $  is a  {\sl truncated polynomial algebra},  namely of type 
$ \, F[G] = \Bbbk \big[ {\{x_i\}}_{i \in \mathcal{I}} \big] \Big/
\big( {\{x_i^{\,p}\}}_{i \in \mathcal{I}} \big) \, $  for some
index set  $ \mathcal{I} \, $.   

\vskip7pt

  {\bf 3.5 Special fibers of  $ {H_\h}' $  and  $ {H_\h}^{\!\vee} $
and deformations.}  \, Given  $ \, H \in \HA_\Bbbk \, $,  \, consider
$ \, H_\h \, $:  \, our goal is to study  $ \, {H_\h}^{\!\vee} \, $
and  $ \, {H_\h}' \, $.
                                            \par
   As for  $ {H_\h}^{\!\vee} $,  the natural map from  $ H $
to  $ \, \widehat{H} := G_{\underline{J}}(H) = {H_\h}^{\!\vee}
\Big/ \h \, {H_\h}^{\!\vee} =: {H_\h}^{\!\vee}{\Big|}_{\h=0} \, $
sends  $ \, {J\phantom{|}}^{\hskip-1pt \infty} := \bigcap_{n \geq 0}
{J\phantom{|}}^{\hskip-1pt n} \, $  to zero, by definition; also,
letting  $ \, H^\vee := H \Big/ {J\phantom{|}}^{\hskip-1pt \infty}
\, $  (a Hopf algebra quotient of  $ H $,  for  $ \underline{J} $
is a Hopf algebra filtration), we have  $ \, \widehat{H} =
\widehat{H^\vee} \, $.  Thus  $ \, {{(H^\vee)}_\h}^{\!\!\vee}
{\Big|}_{\h=0} = \widehat{H^\vee} = \widehat{H} = \U(\gerg_-) \, $
for some graded restricted Lie bialgebra  $ \gerg_- \, $;  \, also,
$ \, {{(H^\vee)}_\h}^{\!\!\vee}{\Big|}_{\h=1} := {{(H^\vee)}_\h}^{\!
\!\vee} \Big/ (\h-1) \, {{(H^\vee)}_\h}^{\!\!\vee} = \sum_{n \geq 0}
\overline{J}^{\,n} = H^\vee \, $  (see \S 3.3).  Thus we can see  $ \,
{{(H^\vee)}_\h}^{\!\!\vee} = \mathcal{R}^\h_{\underline{J}} (H^\vee)
\, $  as a 1-parameter family inside  $ \HA_\Bbbk $  with  {\sl
regular\/}  fibers (that is, they are isomorphic to each other
as  $ \Bbbk $--vector  spaces; indeed, we switch from  $ H $
to  $ H^\vee $  just to achieve this regularity) which links
$ \widehat{H^\vee} $  and  $ H^\vee $  as (polynomial)
deformations of each other, namely
  $$  \U(\gerg_-) = \widehat{H^\vee} = {{(H^\vee)}_\h}^{\!\!\vee}
{\Big|}_{\h=0}  \hskip5pt  \underset{\;{{(H^\vee)}_\h}^{\!\!\vee}}
 {\overset{0 \,\leftarrow\, \h \,\rightarrow\, 1}
{\longleftarrow\joinrel\relbar\joinrel%
\relbar\joinrel\relbar\joinrel\relbar\joinrel\llongrightarrow}}
\hskip6pt  {{(H^\vee)}_\h}^{\!\!\vee}{\Big|}_{\h=1} = H^\vee \; .  $$
   \indent   Now look at  $ \, {\big( {{(H^\vee)}_\h}^{\!\!\vee}
\big)}' \, $.  By construction,  $ \, {\big( {{(H^\vee)}_\h}^{\!
\!\vee} \big)}'{\Big|}_{\h=1} \! = {{( H^\vee )}_\h}^{\!\!\vee}
{\Big|}_{\h=1} \! = H^\vee \, $,
  \, whereas\break
   $ \, {\big(
{{(H^\vee)}_\h}^{\!\!\vee} \big)}'{\Big|}_{\h=1} \! = F[K_-] \, $
for some connected algebraic Poisson group  $ K_- \, $:  in addition,
if  $ \, \text{\it Char}\,(\Bbbk) = 0 \, $  then  $ \, K_- = G_-^\star
\, $  by  Theorem 2.2{\it (c)}.  So  $ \, {\big( {{(H^\vee)}_\h}^{\!
\!\vee} \big)}' \, $  can be thought of as a 1-parameter family inside
$ \HA_\Bbbk \, $,  with regular fibers, linking  $ H^\vee $  and
$ F[G_-^\star] $  as (polynomial) deformations of each other, namely   
  $$  H^\vee = {\big( {{(H^\vee)}_\h}^{\!\!\vee} \big)}'{\Big|}_{\h=1}
\hskip-4pt  \underset{{({{(H^\vee)}_\h}^{\!\!\vee})}'}
{\overset{1 \,\leftarrow\, \h \,\rightarrow\, 0}
{\longleftarrow\joinrel\relbar\joinrel%
\relbar\joinrel\relbar\joinrel\relbar\joinrel\llongrightarrow}}
\hskip0pt  {\big( {{(H^\vee)}_\h}^{\!\!\vee} \big)}'{\Big|}_{\h=0}
\hskip-5pt  = F[K_-]
\hskip9pt  \Big( = F[G_-^\star]  \,\text{\ if \ } 
\text{\it Char}\,(\Bbbk) = 0 \Big) \, .  $$
Therefore  $ H^\vee $  is  {\sl both\/}  a deformation of an enveloping
algebra {\sl and\/}  a deformation of a function algebra, via two
different 1-parameter families (with regular fibers) in  $ \HA_\Bbbk $ 
which match at the value  $ \, \h = 1 \, $,  \, corresponding to the
common element  $ \, H^\vee \, $.  At a glance,
  $$  \U(\gerg_-)  \hskip2pt
\underset{{{(H^\vee)}_\h}^{\!\!\vee}}
{\overset{0 \,\leftarrow\, \h \,\rightarrow\, 1}
{\longleftarrow\joinrel\relbar\joinrel%
\relbar\joinrel\relbar\joinrel\relbar\joinrel\llongrightarrow}}
\hskip2pt  H^\vee  \hskip1pt
\underset{\;{({{(H^\vee)}_\h}^{\!\!\vee})}'}
{\overset{1 \,\leftarrow\, \h \,\rightarrow\, 0}
{\longleftarrow\joinrel\relbar\joinrel%
\relbar\joinrel\relbar\joinrel\relbar\joinrel\llongrightarrow}}
\hskip2pt  F[K_-]
\hskip9pt  \Big( = F[G_-^\star]  \,\text{\ if \ } 
\!\text{\it Char}\,(\Bbbk) = 0 \Big) \, .   \eqno (3.3)  $$
   \indent   Now consider  $ {H_\h}' $.  We have  $ \,
{H_\h}'{\Big|}_{\h=0} := {H_\h}' \Big/ \h \, {H_\h}' =
G_{\underline{D}}(H) =: \widetilde{H} \, $,  \, and  $ \,
\widetilde{H} = F[G_+] \, $  for some connected algebraic
Poisson group  $ G_+ \, $.  On the other hand, we have also
$ \, {H_\h}'{\Big|}_{\h=1} := {H_\h}' \Big/ (\h-1) \, {H_\h}' =
\sum_{n \geq 0} D_n =: H' \, $;  \, note that the latter is a
Hopf subalgebra of  $ H $,  because  $ \underline{D} $  is a
Hopf algebra filtration; moreover we have  $ \, \widetilde{H}
= \widetilde{H'} \, $,  by the very definitions.  Therefore we can
think of  $ \, {H_\h}' = \mathcal{R}^\h_{\underline{D}}(H') \, $  as
a 1-parameter family inside  $ \HA_\Bbbk $  with regular fibers
which links  $ \widetilde{H} $  and  $ H' $  as (polynomial)
deformations of each other, namely
  $$  F[G_+] = \widetilde{H} = {H_\h}'{\Big|}_{\h=0}  \hskip5pt
\underset{\;{H_\h}^{\!\prime}}
{\overset{0 \,\leftarrow\, \h \,\rightarrow\, 1}
{\longleftarrow\joinrel\relbar\joinrel%
\relbar\joinrel\relbar\joinrel\relbar\joinrel\llongrightarrow}}
\hskip6pt  {H_\h}'{\Big|}_{\h=1} = H' \; .  $$
   \indent   Consider also  $ \, {\big({H_\h}'\big)}^\vee \, $: 
by construction  $ \, {\big({H_\h}'\big)}^\vee{\Big|}_{\h=1} \!
= {H_\h}'{\Big|}_{\h=1} \! = H' \, $,  \, whereas  $ \, {\big(
{H_\h}' \big)}^\vee{\Big|}_{\h=0} \! = \U(\gerk_+) \, $  for some
restricted Lie bialgebra  $ \gerk_+ \, $:  in addition, if  $ \,
\text{\it Char}\,(\Bbbk) = 0 \, $  then  $ \, \gerk_+ = \gerg_+^{\,
\times} \, $  by  Theorem 2.2{\it (c)}.  Thus  $ \, {\big( {H_\h}'
\big)}^\vee \, $  can be seen as a 1-parameter family with regular
fibers, inside  $ \HA_\Bbbk \, $,  which links  $ \, \U(\gerk_+) $ 
and  $ H' $  as (polynomial) deformations of each other, namely
  $$  H' = {\big({H_\h}'\big)}^\vee{\Big|}_{\h=1}
\hskip3pt  \underset{\;{({H_\h}^{\!\prime})}^\vee}
{\overset{1 \,\leftarrow\, \h \,\rightarrow\, 0}
{\longleftarrow\joinrel\relbar\joinrel%
\relbar\joinrel\relbar\joinrel\relbar\joinrel\llongrightarrow}}
\hskip4pt  {\big({H_\h}'\big)}^\vee{\Big|}_{\h=0} = \U(\gerk_+)
\hskip9pt  \Big( = U(\gerg_+^{\,\times})  \,\text{\ if \ }
\text{\it Char}\,(\Bbbk) = 0 \Big) \, .  $$
Therefore,  $ H' $  is  {\sl at the same time\/}  a deformation of a
function algebra  {\sl and\/}  a deformation of an enveloping algebra,
via two different 1-parameter families inside  $ \HA_\Bbbk $  (with
regular fibers) which match at the value  $ \, \h = 1 \, $,  \,
corresponding (in both families) to  $ \, H' \, $.  In short,
  $$  F[G_+]  \hskip3pt
\underset{\;{H_\h}^{\!\prime}}
{\overset{0 \,\leftarrow\, \h \,\rightarrow\, 1}
{\longleftarrow\joinrel\relbar\joinrel%
\relbar\joinrel\relbar\joinrel\relbar\joinrel\llongrightarrow}}
\hskip2pt  H'  \hskip1pt
\underset{\;{({H_\h}^{\!\prime})}^\vee}
{\overset{1 \,\leftarrow\, \h \,\rightarrow\, 0}
{\longleftarrow\joinrel\relbar\joinrel%
\relbar\joinrel\relbar\joinrel\relbar\joinrel\llongrightarrow}}
\hskip3pt  \U(\gerk_+)
\hskip9pt  \Big( = U(\gerg_+^{\,\times})  \,\text{\ if \ }
\text{\it Char}\,(\Bbbk) = 0 \Big) \, .   \eqno (3.4)  $$
   \indent   Finally, it is worth noticing that when
$ \, H' = H = H^\vee \, $  formulas (3.3--4) give
%
%
%
%
\begin{align*}
 \hskip-17pt   F[G_+]  \hskip3pt
\underset{\;{H_\h}^{\!\prime}}
{\overset{0 \,\leftarrow\, \h \,\rightarrow\, 1}
{\longleftarrow\joinrel\relbar\joinrel%
\relbar\joinrel\relbar\joinrel\relbar\joinrel\llongrightarrow}}
\hskip2pt {}  &  H'  \hskip1pt
\underset{\;{({H_\h}^{\!\prime})}^\vee}
{\overset{1 \,\leftarrow\, \h \,\rightarrow\, 0}
{\longleftarrow\joinrel\relbar\joinrel%
\relbar\joinrel\relbar\joinrel\relbar\joinrel\llongrightarrow}}
\hskip3pt  \U(\gerk_+)
\hskip9pt  \Big( = U(\gerg_+^{\,\times})  \,\text{\ if \ }
\text{\it Char}\,(\Bbbk) = 0 \Big)   &  \\
\end{align*}
%
 \vskip-53pt  
  $$  || \hskip201pt  $$   
%
%
  $$  H \hskip202pt  \eqno (3.5)  $$   
 \vskip-26pt  
%
\begin{align*}
   &  \hskip4pt   ||   &  \\
   \hskip-17pt   \U(\gerg_-)  \hskip2pt
\underset{{{(H^\vee)}_\h}^{\!\!\vee}}
{\overset{0 \,\leftarrow\, \h \,\rightarrow\, 1}
{\longleftarrow\joinrel\relbar\joinrel%
\relbar\joinrel\relbar\joinrel\relbar\joinrel\llongrightarrow}}
\hskip2pt {}  &  H^\vee  \hskip1pt
\underset{\;{({{(H^\vee)}_\h}^{\!\!\vee})}'}
{\overset{1 \,\leftarrow\, \h \,\rightarrow\, 0}
{\longleftarrow\joinrel\relbar\joinrel%
\relbar\joinrel\relbar\joinrel\relbar\joinrel\llongrightarrow}}
\hskip2pt  F[K_-]
\hskip9pt  \Big( = F[G_-^\star]  \,\text{\ if \ }  \!\text{\it Char}\,
(\Bbbk) = 0 \Big)   &  
\end{align*}
which provides  {\sl four\/}  different regular 1-parameter (polynomial)
deformations from  $ H $  to Hopf algebras encoding geometrical objects
of Poisson type, i.e.~Lie bialgebras or Poisson algebraic groups.

\vskip7pt

  {\bf 3.6 The function algebra case.} \, Let  $ G $  be any algebraic
group over the field  $ \Bbbk $.  Let  $ \, R := \Bbbk[\h] \, $  be
as in \S 3.1, and set  $ \, F_\h[G] := {\big( F[G] \big)}_\h = R
\otimes_\Bbbk F[G] \, $:  \, this is trivially a QFA at  $ \h \, $, 
because  $ \, F_\h[G] \big/ \h \, F_\h[G] = F[G] \, $,  inducing on 
$ G $  the trivial Poisson structure, so that its cotangent Lie bialgebra
is simply  $ \gerg^\times $  with trivial Lie bracket and Lie cobracket
dual to the Lie bracket of  $ \gerg \, $.  In the sequel we identify
\hbox{$ F[G] $  with  $ \, 1 \! \otimes \! F[G] \subset F_\h[G] \, $.}
                                          \par
   We begin by computing  $ \, {F_\h[G]}^\vee \, $  (w.r.t.~$ \h \, $)
and  $ \, {F_\h[G]}^\vee{\Big|}_{\h=0} \! = \widehat{F[G]} =
G_{\underline{J}}\big(F[G]\big) \, $.
                                          \par
   Let  $ \, J := J_{F[G]} \equiv \text{\sl Ker}\, \big( \epsilon_{F[G]}
\big) \, $,  let  $ \, {\{j_b\}}_{b \in \mathcal{S}} \, (\subseteq J \,)
\, $  be a system of parameters of  $ F[G] $,  i.e.~$ \, \{ y_b :=
j_b \mod J^2 \}_{b \in \mathcal{S}} \, $  is a  $ \Bbbk $--basis  of
$ \, Q\big(F[G]\big) = J \Big/ J^2 = \gerg^\times \, $.  Then  $ \,
J^n \big/ J^{n+1} \, $  is  $ \Bbbk $--spanned by  $ \, \big\{\,
j^{\,\underline{e}} \mod J^{n+1} \;\big|\; \underline{e} \in
\N^{\mathcal{S}}_f \, , \; |\underline{e}| = n \,\big\} \, $  for
all  $ n \, $,  where  $ \, \N^{\mathcal{S}}_f := \big\{\, \sigma \in
\N^{\mathcal{S}} \,\big\vert\, \sigma(b) = 0 \;\; \text{\sl for almost
all} \;\; b \in \mathcal{S} \,\big\} \, $  (hereafter, monomials like
the previous ones are  {\sl ordered\/}  w.r.t.~some fixed order of the
index set  $ \mathcal{S} \, $)  and  $ \, |\underline{e}| := \sum_{b
\in \mathcal{S}} \underline{e}(b) \, $.  This implies that  
 \vskip3pt
   \centerline{ $ {F[G]}^\vee \, = \, \sum_{\underline{e}
\in \N^{\mathcal{S}}_f} \Bbbk[\h\,] \cdot \h^{-|\underline{e}|}
j^{\,\underline{e}} \bigoplus \Bbbk[\h\,] \big[\h^{-1}\big]
\, J^\infty \, = \, \sum_{\underline{e} \in \N^{\mathcal{S}}_f}
\Bbbk[\h\,] \cdot {(j^\vee)}^{\,\underline{e}} \bigoplus
\Bbbk[\h\,] \big[\h^{-1}\big] \, J^\infty $ }
 \vskip1pt
\noindent   where  $ \, J^\infty := \bigcap_{n \in \N} J^n \, $  and
$ \, j^\vee_s := \h^{-1} j_s \, $  for all  $ \, s \in \mathcal{S} \, $.
We also get that  $ \, \widehat{F[G]} = G_{\underline{J}}\big(F[G]\big)
\, $  is  $ \Bbbk $--spanned by  $ \, \big\{\, j^{\,\underline{e}}
\mod J^{n+1} \;\big|\; \underline{e} \in \N^{\mathcal{S}}_f \,\big\} \, $,
\, so  $ \, \widehat{F[G]} = G_{\underline{J}}\big(F[G]\big) \, $
is a quotient of  $ \, S(\gerg^\times) \, $.
                                          \par
   Now we distinguish various cases.  First assume  $ \, G \, $  is
{\sl smooth}, i.e.~$ \, \Bbbk^a \otimes_\Bbbk F[G] \, $  is  {\sl
reduced\/}  (where  $ \Bbbk^a $  is the algebraic closure of
$ \Bbbk $),  which is always the case if  $ \hbox{\it Char}\,(\Bbbk)
= 0 \, $.  Then (by standard results on algebraic groups) the above
set spanning  $ \widehat{F[G]} $  is a  $ \Bbbk $--basis:  thus  $ \,
{F_\h[G]}^\vee{\Big|}_{\h=0} \! = \widehat{F[G]} = G_{\underline{J}}
\big(F[G]\big) \cong S(\gerg^\times) \, $  as  $ \Bbbk $--algebras.
In addition, tracking the construction of the co-Poisson Hopf
structure onto  $ \widehat{F[G]} $  we see at once that  $ \,
\widehat{F[G]} \cong S(\gerg^\times) \, $  {\sl as co-Poisson Hopf
algebras too},  where the Hopf structure on  $ S(\gerg^\times) $
is the standard one and the co-Poisson structure is the one induced
from the Lie cobracket of  $ \gerg^\times $  (cf.~[Ga5] for details). 
Note also that  $ \, S(\gerg^\times) = U(\gerg^\times) \, $  because 
$ \gerg^\times $  is Abelian.
                                             \par
   Another ``extreme'' case is when  $ G $  is a  {\sl finite
connected group scheme\/}:  then, assuming for simplicity that
$ \, \Bbbk \, $  be perfect, we have  $ \, F[G] = \Bbbk[x_1,
\dots, x_n] \Big/ \big( x_1^{p^{e_1}}, \dots, x_n^{p^{e_n}}
\big) \, $  for some  $ \, n , e_1, \dots, e_n \in \N \, $.
Modifying a bit the analysis of the smooth case one gets
 \vskip3pt
   \centerline{ $ {F[G]}^\vee \, = \, \sum_{\underline{e} \in \N^n}
\Bbbk[\h\,] \cdot \h^{-|\underline{e}|} x^{\,\underline{e}} \, =
\, \sum_{\underline{e} \in \N^n} \Bbbk[\h\,] \cdot {(x^\vee)}^{\,
\underline{e}} $ }
 \vskip1pt
\noindent   (now  $ \, J^\infty = \{0\} \, $),  and  $ \, {F_\h[G]}^\vee
{\Big|}_{\h=0} \! = \widehat{F[G]} = G_{\underline{J}}\big(F[G]\big)
\cong S(\gerg^\times) \Big/ \big( \bar{x}_1^{p^{e_1}}, \dots,
\bar{x}_n^{p^{e_n}} \big) \, $,  \, where  $ \, \bar{x}_i :=
x_i \mod J^2 \in \gerg^\times \, $.  Now, recall that for any Lie
algebra  $ \gerh $  there is  $ \, \gerh^{{[p\hskip0,7pt]}^\infty}
:= \Big\{\, x^{{[p\hskip0,7pt]}^n} \! := x^{p^n} \,\Big\vert\, x \in
\gerh \, , n \in \N \,\Big\} \, $,  \, the  {\sl restricted Lie algebra
generated by  $ \gerh $}  inside  $ U(\gerh) $,  with  $ p $--operation 
given by  $ \, x^{[p\hskip0,7pt]} := x^p \, $;  \, then one always has 
$ \, U(\gerh) = \u\big( \gerh^{{[p\hskip0.7pt]}^\infty} \big) \, $.  In
our case  $ \, \Big\{ \bar{x}_1^{\,p^{e_1}}, \, \ldots, \, \bar{x}_n^{\,
p^{e_n}} \Big\} \, $  generates a  $ p $--ideal  $ \mathcal{I} $ 
of  $ (\gerg^\times)^{{[p\hskip0,7pt]}^\infty} $,  hence  $ \,
\gerg^\times_{\text{\it res}} := \gerg^{{[p\hskip0,7pt]}^\infty}
\! \Big/ \mathcal{I} \, $  is a restricted Lie algebra too, with 
$ \Bbbk $--basis  $ \, \Big\{ \bar{x}_1^{\,p^{a_1}}, \, \ldots, \,
\bar{x}_n^{\,p^{a_n}} \,\Big|\; a_1 < e_1, \dots, a_n < e_n \Big\} \, $. 
Then the previous analysis gives  $ \; {F_\h[G]}^\vee {\Big|}_{\h=0} =
\u \left( \gerg^\times_{\text{\it res}} \right) \equiv S(\gerg^\times)
\Big/ \Big( \Big\{ \bar{x}_1^{\,p^{e_1}}, \, \ldots, \, \bar{x}_n^{\,
p^{e_n}} \Big\} \Big) \; $  as co-Poisson Hopf algebras.   
                                             \par
   The general case is intermediate.  Assume again for simplicity that
$ \Bbbk $  be perfect.  Let  $ \, F[[G]] \, $  be the  $ J $--adic
completion of  $ \, H = F[G] \, $.  By standard results on algebraic
groups (cf.~[DG]) there is a subset  $ \, {\{x_i\}}_{i \in \mathcal{I}}
\, $  of  $ \, J \, $  such that  $ \, {\big\{\, \overline{x}_i :=
x_i \! \mod J^2 \,\big\}}_{i \in \mathcal{I}} \, $  is a basis of  $ \,
\gerg^\times = J \big/ J^2 \, $  and  $ \; F[[G]] \, \cong \,
\Bbbk\big[\big[ {\{x_i\}}_{i \in \mathcal{I}} \big]\big] \Big/ \Big(
\Big\{ x_i^{\,p^{n(x_i)}} \Big\}_{i \in \mathcal{I}_0} \Big) \; $  (the
algebra of truncated formal power series), for some  $ \, \mathcal{I}_0
\subset \mathcal{I} \, $  and  $ \, {\big( n(x_i) \big)}_{i \in \mathcal{I}_0}
\in \N^{\,\mathcal{I}_0} $.  Since  $ \, G_{\underline{J}}\,\big(F[G]\big)
= G_{\underline{J}}\,\big(F[[G]]\big) \, $,  \, we argue that  $ \;
G_{\underline{J}}\,\big(F[G]\big) \, \cong \, \Bbbk\big[ {\{
\overline{x}_i \}}_{i \in \mathcal{I}} \big] \Big/ \Big( \Big\{
\overline{x}_i^{\,p^{n(x_i)}} \Big\}_{i \in \mathcal{I}_0} \,\Big)
\; $;  \; finally, since  $ \; \Bbbk\big[ {\{ \overline{x}_i
\}}_{i \in \mathcal{I}} \big] \cong S(\gerg^\times) \; $  we get
  $$  G_{\underline{J}}\,\big(F[G]\big) \cong S(\gerg^\times)
\bigg/ \! \Big( \Big\{\, \overline{x}^{\,p^{n(x)}} \Big\}_{x
\in \mathcal{N}(F[G])} \,\Big)  $$  
as algebras,  $ \mathcal{N} \big( F[G] \big) $  being the nilradical
of  $ F[G] $  and  $ \, p^{n(x)} $  is the nilpotency order of  $ \,
x \in \mathcal{N}\big(F[G]\big) \, $.
                                        \par
   Finally, noting that  $ \, \Big( \Big\{\, \overline{x}^{\,p^{n(x)}}
\Big\}_{x \in \mathcal{N}(F[G])} \,\Big) \, $  is a co-Poisson Hopf
ideal of  $ S(\gerg^\times) \, $,  \, like in the smooth case we argue
that the above isomorphism is one of  {\sl co-Poisson Hopf\/}  algebras.
                                            \par
   If  $ \Bbbk $  is not perfect the same analysis applies, but
modifying a bit the previous arguments.
                                             \par
  As for  $ \, F[G]^\vee := F[G] \big/ J^\infty \, $,  one has
(cf.~[Ab], Lemma 4.6.4)  $ \, F[G]^\vee =F[G] \, $  whenever
$ G $  is finite dimensional and there exists no  $ \, f \in
F[G] \setminus \Bbbk \, $  which is separable algebraic over
$ \Bbbk \, $.   
                                             \par
   It is also interesting to consider  $ {\big( {F_\h[G]}^\vee
\big)}' $.  If  $ \, \text{\it Char}\,(\Bbbk) = 0 \, $  Theorem
2.2{\it (c)\/}  gives  $ \, {\big( {F_\h[G]}^\vee \big)}' = F_\h[G]
\, $.  If instead  $ \, \text{\it Char}\,(\Bbbk) = p > 0 \, $,  then
the situation might change dramatically.  Indeed, if  $ G $  has
dimension 0 and eight 1 then   --- i.e., if  $ \, F[G] = \Bbbk\big[
{\{x_i\}}_{i \in \mathcal{I}} \big] \big] \Big/ \big( \big\{ x_i^p
\,\big|\, i \in \mathcal{I} \big\} \big) \, $  as a  $ \Bbbk $--algebra
---   the same analysis as in the zero characteristic case applies, with
a few minor changes, whence one gets again  $ \, {\big( {F_\h[G]}^\vee
\big)}' = F_\h[G] \, $.  Otherwise, let  $ \, y \in J \setminus \{0\}
\, $  be primitive and such that  $ \, y^p \not= 0 \, $  (for instance,
this occurs for  $ \, G \cong \Bbb{G}_a \, $).  Then  $ \, y^p \, $  is
primitive as well, hence  $ \, \delta_n(y^p) = 0 \, $  for each  $ \,
n > 1 \, $.  It follows that  $ \, 0 \not= \h \, {(y^\vee)}^p \in {\big(
{F_\h[G]}^\vee \big)}' \, $,  \, whereas  $ \, \h \, {(y^\vee)}^p \not\in
F_\h[G] \, $,  due to our previous description of  $ {F_\h[G]}^\vee $. 
Thus  $ \, {\big( {F_\h[G]}^\vee \big)}' \supsetneqq {F_\h[G]}^\vee \, $, 
\, a  {\sl counterexample to the first part of  Theorem 2.2{\it (c)}}.   
                                             \par
   What for  $ \, {F[G]}' \, $  and  $ \, \widetilde{F[G]} \, $?  Again,
this depends on the group  $ G $  under consideration.  We provide two
simple examples, both ``extreme'', in a sense, and opposite to each
other.
                                             \par
   Let  $ \, G := \Bbb{G}_a = \hbox{\it Spec}\big(\Bbbk[x]\big)
\, $,  \, so  $ \, F[G] = F[\Bbb{G}_a] = \Bbbk[x] \, $  and  $ \,
F_\h[\Bbb{G}_a] := R \otimes_\Bbbk \Bbbk[x] = R[x] \, $.  Then
since  $ \, \Delta(x) := x \otimes 1 + 1 \otimes x \, $  and  $ \,
\epsilon(x) = 0 \, $  we find  $ \, {F_\h[\Bbb{G}_a]}' = R[\h{}x]
\, $  (like in \S 3.7 below: indeed, this is just a special instance,
for  $ \, F[\Bbb{G}_a] = U(\gerg) \, $  where  $ \gerg $  is the
1-dimensional Lie algebra).  Moreover, iterating one gets easily
$ \, {\big( {F_\h[\Bbb{G}_a]} \big)}' = R\big[\h^2{}x\big] \, $,
$ \, {\Big( {\big( {F_\h [\Bbb{G}_a]}' \big)}' \Big)}' = R \big[
\h^3{}x \big] \, $,  \, and in general  $ \; \Big( \Big( \big(
F_\h [\Bbb{G}_a \underbrace{]'{\big)}'{\Big)}' \cdots
{\Big)}'}_{n} = R\big[\h^n{}x\big] \cong R[x] = F_\h
[\Bbb{G}_a] \; $  for all  $ \, n \in \N \, $.
                                             \par
   Second, let  $ \, G := \Bbb{G}_m = \hbox{\it Spec} \, \Big(
\Bbbk \big[ z^{+1}, z^{-1} \big] \Big) \, $,  \, that is  $ \,
F[G] = F[\Bbb{G}_m] = \Bbbk \big[ z^{+1}, z^{-1} \big] \, $  so
that  $ \, F_\h[\Bbb{G}_m] := R \otimes_\Bbbk \Bbbk \big[ z^{+1},
z^{-1} \big] = R \big[ z^{+1}, z^{-1} \big] \, $.  Then since
$ \, \Delta \big( z^{\pm 1} \big) := z^{\pm 1} \otimes z^{\pm 1}
\, $  and  $ \, \epsilon \big( z^{\pm 1} \big) = 1 \, $  we
find  $ \, \Delta^n \big( z^{\pm 1} \big) = {\big( z^{\pm 1}
\big)}^{\otimes n} \, $  and  $ \, \delta_n \big( z^{\pm 1}
\big) = {\big( z^{\pm 1} - 1 \big)}^{\otimes n} \, $  for
all  $ \, n \in \N \, $.  It follows easily from that
$ \, {F_\h[\Bbb{G}_m]}' = R \cdot 1 \, $,  \, the
trivial possibility (see also \S 3.12 later on).

\vskip7pt

  {\bf 3.7 The enveloping algebra case.} \, Let  $ \gerg $  be any
Lie algebra over the field  $ \Bbbk $,  and  $ U(\gerg) $  its
universal enveloping algebra with its standard Hopf structure.
Assume  $ \, \Char(\Bbbk) = 0 \, $,  and let  $ \, R = \Bbbk[\h] \, $,
as in \S 3.1,  and set  $ \, U_\h(\gerg) := R \otimes_\Bbbk U(\gerg)
= {\big( U(\gerg) \big)}_\h \, $.  Then  $ U_\h(\gerg) $  is trivially
a QrUEA at  $ \h $,  for  $ \, U_\h(\gerg) \big/ \h \, U_\h(\gerg) =
U(\gerg) \, $,  inducing on  $ \gerg $  the trivial Lie cobracket.
Thus the dual Poisson group is just  $ \gerg^\star $  (the
topological dual of  $ \gerg $  w.r.t.~the weak topology)
w.r.t.~addition, with  $ \gerg $  as
cotangent Lie bialgebra and function algebra  $ \, F[\gerg^\star]
= S(\gerg) \, $:  \, the Hopf structure is the standard one, and the
Poisson structure is the one induced by  $ \, \{x,y\} := [x,y] \, $
for all  $ \, x $,  $ y \in \gerg \, $  (it is the  {\sl
Kostant-Kirillov structure\/}  on  $ \gerg^\star \, $).
                                          \par
   Similarly, if  $ \, \text{\it Char}\,(\Bbbk) = p > 0 \, $  and
$ \gerg $  is any restricted Lie algebra over  $ \Bbbk $,  let
$ \, \u(\gerg) \, $  be its restricted universal enveloping algebra,
with its standard Hopf structure.  Then if  $ \, R = \Bbbk[\h] \, $
the Hopf  $ R $--algebra  $ \, U_\h(\gerg) := R \otimes_\Bbbk \u(\gerg)
= {\big(\u(\gerg)\big)}_\h \, $  is a QrUEA at  $ \h $,  because  $ \,
\u_\h(\gerg) \big/ \h \, \u_\h(\gerg) = \u(\gerg) \, $,  inducing on
$ \gerg $  the trivial Lie cobracket: then the dual Poisson group is
again  $ \gerg^\star $, with cotangent Lie bialgebra  $ \gerg $  and
function algebra  $ \, F[\gerg^\star] = S(\gerg) \, $  (the Poisson
Hopf structure being as above).  Recall also that  $ \, U(\gerg) =
\u\big( \gerg^{{[p\hskip0.7pt]}^\infty} \big) \, $  (cf.~\S 3.6).
                                          \par
   First we compute  $ \, {\u_\h(\gerg)}' $  (w.r.t.~the prime  $ \h \, $) 
using (3.2), i.e.~computing the filtration  $ \underline{D} \, $.   
                                          \par   
   By the PBW theorem, once an ordered basis  $ B $  of  $ \gerg $  is
fixed  $ \u(\gerg) $  admits as basis the set of ordered monomials
in the elements of  $ B $  whose degree (w.r.t.~each element of  $ B $)
is less than  $ p \, $;  this yields a Hopf algebra filtration of
$ \u(\gerg) $  by the total degree, which we refer to as  {\sl the
standard filtration}.  Then a straightforward calculation shows that 
$ \underline{D} $  coincides with the standard filtration.  This and
(3.2) imply  $ \, {\u_\h(\gerg)}' = \langle \tilde{\gerg} \rangle =
\langle \h \, \gerg \rangle \, $:  \, hereafter  $ \, \tilde{\gerg}
:= \h \, \gerg \, $,  \, and similarly  $ \, \tilde{x} := \h \, x \, $ 
for all  $ \, x \in \gerg \, $.  Then the relations  $ \; x \, y - y \,
x = [x,y] \; $  and  $ \; z^p = z^{[p\hskip0,7pt]} \, $  in  $ \u(\gerg) $ 
yield  $ \; \tilde{x} \, \tilde{y} - \tilde{y} \, \tilde{x} = \h \,
\widetilde{[x,y]} \equiv 0 \mod \h \, {\u_\h(\gerg)}' \; $  and also 
$ \; \tilde{z}^p = \h^{p-1} \widetilde{z^{[p\hskip0,7pt]}} \equiv 0
\mod \h \, {\u_\h(\gerg)}' \; $;  \; therefore, from
%
%
  $ \; \u_\h(\gerg) = T_R(\gerg) \Big/ \big( \big\{\,
x \, y - y \, x - [x,y] \, , \, z^p - z^{[p\hskip0,7pt]} \;\big\vert\;
x, y, z \in \gerg \,\big\} \big) \; $  we get
 \vskip-10pt
  $$  \displaylines{
   {\u_\h(\gerg)}' = \langle \tilde{\gerg} \rangle
\,\;{\buildrel {\h \rightarrow 0} \over
{\relbar\joinrel\relbar\joinrel\llongrightarrow}}\;\,
\widetilde{\u(\gerg)} \, = \, T_\Bbbk(\tilde{\gerg})
\bigg/ \Big( \Big\{\, \tilde{x} \, \tilde{y} -
\tilde{y} \, \tilde{x} \, , \, \tilde{z}^p \;\Big\vert\;
\tilde{x}, \tilde{y}, \tilde{z} \in \tilde{\gerg} \,\Big\} \Big)
=   \hfill  \cr
   = T_\Bbbk (\gerg) \! \Big/ \! \big( \big\{\, x \, y - y \,
x \, , \, z^p \;\big\vert\; x, y, z \in \gerg \,\big\} \big)
= S_\Bbbk (\gerg) \! \Big/ \! \big( \big\{\, z^p \,\big\vert\;
z \in \gerg \,\big\} \big) = F[\gerg^\star] \! \Big/ \! \big(
\big\{\, z^p \,\big\vert\; z \in \gerg \,\big\} \big)  \cr }  $$
that is  $ \; \widetilde{\u(\gerg)} := G_{\underline{D}} \big(
\u(\gerg) \big) = {\u_\h(\gerg)}' \big/ \h \, {\u_\h(\gerg)}' \cong
F[\gerg^\star] \Big/ \big( \big\{\, z^p \,\big\vert\; z \in \gerg
\,\big\} \big) \; $  {\sl as Poisson Hopf algebras}.  In particular,
{\it this means that  $ \widetilde{\u(\gerg)} $  is the function
algebra of, and  $ {\u_\h(\gerg)}' $  is a QFA (at $ \h \, $)  for,
a non-reduced algebraic Poisson group of dimension 0 and height 1,
whose cotangent Lie bialgebra is  $ \gerg \, $,  hence which is dual
to  $ \, \gerg \, $};  \, thus, in a sense, part  {\it (c)\/}  of
Theorem 2.2 is still valid in this  case too.

\vskip5pt

   {\sl $ \underline{\hbox{\it Remark}} $:}  \, Note that this last
result reminds the classical formulation of the analogue of Lie's
Third Theorem in the context of group-schemes:  {\it Given a
restricted Lie algebra  $ \gerg $,  there exists a group-scheme
$ G $  {\sl of dimension 0 and height 1\/}  whose tangent Lie algebra
is  $ \, \gerg \, $}  (see e.g.~[DG]).  Here we have just given sort
of a ``dual Poisson-theoretic version'' of this fact, in that our
result sounds as follows:  {\it Given a restricted Lie algebra
$ \gerg $,  there exists a Poisson group-scheme  $ G $  {\sl of
dimension 0 and height 1\/}  whose  {\sl cotangent}  Lie algebra
is  $ \, \gerg \, $}.   

\vskip5pt

   As a byproduct, since  $ \, U_\h(\gerg) = \u_\h\big( \gerg^{{[p
\hskip0,7pt]}^\infty} \big) \, $  we have also  $ \; {U_\h(\gerg)}' =
{\u_\h\big( \gerg^{{[p\hskip0,7pt]}^\infty} \big)}' \, $,  \, whence
  $$  {U_\h(\gerg)}' = {\u_\h\big( \gerg^{{[p\hskip0,7pt]}^\infty}
\big)}' \,{\buildrel {\h \rightarrow 0} \over
{\relbar\joinrel\relbar\joinrel\loongrightarrow}}\; S_\Bbbk \Big(
\gerg^{{[p\hskip0,7pt]}^\infty} \Big) \hskip-2pt \bigg/ \hskip-2pt \Big(
{\big\{\, z^p \,\big\}}_{z \in \gerg^{{[p\hskip0,7pt]}^\infty}} \Big)
\, = \, F \! \left[ {\big( \gerg^{{[p\hskip0,7pt]}^\infty} \big)}^{\!
\star} \right] \hskip-2pt \bigg/ \hskip-2pt \Big( {\big\{\, z^p \,
\big\}}_{z \in \gerg^{{[p\hskip0,7pt]}^\infty}} \Big) \, .  $$   
   \indent   Furthermore,  $ \, {\u_\h(\gerg)}' = \langle \tilde{\gerg}
\rangle \, $  implies that  $ \, I_{{\u_\h(\gerg)}'} \, $  is generated
(as an ideal) by  $ \, \h \, R \cdot 1_{\u_\h(\gerg)} + R \, \tilde{\gerg}
\, $,  \, hence  $ \, \h^{-1} I_{{\u_\h(\gerg)}'} \, $  is generated by 
$ \, R \cdot 1 + R \, \gerg \, $,  \, therefore   
  $$  {\big( {\u_\h(\gerg)}' \big)}^\vee := \, {\textstyle
\bigcup\nolimits_{n \geq 0}} {\big( \h^{-1} I_{{\u_\h(\gerg)}'} \big)}^n
= \, {\textstyle \bigcup\nolimits_{n \geq 0}} {\big( R \cdot 1 + R \,
\gerg \big)}^n = \, \u_\h(\gerg) \; .  $$   
This means that also part  {\it (b)}  of Theorem 2.2 is still valid,
though now  $ \, \text{\it Char}\,(\Bbbk) > 0 \, $.
                                          \par
   When  $ \, \text{\it Char}\,(\Bbbk) = 0 \, $  and we look at 
$ U(\gerg) $,  the like argument applies:  $ \underline{D} $  coincides
with the standard filtration of  $ U(\gerg) $  provided by the total
degree, via the PBW theorem.  This and (3.2) imply  $ \, {U(\gerg)}' =
\langle \tilde{\gerg} \rangle = \langle \h \, \gerg \rangle\, $,  so that
from the presentation  $ \; U_\h(\gerg) = T_R(\gerg) \! \Big/ \! \big(
{\big\{\, x \, y - y \, x - [x,y] \,\big\}}_{x, y, z \in \gerg\,} \big)
\, $ we get  $ \, {U_\h(\gerg)}' \! = T_{\!R}(\tilde{\gerg}) \! \Big/
\! \Big( {\Big\{ \tilde{x} \, \tilde{y} - \tilde{y} \, \tilde{x} -
\h \cdot \widetilde{[x,y]} \Big\}}_{\tilde{x}, \tilde{y} \in
\tilde{\gerg}\,} \Big) \, $,  \, whence we get at once
  $$  {U_\h(\gerg)}' \;{\buildrel {\h \rightarrow 0}
\over {\relbar\joinrel\relbar\joinrel\llongrightarrow}}\;\,
\widetilde{U(\gerg)} \; \cong \; T_\Bbbk(\tilde{\gerg})\Big/
\big( \big\{\, \tilde{x} \, \tilde{y} - \tilde{y} \, \tilde{x}
\;\big\vert\; \tilde{x}, \tilde{y} \in \tilde{\gerg} \,\big\} \big)
\; \cong \; S_\Bbbk(\gerg) \; = \; F[\gerg^\star]  $$
i.e.~$ \; \widetilde{U(\gerg)} := G_{\underline{D}} \big( U(\gerg)
\big) = {U_\h(\gerg)}' \big/ \h \, {U_\h(\gerg)}' \cong F[\gerg^\star]
\; $  as  {\sl Poisson Hopf\/}  algebras, as predicted by  Theorem
2.2{\it (c)}.  Moreover,
      \hbox{$ \, {U_\h(\gerg)}' = \langle \tilde{\gerg}
\rangle = T(\tilde{\gerg}) \Big/ \big( \big\{\, \tilde{x} \, \tilde{y}
- \tilde{y} \, \tilde{x} = \h \cdot \widetilde{[x,y]} \;\big\vert\;
\tilde{x}, \tilde{y} \in \tilde{\gerg} \,\big\} \big) \, $}
  implies
that  $ \, I_{{U_\h(\gerg)}'} \, $  is generated by  $ \, \h \, R
\cdot 1_{U_\h(\gerg)} + R \, \tilde{\gerg} \, $:  \, thus  $ \,
\h^{-1} I_{{U_\h(\gerg)}'} \, $  is generated by  $ \, R \cdot
1_{U_\h(\gerg)} + R \, \gerg \, $,  \, so  $ \, {\big( {U_\h
(\gerg)}' \big)}^\vee := \bigcup\limits_{n \geq 0} {\big( \h^{-1}
I_{{U_\h(\gerg)}'} \big)}^n = \bigcup\limits_{n \geq 0} {\big(
R \cdot 1_{U_\h(\gerg)} + R \, \gerg \big)}^n = U_\h(\gerg) \, $,
\, agreeing     
       \hbox{with  Theorem 2.2{\it (b)}.}   
                                          \par
   What for the functor  $ \, {(\ )}^\vee \, $?  This heavily
depends on the  $ \gerg $  we start from!
                                          \par
   First assume  $ \, \hbox{\it Char}\,(\Bbbk) = 0 \, $.  Let
$ \, \gerg_{(1)} := \gerg \, $,  $ \, \gerg_{(k)} := \big[ \gerg,
\gerg_{(k-1)} \big] \, $  ($ k \in \N_+ $),  be the  {\sl lower
central series\/}  of  $ \gerg \, $.  Pick subsets  $ \, B_1 \, $,
$ B_2 \, $,  $ \dots \, $,  $ B_k \, $,  $ \dots \, $  ($ \subseteq
\gerg \, $)  such that  $ \, B_k \! \mod \gerg_{(k+1)} \, $  be a
$ \Bbbk $--basis  of  $ \, \gerg_{(k)} \big/ \gerg_{(k+1)} \, $
(for all  $ \, k \in \N_+ \, $),  \, pick also a  $ \Bbbk $--basis
$ B_{\infty} $  of  $ \, \gerg_{(\infty)} := \bigcap_{\, k \in \N_+}
\, $,  \, and set  $ \, \partial(b) := k \, $  for any  $ \, b \in B_k
\, $  and each  $ \, k \in \N_+ \cup \{\infty\} \, $.  Then  $ \, B :=
\left( \bigcup_{k \in \N_+} B_k \right) \cup B_{\infty} \, $  is a
$ \Bbbk $--basis  of  $ \gerg \, $;  \, we fix a total order on it.
Applying the PBW theorem to this ordered basis of  $ \gerg $  we
get that  $ \, J^n \, $  has basis the set of ordered monomials
$ \, \big\{\, b_1^{e_1} b_2^{e_2} \cdots b_s^{e_s} \;\big|\; s \in
\N_+ \, , b_r \in B \, , \, \sum_{r=1}^s b_r \, \partial(b_r) \geq
n \,\big\} \, $.  Then one finds that  $ \, {U_\h(\gerg)}^\vee \, $ 
is generated by  $ \, \big\{\, \h^{-1} b \;\big|\; b \in B_1 \setminus
B_2 \,\big\} \, $  (as a unital  $ R $--algebra)  and it is the direct
sum
 \vskip2pt
 \centerline{ $  {U_\h(\gerg)}^\vee = \, \bigg( \hskip-3pt
\oplus_{\substack{
\hskip-15pt  s \in \N_+  \\   b_r \in B \setminus B_\infty  }}
\hskip-15pt  R \, \big( \h^{-\partial(b_1)} b_1 \big)^{e_1} \cdots
\big( \h^{-\partial(b_s)} b_s \big)^{e_s} \hskip-2pt \bigg) \bigoplus
\bigg( \hskip-3pt \oplus_{\substack{s \in \N_+, \, b_r \in B  \\  
\exists \, \bar{r} \, : \, b_{\bar{r}} \in B_\infty}} \; R\big[
\h^{-1} \big] \, b_1^{e_1} \cdots b_s^{e_s} \hskip-2pt \bigg) $ }
   From this it follows at once that  $ \; {U_\h(\gerg)}^\vee \Big/ \h \,
{U_\h(\gerg)}^\vee \cong \, U\left( \gerg \big/ \gerg_{(\infty)} \right)
\; $  via an isomorphism which maps  $ \, \h^{-\partial(b)} b \mod
\h \, {U_\h(\gerg)}^\vee \, $  to  $ \; b \! \mod \gerg_{(\infty)}
\in \gerg \big/ \gerg_{(\infty)} \subset U \left( \gerg \big/
\gerg_{(\infty)} \right) \, $  for all  $ \, b \in B \setminus
B_{\infty} \, $  and maps  $ \, \h^{-n} b \! \mod \h \,
{U_\h(\gerg)}^\vee \, $  to  $ \, 0 \, $  for all  $ \, b
\in B \setminus B_{\infty} \, $  and all  $ \, n \in \N \, $.
                                          \par
   Now assume  $ \, \hbox{\it Char}\,(\Bbbk) = p > 0 \, $.  Then
in addition to the previous considerations one has to take into
account the filtration of  $ \, \u(\gerg) \, $  induced by both
the lower central series of  $ \gerg $  {\sl and\/}  the
$ p $--filtration  of  $ \gerg \, $,  that is  $ \, \gerg \supseteq
\gerg^{[p\hskip0,7pt]} \supseteq \gerg^{{[p\hskip0,7pt]}^2} \supseteq
\cdots \supseteq \gerg^{{[p\hskip0,7pt]}^n} \supseteq \cdots \, $,
\, where  $ \, \gerg^{{[p\hskip0,7pt]}^n} \, $  is the restricted
Lie subalgebra generated by  $ \, \big\{\, x^{{[p\hskip0,7pt]}^n}
\,\big|\, x \in \gerg \,\big\} \, $  and  $ \, x \mapsto
x^{[p\hskip0,5pt]} \, $  is the  $ p \, $--operation  in
$ \gerg \, $:  these encode the  $ J $--filtration  of
$ \u(\gerg) \, $,  hence of  $ \, \u_\h(\gerg) \, $, 
\, so permit to describe  $ {\u_\h(\gerg)}^\vee $.   
                                          \par
   In detail, for any restricted Lie algebra  $ \gerh $,
let  $ \, \gerh_n := \Big\langle \bigcup_{(m \, p^k \geq n}
{(\gerh_{(m)})}^{[p^k]} \Big\rangle \, $  for all  $ \, n \in
\N_+ \, $  \, (where  $ \langle X \rangle $  denotes the Lie
subalgebra of  $ \gerh $  generated by  $ X \, $)  and  $ \,
\gerh_\infty := \bigcap_{\,n \in \N_+} \gerh_n \, $:  \, we
call  $ \, {\big\{ \gerh_n \big\}}_{n \in \N_+} $  {\sl the
$ p $--lower  central series of\/}  $ \gerh \, $.  It is a  {\sl
strongly central series\/}  of  $ \gerh $,  i.e.~a central series
of  $ \gerh $  such that  $ \, [\gerh_m,\gerh_n] \leq \gerh_{m+n}
\, $  for all  $ m $,  $ n \, $,  and  $ \, \gerh_n^{\,[p\,]} \leq
\gerh_{n+1} \, $  for all  $ n \, $.
                                          \par
   Applying these tools to  $ \, \gerg \subseteq \u(\gerg) \, $
the very definitions give  $ \, \gerg_n \subseteq J^n \, $  (for
all  $ \, n \in \N \, $)  where  $ \, J := J_{\u(\gerg)} \, $:
more precisely, if  $ B $  is an ordered basis of  $ \gerg $
then the (restricted) PBW theorem for  $ \, \u(\gerg) \, $
implies that  $ \, J^n \big/ J^{n+1} \, $  admits as
$ \Bbbk $--basis  the set of ordered monomials of the
form  $ \, x_{i_1}^{e_1} x_{i_2}^{e_2} \cdots x_{i_s}^{e_s}
\, $  such that  $ \, \sum_{r=1}^s e_r \partial(x_{i_r}) = n \, $
where  $ \, \partial(x_{i_r}) \in \N \, $  is uniquely determined
by the condition  $ \, x_{i_r} \in \gerg_{\partial(x_{i_r})}
\setminus \gerg_{\partial(x_{i_r})+1} \, $  and each  $ \,
x_{i_k} \, $  is a fixed lift in  $ \gerg $  of an element of
a fixed ordered basis of  $ \, \gerg_{\partial(x_{i_k})} \Big/
\gerg_{\partial(x_{i_k})+1} \, $.  This yields an explicit
description of  $ \underline{J} \, $,  hence of  $ {\u(\gerg)}^\vee $
and  $ {\u_\h(\gerg)}^\vee $,  like before: in particular  
  $$  \widehat{\u_\h(\gerg)} \; := \; {\u_\h(\gerg)}^\vee \Big/
\h \, {\u_\h(\gerg)}^\vee \, \cong \; \u\left( \gerg \big/
\gerg_{\infty} \right) \; .  $$   

\noindent
{\bf Definition 3.8.}{\em \, We call  {\sl pre-restricted universal
enveloping algebra}  ({\sl =PrUEA})  any  $ \, H \in \HA_\Bbbk \, $
which is down-filtered by  $ \underline{J} $  (i.e.,  $ \, \bigcap_{n
\in \N} J^n = \{0\} \, $),  \, and  $ \, \PrUEA \, $  the full
subcategory of\/  $ \HA_\Bbbk $  of all PrUEAs.  We call  {\sl
pre-function algebra}  ({\sl =PFA})  any  $ \, H \in \HA_\Bbbk \, $
which is up-filtered by  $ \underline{D} $  (i.e.,  $ \, \bigcup_{n
\in \N} D_n = H \, $),  \, and  $ \, \PFA \, $  the full subcategory
of\/  $ \HA_\Bbbk $  of all PFAs.
}

\vskip7pt
\noindent
{\bf Theorem 3.9}{\em \, ({\sl ``The Crystal Duality Principle''})
                                       \hfill\break
  \indent   (a) \, $ H \mapsto H^\vee := H \big/ J_{{}_H}^{\, \infty}
\, $  and  $ \, H \mapsto H' := \bigcup_{n \in \N} D_n \, $  define
functors  $ \; {(\ )}^\vee \colon \, \HA_\Bbbk \longrightarrow
\HA_\Bbbk \, $  and  $ \; {(\ )}' \colon \, \HA_\Bbbk \longrightarrow
\HA_\Bbbk \, $  respectively whose image are  $ \PrUEA $  and  $ \PFA $
respectively.
                                          \hfill\break
  \indent   (b) \, Let  $ \, H \in \HA_\Bbbk \, $.  Then  $ \,
\widehat{H} := G_{\underline{J}}(H) \cong \U(\gerg) \, $  as graded
co-Poisson Hopf algebras, for some restricted Lie bialgebra  $ \gerg $
which is graded as a Lie algebra.  In particular, if  $ \, \Char(\Bbbk)
= 0 \, $  and  $ \, \dim(H) \in \N \, $  then  $ \, \widehat{H}
= \Bbbk \!\cdot\! 1 \, $  and  $ \, \gerg = \{0\} \, $.
                                        \hfill\break
   \indent   More in general, the same holds if  $ \, H = B \, $
is a\/  $ \Bbbk $--bialgebra.
                                        \hfill\break
  \indent   (c) \, Let  $ \, H \in \HA_\Bbbk \, $.  Then  $ \,
\widetilde{H} := G_{\underline{D}}(H) \cong F[G] \, $,  as graded
Poisson Hopf algebras, for some connected algebraic Poisson group
$ G $  whose variety of closed points form a (pro)affine space.  If
$ \, \Char(\Bbbk) = 0 \, $  then  $ \, F[G] = \widetilde{H} \, $  is
a polynomial algebra, i.e.~$ \, F[G] = \Bbbk \big[ {\{x_i \}}_{i \in
\mathcal{I}} \big] \, $  (for some set  $ \mathcal{I} $);  in particular, if
$ \, \dim(H) \in \N \, $  then  $ \, \widetilde{H} = \Bbbk \cdot 1 \, $
and  $ \, G = \{1\} \, $.  If  $ \, p := \Char(\Bbbk) > 0 \, $  then
$ G $  has dimension 0 and height 1, and if\/  $ \Bbbk $  is perfect
then  $ \, F[G] = \widetilde{H} \, $  is a truncated polynomial
algebra, i.e.~$ \, F[G] = \Bbbk \big[ {\{x_i\}}_{i \in \mathcal{I}}
\big] \Big/ \big( {\{x_i^{\,p}\}}_{i \in \mathcal{I}} \big) \, $
(for some set  $ \mathcal{I} $).
                                          \hfill\break
   \indent   More in general, the same holds if  $ \, H = B \, $
is a\/  $ \Bbbk $--bialgebra.
                                          \hfill\break
   \indent   (d) \, For every  $ \, H \in \HA_\Bbbk \, $,  \, there
exist two 1-parameter families  $ \, {{(H^\vee)}_\h}^{\!\!\vee}
= \mathcal{R}^\h_{\underline{J}}(H^\vee) \, $  and  $ \, {\big(
{{(H^\vee)}_\h}^{\!\!\vee} \big)}' \, $  in  $ \, \HA_\Bbbk \, $
giving deformations of  $ H^\vee $  with regular fibers
 \vskip-9pt
  $$  \left.  \hbox{$ \begin{matrix}
 \text{\ if \ }  \!\text{\it Char}\,(\Bbbk) = 0 \, ,
\quad\;  U(\gerg_-)  \\
 \text{\ if \ }  \!\text{\it Char}\,(\Bbbk) > 0 \, ,
\quad\;  \u(\gerg_-)
                      \end{matrix} $}  \right\}
 = \widehat{H}  \hskip3pt
\underset{{{(H^\vee)}_\h}^{\!\!\vee}} 
{\overset{0 \,\leftarrow\, \h \,\rightarrow\, 1}
{\longleftarrow\joinrel\relbar\joinrel%
\relbar\joinrel\relbar\joinrel\llongrightarrow}}
\hskip2pt  H^\vee  \hskip0pt
\underset{\;{({{(H^\vee)}_\h}^{\!\!\vee})}'}  
{\overset{1 \,\leftarrow\, \h \,\rightarrow\, 0}
{\longleftarrow\joinrel\relbar\joinrel%
\relbar\joinrel\relbar\joinrel\llongrightarrow}}
\hskip3pt  \hbox{$ \begin{cases}
         F[K_-] = F[G_-^\star]  \\
         F[K_-]
                   \end{cases} $}  $$
 \vskip-3pt
 \noindent
and two 1-parameter families  $ \, {H_\h}^{\!\prime}
= \mathcal{R}^\h_{\underline{D}}(H') \, $  and  $ \,
{({H_\h}^{\!\prime})}^\vee \, $  in  $ \, \HA_\Bbbk \, $
giving deformations
 \vskip-9pt
  $$  F[G_+] = \widetilde{H}  \hskip3pt
\underset{\;{H_\h}^{\!\prime}}  
{\overset{0 \,\leftarrow\, \h \,\rightarrow\, 1}
{\longleftarrow\joinrel\relbar\joinrel%
\relbar\joinrel\relbar\joinrel\llongrightarrow}}
\hskip2pt  H'  \hskip0pt
\underset{\;{({H_\h}^{\!\prime})}^\vee}  
{\overset{1 \,\leftarrow\, \h \,\rightarrow\, 0}
{\longleftarrow\joinrel\relbar\joinrel%
\relbar\joinrel\relbar\joinrel\llongrightarrow}}
\hskip3pt  \hbox{$ \begin{cases}
      U(\gerk_+) = U(\gerg_+^{\,\times})
&   \quad  \text{\ if \ }  \text{\it Char}\,(\Bbbk) = 0  \\
      \u(\gerk_+)
&   \quad  \text{\ if \ }  \text{\it Char}\,(\Bbbk) > 0  \\
                   \end{cases} $}  $$
 \vskip-3pt
%
%
%
%
%
 \noindent
of  $ H' $  with regular fibers, where  $ G_+ $  is like  $ G $  in (c),
$ K_- $  is a connected algebraic Poisson group,  $ \gerg_- $  is like 
$ \gerg $  in (b),  $ \gerk_+ $  is a (restricted, if  $ \, \text{\it
Char}\,(\Bbbk) >  0 \, $)  Lie bialgebra,  $ \gerg_+^{\,\times} $  is
the cotangent Lie bialgebra to  $ G_+ $  and  $ G_-^\star $  is a
connected algebraic Poisson group with cotangent Lie bialgebra 
$ \gerg_- \, $.
                                        \hfill\break
   \indent   (e) \, If  $ \, H = F[G] \, $  is the function algebra
of an algebraic Poisson group  $ G $,  then  $ \widehat{F[G]} $  is
a bi-Poisson Hopf algebra (cf.~[KT], \S 1), namely 
  $$ \; \widehat{F[G]} \; \cong \; S(\gerg^\times) \bigg/ \! \Big(
\Big\{ \overline{x}^{\,p^{n(x)}} \Big\}_{x \in\mathcal{N}_{F[G]}} \Big)
\; \cong \; U(\gerg^\times) \bigg/ \! \Big( \Big\{ \overline{x}^{\,
p^{n(x)}} \Big\}_{x \in \mathcal{N}_{F[G]}} \Big)  $$  
where  $ \mathcal{N}_{F[G]} $  is the nilradical of  $ F[G] $, 
$ \, p^{n(x)} \, $  is the order of nilpotency of  $ \, x \in
\mathcal{N}_{F[G]} $  and the bi-Poisson Hopf structure of 
$ \; S(\gerg^\times) \bigg/ \! \Big( \Big\{ \overline{x}^{\,
p^{n(x)}} \Big\}_{x \in \mathcal{N}_{F[G]}} \Big) \; $  is the
quotient one from  $ \, S(\gerg^\times) \, $;  in particular if
the group  $ G $  is reduced then  $ \; \widehat{F[G]} \cong
S(\gerg^\times) \cong U(\gerg^\times) \; $.
                                      \hfill\break
   \indent   (f) \, If  $ \, \text{Char}\,(\Bbbk) = 0 \, $  and  $ \,
H = U(\gerg) \, $  is the universal enveloping algebra of some Lie
bialgebra  $ \gerg \, $,  then  $ \, \widetilde{U(\gerg)} \, $  is a
bi-Poisson Hopf algebra, namely  $ \; \widetilde{U(\gerg)} \, \cong
\, S(\gerg) \, = \, F[\gerg^\star] \; $  where the bi-Poisson Hopf
structure on  $ S(\gerg) $  is the canonical one.
                                      \hfill\break
   \indent   If  $ \, \text{Char}\,(\Bbbk) = p > 0 \, $  and  $ \,
H = \u(\gerg) \, $  is the restricted universal enveloping algebra
of some restricted Lie bialgebra  $ \gerg \, $,  then  $ \,
\widetilde{\u(\gerg)} \, $  is a bi-Poisson Hopf algebra, namely
we have  $ \; \widetilde{\u(\gerg)} \, \cong \, S(\gerg) \Big/
\big( \big\{ x^p \,\big|\, x \in \gerg \big\} \big) \, = \,
F[G^\star] \; $  where the bi-Poisson Hopf structure on
$ \, S(\gerg) \Big/ \big( \big\{ x^p \,\big|\, x \in \gerg
\big\} \big) \, $  is induced by the canonical one on
$ S(\gerg) $,  and  $ G^\star $  is a connected algebraic
Poisson group of dimension 0 and height 1 whose cotangent
Lie bialgebra is  $ \gerg \, $.
                                      \hfill\break
  \indent   (g) \, Let  $ \, H $,  $ K \in \HA_\Bbbk \, $  and let
$ \, \pi \colon \, H \times K \loongrightarrow \Bbbk \, $  be a Hopf
pairing.  Then  $ \pi $  induce a  {\sl filtered}  Hopf pairing  $ \,
\pi_f \, \colon \, H^\vee \times K' \loongrightarrow \Bbbk \, $,  \,
a  {\sl graded}  Hopf pairing  $ \, \pi_{{}_G} \, \colon \, \widehat{H}
\times \widetilde{K} \loongrightarrow \Bbbk \, $,  both perfect on the
right, and Hopf pairings over  $ \Bbbk[\h\,] $  (notation of\/ \S 3.1)
$ \, H_\h \times K_\h \loongrightarrow \Bbbk[\h\,] \, $  and  $ \,
{H_\h}^{\!\vee} \times {K_\h}^{\!\prime} \loongrightarrow \Bbbk[\h\,]
\, $,  \, the latter being perfect on the right.  If in addition the
pairing  $ \, \pi_f \, \colon \, H^\vee \times K' \loongrightarrow
\Bbbk \, $  is perfect, then all other induced pairings are perfect
as well, and  $ {H_\h}^{\!\vee} $  and  $ {K_\h}^{\!\prime} $  are
dual to each other.  
                                      \hfill\break
   \indent   The left-right symmetrical results hold too.
}\\

%
 \vskip-15pt   
\noindent
{\it Proof.} \, Everything follows from the previous analysis, but
for  {\it (g)},  to be found in [Ga5] or [Ga6].   \qed

\vskip4pt

  {\bf Remarks 3.10.}  \; {\it (a)} \, Though usually introduced
in a different way,  $ H' $  is an object pretty familiar to Hopf
algebraists: it is the  {\sl connected component\/}  of  $ H \, $ 
(see [Ga6] for a proof); in particular,  $ H $  is a PFA iff it
is connected.  Nevertheless, the remarkable properties of  $ \,
\widetilde{H} = G_{\underline{D}}(H) \, $  in  Theorem 3.9{\it
(c)\/}  seems to have been unknown so far.  Similarly, the ``dual''
construction of  $ H^\vee $  and the important properties of  $ \,
\widehat{H} = G_{\underline{J}}(H) \, $  in  Theorem 3.9{\it (b)\/}
seem to be new.
                                            \par
   {\it (b)} \, Theorem 3.9{\it (f)\/}  reminds the classical
formulation of the analogue of Lie's Third Theorem for group-schemes,
i.e.:  {\it Given a restricted Lie algebra  $ \gerg $, there exists
a group-scheme  $ G $  {\sl of dimension 0 and height 1\/}  whose
tangent Lie algebra is  $ \, \gerg \, $}  (see e.g.~[DG]).  Our
result gives just sort of a ``dual Pois\-son-theoretic version'' of
this fact, in that it sounds as follows:  {\it Given a restricted
Lie algebra  $ \gerg \, $,  there exists a Poisson group-scheme 
$ G $  {\sl of dimension 0 and height 1\/}  
   \hbox{whose  {\sl cotangent}  Lie
algebra is  $ \, \gerg \, $}.}   
                                    \par
   {\it (c)} \, Part  {\it (d)}  of Theorem 3.9 is quite interesting
for applications in physics.  In fact, let  $ H $  be a Hopf algebra
which describes the symmetries of some physically meaningful system,
but has no geometrical meaning, and assume also  $ \, H' = H = H^\vee
\, $.  Then  Theorem 3.9{\it (d)}  yields a recipe to deform  $ H $ 
to four Hopf algebras with geometrical content, which means having
two Poisson groups and two Lie bialgebras attached to  $ H $,  hence
a rich ``Poisson geometrical symmetry'' underlying the physical
system.  As  $ \Bbb{R} $  (the typical ground field) has zero
characteristic, we have in fact two pairs of mutually dual Poisson
groups along with their tangent Lie bialgebras.  A nice application
is in [Ga7].

\vskip7pt

  {\bf 3.11. The hyperalgebra case.} \, Let  $ G $  be an algebraic
group, which for simplicity we assume to be finite-dimensional.  Let 
$ \hyp(G) $  be the hyperalgebra of  $ G $  (cf.~\S 1.1), which is
connected cocommutative.  Recall also the Hopf algebra morphism  $ \,
\Phi : U(\gerg) \longrightarrow \hyp(G) \, $;  \, if  $ \,
\text{\it Char}\,(\Bbbk) = 0 \, $  then  $ \Phi $  is an isomorphism,
so  $ \hyp(G) $  identifies to  $ U(\gerg) $;  \, if  $ \, \text{\it
Char}\,(\Bbbk) > 0 \, $  then  $ \Phi $  factors through  $ \u(\gerg) $
and the induced morphism  $ \, \overline{\Phi} : \u(\gerg)
\longrightarrow \hyp(G) \, $  is injective, so that  $ \u(\gerg) $
identifies with a Hopf subalgebra of  $ \hyp(G) $.  Now we study
$ {\hyp(G)}' $,  $ {\hyp(G)}^\vee $,  $ \widetilde{\hyp(G)} $,
$ \widehat{\hyp(G)} $,  the key tool being the existence of a
perfect (= non-degenerate) Hopf pairing between  $ F[G] $  and 
$ \hyp(G) $.
                                               \par
   One can prove (see [Ga6]) that a Hopf  $ \Bbbk $--algebra  $ H $
is connected iff  $ \, H = H' $.  As  $ \hyp(G) $  is connected, we
have  $ \, \hyp(G) = {\hyp(G)}' \; $.  Now,  Theorem 3.9{\it (c)\/}
gives  $ \, \widetilde{\hyp(G)} := G_{\underline{D}}\big(\hyp(G)\big)
= F[\varGamma\,] \, $  for some connected algebraic Poisson group
$ \varGamma \, $;  Theorem 3.9{\it (e)}  yields  
 \vskip-17pt   
  $$  \widehat{F[G]} \, \cong \, S(\gerg^*) \bigg/ \! \Big( \Big\{\,
\overline{x}^{\,p^{n(x)}} \Big\}_{x \in \mathcal{N}_{F[G]}} \Big)
\, = \, \u \bigg( P \bigg( S(\gerg^*) \bigg/ \! \Big( \Big\{\,
\overline{x}^{\,p^{n(x)}} \Big\}_{x \in \mathcal{N}_{F[G]}} \Big)
\bigg) \bigg) \, = \, \u \Big( \big( \gerg^* \big)^{p^\infty} \Big)  $$    
 \vskip-1pt   
\noindent   
with  $ \, \big( \gerg^* \big)^{p^\infty} := \text{\sl Span}\, \Big(
\Big\{\, x^{p^n} \,\Big\vert\; x \in \gerg^* \, , n \in \N \,\Big\} \Big)
\subseteq \widehat{F[G]} \, $,  \, and noting that  $ \, \gerg^\times
= \gerg^* \, $.  On the other hand, exactly like for  $ U(\gerg) $
and  $ \u(\gerg) $  respectively in case  $ \, \Char(\Bbbk) = 0 \, $
and  $ \, \Char(\Bbbk) > 0 \, $,  \, the filtration  $ \underline{D} $
of  $ \hyp(G) $  is nothing but the natural filtration given by the
order of differential operators: this implies immediately  $ \;
{\hyp(G)_\h}' := {\big( \Bbbk[\h\,] \otimes_\Bbbk \hyp(G) \big)}'
\, = \, \big\langle \big\{\, \h^n x^{(n)} \,\big|\; x \in \gerg \, ,
n \in \N \,\big\} \big\rangle \, $,  \, where  $ x^{(n)} $  denotes
the  $ n $--th  divided power of  $ \, x \in \gerg \, $  (recall that
$ \hyp(G) $  is generated as an algebra by all the  $ x^{(n)} $'s,
some of which might be zero).  It is then immediate to check that
the graded Hopf pairing between  $ \, {\hyp(G)_\h}' \Big/ \h \,
{\hyp(G)_\h}' = \widetilde{\hyp(G)} = F[\varGamma] \, $  and  $ \,
\widehat{F[G]} \, $  from  Theorem 3.9{\it (f)\/}  is perfect.
From this one argues that the cotangent Lie bialgebra of
$ \varGamma $  is isomorphic to  $ \, \Big( \big( \gerg^*
\big)^{p^\infty} \Big)^* \, $.
                                               \par
   As for  $ {\hyp(G)}^\vee $  and  $ \widehat{\hyp(G)} $,  the
situation is much like for  $ U(\gerg) $  and  $ \u(\gerg) $,
in that it strongly depends on the algebraic nature of  $ G $
(cf.~\S 3.7).

\vskip7pt

  {\bf 3.12 The CDP on group algebras and their duals.} \, In
this section,  $ G $  is any abstract group.  We divide the
subsequent material in several subsections.

\vskip5pt

\noindent   {\it  $ \underline{\hbox{\sl Group-related algebras}} $.}
\, For any commutative unital ring  $ \A \, $,  by  $ \, \A[G] \, $
we mean the group algebra of  $ G $  over  $ \A \, $;  \, when  $ G $
is  {\sl finite},  we denote by  $ \, A_\A(G) := {\A[G]}^* \, $  (the
linear dual of  $ \A[G] \, $)  the function algebra of  $ G $  over
$ \A \, $.  Our aim is to apply the Crystal Duality Principle to
$ \Bbbk[G] $  and  $ A_\Bbbk(G) \, $  with their standard Hopf algebra
structure: hereafter  $ \Bbbk $  is a field and  $ \, R := \Bbbk[\h]
\, $  as in \S 5.1, with  $ \, p := \text{\it Char}\,(\Bbbk) \, $.
                                          \par
   Recall that  $ \, H := \A[G] \, $  admits  $ G $  itself
as a distinguished basis, with Hopf algebra structure given by
$ \, g \cdot_{{}_H} \gamma := g \cdot_{{}_G} \gamma \, $,  $ \,
1_{{}_H} := 1_{{}_G} \, $,  $ \, \Delta(g) := g \otimes g \, $,
$ \, \epsilon(g) := 1 \, $,  $ \, S(g) := g^{-1} \, $,  \, for all
$ \, g, \gamma \in G \, $.  Dually,  $ \, H := A_\A(G) \, $  has
basis  $ \, \big\{ \varphi_g \,\big|\, g \! \in \! G \big\} \, $
dual to the basis  $ G $  of  $ \A[G] \, $,  \, with  $ \, \varphi_g
(\gamma) := \delta_{g,\gamma} \, $  for all  $ \, g, \gamma \in G
\, $;  its Hopf algebra structure is given by  $ \, \varphi_g \cdot
\varphi_\gamma := \delta_{g,\gamma} \varphi_g \, $,  $ \, 1_{{}_H}
:= \sum_{g \in G} \varphi_g \, $,  $ \, \Delta(\varphi_g) :=
\sum_{\gamma \cdot \ell = g} \varphi_\gamma \otimes \varphi_\ell
\, $,  $ \, \epsilon(\varphi_g) := \delta_{g,1_G} \, $,  $ \,
S(\varphi_g) := \varphi_{g^{-1}} \, $,  \, for all  $ \, g, \gamma
\in G \, $.  In particular,  $ \, R[G] = R \otimes_\Bbbk \Bbbk[G]
\, $  and  $ \, A_R[G] = R \otimes_\Bbbk A_\Bbbk[G] \, $.  Our
first result is

\vskip5pt

\noindent   {\it  $ \underline{\hbox{\sl Theorem A}} $:}  $ \,
{{\big(\Bbbk[G]\big)}_\h}^{\!\prime} = R \cdot 1 \, $,  $ \;
{\Bbbk[G]}' = \Bbbk \cdot 1 \; $  {\it and}  $ \; \widetilde{\Bbbk[G]}
= \Bbbk \cdot 1 = F\big[\{*\}\big] \, $.

%
%
\noindent
 {\it Proof.} \, 
The claim follows easily from the formula  $ \, \delta_n(g)
= {(g - 1)}^{\otimes n} \, $,  \, for  $ \, g \in G \, $, 
$ \, n \in \N \, $.   \qed
%
%

\vskip5pt

\noindent   {\it  $ \underline{\hbox{\sl  $ {R[G]}^\vee $,  $ \,
{\Bbbk[G]}^\vee $,  $ \widehat{\Bbbk[G]} $  \ and the dimension
subgroup problem}} $.} \, In contrast with the triviality result in
Theorem A above, things are more interesting for  $ \, {R[G]}^\vee \!
= {{\big( \Bbbk[G] \big)}_\h}^{\!\!\vee} \, $,  $ \, {\Bbbk[G]}^\vee $
and  $ \, \widehat{\Bbbk[G]} \, $.  Note however that since  $ \,
\Bbbk[G] \, $  is cocommutative the induced Poisson cobracket on
$ \, \widehat{\Bbbk[G]} \, $  is trivial, hence the Lie cobracket
of  $ \, \gerk_G := P \Big( \widehat{\Bbbk[G]} \Big) \, $  is
trivial as well.
                                          \par
   Studying  $ {\Bbbk[G]}^\vee $  and  $ \widehat{\Bbbk[G]} $
amounts to study the filtration  $ {\big\{ J^n \big\}}_{n \in \N}
\, $,  \, with  $ J := \text{\sl Ker}\,(\epsilon_{{}_{\Bbbk[G]}}) $,
\, which is a classical topic.  Indeed, for  $ \, n \! \in \! \N \, $
let  $ \, D_n(G) := \big\{\, g \in G \,\big|\, (g \! - \! 1) \in J^n
\,\big\} \, $: \, this is a characteristic subgroup of  $ G $,  called
{\sl the  $ n^{\text{th}} $  dimension subgroup of  $ G \, $}.  All these
form a filtration inside  $ G \, $:  \, characterizing it in terms of
$ G $  is the  {\sl dimension subgroup problem},  which (for group
algebras over fields) is completely solved (see [Pa], Ch.~11, \S 1,
and [HB], and references therein); this also gives a description of
$ \, \big\{ J^n \big\}_{n \in \N_+} \, $.  Thus we find ourselves
within the domain of classical group theory: now we use the results
which solve the dimension subgroup problem to argue a description of
$ {\Bbbk[G]}^\vee $,  $ \widehat{\Bbbk[G]} $  and  $ {R[G]}^\vee $,
and later on we'll get from this a description
      \hbox{of  $ \big( {R[G]}^\vee \big)' $  and its semiclassical
limit too.}
                                          \par
   By construction,  $ J $  has  $ \Bbbk $--basis  $ \, \big\{ \eta_g
\,\big|\; g \! \in \! G \setminus \{1_{{}_G}\} \big\} \, $,  \, where
$ \, \eta_g := (g\!-\!1) \, $.  Then  $ \, {\Bbbk[G]}^\vee \, $  is
generated by  $ \, \big\{\, \eta_g \! \mod J^\infty \,\big|\; g \in
G \setminus \{1_{{}_G}\} \big\} \, $,  \, and  $ \, \widehat{\Bbbk[G]}
\, $  by  $ \, \big\{\, \overline{\,\eta_g} \;\big|\; g \! \in \! G
\setminus \{1_{{}_G}\} \big\} \, $:  \, hereafter  $ \, \overline{x}
:= x \mod J^{n+1} \, $  for all  $ \, x \in J^n \, $,  \, that is
$ \overline{x} $  is the element in  $ \widehat{\Bbbk[G]} $  which
corresponds to  $ \, x \in \Bbbk[G] \, $.  Moreover,  $ \, \overline{g}
= \overline{\,1 + \eta_g} = \overline{1} \, $  for all  $ \, g \in G \, $;
\, also,  $ \; \Delta \big(\overline{\,\eta_g}\big) = \overline{\,\eta_g}
\otimes \overline{g} + 1 \otimes \overline{\,\eta_g} = \overline{\,\eta_g}
\otimes 1 + 1 \otimes \overline{\,\eta_g} \; $:  \; thus  $ \overline{\,
\eta_g} $  is primitive, so  $ \, \big\{\,\overline{\,\eta_g} \;\big|\;
g \! \in \! G \setminus \{1_{{}_G}\} \big\} \, $  generates  $ \,
\gerk_G := P \Big( \widehat{\Bbbk[G]} \Big) \, $.

\vskip5pt

\noindent   {\it  $ \underline{\hbox{\sl The Jennings-Hall theorem}} $.}
\, The description of  $ D_n(G) $  is given by the Jennings-Hall theorem,
which we now recall.  The construction involved strongly depends on
whether  $ \, p := \text{\it Char}\,(\Bbbk) \, $  is zero or not, so
we shall distinguish these two cases.
                                          \par
   First assume  $ \, p = 0 \, $.  Let  $ \, G_{(1)} := G \, $,  $ \,
G_{(k)} := (G,G_{(k-1)}) \, $  ($ k \in \N_+ $),  form the  {\sl lower
central series\/}  of  $ G \, $;  hereafter  $ (X,Y) $  is the
commutator subgroup of  $ G $  generated by the set of commutators
$ \, \big\{ (x,y) := x \, y \, x^{-1} y^{-1} \,\big|\, x \in X, y
\in Y \big\} \, $:  \, this is a  {\sl strongly central series\/}
in  $ G $,  which means a central series  $ \, {\big\{ G_k \big\}}_{k
\in \N_+} \, $  (\,= decreasing filtration of normal subgroups, each
one centralizing the previous one) of  $ G $  such that  $ \, (G_m,
G_n) \leq G_{m+n} \, $  for all  $ m \, $,  $ n \, $.  Then let  $ \,
\sqrt{G_{(n)}} := \big\{ x \in G \,\big|\, \exists \, s \in \N_+ : x^s
\in G_{(n)} \big\} \, $  for all  $ n \in \N_+ \, $:  these form a
descending series of characteristic subgroups in  $ G $,  such that
each composition factor  $ \, A^G_{(n)} := \sqrt{G_{(n)}} \Big/ \!
\sqrt{G_{{(n+1)}}} \, $  is torsion-free Abelian.  Therefore  $ \,
\L_0(G) := \bigoplus_{n \in \N_+} A^G_{(n)} \, $  is a graded Lie ring,
with Lie bracket  $ \, \big[ \overline{g}, \overline{\ell} \,\big] :=
\overline{(g,\ell\,)} \, $  for all  {\sl homogeneous}  $ \overline{g} $,
$ \overline{\ell} \in \L_0(G) \, $,  \, with obvious notation.  It is
easy to see that the map  $ \; \Bbbk \otimes_\Z \L_0(G) \longrightarrow
\gerk_G \, $,  $ \, \overline{g} \mapsto \overline{\eta_g} \, $,  \;
is an epimorphism  {\sl of graded Lie rings\/}:  \, thus  {\sl
the Lie algebra  $ \, \gerk_G \, $  is a quotient of  $ \; \Bbbk
\otimes_\Z \L_0(G) \, $};  in fact, the above is an isomorphism,
see below.  We use notation  $ \, \partial(g) := n \, $  for
all  $ \, g \in \sqrt{G_{(n)}} \, \setminus \sqrt{G_{(n+1)}} \; $.
                                     \par
   For each $ \, k \in \N_+ \, $  pick in  $ A^G_{(k)} $  a subset
$ \overline{B}_k $  which is a  $ \Bbb{Q} $--basis  of  $ \, \Bbb{Q}
\otimes_\Z A^G_{(k)} \, $;  \, for each  $ \, \overline{b} \in
\overline{B}_k \, $,  \, choose a fixed  $ \, b \in \sqrt{G_{(k)}}
\, $  such that its coset in  $ A^G_{(k)} $  be  $ \overline{b} $,  \,
and denote by  $ \, B_k \, $  the set of all such elements  $ b \, $.
Let  $ \, B := \bigcup_{k \in \N_+} B_k \, $:  \, we call such a set
{\sl t.f.l.c.s.-net\/}  (\,=\,``torsion-free-lower-central-series-net'')
on  $ G $.  Clearly  $ \, B_k = \Big( B \cap \sqrt{G_{(k)}} \,\Big)
\setminus \Big( B \cap \sqrt{G_{(k+1)}} \,\Big) \, $  for all
$ k \, $.  By an  {\sl ordered t.f.l.c.s.-net\/}  is meant a
t.f.l.c.s.-net  $ B $  which is totally ordered in such a way that:
{\it (i)\/}  if  $ \, a \in B_m \, $,  $ \, b \in B_n \, $,  $ \, m
< n \, $,  \, then  $ \, a \preceq b \, $;  \, {\it (ii)\/}  for each
$ k $,  every non-empty subset of  $ B_k $  has a greatest element.
As a matter of fact, an ordered t.f.l.c.s.-net always exists.
                                     \par
   Now assume instead  $ \, p > 0 \, $.  The situation is similar, but
we must also consider the  $ p $--power  operation in the group  $ G $
and in the restricted Lie algebra  $ \gerk_G \, $.  Starting from the
lower central series  $ \, {\big\{G_{(k)}\big\}}_{k \in \N_+} $,  define
$ \, G_{[n]} := \prod_{k p^\ell \geq n} {(G_{(k)})}^{p^\ell} \; $  for
all  $ \, n \in \N_+ \, $  (hereafter, for any group  $ \varGamma $  we
denote  $ \varGamma^{p^e} $  the subgroup generated by  $ \, \big\{
\gamma^{p^e} \,\big|\, \gamma \! \in \! \varGamma \,\big\} \, $):  \,
this gives another strongly central series  $ \, {\big\{G_{[n]}\big\}}_{n
\in \N_+} $  in  $ G $,  \, with the additional property that  $ \,
{(G_{[n]})}^p \leq G_{[n+1]} \, $  for all  $ n \, $,  \, called
{\sl the  $ p $--lower  central series of\/}  $ G \, $.  Then  $ \,
\mathcal{L}_p(G) := \oplus_{n \in \N_+} G_{[n]} \big/ G_{[n+1]} \, $  is
a graded restricted Lie algebra over  $ \, \Z_p := \Z \big/ p \, \Z \, $,
\, with operations  $ \, \overline{g} + \overline{\ell} := \overline{g
\cdot \ell} \, $,  $ \, \big[\overline{g},\overline{\ell}\,\big] :=
\overline{(g,\ell\,)} \, $,  $ \, \overline{g}^{\,[p\,]} := \overline{g^p}
\, $,  \, for all  $ \, g $,  $ \ell \in G \, $  (cf.~[HB], Ch.~VIII,
\S 9).  Like before, we consider the map  $ \; \Bbbk \otimes_{\Z_p}
\mathcal{L}_p(G) \longrightarrow \gerk_G \, $,  $ \, \overline{g} \mapsto
\overline{\eta_g} \, $,  \; which now is an epimorphism  {\sl of graded
restricted Lie  $ \Z_p $--algebras},  whose image spans  $ \gerk_G $
over  $ \Bbbk \, $:  \, therefore  {\sl  $ \, \gerk_G \, $  is a quotient
of  $ \; \Bbbk \otimes_{\Z_p} \mathcal{L}_p(G) \, $};  \, in fact, the
above is an isomorphism, see below.  Finally, we introduce also the 
notation  $ \, d(g) := n \, $  for all  $ \, g \in G_{[n]} \setminus
G_{[n+1]} \, $.  
                                      \par
   For each  $ \, k \in \N_+ \, $  choose a  $ \Z_p $--basis
$ \overline{B}_k $  of the  $ \Z_p $--vector  space  $ \, G_{[k]} \big/
G_{[k+1]} \, $;  \, for each  $ \, \overline{b} \in \overline{B}_k \, $,
\, fix  $ \, b \in G_{[k]} \, $  such that  $ \, \overline{b} = b \,
G_{[k+1]} \, $,  \, and let  $ \, B_k \, $  be the set of all such
elements  $ b \, $.  Let  $ \, B := \bigcup_{k \in \N_+} B_k \, $:
\, such a set will be called a  {\sl  $ p $-l.c.s.-net\/}  (=
``$ p $-lower-central-series-net''; the terminology in [HB] is
``$ \kappa $-net'') on  $ G $.  Of course  $ \, B_k = \big( B \cap
G_{[k]} \big) \setminus \big( B \cap G_{[k+1]} \big) \, $  for all
$ k \, $.  By an  {\sl ordered  $ p $-l.c.s.-net\/}  we mean a
$ p $-l.c.s.-net  $ B $  which is totally ordered in such a way
that:  {\it (i)\/}  if  $ \, a \in B_m \, $,  $ \, b \in B_n \, $,
$ \, m < n \, $,  \, then  $ \, a \preceq b \, $;  \, {\it (ii)\/}
for each  $ k $,  every non-empty subset of  $ B_k $  has a greatest
element (like for  $ \, p = 0 \, $).  Again, it is known that 
$ p $-l.c.s.-nets  always do exist.
                                      \par
   We can now describe each  $ D_n(G) $,  hence also each graded summand
$ \, J^n \big/ J^{n+1} \, $  of  $ \, \widehat{\Bbbk[G]} $,  in terms
of the lower central series or the  $ p $--lower  central series of
$ G \, $,  more precisely in terms of a fixed ordered t.f.l.c.s.-net
or  $ p $-l.c.s.-net.  To unify notations, set  $ \, G_n := G_{(n)}
\, $,  $ \, \theta(g) := \partial(g) \, $  if  $ \, p \! = \! 0 \, $,
\, and  $ \, G_n := G_{[n]} \, $,  $ \, \theta(g) := d(g) \, $  if
$ \, p \! > \! 0 \, $,  \, set  $ \, G_\infty := \bigcap_{n \in N_+}
\!\! G_n \, $,  \, let  $ \, B := \bigcup_{k \in \N_+} B_k \, $  be
an ordered t.f.l.c.s.-net or  $ p $-l.c.s.-net according to whether
$ \, p \! = \! 0 \, $  or  $ \, p \! > \! 0 \, $,  \, and set  $ \,
\ell(0) := + \infty \, $  and  $ \, \ell(p) := p \, $  for  $ \, p
> 0 \, $.  The key result we need is

\vskip4pt

\noindent   {\it  $ \underline{\hbox{\sl Jennings-Hall theorem}} $
(cf.~[HB], [Pa] and references therein).  Let  $ \, p:= \text{\it
Char}\,(\Bbbk) \, $.
                                       \par
   (a) \, For all  $ \, g \in G \, $,  $ \; \eta_g
\in J^n \Longleftrightarrow g \in \! G_n \, $.  Therefore
$ \, D_n(G) = G_n \; $  for all  $ \, n \in \N_+ \, $.
                                       \par
   (b) \, For any  $ \, n \in \N_+ \, $,  the set of ordered monomials
 \vskip-23pt
  $$  \Bbb{B}_n \, := \, \Big\{\, {\overline{\,\eta_{b_1}}}^{\;e_1}
\cdots {\overline{\,\eta_{b_r}}}^{\;e_r} \;\Big|\; b_i \in B_{d_i}
\, , \; e_i \in \N_+ \, , \; e_i < \ell(p) \, , \; b_1 \precneqq
\cdots \precneqq b_r \, , \; {\textstyle \sum}_{i=1}^r e_i \, d_i
= n \,\Big\}  $$
 \vskip-9pt
\noindent   is a\/  $ \Bbbk $--basis  of  $ \, J^n \big/ J^{n+1} \, $,
\, and  $ \; \Bbb{B} \, := \, \{1\} \cup \bigcup_{n \in \N} \Bbb{B}_n
\; $  is a\/  $ \Bbbk $--basis  of  $ \; \widehat{\Bbbk[G]} \, $.
                                          \par
   (c) \;  $ \big\{\, \overline{\,\eta_b} \,\;\big|\;
b \in B_n \big\} \; $  is a  $ \Bbbk $--basis  of the  $ n $--th
graded summand  $ \, \gerk_G \cap \big( J^n \big/ J^{n+1} \big) \, $
of the graded restricted Lie algebra\/  $ \gerk_G \, $,  \, and
$ \, \big\{\, \overline{\,\eta_b} \,\;\big|\; b \in B \,\big\} \; $
is a\/  $ \Bbbk $--basis  of\/  $ \gerk_G \, $.
                                          \par
   (d) \;  $ \big\{\, \overline{\,\eta_b} \,\;\big|\;
b \in B_1 \big\} \; $  is a minimal set of generators of
the (restricted) Lie algebra\/  $ \gerk_G \, $.
                                          \par
   (e) \; The map  $ \; \Bbbk \otimes_\Z \L_p(G) \longrightarrow
\gerk_G \, $,  $ \, \overline{g} \mapsto \overline{\,\eta_g} \, $,
\, is an isomorphism of graded restricted Lie algebras.  Therefore
$ \; \widehat{\Bbbk[G]} \, \cong \, \U \big( \Bbbk \otimes_\Z \L_p(G)
\big) \; $  as Hopf algebras.
                                          \par
   (f) \;  $ J^\infty \, = \, \hbox{\sl Span} \big( \big\{\,
\eta_g \,\big|\, g \in G_\infty \,\big\} \big) \, $,  \, whence\/
$ \; {\Bbbk[G]}^\vee \cong \, \bigoplus_{\overline{g} \in G/G_\infty}
\! \Bbbk \cdot \overline{g} \; \cong \, \Bbbk \big[ G \big/ G_\infty
\big] \; $.   \qed}

\vskip8pt

   Recall that  $ A\big[x,x^{-1}\big] $  (for any  $ A $)  has
$ A $--basis  $ \, \big\{ {(x \!-\! 1)}^n x^{-[n/2]} \,\big|\, n \in
\N \big\} \, $,  \, where  $ [q] $  is the integer part of  $ \, q \in
\Bbb{Q} \, $.  Then from Jennings-Hall theorem and (5.2) we argue

\vskip8pt

\noindent   {\it  $ \underline{\hbox{\sl Proposition B}} $.
\, Let  $ \; \chi_g := \h^{-\theta(g)} \eta_g \, $,  \, for
all  $ \, g \in \{G\} \setminus \{1\} \, $.  Then
  $$  \displaylines{
   \hskip25pt   {R[G]}^\vee  \, = \;  \Big( {\textstyle
\bigoplus_{\substack{ \, b_i \in B, \; 0 < e_i < \ell(p)  \\
                r \in \N, \; b_1 \precneqq \cdots \precneqq b_r }}} \,
R \cdot \chi_{b_1}^{\;e_1} \, b_1^{\,-[e_1\!/2]} \cdots \chi_{b_r}^{\;
e_r} \, b_r^{\,-[e_r/2]} \Big) \,{\textstyle
\bigoplus}\; R\big[\h^{-1}\big] \cdot J^\infty  \; =   \hfill \qquad  \cr
   {} \; \hfill   = \;  \Big( {\textstyle
\bigoplus_{\substack{ \, b_i \in B, \; 0 < e_i < \ell(p)  \\
                r \in \N, \; b_1 \precneqq \cdots \precneqq b_r }}} \,
R \cdot \chi_{b_1}^{\,e_1} \, b_1^{\,-[e_1\!/2]} \cdots \chi_{b_r}^{\,
e_r} \, b_r^{\,-[e_r/2]} \Big) \,{\textstyle \bigoplus}\,
\Big( \hskip1pt {\textstyle \sum_{\gamma \in G_\infty}}
\hskip-0pt R\big[\h^{-1}\big] \cdot \eta_\gamma \Big) \; ;
\hskip35pt  \cr }  $$
If  $ \, J^\infty \! = \! J^n $  for some  $ \, n \! \in \! \N $
(iff  $ \, G_\infty \! = G_n $)  we can drop the factors  $ \,
b_1^{-[e_1\!/2]}, \dots, b_r^{-[e_r/2]} \, $.   \qed }  

\vskip5pt

\noindent   {\it  $ \underline{\hbox{\sl Poisson groups from
$ \Bbbk[G] $}} $.} \, The previous discussion attached to the
abstract group  $ G $  the (maybe restricted) Lie algebra  $ \gerk_G $
which, by Jennings-Hall theorem, is just the scalar extension of the Lie
ring  $ \L_{\text{\it Char}(\Bbbk)} $  associated to  $ G $  via the
central series of the  $ G_n $'s;  in particular the functor  $ \, G
\mapsto \gerk_G \, $  is one considered since long in group theory. 
Now, by  Theorem 5.8{\it (d)\/}  we know that  $ \, {\big( {R[G]}^\vee
\big)}' \, $  is a QFA, with  $ \, {\big({R[G]}^\vee \big)}'{\Big|}_{\h=0}
= F\big[\varGamma_G\big] \, $  for some connected Poisson group
$ \varGamma_G \, $.  This defines a functor  $ \, G \mapsto \varGamma_G
\, $  from abstract groups to connected Poisson groups, of dimension
zero and height 1 if  $ \, p > 0 \, $;  \, in particular, this
$ \varGamma_G $  {\sl is a new invariant for abstract groups}.
                                            \par
   The description of  $ {R[G]}^\vee $  in Proposition B above leads
us to an explicit description of  $ {\big({R[G]}^\vee \big)}' \, $, 
hence of  $ \, {\big({R[G]}^\vee \big)}'{\Big|}_{\h=0} \! = F \big[
\varGamma_G\big] \, $  and of  $ \varGamma_G \, $.  Indeed direct
inspection gives  $ \, \delta_n\big(\chi_g\big) = \h^{(n-1) \theta(g)}
\chi_g^{\;\otimes n} \, $,  \, so  $ \, \psi_g := \h \, \chi_g = \h^{1
- \theta(g)} \eta_g \in {\big( {R[G]}^\vee \big)}' \setminus \h \,
{\big( {R[G]}^\vee \big)}' \, $  for each  $ \, g \in G \setminus
G_\infty \, $,  \, whilst for  $ \, \gamma \in G_\infty $  we have
$ \, \eta_\gamma \in J^\infty \, $  which implies  $ \, \eta_\gamma
\in {\big({R[G]}^\vee \big)}' \, $  and even  $ \, \eta_\gamma \in
\bigcap_{n \in \N} \h^n {\big({R[G]}^\vee \big)}' \, $.  Thus
$ {\big({R[G]}^\vee\big)}' $  is generated by  $ \, \big\{\, \psi_g
\;\big|\; g \in G \setminus \{1\} \big\} \cup \big\{\, \eta_\gamma
\,\big|\, \gamma \in G_\infty \big\} \, $.  Moreover,  $ \, g = 1 +
\h^{\theta(g)-1} \psi_g \in {\big( {R[G]}^\vee \big)}' \, $  for every
$ \, g \in G \setminus G_\infty \, $,  \, and  $ \, \gamma = 1 + (\gamma
- 1) \in 1 + J^\infty \subseteq {\big({R[G]}^\vee \big)}' \, $  for  $ \,
\gamma \in G_\infty \, $.  This and the previous analysis along with
Proposition B prove next result, which in turn is the basis for
Theorem D below.

\vskip8pt

\noindent   {\it  $ \underline{\hbox{\sl Proposition C}} $.
 \vskip-19pt
  $$  \displaylines{
   \hskip15pt   {\big({R[G]}^\vee\big)}' \; = \; \Big( {\textstyle
\bigoplus_{\substack{ \, b_i \in B, \; 0 < e_i < \ell(p)  \\
                r \in \N, \; b_1 \precneqq \cdots \precneqq b_r }}} \,
R \cdot \psi_{b_1}^{\;e_1} \, b_1^{\,-[e_1\!/2]} \cdots \psi_{b_r}^{\;
e_r} \, b_r^{\,-[e_r/2]} \Big) \,{\textstyle \bigoplus}\;
R\big[\h^{-1}\big] \cdot J^\infty  \; =   \hfill  \cr
   {} \hskip32pt   \; = \;  \Big( {\textstyle
\bigoplus_{\substack{ \, b_i \in B, \; 0 < e_i < \ell(p)  \\
                r \in \N, \; b_1 \precneqq \cdots \precneqq b_r }}} \,
R \cdot \psi_{b_1}^{\,e_1} \, b_1^{\,-[e_1\!/2]} \cdots \psi_{b_r}^{\,
e_r} \, b_r^{\,-[e_r/2]} \Big) \,{\textstyle \bigoplus}\,
\Big( \hskip1pt {\textstyle \sum_{\gamma \in G_\infty}}
\hskip-0pt R\big[\h^{-1}\big] \cdot \eta_\gamma \Big) \; .  \cr }  $$
 \vskip-3pt
In particular,  $ \; {\big({R[G]}^\vee\big)}' = R[G] \; $  {\sl if
and only if  $ \; G_2 = \{1\} = G_\infty \; $.}  If in addition  $ \,
J^\infty \! = \! J^n \, $  for some  $ \, n \! \in \! \N \, $  (iff
$ \, G_\infty = G_n $)  then we can drop the factors  $ \, b_1^{\,
-[e_1\!/2]}, \dots, b_r^{\,-[e_r/2]} \, $.   \qed}

\vskip8pt

\noindent   {\it  $ \underline{\hbox{\sl Theorem D}} $.  \, Let
$ \; x_g := \psi_g \mod \h \; {\big( {R[G]}^\vee \big)}' \, $,
$ \, z_g := g \mod \h \; {\big( {R[G]}^\vee \big)}' \, $  for
all  $ \, g \not= 1 \, $,  \, and  $ \, B_1 := \big\{\, b \in
\! B \,\big\vert\, \theta(b) = 1 \big\} \, $,  $ \, B_> :=
\big\{\, b \in \! B \,\big\vert\, \theta(b) > 1 \big\} \, $.
                                        \hfill\break
   \indent   (a) \, If  $ \, p = 0 \, $,  \, then  $ \, F \big[
\varGamma_G \big] = {\big( {R[G]}^\vee \big)}'{\Big|}_{\h=0} $  is
{\sl a polynomial/Laurent polynomial algebra}, namely  $ \, F \big[
\varGamma_G \big] = \Bbbk \big[ {\{x_b\}}_{b \in B_>} \cup {\big\{
{z_b}^{\!\pm 1} \big\}}_{b \in B_1} \big] \, $,  \; the  $ x_b $'s 
being primitive and the  $ z_b $'s  being group-like.  In particular 
$ \, \varGamma_G \cong \big( \Bbb{G}_a^{\, \times B_>} \big) \times
\big( \Bbb{G}_m^{\,\times B_1} \big) \, $ as algebraic groups,
i.e.~$ \varGamma_G $    
   \hbox{is a (pro)affine space times a torus.}
                                        \hfill\break
   \indent   (b) \, If  $ \, p > 0 \, $,  \, then  $ \, F \big[
\varGamma_G \big] = {\big( {R[G]}^\vee \big)}'{\Big|}_{\h=0} $  is
{\sl a truncated polynomial/Laurent polynomial algebra}, namely  $ \,
F \big[ \varGamma_G \big] = \, \Bbbk \big[ {\{x_b\}}_{b \in B_>} \cup
{\big\{ {z_b}^{\!\pm 1} \big\}}_{b \in B_1} \big] \Big/ \big( \{ x_b^{\,p}
\}_{b \in B_>} \cup \{ z_b^{\,p} - 1 \} \big) \, $,  \; the  $ x_b $'s
being primitive and the  $ z_b $'s  being group-like.  In particular
$ \, \varGamma_G \cong \big( {{\boldsymbol\alpha}_p}^{\!\times B_>}
\big) \times \big( {{\boldsymbol\mu}_p}^{\! \times B_1} \big)
\, $  as algebraic groups of dimension zero and height 1.
                                        \hfill\break
   \indent   (c) \, The Poisson group  $ \varGamma_G $  has cotangent
Lie bialgebra  $ \gerk_G \, $,  that is  $ \, \text{\it coLie}\,
(\varGamma_G) = \gerk_G \, $.}

%
%
 \vskip4pt
 \noindent   {\it Proof.}
{\it (a)} \, The very definitions give  $ \, \partial(g\,\ell\,)
\geq \partial(g) + \partial(\ell\,) \, $  for all  $ \, g, \ell
\in G \, $,  \, so that  $ \; [\psi_g,\psi_\ell\big] =
%
%
\h^{1 - \partial(g)
- \partial(\ell) + \partial((g,\ell))} \, \psi_{(g,\ell)} \, g \, \ell
\in \h \cdot {\big({R[G]}^\vee\big)}' \, $,  \; which proves (directly)
that  $ \, {\big({R[G]}^\vee\big)}'{\Big|}_{\h=0} \, $  is commutative!
Moreover, the relation  $ \, 1 = g^{-1} \, g =
g^{-1} \, \big(1 + \h^{\partial(g)-1} \psi_g \big) \, $  (for any
$ \, g \in G \, $)  yields  $ \, z_{g^{-1}} = {z_g}^{\!-1} \, $  iff
$ \, \partial(g) = 1 \, $  and  $ \, z_{g^{-1}} = 1 \, $  iff  $ \,
\partial(g) > 1 \, $.  Noting also that  $ \, J^\infty \equiv 0 \!\!
\mod \h \, {\big({R[G]}^\vee\big)}' \, $  and  $ \, g = 1 + \h^{\partial(g)
- 1} \psi_g \equiv 1 \mod \h \, {\big({R[G]}^\vee \big)}' \, $  for  $ \,
g \in G \setminus G_\infty \, $,  \, and also  $ \, \gamma = 1 +
(\gamma - 1) \in 1 + J^\infty \equiv 1 \mod \h \, {\big({R[G]}^\vee
\big)}' \, $  for  $ \, \gamma \in G_\infty \, $,  \, Proposition C
gives
  $$  F\big[\varGamma_G\big]  \; = \;
{\big({R[G]}^\vee\big)}'{\Big|}_{\h=0}
\, = \;  \Big( {\textstyle
\bigoplus_{\hskip-1pt   \substack{ \, b_i \in B_>, \; e_i \in \N_+  \\
                r \in \N, \; b_1 \precneqq \cdots \precneqq b_r }}} \,
\Bbbk \cdot x_{b_1}^{\,e_1} \cdots x_{b_r}^{\,e_r} \Big)
\,{\textstyle \bigoplus}\, \Big( {\textstyle
\bigoplus_{\hskip-3pt   \substack{ \, b_i \in B_1, \; a_i \in \Z  \\
                s \in \N, \; b_1 \precneqq \cdots \precneqq b_s }}}
\, \Bbbk \cdot z_{b_1}^{\,a_1} \cdots z_{b_s}^{\,a_s} \Big)  $$
so  $ F\big[\varGamma_G\big] $  is a polynomial-Laurent polynomial
algebra as claimed.  Similarly  $ \, \Delta(z_g) = z_g \otimes z_g \, $ 
for all  $ \, g \in G \, $  and  $ \, \Delta(x_g) = x_g \otimes 1 +
1 \otimes x_g \, $  iff  $ \, \partial(g) > 1 \, $;  \, so the 
$ z_b $'s  are group-like and the  $ x_b $'s  primitive.   
                                          \par
   {\it (b)} \, The definition of  $ d $  implies  $ \, d(g\,\ell\,)
\geq d(g) + d(\ell\,) \, $  ($ g, \ell \in G $),  \, whence we get
$ \; [\psi_g,\psi_\ell] \, = \, \h^{1 - d(g) - d(\ell) + d((g,\ell))}
\, \psi_{(g,\ell)} \, g \, \ell \, \in \, \h \cdot {\big( {R[G]}^\vee
\big)}' \, $,  \, proving that  $ \, {\big( {R[G]}^\vee \big)}'
{\Big|}_{\h=0} \, $  is commutative.  In addition  $ \; d(g^p)
\geq p \; d(g) \, $,  \, so  $ \; \psi_g^{\;p} = \h^{\, p \,
(1 - d(g))} \, \eta_g^{\;p} =
%
%
\h^{\, p - 1 + d(g^p) - p\,d(g)} \,
\psi_{g^p} \in \h \cdot {\big({R[G]}^\vee\big)}' \, $,  \, whence
$ \, {\big( \psi_g^{\;p}{\big|}_{\h=0} \big)}^p = 0 \, $  inside  $ \,
{\big( {R[G]}^\vee \big)}'{\Big|}_{\h=0} \! = F \big[ \varGamma_G \big]
\, $,  \, which proves that  $ \varGamma_G $  has dimension 0 and height
1.  Finally  $ \; b^p = {(1 + \psi_b)}^p = 1 + {\psi_b}^{\!p} \equiv 1
\mod \h \, {\big( {R[G]}^\vee \big)}' \, $  for all  $ \, b \in B_1 \, $,
\, so  $ \, b^{-1} \equiv b^{p-1} \!\! \mod \h \, {\big( {R[G]}^\vee
\big)}' \, $.
   \hbox{Thus letting  $ \, x_g := \psi_g \!\! \mod \h \; {\big(
{R[G]}^\vee \big)}' $  (for  $ \, g \! \not= \! 1 $)  we get}
  $$  F\big[\varGamma_G\big] \; = \; {\big({R[G]}^\vee\big)}'{\Big|}_{\h=0}
\, = \; \Big( {\textstyle \bigoplus_{\substack{
                \, b_i \in B_>, \; 0 < e_i < p  \\
                r \in \N, \; b_1 \precneqq \cdots \precneqq b_r }}}
\; \Bbbk \cdot x_{b_1}^{\,e_1} \cdots x_{b_r}^{\,e_r} \Big) 
\, {\textstyle \bigoplus} \,
\Big( {\textstyle \bigoplus_{\substack{ \, b_i \in B_1, \; 0 < e_i < p  \\
                s \in \N, \; b_1 \precneqq \cdots \precneqq b_s }}}
\; \Bbbk \cdot z_{b_1}^{\,e_1} \cdots z_{b_s}^{\,e_s} \Big)  $$
just like for  {\it (a)\/}  and also taking care that  $ \, z_b = x_b
+ 1 \, $  and  $ \, z_b^{\,p} = 1 \, $  for  $ \, b \in B_1 \, $.
Therefore  $ \, {\big({R[G]}^\vee\big)}'{\Big|}_{\h=0} \, $  is a
truncated polynomial/Laurent polynomial algebra as claimed.  The
properties of the  $ x_b $'s  and the  $ z_b $'s  w.r.t.~the Hopf
     \hbox{structure are then proved like for \!  {\it (a)\/}  again.}
                                                 \par
   {\it (c)} \, The augmentation ideal  $ \, \germ_e \, $  of  $ \,
{\big( {R[G]}^\vee \big)}'{\Big|}_{\h=0} = F\big[\varGamma_G\big] \, $
is generated by  $ \, {\{ x_b \}}_{b \in B} \, $;  \, then  $ \; \h^{-1}
\, [\psi_g,\psi_\ell\big] \, =
%
%
\, \h^{\,
\theta((g,\ell)) - \theta(g) - \theta(\ell)} \, \psi_{(g,\ell)}
\, \big( 1 + \h^{\, \theta(g) - 1} \psi_g \,\big) \, \big( 1 + \h^{\,
\theta(\ell) - 1} \psi_\ell \,\big) \, $
by the
previous computation, whence at  $ \, \h = 0 \, $
one has  $ \; \big\{ x_g \, , x_\ell \big\} \, \equiv \,
x_{(g,\ell)} \mod \, \germ_e^{\,2} \, $  if  $ \, \theta
\big( (g,\ell\,) \big) = \theta(g) + \theta(\ell\,) \, $,  and
$ \; \big\{ x_g \, , x_\ell \big\} \, \equiv \, 0 \mod \, \germ_e^{\,2}
\, $  if  $ \, \theta \big( (g,\ell\,) \big) > \theta(g) + \theta(\ell\,)
\, $.  This means that the cotangent Lie bialgebra  $ \, \germ_e \Big/
\germ_e^{\,2} \, $  of  $ \varGamma_G \, $  is isomorphic to  $ \gerk_G
\, $,  \, as claimed.   \qed
%
%

\vskip4pt

\noindent   {\it  $ \underline{\hbox{\sl Remarks}} $:}  {\it (a)}
\, Theorem D claims that the connected Poisson group  $ \, K_G^\star
:= \varGamma_G \, $  is  {\sl dual\/}  to  $ \gerk_G $  in the sense
of \S 1.1.  Since  $ \; {R[G]}^\vee{\Big|}_{\h=0} \! = \U(\gerk_G)
\; $  and  $ \; {\big( {R[G]}^\vee \big)}'{\Big|}_{\h=0} \! = F \big[
K_G^\star \big] \, $,  \,  {\sl this gives a close analogue, in positive
characteristic, of the second half of Theorem 2.2{\it (c)}.}
                                                \par
  {\it (b)} \, Theorem D provides functorial recipes to attach to
each abstract group  $ G $  and each field  $ \Bbbk $  a connected
Abelian algebraic Poisson group over\/  $ \Bbbk $,  namely  $ \;
G \mapsto \varGamma_G \; $,  \; explicitly described as algebraic
group and such that  $ \, \text{\it coLie}\,(K_G^\star) = \gerk_G
\, $.  {\sl Every such  $ \varGamma_G $  (for given  $ \Bbbk $)
is then an invariant of  $ G \, $,  a new one to the author's
knowledge.  Indeed, it is perfectly equivalent to the well-known
invariant\/  $ \gerk_G $}  (over the same  $ \Bbbk $),  because
clearly  $ \, G_1 \! \cong G_2 \, $  implies  $ \, \gerk_{G_1} \!
\cong \gerk_{G_2} \, $,  \, whereas  $ \, \gerk_{G_1} \! \cong
\gerk_{G_2} \, $  implies that  $ G_1 $  and  $ G_2 $  are
isomorphic as algebraic groups   --- by  Theorem D{\it
(a--b)}  ---   and bear isomorphic Poisson structures
--- by  part  {\it (c)}  of Theorem D  ---   whence
$ \, G_1 \! \cong G_2 \, $  as algebraic Poisson groups.

\vskip5pt

\noindent   {\it  $ \underline{\hbox{\sl The case of  $ A_\Bbbk(G)
\, $}} $.} \, Let's now dwell upon  $ \, H := A_\Bbbk(G) \, $,  \,
for a  {\sl finite\/}  group  $ G \, $.
                                          \par
   Let  $ \A $  be a commutative unital ring, and  $ \Bbbk $,
$ \, R := \Bbbk[\h] \, $  be as before.  Since  $ \, A_\A(G) :=
{\A[G]}^* \, $,  \, we have  $ \, \A[G] = {A_\A(G)}^* \, $,  \,
so there is a perfect Hopf pairing  $ \, A_\A(G) \times \A[G]
\longrightarrow \A \, $.
%
%
Our first result is

\vskip5pt

\noindent   {\it  $ \underline{\hbox{\sl Theorem E}} $.  $ \,
{A_R(G)}^\vee \! = R \cdot 1 \oplus R \big[ \h^{-1} \big] \,
J = {\big({A_R(G)}^\vee\big)}' \, $,  $ \, {A_\Bbbk(G)}^\vee
= \Bbbk \!\cdot\! 1 \, $,  $ \, \widehat{A_\Bbbk(G)} =
{A_R(G)}^\vee{\Big|}_{\h=0} \! = \, \Bbbk \cdot 1 = \,
\U(\mathbf{0}) \; $  and  $ \; {\big( {A_R(G)}^\vee
\big)}'{\Big|}_{\h=0} \! = \, \Bbbk \cdot 1 =
F\big[\{*\}\big] \; $.}

%
%
 \vskip4pt
 \noindent   {\it Proof.}
 By construction  $ \, J := \text{\sl Ker}\,(\epsilon_{\!
{}_{A_\Bbbk(G\,)}}) \, $  has  $ \Bbbk $--basis  $ \, \big\{ \varphi_g
\big\}_{g \in G \setminus \{1_{{}_G}\}} \cup \big\{ \varphi_{1_G} \! -
1_{\!{}_{A_\Bbbk(G\,)}} \big\} \, $,  \, and since  $ \, \varphi_g =
{\varphi_g}^{\!2} \, $  for all  $ g $  and  $ \, {(\varphi_{1_G} \!
- \! 1)}^2 = -(\varphi_{1_G} \!-\! 1) \, $  we have  $ \, J = J^\infty
\, $,  \, so  $ \, {A_\Bbbk(G)}^\vee = \Bbbk \!\cdot\! 1 \, $  and
$ \, \widehat{A_\Bbbk(G)} = \Bbbk \!\cdot\! 1 \, $.  Similarly,  $ \,
{A_R(G)}^\vee \, $  is generated by  $ \, \big\{ \h^{-1} \varphi_g
\big\}_{g \in G \setminus \{1_{{}_G}\}} \cup \big\{ \h^{-1} (\varphi_{1_G}
\! - \! 1_{\!{}_{A_R(G\,)}}) \big\} \, $;  \, moreover,  $ \, J = J^\infty
\, $  implies  $ \, \h^n J \subseteq {A_R(G)}^\vee \, $  for all  $ \,
n \, $,  \, whence  $ \, {A_R(G)}^\vee = R \, 1 \oplus R[\h^{-1}]
J \, $.  Then  $ \, J_{{A_R(G)}^\vee} = R\big[\h^{-1}\big] J \subseteq
\h \, A_R(G) \, $,  \, which implies  $ \, {\big({A_R(G)}^\vee\big)}'
= {A_R(G)}^\vee \, $:  \, in particular,  $ \, {\big( {A_R(G)}^\vee
\big)}'{\Big|}_{\h=0} = {A_R(G)}^\vee{\Big|}_{\h=0} \! = \Bbbk \cdot
1 \, $,  \, as claimed.   \qed
%
%

\vskip5pt

\noindent   {\it  $ \underline{\hbox{\sl Poisson groups from
$ A_\Bbbk(G) \, $}} $.} \, Now we look at  $ {A_R(G)}' $,
$ {A_\Bbbk(G)}' $  and  $ \widetilde{A_\Bbbk(G)} \, $.  By
construction  $ A_R(G) $  and  $ R[G] $  are in perfect Hopf pairing,
and are free  $ R $--modules  of  {\sl finite rank}.  In this case,
using a general result about the relation between Drinfeld's functors
and Hopf pairings (namely, Proposition 4.4 in [Ga5]) one finds  $ \,
{A_R(G)}' = {\big({R[G]}^\vee\big)}^\bullet = {\big( {R[G]}^\vee \big)}^*
\, $:  thus  $ \, {A_R(G)}' \, $  is the  {\sl dual\/}  Hopf algebra to 
$ {R[G]}^\vee $.  Then from Proposition B we can argue an explicit
description of  $ {A_R(G)}' $,  whence also of  $ \big({A_R(G)}'\big)^{\!
\vee} $.  Now, in proving  Theorem 3.9{\it (g)\/}  one also shows that 
$ \, {A_\Bbbk(G)}' = \big( J_{\Bbbk[G]}^{\;\infty} \big)^\perp \, $; 
therefore there is a perfect filtered Hopf pairing  $ \, {\Bbbk[G]}^\vee
\times {A_\Bbbk(G)}' \longrightarrow \Bbbk \, $  and a perfect graded
Hopf pairing  $ \, \widetilde{A_\Bbbk(G)} \times \widehat{\Bbbk[G]} \!
\longrightarrow \Bbbk \, $.  Thus  $ \, {A_\Bbbk(G)}' \! \cong \! {\big(
{\Bbbk[G]}^\vee \big)}^* \, $  as filtered Hopf algebras and  $ \,
\widetilde{A_\Bbbk(G)} \cong {\big( \widehat{\Bbbk[G]}\big)}^* \, $  as
graded Hopf algebras.  If  $ p = 0 \, $  then  $ J = J^\infty $,  \, as
each  $ \, g \in G \, $  has finite order and  $ \, g^n = 1 \, $  implies 
$ \, g \in G_\infty \, $:  \, then  $ \, {\Bbbk[G]}^\vee \! = \Bbbk \cdot
1 = \widehat{\Bbbk[G]} \, $,  \, so  $ {A_\Bbbk(G)}' = \Bbbk \cdot 1 =
\widetilde{A_\Bbbk(G)} \, $.  If  $ \, p > 0 \, $  instead, this analysis
gives  $ \, \widetilde{A_\Bbbk(G)} = {\big( \widehat{\Bbbk[G]} \big)}^*
= {\big( \u(\gerk_G) \big)}^* = F[K_G] \, $,  \, where  $ \, K_G  \, $ 
is a connected Poisson group of dimension 0, height 1 and tangent Lie
bialgebra  $ \gerk_G \, $.  Thus   

\vskip7pt

\noindent   {\it  $ \underline{\hbox{\sl Theorem F}} $.
                                                \par
   (a) \, There is a second functorial recipe to attach to each
finite abstract group a connected algebraic Poisson group of
dimension zero and height 1 over any field\/  $ \Bbbk $  with
$ \, \text{\it Char}\,(\Bbbk) > 0 \, $,  namely  $ \; G \mapsto
K_G := \text{\it Spec}\,\Big( \widetilde{A_\Bbbk(G)} \Big) $.
This  $ K_G $  is Poisson dual to  $ \varGamma_G $  of Theorem
D in the sense of \S 1.1, in that  $ \; \text{\it Lie}\,(K_G)
= \, \gerk_G = \text{\it coLie}\,(\varGamma_G) \; $.
                                                \par
   (b) \, If  $ \, p:= \text{\it Char}\,(\Bbbk) > 0 \, $,  then
$ \; \big({A_R(G)}'\big)^{\!\vee}\Big|_{\h=0} \! = \, \u \big(
\gerk_G^{\,\times} \big) = S\big(\gerk_G^{\,\times}\big) \Big/ \big(
\big\{ x^p \,\big|\, x \in \gerk_G^{\,\times} \big\} \big) \; $.}

%
%
 \vskip4pt
 \noindent   {\it Proof.}
 Claim  {\it (a)\/}  is the outcome of the discussion
above.  Part  {\it (b)\/}  instead requires an explicit description
of  $ \, \big({A_R(G)}'\big)^{\!\vee} $.  Since  $ \, {A_R(G)}'
\cong \big({R[G]}^\vee\big)^* \, $,  \, from Proposition B we get
$ \; {A_R(G)}' = \, \Big(\! \bigoplus_{\substack{ \, b_i \in B,
\; 0 < e_i < p  \\
                 r \in \N, \; b_1 \precneqq \cdots \precneqq b_r }}
\, R \cdot \rho_{b_1, \dots, b_r}^{e_1, \dots, e_r} \Big) \; $
where each  $ \, \rho_{b_1, \dots, b_r}^{e_1, \dots, e_r} \, $  is
defined by
                      \hfill\break
 \vskip-5pt
   \centerline{$ \Big\langle \, \rho_{b_1, \dots, b_r}^{e_1, \dots,
e_r} \; , \; \chi_{\beta_1}^{\;\varepsilon_1} \, \beta_1^{\,
-[\varepsilon_1\!/2]} \cdots \chi_{\beta_s}^{\;\varepsilon_s} \,
\beta_s^{\,-[\varepsilon_s/2]} \,\Big\rangle \; = \; \delta_{r,s}
\, \prod_{i=1}^r \delta_{b_i,\beta_i} \delta_{e_i, \varepsilon_i} $}
 \vskip9pt
\noindent   (for all  $ \, b_i, \beta_j \in B \, $  and  $ \, 0 <
e_i, \varepsilon_j < p \, $).  Now, using notation of \S 1.3,  $ \,
K_\infty \subseteq K' \, $  for any  $ \, K \in \HA \, $,  \, whence
$ \, K' = \pi^{-1} \big(\, \overline{K}^{\;\prime} \,\big) \, $  where
$ \; \pi \, \colon \, K \relbar\joinrel\twoheadrightarrow K \big/
K_\infty =: \overline{K} \; $  is the canonical projection.  So let
$ \, K := {R[G]}^\vee \, $,  $ \, H := {A_R(G)}' \, $;  \, Proposition
B gives  $ \, K_\infty = R\big[\h^{-1}\big] \cdot J^\infty \, $  and
provides at once a description of  $ \, \overline{K} \, $;  \, from this
and the previous description of  $ H $  one sees also that in the present
case  $ K_\infty $  is exactly the right kernel of the natural pairing
$ \, H \times K \longrightarrow R \, $,  \, which is perfect on the left,
so that the induced pairing  $ \, H \times \overline{K} \longrightarrow
R \, $  is perfect.  By construction its specialization at  $ \, \h
= 0 \, $  is the natural pairing  $ \, F[K_G] \times \u(\gerk_G)
\longrightarrow \Bbbk \, $,  \, which is perfect too.  Then one applies
Proposition 4.4{\it (c)\/}  of [Ga5] (with  $ \overline{K} $  playing
the r\^{o}le of  $ K $  therein), which yields  $ \, \overline{K}^{\;
\prime} \! = \big(H^\vee\big)^\bullet = \Big(\!\! \big({A_R(G)}'\big)^{\!
\vee} \Big)^{\!\bullet} \, $.  By construction  $ \, \overline{K}^{\;
\prime} \! = \big({R[G]}^\vee\big)' \! \Big/ \! \big( R\big[\h^{-1}\big]
\cdot J^\infty \big) \, $,  \, and Proposition C describes the latter as
$ \; \overline{K}^{\;\prime} = \, \Big(\!
          \bigoplus_{\substack{ \, b_i \in B, \; 0 < e_i < p  \\
                        r \in \N, \; b_1 \precneqq \cdots \precneqq b_r }}
R \cdot \overline{\psi}_{b_1}^{\;e_1} \cdots \overline{\psi}_{b_r}^{\;e_r}
\Big) \, $,   
\; where  $ \, \overline{\psi}_{b_i} := \psi_{b_i}
\mod R\big[\h^{-1}\big] \cdot J^\infty \, $  for all  $ i \, $;  \,
since  $ \, \overline{K}^{\;\prime} \! = \Big(\!\! \big( {A_R(G)}'
\big)^{\!\vee} \Big)^{\!\bullet}  \, $  and  $ \, \psi_g = \h^{+1}
\chi_g \, $,  \, this yields
$ \; \big({A_R(G)}'\big)^{\!\vee} = \, \Big(\!
             \bigoplus_{\substack{ \, b_i \in B, \; 0 < e_i < p  \\
                r \in \N, \; b_1 \precneqq \cdots \precneqq b_r }}
R \cdot \h^{- \sum_i e_i d(b_i)} \rho_{b_1, \dots, b_r}^{e_1,
\dots, e_r} \Big) \cong \big(\overline{K}^{\;\prime}\big)^* \, $,
\; whence we get  $ \; \big({A_R(G)}'\big)^{\!\vee}\Big|_{\h=0} \! \cong
\big(\overline{K}^{\;\prime}\big)^*\Big|_{\h=0} \! = \big(K'\big|_{\h=0}
\big)^* \! = \Big(\! \big({R[G]}^\vee\big)'\big|_{\h=0}\Big)^* \! \cong
{F\big[\varGamma_G\big]}^* = \u \big( \gerk_G^{\,\times} \big) =
S\big(\gerk_G^{\,\times}\big) \Big/ \big( \big\{ x^p \,\big|\, x \in
\gerk_G^{\,\times} \big\} \big) \; $  as claimed, the latter identity
being trivial (for  $ \, \gerk_G^{\,\times} \, $  is Abelian).   \qed
%
%

\vskip11pt

\noindent   {\it  $ \underline{\hbox{\sl Remarks}} $:}  {\it (a)} \,
This  $ K_G $  is another invariant for  $ G $,  but again equivalent
to  $ \gerk_G \, $.
                                                \par
   {\it (b)} \, {\sl Theorem F\,{\it (b)\/}  is a positive
characteristic analogue for  $ \, F_\h[G] = {A_R(G)}' \, $
of the first half of  Theorem 2.2{\it (c)}.}

\vskip11pt

\noindent   {\it  $ \underline{\hbox{\sl Examples}} $:}
 \vskip5pt
   {\sl (1) \, Finite Abelian  $ p \, $--groups.} \, Let  $ p \, $
be a prime number and  $ \, G := \Z_{p^{e_1}} \times \Z_{p^{e_2}}
\times \cdots \times \Z_{p^{e_k}} \, $  ($ k, e_1, \ldots, e_k \in
\N \, $),  \, with  $ \, e_1 \geq e_2 \geq \cdots \geq e_k \, $.
Let  $ \Bbbk $  be a field with  $ \, \Char(\Bbbk) = p > 0 \, $,
\, and  $ \, R := \Bbbk[\h] \, $  as above, so that  $ \,
\Bbbk[G]_\h = R[G] \, $.
                                              \par
   First,  $ \, \gerk_G \, $  is Abelian, because  $ G $  is.
Let  $ g_i $  be a generator of  $ \, \Z_{p^{e_i}} \, $  (for
all  $ i \, $),  identified with its image in  $ G \, $.  Since
$ G $  is Abelian we have  $ \, G_{[n]} = G^{p^n} \, $  (for all
$ n \, $),  and an ordered  $ p $-l.c.s.-net  is  $ \, B :=
\bigcup_{r \in \N_+} B_r \, $  with  $ \, B_r := \Big\{\, g_1^{\,p^r},
\, g_2^{\,p^r}, \, \ldots , \, g_{j_r}^{\,p^r} \Big\} \, $  where
$ j_r $  is uniquely defined by  $ \; e_{j_r} > r \, $,  $ \; e_{j_r
+ 1} \leq r \, $.  Then  $ \, \gerk_G \, $  has  $ \Bbbk $--basis
$ \, {\big\{\, \overline{\,\eta_{g_i^{p^{s_i}}}} \,\big\}}_{1 \leq i
\leq k; \; 0 \leq s_i < e_i} \, $,  \, and minimal set of generators
(as a restricted Lie algebra)  $ \, \big\{\,\overline{\,\eta_{g_1}}
\, , \, \overline{\,\eta_{g_2}} \, , \, \ldots, \, \overline{\,
\eta_{g_k}} \,\big\} \, $,  \, for the  $ p $--operation  of
$ \gerk_G $  is  $ \, {\big(\overline{\,\eta_{g_i^{p^s}}}
\,\big)}^{[p\hskip0,5pt]} = \overline{\,\eta_{g_i^{p^{s+1}}}} \, $,
\, and the order of nilpotency of each  $ \, \overline{\,\eta_{g_i}}
\, $  is exactly  $ p^{e_i} $,  i.e.~the order of  $ g_i \, $.  In
addition  $ \, J^\infty = \{0\} \, $  so  $ \, {\Bbbk[G]}^\vee \!
= \Bbbk[G] \, $.  The outcome is  $ \; {\Bbbk[G]}^\vee \! = \,
\Bbbk[G] \; $  and
 \vskip3pt
   \centerline{ $  \widehat{\Bbbk[G]}  \, = \,  \u(\gerk_G)  \; = \;
U(\gerk_G) \bigg/ \Big( {\Big\{ {\big(\overline{\,\eta_{g_i^{p^s}}}
\,\big)}^p - \overline{\,\eta_{g_i^{p^{s+1}}}} \,\Big\}}_{1 \leq i
\leq k}^{0 \leq s < e_i} \bigcup \; {\Big\{ {\big(\overline{\,
\eta_{g_i^{p^{e_i-1}}}}\,\big)}^p \,\Big\}}_{1 \leq i \leq k}
\, \Big) $ }
 \vskip3pt
\noindent
whence  $ \; \widehat{\Bbbk[G]} \, \cong \, \Bbbk[x_1,\dots,x_k] \bigg/
\! \Big( \Big\{\, x_i^{p^{e_i}} \;\Big|\; 1 \leq i \leq k \,\Big\} \Big)
\, $,  \; via  $ \; \overline{\,\eta_{g_i^{p^s}}} \mapsto x_i^{\,p^s}
\, $  (for all  $ i $,  $ s \, $).
                                              \par
   As for  $ {\Bbbk[G]}_\h^{\,\vee} $,  for all  $ \, r < e_i \, $  we
have  $ \, d\big(g_i^{p^r}\big) = p^r \, $  and so  $ \, \chi_{g_i^{p^r}}
= \h^{-p^r} \big( g_i^{p^r} \!-\! 1 \big) \, $  and  $ \, \psi_{g_i^{p^r}}
= \h^{1-p^r} \big( g_i^{p^r} \!-\! 1 \big) \, $;  \, since  $ \,
G_{[\infty]} = \{1\} \, $  (or, equivalently,  $ \, J^\infty = \{0\}
\, $)  and everything is Abelian, from the general theory we conlude
that both  $ {\Bbbk[G]}_\h^{\,\vee} $  and  $ \Big( {\Bbbk[G]}_\h^{\,
\vee} \Big)' $  are truncated-polynomial algebras, in the
$ \chi_{g_i^{p^r}} $'s  and in the  $ \psi_{g_i^{p^r}} $'s
respectively, namely
 \vskip3pt
   \centerline{ $ \begin{array}{rl}
   {\Bbbk[G]_\h}^{\!\vee} \!  &  \! = \; \Bbbk[\h] \Big[ \big\{\,
\chi_{g_i^{p^s}} \big\}_{1 \leq i \leq k \, ; \; 0 \leq s < e_i} \Big]
\; \cong \;  \Bbbk[\h]\big[\,y_1,\dots,y_k\big] \bigg/ \! \Big(
\Big\{\, y_i^{p^{e_i}} \;\Big|\; 1 \leq i \leq k \,\Big\} \Big)  \\
   \Big( {\Bbbk[G]}_\h^{\,\vee} \Big)' \!  &  \! = \, \Bbbk[\h] \Big[
\big\{\, \psi_{g_i^{p^s}} \big\}_{1 \leq i \leq k \, ; \; 0 \leq s < e_i}
\Big] \, \cong \, \Bbbk[\h] \Big[ \big\{\, z_{i,s} \big\}_{1 \leq i \leq
k \, ; \; 0 \leq s < e_i} \Big] \bigg/ \! \Big( \Big\{\, {z_{i,s}}^{\!p}
\;\Big|\; 1 \leq i \leq k \,\Big\} \Big)  \end{array} $ }
 \vskip4pt
\noindent
via the isomorphisms given by  $ \, \overline{\,\chi_{g_i^{p^s}}}
\mapsto y_i^{\,p^s} \, $  and  $ \, \overline{\,\psi_{g_i^{p^s}}}
\mapsto z_{i,s} \, $  (for all  $ i $,  $ s \, $).  When  $ \, e_1
> 1 \, $  this implies  $ \, {\big( {\Bbbk[G]}_\h^{\,\vee\,} \big)}'
\supsetneqq {\Bbbk[G]}_\h \; $,  \, that is a  {\sl counterexample\/} 
to Theorem 2.2{\it (b)}.  Setting  $ \; \overline{\psi_{g_i^{p^s}}} :=
\psi_{g_i^{p^s}} \mod \h \; {\big( {\Bbbk[G]}_\h^{\,\vee\,} \big)}' \; $ 
(for all  $ \, 1 \leq i \leq k \, $,  $ \, 0 \leq s < e_i \, $)  we have
                              \hfill\break   
   \centerline{ $ F\big[K_G^\star\big] \, = \, {\big(
{\Bbbk[G]}_\h^{\,\vee\,} \big)}'{\Big|}_{\h=0} = \,
\Bbbk \Big[ \big\{ \overline{\psi_{g_i^{p^s}}}
\,\big\}_{1 \leq i \leq k}^{0 \leq s < e_i} \Big] \, \cong
\, \Bbbk\,\Big[ \big\{\, w_{i,s} \big\}_{1 \leq i \leq k}^{0
\leq s < e_i} \Big] \bigg/ \! \Big( \Big\{\, w_{i,s}^{\!p}
\Big|\; 1 \!\leq\! i \!\leq\! k \,\Big\} \Big) $ }
(via  $ \, \overline{\psi_{g_i^{p^s}}} \mapsto w_{i,s} \, $)
as a  $ \Bbbk $--algebra.  The Poisson bracket trivial, and
the  $ w_{i,s} $'s  are primitive for  $ \, s > 1 \, $  and
$ \, \Delta(w_{i,1}) = w_{i,1} \otimes 1 + 1 \otimes w_{i,1}
+ w_{i,1} \otimes w_{i,1} \, $  for all  $ 1 \leq i \leq k \, $.
%
%
   If instead  $ \; e_1 = \cdots = e_k = 1 \, $,  \, then  $ \,
{\big( {\Bbbk[G]}_\h^{\,\vee\,} \big)}' = {\Bbbk[G]}_\h \, $;  \,
this is an analogue of  Theorem 2.2{\it (b)},  though now  $ \,
\Char(\Bbbk) > 0 \, $,  \, in that  $ \, {\Bbbk[G]}_\h \, $  is
a QFA, with  $ \, {\Bbbk[G]}_\h{\Big|}_{\h=0} \! = \Bbbk[G] = F
\big[\,\widehat{G}\,\big] \, $  where  $ \, \widehat{G} \, $  is the 
{\sl group of characters\/}  of  $ G \, $.  But then  $ \, F \big[
\widehat{G} \,\big] = \Bbbk[G] = {\Bbbk[G]}_\h{\Big|}_{\h=0} \! =
{\big( {\Bbbk[G]}_\h^{\,\vee\,} \big)}'{\Big|}_{\h=0} \! = F \big[
K_G^\star \big] \, $  (by our general analysis) so 
$ \widehat{G} $  can be realized as a finite, connected, Poisson
group-scheme of dimension 0 
   \hbox{and height 1 dual to  $ \gerk_G \, $, 
namely  $ K_G^\star \, $.}
                                                    \par
   Finally, a direct easy calculation shows that   --- letting  $ \,
\chi^*_g := \h^{\,d(g)} \, (\varphi_g - \varphi_1) \in {A_\Bbbk(G)}'_\h
\, $  and  $ \, \psi^*_g := \h^{\,d(g)-1} \, (\varphi_g - \varphi_1)
\in {\big({A_\Bbbk(G)}'\big)}^\vee_\h \, $  (for all  $ \, g \in G
\setminus \{1\} \, $)  ---   we have also
 \vskip5pt
   \centerline{ $  \begin{array}{rl}
   {A_\Bbbk(G)}_\h^{\,\prime} \!  &  \!\! = \; \Bbbk[\h] \Big[ \big\{\,
\chi^*_{g_i^{p^s}} \big\}_{1 \leq i \leq k}^{0 \leq s < e_i}
\Big] \, \cong \; \Bbbk[\h] \Big[ \big\{ Y_{i,j} \big\}_{1
\leq i \leq k}^{0 \leq s < e_i} \Big] \bigg/ \! \Big( \big\{
Y_{i,j}^{\;p} \big\}_{1 \leq i \leq k}^{0 \leq s < e_i} \Big)  \\
   {\big({A_\Bbbk(G)}_\h^{\,\prime}\,\big)}^\vee \!  &  \!\! = \;
\Bbbk[\h] \Big[ \big\{\, \psi^*_{g_i^{p^s}} \big\}_{1 \leq i
\leq k}^{0 \leq s < e_i} \Big] \, \cong \; \Bbbk[\h] \Big[ \big\{
Z_{i,s} \big\}_{1 \leq i \leq k}^{0 \leq s < e_i} \Big] \bigg/ \!
\Big( \big\{ Z_{i,s}^{\;p} - Z_{i,s} \big\}_{1 \leq i \leq k}^{0 \leq
s < e_i} \Big)  \end{array}$ }
 \vskip5pt
 \noindent
via the isomorphisms given by  $ \, \chi^*_{g_i^{p^s}} \mapsto
Y_{i,s} \, $  and  $ \, \psi^*_{g_i^{p^s}} \mapsto Z_{i,s} \, $,
\, from which one also gets the analogous descriptions of  $ \,
{A_\Bbbk(G)}_\h^{\,\prime}{\Big|}_{\h=0} \! = \widetilde{A_\Bbbk(G)}
= F[K_G] \, $  and of  $ \, {\big({A_\Bbbk(G)}_\h^{\,\prime}\,
\big)}^{\!\vee}{\Big|}_{\h=0} \! = \u(\gerk_G^\times) \; $.
 \vskip7pt
   {\sl (2) \, A non-Abelian  $ p \, $--group.} \, Let  $ p \, $
be a prime number,  $ \Bbbk $  be a field with  $ \, \Char(\Bbbk)
= p > 0 \, $,  \, and  $ \, R := \Bbbk[\h] \, $  as above, so
that  $ \, \Bbbk[G]_\h = R[G] \, $.
                                              \par
   Let  $ \, G := \Z_p \ltimes \Z_{p^{\,2}} \, $,  \, that is the group
with generators  $ \, \nu $,  $ \tau \, $  and relations  $ \, \nu^p
= 1 \, $,  $ \, \tau^{p^2} = 1 \, $,  $ \, \nu \, \tau \, \nu^{-1}
= \tau^{1+p} \, $.  In this case,  $ \, G_{[2]} = \cdots = G_{[p\,]}
= \big\{ 1, \tau^p \,\big\} \, $,  $ \, G_{[p+1]} = \{1\} \, $,  \,
so we can take  $ \, B_1 = \{\nu \, , \tau \,\} \, $  and  $ \, B_p
= \big\{ \tau^p \,\big\} \, $  to form an ordered  $ p $-l.c.s.-net
$ \, B := B_1 \cup B_p \, $  w.r.t.~the ordering  $ \, \nu \preceq
\tau \preceq \tau^p \, $.  Noting also that  $ \, J^\infty = \{0\}
\, $  (for  $ \, G_{[\infty]} = \{1\} \, $),  we have
 \vskip3pt
   \centerline{ $ {{\Bbbk[G]}_\h}^{\!\vee} \; = \; {\textstyle
\bigoplus_{a,b,c=0}^{p-1}} \, \Bbbk[\h] \cdot \chi_\nu^{\,a}
\, \chi_\tau^{\,b} \, \chi_{\tau^p}^{\,c} \; = \; {\textstyle
\bigoplus_{a,b,c=0}^{p-1}} \, \Bbbk[\h] \, \h^{- a - b - c \, p}
\cdot {(\nu-1)}^a \, {(\tau-1)}^b \, {\big( \tau^p - 1 \big)}^c $ }
 \vskip3pt
\noindent
as  $ \Bbbk[\h] $--modules,  since  $ \, d(\nu) = 1 = d(\tau) \, $
and  $ \, d\big(\tau^p\big)) = p \, $,  \, with  $ \, \Delta(\chi_g)
= \chi_g \otimes 1 + 1 \otimes \chi_g + \h^{d(g)} \, \chi_g \otimes
\chi_g \, $  for all  $ \, g \in B \, $.  As a direct consequence
we have also
 \vskip3pt
   \centerline{ $ {\textstyle \bigoplus_{a,b,c=0}^{p-1}} \, \Bbbk
\cdot \overline{\chi_\nu}^{\;a} \, \overline{\chi_\tau}^{\;b} \,
\overline{\chi_{\tau^p}}^{\;c} \; = \; {{\Bbbk[G]}_\h}^{\!\vee}
{\Big|}_{\h=0} \; \cong \; \widehat{\Bbbk[G]} \; = \; {\textstyle
\bigoplus_{a,b,c=0}^{p-1}} \, \Bbbk \cdot \overline{\eta_\nu}^{\;a} \,
\overline{\eta_\tau}^{\;b} \, \overline{\eta_{\tau^p}}^{\;c} \, . $ }
 \vskip3pt
   The two relations  $ \, \nu^p = 1 \, $  and  $ \, \tau^{p^2} =
1 \, $  within  $ G $  yield trivial relations inside  $ \Bbbk[G] $
and  $ {\Bbbk[G]}_\h \, $;  instead, the relation  $ \, \nu \, \tau
\, \nu^{-1} = \tau^{1+p} \, $  turns into  $ \, [\eta_\nu,\eta_\tau]
= \eta_{\tau^p} \cdot \tau \, \nu \, $,  \, which gives  $ \, [\chi_\nu,
\chi_\tau] = \h^{p-2} \, \chi_{\tau^p} \cdot \tau \, \nu \, $  in
$ {{\Bbbk[G]}_\h}^{\!\vee} $.  Therefore  $ \; [\, \overline{\chi_\nu}
\, , \, \overline{\chi_\tau}\,] = \delta_{p,2} \,
\overline{\chi_{\tau^p}} \, $.  Since  $ \, [\,
\overline{\chi_\tau} \, , \, \overline{\chi_{\tau^p}}\,] = 0 =
[\,\overline{\chi_\nu} \, , \, \overline{\chi_{\tau^p}} \,] \, $
(because  $ \, \nu \, \tau^p \, \nu^{-1} = {\big( \tau^{1+p}
\big)}^p = \tau^{p + p^2} = \tau^p \, $)  and  $ \, \{
\overline{\chi_\nu} \, , \, \overline{\chi_\tau} \, , \,
\overline{\chi_{\tau^p}} \,\} \, $  is a  $ \Bbbk $--basis  of
$ \, \gerk_G = \mathcal{L}_p(G) \, $,  \, we conclude that the latter
has trivial or non-trivial Lie bracket according to whether  $ \,
p \not= 2 \, $  or  $ \, p = 2 \, $.  In addition, we have the
relations  $ \, \chi_\nu^{\;p} = 0 \, $,  $ \, \chi_{\tau^p}^{\;p}
= 0 \, $  and  $ \, \chi_\tau^{\;p} = \chi_{\tau^p} \, $:  \, these
give analogous relations in  $ \, {{\Bbbk[G]}_\h}^{\!\vee}{\Big|}_{\h=0}
\, $,  \, which define the  $ p $--operation  of  $ \gerk_G \, $, 
namely  $ \; \overline{\chi_\nu}^{\;[p\,]} = 0 \, $,  $ \;
\overline{\chi_{\tau^p}}^{\;[p\,]} = 0 \, $,  $ \;
\overline{\chi_\tau}^{\;[p\,]} = \chi_{\tau^p} \; $.
                                              \par
   To sum up, we have a complete presentation for  $ {R[G]}^\vee $
by generators and relations, that is   
 \vskip3pt
  \centerline{ $ {{\Bbbk[G]}_\h}^{\!\vee}  \; \cong \;\;
\Bbbk[\h] \big\langle v_1, v_2, v_3 \big\rangle \bigg/ \!
\bigg( \hskip-9pt
%
%
\hbox{ $ \begin{matrix}
 &  v_1 \, v_2 - v_2 \, v_1 - \h^{p-2} \, v_3 \,
(1 + \h \, v_2) \, (1 + \h \, v_1)  \\
 &  v_1 \, v_3 - v_3 \, v_1 \, ,  \quad  v_1^{\,p} \, ,
\quad  v_2^{\,p} - v_3 \, , \quad  v_3^{\,p} \, , \quad
v_2 \, v_3 - v_3 \, v_2
       \end{matrix} $ }  \bigg) $ }
 \vskip3pt
\noindent
via  $ \, \chi_\nu \mapsto v_1 \, $,  $ \, \chi_\tau \mapsto v_2 \, $,
$ \, \chi_{\tau^p} \mapsto v_3 \, $.  Similarly (as a consequence) we
have the presentation
 \vskip3pt
  \centerline{ $ \widehat{\Bbbk[G]}  \; = \; 
{{\Bbbk[G]}_\h}^{\!\vee}{\Big|}_{\h=0} \; \cong \;\; 
\Bbbk \, \big\langle y_1, y_2, y_3 \big\rangle \bigg/ \!
%
%
\bigg( \hbox{ $ \begin{matrix}
     y_1 \, y_2 - y_2 \, y_1 - \delta_{p,2} \, y_3 \, ,
\qquad  y_2^{\,p} - y_3  \\
     y_1 \, y_3 - y_3 \, y_1 \, ,  \quad  y_1^{\,p} \, ,
\quad  y_3^{\,p} \, ,  \quad  y_2 \, y_3 - y_3 \, y_2
       \end{matrix} $ }  \bigg) $ }
 \vskip3pt
\noindent
via  $ \, \overline{\chi_\nu} \mapsto y_1 \, $,  $ \,
\overline{\chi_\tau} \mapsto y_2 \, $,  $ \, \overline{\chi_{\tau^p}}
\mapsto y_3 \, $,  \, with  $ p $--operation  as above and the
$ y_i $'s  being primitive.
                                          \par
   {\it  $ \underline{\text{Remark}} $:}  \, if  $ \, p \not= 2 \, $
exactly the same result holds for  $ \, G = \Z_p \times \Z_{p^2} \, $,
\, i.e.~$ \; \gerk_{\, \Z_p \hskip-0,5pt \ltimes \Z_{p^2}} =
\gerk_{\, \Z_p \hskip-0,5pt \times \Z_{p^2}} \; $:  \; this
shows that the restricted Lie bialgebra  $ \gerk_G $  may be
     \hbox{not enough to recover the group  $ G \, $.}
                                          \par
   As for  $ \, {\big( {{\Bbbk[G]}_\h}^{\!\vee} \big)}' $,  \, it is
generated by  $ \, \psi_\nu = \nu - 1 $,  $ \, \psi_\tau = \tau - 1 $,
$ \, \psi_{\tau^p} = \h^{1-p} \big(\tau^p - 1 \big) \, $,  \, with
relations  $ \; \psi_\nu^{\;p} = 0 \, $,  $ \; \psi_\tau^{\;p}
= \h^{p-1} \psi_{\tau^p} \, $,  $ \; \psi_{\tau^p}^{\;\;p} =
0 \, $,  $ \; \psi_\nu \, \psi_\tau - \psi_\tau \, \psi_\nu =
\h^{\,p-1} \psi_{\tau^p} \, (1 + \psi_\tau) \, (1 + \psi_\nu) \, $,
$ \; \psi_\tau \, \psi_{\tau^p} - \psi_{\tau^p} \, \psi_\tau = 0 \, $,
and  $ \; \psi_\nu \, \psi_{\tau^p} - \psi_{\tau^p} \, \psi_\nu = 0
\, $.  In particular  $ \; {\big( {{\Bbbk[G]}_\h}^{\!\vee} \big)}'
\supsetneqq {\Bbbk[G]}_\h \, $,  \, and

 \vskip3pt
  \centerline{ $ {\big( {{\Bbbk[G]}_\h}^{\!\vee} \big)}'  \,\; \cong
\;\;  \Bbbk[\h] \, \big\langle u_1, u_2, u_3 \big\rangle \bigg/
\! \bigg( \hskip-12pt
%
%
\hbox{ $ \begin{matrix}
   u_1 \, u_3 - u_3 \, u_1 \, ,
\quad  u_2^{\;p} - \h^{1-p} \, u_3  \, ,  \quad
u_2 \, u_3 - u_3 \, u_2  \\
   \quad  u_1^{\;p} \, ,  \quad  u_1 \, u_2 - u_2 \, u_1 -
\h^{\,p-1} \, u_3 \, (1 + u_2) \, (1 + u_1) \, , \quad  u_3^{\;\;p}
       \end{matrix} $ }  \hskip-1pt \bigg) $ }
 \vskip3pt
\noindent
via  $ \, \psi_\nu \mapsto u_1 \, $,  $ \, \psi_\tau \mapsto u_2
\, $,  $ \, \psi_{\tau^p} \mapsto u_3 \, $.  Letting  $ \; z_1 :=
\psi_\nu{\big|}_{\h=0} \! + 1 \, $,  $ \; z_2 := \psi_\tau{\big|}_{\h=0}
\! + 1 \; $  and  $ \; x_3 := \psi_{\tau^p}{\big|}_{\h=0} \; $  this gives
$ \, {\big( {{\Bbbk[G]}_\h}^{\!\vee} \big)}'{\Big|}_{\h=0} \!\! = \Bbbk
\big[z_1,z_2,x_3\big] \Big/ \big( z_1^{\,p} \! - \! 1, z_2^{\,p} \! -
\! 1, x_3^{\,p} \,\big) \, $  as a  $ \Bbbk $--algebra,  with the
$ z_i $'s  group-like,  $ x_3 $  primitive  (cf.~Theorem D\,{\it
(b)\/}),  and Poisson bracket given by  $ \, \big\{ z_1, z_2 \big\}
= \delta_{p,2} \, z_1 \, z_2 \, x_3 \, $,  $ \, \big\{ z_2, x_3 \big\}
= 0 \, $  and  $ \, \big\{z_1, x_3\big\} = 0 \, $.  Thus  $ \, {\big(
{{\Bbbk[G]}_\h}^{\!\vee} \big)}'{\Big|}_{\h=0} \! = F[\varGamma_G]
\, $  with  $ \, \varGamma_G \cong {\boldsymbol\mu}_p \times
{\boldsymbol\mu}_p \times {\boldsymbol\alpha}_p \, $  as algebraic
groups, with Poisson structure such that  $ \,\text{\it coLie}\,
(\varGamma_G) \cong \gerk_G \; $.
                                         \par
   Since  $ \, G_\infty = \{1\} \, $  the general theory ensures that
$ \, {A_\Bbbk(G)}' = A_\Bbbk(G) \, $.  We leave to the interested
reader the task of computing the filtration  $ \underline{D} $  of
$ A_\Bbbk(G) $,  and consequently describe  $ \, {A_R(G)}' \, $,
$ \, \big({A_R(G)}'\big)^{\!\vee} \, $,  $ \, \widetilde{A_\Bbbk(G)}
\, $  and the connected Poisson group  $ \, K_G := \text{\it Spec}\,
\big( \widetilde{A_\Bbbk(G)} \big) \; $.
 \vskip5pt
   {\sl (3) \, An Abelian infinite group.} \,  Let  $ \, G = \Z^n
\, $  (written multiplicatively with generators  $ \, e_1, \dots,
e_n \, $),  then  $ \, \Bbbk[G] = \Bbbk[\Z^n] = \Bbbk \big[ e_1^{\pm 1},
\dots, e_n^{\pm 1} \big] \, $  (the ring of Laurent polynomials).
This is the function algebra of the algebraic group
$ {\Bbb{G}_m}^{\hskip-3pt n} $   --- the  $ n $--dimensional
torus on  $ \Bbbk $  ---   which is exactly the character group
of  $ \Z^n \, $,  \, thus we get back to the function algebra case.

 \vskip1,1truecm

\centerline {\bf \S \; 4 \  First example: the Kostant-Kirillov
structure }

\vskip10pt

  {\bf 4.1 Classical and quantum setting.} \, Let  $ \gerg $  and 
$ \gerg^\star $  be as in \S 3.7, consider  $ \gerg $  as a Lie
bialgebra with trivial Lie cobracket and look at  $ \gerg^\star $
as its dual Poisson group, whose Poisson structure then is exactly
the Kostant-Kirillov one.  Take as ground ring  $ \, R := \Bbbk[\nu]
\, $  (a PID, hence a 1dD): we shall consider the primes  $ \, \h =
\nu \, $  and  $ \, \h = \nu - 1 \, $,  \, and we'll find quantum
groups at either of them for both  $ \gerg $  and  $ \gerg^\star \, $.
                                          \par
   To begin with, we assume  $ \, \hbox{\it Char}\,(\Bbbk) = 0 \, $,  \,
and postpone to \S 4.4 the case  $ \, \hbox{\it Char}\,(\Bbbk) > 0 \, $.
                                  \par
   Let  $ \, \gerg_\nu := \gerg[\nu] = \Bbbk[\nu] \otimes_\Bbbk
\gerg \, $,  \, endow it with the unique  $ \Bbbk[\nu] $--linear
Lie bracket  $ \, {[\ ,\ ]}_\nu \, $  given by  $ \, {[x,y]}_\nu
:= \nu \, [x,y] \, $  for all  $ \, x $,  $ y \in \gerg \, $,
\, and define
  $$  H := U_{\Bbbk[\nu]}(\gerg_\nu) = T_{\Bbbk[\nu]}(\gerg_\nu)
\Big/ \big(\big\{\, x \cdot y - y \cdot x - \nu \, [x,y] \;\big|\;
x, y \in \gerg \,\big\} \big)  $$  
the universal enveloping algebra of the Lie  $ \Bbbk[\nu] $--algebra
$ \, \gerg_\nu \, $,  \, endowed with its natural structure of Hopf
algebra.  Then  $ H $  is a free  $ \Bbbk[\nu] $--algebra,
so that  $ \, H \in \HA \, $  and  $ \, H_F := \Bbbk(\nu)
\otimes_{\Bbbk[\nu]} H \in \HA_F \, $  (see \S 1.3); its
specializations at  $ \, \nu = 1 \, $  and at  $ \, \nu =
0 \, $  are  $ \; H \Big/ (\nu\!-\!1) \, H \, = \, U(\gerg)
\; $,  \; as a  {\sl co-Poisson}  Hopf algebra, and  $ \;
H \Big/ \nu \, H \, = \, S(\gerg) \, = \, F[\gerg^\star] \; $, 
\; as a  {\sl Poisson}  Hopf algebra.  In a more suggesting way,
we can also express this with notation like  $ \; H \,{\buildrel
{\, \nu \rightarrow 1 \,} \over \llongrightarrow}\, U(\gerg) \, $, 
$ \; H \,{\buildrel {\, \nu \rightarrow 0 \,} \over \llongrightarrow}\,
F[\gerg^\star] \, $.  So  {\sl  $ H $  is a QrUEA at  $ \, \h := (\nu
\! - \! 1) \, $  and a QFA at  $ \, \h := \nu \, $};  so we'll 
   \hbox{consider Drinfeld's functors for  $ H $  at  $ (\nu\!-\!1) $ 
and at  $ (\nu) \, $.}    

\vskip7pt

  {\bf 4.2 Drinfeld's functors at  $ (\nu) $.} \, 
Let  $ \;
{(\ )}^{\vee_{\!(\nu)}} \, \colon \, \HA \longrightarrow \HA \; $
and  $ \; {(\ )}^{\prime_{(\nu)}} \, \colon \, \HA \longrightarrow
\HA \; $  be the Drinfeld's functors at  $ \, (\nu) \, \big( \in
\text{\it Spec}\big(\Bbbk[\nu]\big) \, \big) \, $.  By definitions
$ \, J := \text{\sl Ker} \, \big( \epsilon \, \colon \, H
\longrightarrow \Bbbk[\nu] \big) \, $  is nothing but the 2-sided
ideal of  $ \, H := U(\gerg_\nu) \, $  generated by  $ \gerg_\nu $
itself; thus  $ \, H^{\vee_{\!(\nu)}} $,  which by definition is
the unital  $ \Bbbk[\nu] $--subalgebra  of  $ H_F $  generated by
$ \, J^{\vee_{\!(\nu)}} := \nu^{-1} J \, $,  \, is just the unital
$ \Bbbk[\nu] $--subalgebra  of  $ H_F $  generated by  $ \;
{\gerg_\nu}^{\!\vee_{\!(\nu)}} := \nu^{-1} \, \gerg_\nu \, $.
Now consider the  $ \Bbbk[\nu] $--module  isomorphism  $ \;
{(\ )}^{\!\vee_{\!(\nu)}} \, \colon \, \gerg_\nu \,{\buildrel
\cong \over \longrightarrow}\, {\gerg_\nu}^{\!\vee_{\!(\nu)}}
:= \nu^{-1} \, \gerg_\nu \; $  given by  $ \, z \mapsto z^\vee
:= \nu^{-1} z \in {\gerg_\nu}^{\!\vee_{\!(\nu)}} \, $  for all
$ z \in \gerg_\nu \, $;  consider on  $ \, \gerg_\nu := \Bbbk[\nu]
\otimes_\Bbbk \gerg \, $  the natural Lie algebra structure (with
trivial Lie cobracket), given by scalar extension from  $ \gerg
\, $,  and push it over  $ {\gerg_\nu}^{\!\vee_{\!(\nu)}} $  via
$ {(\ )}^{\!\vee_{\!(\nu)}} \, $,  so that  $ {\gerg_\nu}^{\!
\vee_{\!(\nu)}} $  is isomorphic to  $ \, \gerg_\nu^{\text{\it
nat}} \, $  (i.e.~$ \gerg_\nu $  carrying the natural Lie bialgebra
structure) as a Lie bialgebra.  Consider  $ \, x^\vee $,  $ y^\vee
\in {\gerg_\nu}^{\!\vee_{\!(\nu)}} $  (with  $ \, x $,  $ y \in
\gerg_\nu \, $):  \, then  $ \; H^{\vee_{\!(\nu)}} \ni \big(
x^\vee \, y^\vee - y^\vee \, x^\vee \big) = \nu^{-2} \big( x
\, y - y \, x \big) = \nu^{-2} \, {[x,y]}_\nu = \nu^{-2} \, \nu
\, [x,y] = \nu^{-1} \, [x,y] = {[x,y]}^\vee =: \big[ x^\vee,
y^\vee \big] \in {\gerg_\nu}^{\!\vee_{\!(\nu)}} \, $.  Therefore
we can conclude at once that  $ \; H^{\vee_{\!(\nu)}} = U
\big( {\gerg_\nu}^{\!\vee_{\!(\nu)}} \big) \cong U \big(
\gerg_\nu^{\text{\it nat}} \big) \; $.
                                        \par
   As a first consequence,  $ \, {\big( H^{\vee_{\!(\nu)}}
\big)}{\Big|}_{\nu=0} \cong U \big( \gerg_\nu^{\text{\it nat}}
\big) \Big/ \nu \, U \big( \gerg_\nu^{\text{\it nat}} \big)
= U \Big( \gerg_\nu^{\text{\it nat}} \big/ \nu \,
\gerg_\nu^{\text{\it nat}} \Big) = U(\gerg) \, $,  \; that is
$ \; H^{\vee_{\!(\nu)}} \,{\buildrel {\, \nu \rightarrow 0 \,}
\over \llongrightarrow}\, U(\gerg) \, $,  \, thus agreeing with
the second half of  Theorem 2.2{\it (c)}.
                                        \par
   Second, look at  $ \, {\big( H^{\vee_{\!(\nu)}}
\big)}^{\prime_{(\nu)}} $.  Since  $ \, H^{\vee_{\!(\nu)}}
= U \big( {\gerg_\nu}^{\!\vee_{\!(\nu)}} \big) \, $,  \, and
$ \, \delta_n(\eta) = 0 \, $  for all  $ \, \eta \in U \big(
{\gerg_\nu}^{\!\vee_{\!(\nu)}} \big) \, $  such that  $ \,
\partial(\eta) < n \, $  (cf.~Lemma 4.2{\it (d)\/}  in [Ga5]),
it is easy to see that   
  $$  {\big( H^{\vee_{\!(\nu)}} \big)}^{\prime_{(\nu)}} \, = \,
\big\langle \nu \, {\gerg_\nu}^{\!\vee_{\!(\nu)}} \big\rangle
\, = \, \big\langle \nu \, \nu^{-1} \gerg_\nu \big\rangle \,
= \, U(\gerg_\nu) \, = \, H  $$   
(hereafter  $ \langle S \, \rangle $  is the subalgebra
generated by  $ S \, $),  so  $ \, {\big( H^{\vee_{\!(\nu)}}
\big)}^{\prime_{(\nu)}} = H \, $,  \, which agrees with
Theorem 2.2{\it (b)}.
%
%
   Finally, proceeding as in \S 3.7 we see that  $ \,
H^{\prime_{(\nu)}} = U(\nu\,\gerg_\nu) \, $,  \, whence  $ \,
{\big( H^{\prime_{(\nu)}} \big)}{\Big|}_{\nu=0} \!\! = {\big( U
(\nu \, \gerg_\nu) \big)}{\Big|}_{\nu=0} \!\! \cong S({\gerg}_{ab})
= F \big[ \gerg^\star_{\delta-{\text{\it ab}}} \big] \, $  where
$ \, {\gerg}_{ab} \, $,  resp.~$ \, \gerg^\star_{\delta-{\text{\it
ab}}} \, $,  is simply  $ \gerg \, $,  resp.~$ \gerg^\star \, $, 
endowed with the trivial Lie bracket, resp.~cobracket,  so that 
$ \, {\big( H^{\prime_{(\nu)}} \big)} {\Big|}_{\nu=0} \! \cong
S({\gerg}_{ab}) = F \big[ \gerg^\star_{\delta-{\text{\it ab}}}
\big] \, $  has trivial Poisson bracket.  Iterating this procedure
one finds that all further images  $ \, \Big( \cdots {\big(
{(H)}^{\prime_{(\nu)}} \big)}^{\prime_{(\nu)}} \cdots
\Big)^{\prime_{(\nu)}} \, $  of the functor 
$ {(\ )}^{\prime_{(\nu)}} $  applied many times
to  $ H $  are pairwise isomorphic; thus in particular
they all have the same specialization at  $ (\nu) $,  namely 
$ \; {\bigg( \! {\Big( \cdots {\big( {(H)}^{\prime_{(\nu)}}
\big)}^{\prime_{(\nu)}} \cdots \Big)}^{\prime_{(\nu)}}
\bigg)}{\bigg|}_{\nu=0} \cong \, S({\gerg}_{ab}) =
F \big[ \gerg^\star_{\delta-{\text{\it ab}}} \big] \; $.  

\vskip7pt

  {\bf 4.3 Drinfeld's functors at  $ (\nu - 1) $.} \, Now
we consider the non-zero prime  $ \, (\nu - 1) \, \big( \! \in
\text{\it Spec}\big(\Bbbk[\nu]\big) \, \big) \, $;  \, let  $ \;
{(\ )}^{\vee_{\!(\nu - 1)}} \, \colon \, \HA \longrightarrow
\HA \; $  and  $ \; {(\ )}^{\prime_{(\nu - 1)}} \, \colon \,
\HA \longrightarrow \HA \; $  be the corresponding Drinfeld's
functors.  Set  $ \, {\gerg_\nu}^{\prime_{(\nu - 1)}} :=
(\nu - 1) \, \gerg_\nu \, $,  \, let  $ \; \, \colon \,
{\gerg_\nu} \,{\buildrel \cong \over \longrightarrow}\,
{\gerg_\nu}^{\!\prime_{(\nu - 1)}} := (\nu - 1)
\, \gerg_\nu \; $  be the  $ \Bbbk[\nu] $--module  isomorphism
given by  $ \, z \mapsto z' := (\nu - 1) \, z \in
{\gerg_\nu}^{\!\prime_{(\nu - 1)}} \, $  for all
$ \, z \in \gerg_\nu \, $,  \, and push over via it the
Lie bialgebra structure of  $ \gerg_\nu $  to an isomorphic
Lie bialgebra structure on  $ \, {\gerg_\nu}^{\!\prime_{(\nu - 1)}}
\, $,  \, whose Lie bracket will be denoted by  $ \,
{[\ ,\ ]}_\ast \, $.  Notice then that we have Lie bialgebra
isomorphisms  $ \; \gerg \, \cong \, \gerg_\nu \big/ (\nu - 1)
\, \gerg_\nu \, \cong \, {\gerg_\nu}^{\!\prime_{(\nu - 1)}}
\big/ (\nu - 1) \, {\gerg_\nu}^{\!\prime_{(\nu - 1)}} \, $.
                                              \par
   Since  $ \, H := U(\gerg_\nu) \, $  it is easy to see by direct
computation that
  $$  H^{\prime_{(\nu - 1)}} \, = \, \big\langle (\nu - 1) \,
\gerg_\nu \big\rangle \, = \, U \big( {\gerg_\nu}^{\!\prime_{(\nu
\! - \! 1)}} \big)   \eqno (4.1)  $$   
where  $ {\gerg_\nu}^{\!\prime_{(\nu - 1)}} $  is seen as a Lie 
$ \Bbbk[\nu] $--subalgebra  of  $ \, \gerg_\nu \, $.  Now, if 
$ \, x' $,  $ y' \in {\gerg_\nu}^{\prime_{(\nu - 1)}} \, $
(with  $ \, x $,  $ y \in \gerg_\nu \, $),  then   
  $$  x' \, y' - y' \, x' = {(\nu - 1)}^2 \big( x \, y -
y \, x \big) = {(\nu - 1)}^2 \, {[x,y]}_\nu = (\nu - 1) \,
{{[x,y]}_\nu}^{\!\prime} = (\nu - 1) \, {\big[x',y'\big]}_\ast
\; .   \eqno (4.2)  $$   
   \indent   This and (4.1) show at once that  $ \, {\big(
H^{\prime_{(\nu - 1)}} \big)}{\Big|}_{(\nu - 1)=0} \! = S
\Big( {\gerg_\nu}^{\prime_{(\nu - 1)}} \big/ (\nu - 1) \,
{\gerg_\nu}^{\prime_{(\nu - 1)}} \Big) \, $  as Hopf algebras,
and also as  {\sl Poisson}  algebras: indeed, the latter holds
because the Poisson bracket  $ \, \{\ ,\ \} \, $  of  $ \, S
\Big( {\gerg_\nu}^{\prime_{(\nu - 1)}} \big/ (\nu - 1) \,
{\gerg_\nu}^{\prime_{(\nu - 1)}} \Big) \, $  inherited from
$ H^{\prime_{(\nu - 1)}} $  (by specialization) is uniquely
determined by its restriction to  $ \, {\gerg_\nu}^{\prime_{(\nu
\! - \! 1)}} \big/ (\nu - 1) \, {\gerg_\nu}^{\prime_{(\nu \! -
\! 1)}} \, $,  \, and on the latter space we have  $ \; \{\ ,\ \}
= {[\ ,\ ]}_\ast \, \mod (\nu - 1) \, {\gerg_\nu}^{\prime_{(\nu
\! - \! 1)}} \; $  (by (4.2)).  Finally, since  $ \,
{\gerg_\nu}^{\prime_{(\nu - 1)}} \big/ (\nu - 1) \,
{\gerg_\nu}^{\prime_{(\nu - 1)}} \cong \gerg \, $  as Lie
algebras we have  $ \, {\big( H^{\prime_{(\nu - 1)}} \big)}
{\Big|}_{(\nu - 1)=0} \! = S(\gerg) = F[\gerg^\star] \, $  as
Poisson Hopf algebras, or, in short,  $ \; H^{\prime_{(\nu
\! - \! 1)}} \,{\buildrel {\, \nu \rightarrow 1 \,} \over
\llongrightarrow}\, F[\gerg^\star] \, $,  \, as prescribed by
the ``first half\/'' of  Theorem 2.2{\it (c)}.
                                        \par
   Second, look at  $ \, {\big( H^{\prime_{(\nu - 1)}} \big)}^{\vee_{\!
(\nu - 1)}} $.  Since  $ \, H^{\prime_{(\nu - 1)}} = U \big(
{\gerg_\nu}^{\!\prime_{(\nu - 1)}} \big) \, $,  \, the kernel 
$ \, \text{\sl Ker} \, \big( \epsilon \, \colon \, H^{\prime_{(\nu - 1)}}
\longrightarrow \Bbbk[\nu] \big) =: J^{\prime_{(\nu - 1)}} \, $  is
just the 2-sided ideal of  $ \, H^{\prime_{(\nu - 1)}} = U \big(
{\gerg_\nu}^{\!\prime_{(\nu - 1)}} \big) \, $  generated by
$ {\gerg_\nu}^{\!\prime_{(\nu - 1)}} \, $.  Therefore  $ \, {\big(
H^{\prime_{(\nu - 1)}} \big)}^{\vee_{\!(\nu - 1)}} $,  \, generated
by  $ \, {\big( J^{\prime_{(\nu - 1)}} \big)}^{\vee_{\!(\nu - 1)}}
:= {(\nu - 1)}^{-1} J^{\prime_{(\nu - 1)}} \, $ as a unital
$ \Bbbk[\nu] $--subalgebra  of  $ \, {\big(
H^{\prime_{(\nu - 1)}} \big)}_F = H_F \, $,  \,
is just the unital  $ \Bbbk[\nu] $--subalgebra
of  $ H_F $  generated by  $ \; {(\nu - 1)}^{-1}
{\gerg_\nu}^{\!\prime_{(\nu - 1)}} = {(\nu - 1)}^{-1}
(\nu - 1) \, \gerg_\nu = \gerg_\nu \, $,  \, that is  $ \,
{\big( H^{\prime_{(\nu - 1)}} \big)}^{\vee_{\!(\nu - 1)}} =
U(\gerg_\nu) = H \, $,  \, confirming  Theorem 2.2{\it (b)}.
                                        \par
   Finally, for  $ H^{\vee_{(\nu - 1)}} $  one has essentially
the same feature as in \S 3.7, and the analysis therein can be
repeated; the final result then will depend on the nature of 
$ \gerg \, $,  in particular on its lower central series.

\vskip7pt

  {\bf 4.4 The case of positive characteristic.} \, Let us
consider now a field  $ \Bbbk $  such that  $ \, \hbox{\it
Char}\,(\Bbbk) = p > 0 \, $.  Starting from  $ \gerg $  and
$ \, R := \Bbbk[\nu] \, $  as in \S 4.1, define  $ \gerg_\nu $
like therein, and consider  $ \, H := U_{\Bbbk[\nu]}(\gerg_\nu)
= U_R(\gerg_\nu) \, $.  Then we have  $ \; H \Big/ (\nu\!-\!1) \, H
\, = \, U(\gerg) \, = \, \u \Big( \gerg^{{[p\hskip0,7pt]}^\infty} \Big)
\; $  as a  {\sl co-Poisson}  Hopf algebra and  $ \; H \Big/ \nu \,
H \, = \, S(\gerg) \, = \, F[\gerg^\star] \; $  as a  {\sl Poisson}
Hopf algebra; therefore  $ H $  is a QrUEA at  $ \, \h := (\nu\!-\!1)
\, $  (for  $ \u \Big( \gerg^{{[p\hskip0,7pt]}^\infty} \Big) \, $)
and is a QFA at  $ \, \h := \nu \, $  (for  $ F[\gerg^\star] \, $). 
Now we go and study Drinfeld's functors for  $ H $  at
$ (\nu\!-\!1) $  and at  $ (\nu) $.
                                          \par
   Exactly the same procedure as before shows again that
$ \; H^{\vee_{\!(\nu)}} = U \big( {\gerg_\nu}^{\!\vee_{\!(\nu)}}
\big) \, $,  \, from which it follows that  $ \, {\big(
H^{\vee_{\!(\nu)}} \big)}{\Big|}_{\nu=0} \! \cong U(\gerg)
\, $,  \, i.e.~in short  $ \; H^{\vee_{\!(\nu)}} \,{\buildrel
{\, \nu \rightarrow 0 \,} \over \llongrightarrow}\, U(\gerg)
\, $,  \, which is a result quite ``parallel'' to the second
half of  Theorem 2.2{\it (c)}.  Changes occur when looking at
$ \, {\big( H^{\vee_{\!(\nu)}} \big)}^{\prime_{(\nu)}} $:  since
$ \, H^{\vee_{\!(\nu)}} = U \big( {\gerg_\nu}^{\!\vee_{\!(\nu)}}
\big) = \u \Big( {\big( {\gerg_\nu}^{\!\vee_{\!(\nu)}} \big)}^{{[p
\hskip0,7pt]}^\infty} \Big) \, $  we have  $ \, \delta_n(\eta) =
0 \, $  for all  $ \, \eta \in \u \Big( \! {\big( {\gerg_\nu}^{\!
\vee_{\!(\nu)}} \big)}^{{[p\hskip0,7pt]}^\infty} \Big) \, $  such
that  $ \, \partial(\eta) < n \, $  w.r.t.~the standard filtration
of  $ \, \u \Big( {\big( {\gerg_\nu}^{\!\vee_{\!(\nu)}} \big)}^{{[p
\hskip0,7pt]}^\infty} \Big) \, $  (cf.~the proof of  Lemma 4.2{\it
(d)\/}  in [Ga5], which clearly adapts to the present situation):
this implies
  $$  {\big( H^{\vee_{\!(\nu)}} \big)}^{\prime_{(\nu)}}
= \Big\langle \nu \cdot {\big( {\gerg_\nu}^{\! \vee_{\!(\nu)}}
\big)}^{{[p\hskip0,7pt]}^\infty} \Big\rangle  \qquad  \Big( \subset
\, \u \Big( \nu \cdot {\big( {\gerg_\nu}^{\!\vee_{\!(\nu)}}
\big)}^{{[p\hskip0,7pt]}^\infty} \Big) \; \Big)  $$
which is strictly  {\sl bigger\/}  than  $ H $,  because we have 
$ \; \Big\langle \nu \cdot {\big( {\gerg_\nu}^{\! \vee_{\!(\nu)}}
\big)}^{{[p\hskip0,7pt]}^\infty} \Big\rangle = \Big\langle \sum\limits_{n
\geq 0} \nu \cdot {\big( {\gerg_\nu}^{\! \vee_{\!(\nu)}}
\big)}^{{[p\hskip0,7pt]}^n}
               \Big\rangle = $\break
    $ = \Big\langle \gerg_\nu + \nu^{1-p} \, \big\{\, x^p
\,\big|\, x \in \gerg_\nu \big\} + \nu^{1-p^2} \Big\{\, x^{p^2}
\,\Big|\, x \in \gerg_\nu \Big\} + \cdots \Big\rangle \,
\supsetneqq \, U(\gerg_\nu) = H \, $.
                                        \par
   Finally, proceeding as above it is easy to see that  $ \,
H^{\prime_{(\nu)}} = \Big\langle \nu \, P \big( U(\gerg_\nu)
\big) \Big\rangle = \Big\langle \nu \, \gerg^{{[p\hskip0,7pt]}^\infty}
\Big\rangle \, $  whence, letting  $ \, \tilde{\gerg} :=
\nu \, \gerg \, $  and  $ \, \tilde{x} := \nu \, x \, $
for all  $ \, x \in \gerg \, $,  we have
  $$  H^{\prime_{(\nu)}} \; = \; T_R (\tilde{\gerg}) \bigg/
\Big( \Big\{\, \tilde{x} \, \tilde{y} - \tilde{y} \, \tilde{x}
- \nu^2 \, \widetilde{[x,y]} \, , \, \tilde{z}^p - \nu^{p-1}
\widetilde{z^{[p\hskip0,7pt]}} \;\Big\vert\;\, x, y, z \in \gerg
\,\Big\} \Big)  $$   
so that  $ \; H^{\prime_{(\nu)}}
\;{\buildrel {\nu \rightarrow 0} \over
{\relbar\joinrel\llongrightarrow}}\;
T_\Bbbk (\tilde{\gerg}) \bigg/ \Big( \Big\{\, \tilde{x} \,
\tilde{y} - \tilde{y} \, \tilde{x} \, , \, \tilde{z}^p
\;\Big\vert\; \tilde{x}, \tilde{y}, \tilde{z} \in
\tilde{\gerg} \,\Big\} \Big) =
S_\Bbbk (\gerg_{\text{\it ab}}) \Big/ \big( \big\{\,
z^p \,\big\vert\; z \in \gerg
              \,\big\} \big) = $\break
      $ = F[\gerg^\star_{\delta-\text{\it ab}}] \Big/ \big(
\big\{\, z^p \,\big\vert\; z \in \gerg \,\big\} \big) \, $,
\, that is  $ \; H^{\prime_{(\nu)}}{\Big|}_{\nu=0} \!\! \cong
F[\gerg^\star_{\delta-\text{ab}}] \Big/ \big( \big\{\, z^p
\,\big\vert\; z \in \gerg \,\big\} \big) \; $  {\sl as Poisson
Hopf algebras},  where  $ \gerg_{\text{\it ab}} $  and
$ \gerg^\star_{\delta-\text{\it ab}} $  are as above.
Therefore  $ H^{\prime_{(\nu)}} $  is a QFA (at  $ \h = \nu
\, $)  for a non-reduced, zero-dimensional algebraic Poisson
group of height 1, whose cotangent Lie bialgebra is the vector
space  $ \gerg $  with trivial Lie bialgebra structure: this
again yields somehow an analogue of part  {\it (c)\/}  of
Theorem 2.2 for the present case.  If we iterate, we find that
all further images  $ \, \Big( \cdots {\big( {(H)}^{\prime_{(\nu)}}
\big)}^{\prime_{(\nu)}} \cdots \Big)^{\prime_{(\nu)}} \, $  of
the functor  $ {(\ )}^{\prime_{(\nu)}} $  applied to  $ H $ 
are pairwise isomorphic, so that   
  $$  {{\Big( \cdots {\big( {(H)}^{\prime_{(\nu)}}
\big)}^{\prime_{(\nu)}} \cdots \Big)}^{\prime_{(\nu)}}}
{\bigg|}_{\nu=0} \, \cong \;\, S({\gerg}_{ab}) \Big/
\big( \big\{\, z^p \,\big\vert\; z \in \gerg \,\big\}
\big) = F \big[ \gerg^\star_{\delta-{\text{\it ab}}}
\big] \Big/ \big( \big\{\, z^p \,\big\vert\; z \in
\gerg \,\big\} \big) \; .  $$
                                        \par
   Now for Drinfeld's functors at  $ (\nu - 1) $.  Up to minor
changes, with the same procedure and notations as in \S 4.3 we
get analogous results.  First of all, a result analogous to (4.1)
holds, namely  $ \; H^{\prime_{(\nu - 1)}} = \Big\langle (\nu - 1)
\cdot P\big( U(\gerg_\nu) \big) \Big\rangle = \Big\langle (\nu - 1)
\, {( \gerg_\nu )}^{{[p\hskip0,7pt]}^\infty} \Big\rangle = \bigg\langle
{\Big( {(\gerg_\nu)}^{{[p\hskip0,7pt]}^\infty} \Big)}^{\!\prime_{(\nu
- 1)}} \bigg\rangle \, $,  \, which yields
  $$  \displaylines{
   H^{\prime_{(\nu - 1)}} \; = \; T_R \left( \! {\Big(
{(\gerg_\nu)}^{{[p\hskip0,7pt]}^\infty} \Big)}^{\!\prime_{(\nu - 1)}}
\right) \Bigg/ \bigg( \Big\{\, x' \, y' - y' \, x' - (\nu-1) \,
{\big[ x', y']}_* \, ,  \, {(x')}^p - {(\nu - 1)}^{p-1} {\big(
x^{[p\hskip0,7pt]} \big)}' \;\Big|\;  \cr
   \hfill   \;\Big|\; x, y \in {(\gerg_\nu)}^{{[p\hskip0,7pt]}^\infty}
\,\Big\} \bigg)  \cr }  $$
and consequently  $ \; H^{\prime_{(\nu - 1)}}{\Big|}_{(\nu - 1) = 0}
\cong \, S_\Bbbk (\gerg) \Big/ \big( \big\{\, x^p \;\big|\; x \in \gerg
\,\big\} \big) \, = \, F[\gerg^\star] \Big/ \big( \big\{\, x^p \;\big|\;
x \in \gerg \,\big\} \big) \; $  as Poisson Hopf algebras: in a nutshell, 
$ \; H^{\prime_{(\nu - 1)}} \,{\buildrel {\, \nu \rightarrow 1 \,} \over
\llongrightarrow}\, F[\gerg^\star] \Big/ \big( \big\{\, x^p \;\big|\;
x \in \gerg \,\big\} \big) \; $.
                                        \par
   Iterating, one finds again that all  $ \, \Big( \cdots {\big(
{(H)}^{\prime_{(\nu)}} \big)}^{\prime_{(\nu - 1)}} \cdots   
\Big)^{\prime_{(\nu)}} \, $   \hbox{are pairwise isomorphic, so}
  $$  {{\Big( \! \cdots {\big( {(H)}^{\prime_{(\nu - 1)}}
\big)}^{\prime_{((\nu - 1)}} \! \cdots \Big)}^{\prime_{(\nu
- 1)}}}{\bigg|}_{(\nu - 1)=0}  \hskip-5pt  \cong \; S({\gerg}_{ab})
\Big/ \big( \big\{ z^p \,\big\vert\; z \! \in \! \gerg \big\} \big)
= F \big[ \gerg^\star_{\delta-{\text{\it ab}}} \big] \Big/ \big(
\big\{ z^p \,\big\vert\; z \! \in \! \gerg \,\big\} \big) \; .  $$   
                                        \par
   Further on, one has  $ \, {\big( H^{\prime_{(\nu - 1)}}
\big)}^{\vee_{(\nu - 1)}} = {\Big\langle (\nu - 1) \, {(
\gerg_\nu )}^{{[p\hskip0,7pt]}^\infty} \Big\rangle}^{\vee_{(\nu
\!-\!1)}} = \big\langle {(\nu - 1)}^{-1} \cdot (\nu - 1)
                        \, \gerg_\nu \big\rangle = $\break
$ = \big\langle \gerg_\nu \big\rangle = U_R(\gerg_\nu) =: H \, $,
\, which perfectly agrees with  Theorem 2.2{\it (b)}.  Finally,
for  $ H^{\vee_{(\nu - 1)}} $  one has again the same feature
as in \S 3.7: one has to apply the analysis therein, however the
$ p $--filtration  in this case is ``harmless'', since it is
``encoded'' in the standard filtration of  $ U(\gerg) $.  In
any case the final result will depend on the lower central
series of  $ \gerg \, $.
                                        \par
   Second, we assume in addition that  $ \gerg $  be a  {\sl
restricted\/}  Lie algebra and consider  $ \, H := \u_{\Bbbk[\nu]}
(\gerg_\nu) = \u_R(\gerg_\nu) \, $.  In this case we have  $ \; H
\Big/ (\nu - 1) \, H \, = \, \u(\gerg) \; $  as a  {\sl co-Poisson\/} 
Hopf algebra,  and  $ \; H \Big/ \nu \, H \, = \, S(\gerg) \Big/
\big( \big\{\, z^p \,\big\vert\; z \in \gerg \,\big\} \big) \, =
\, F[\gerg^\star] \Big/ \big( \big\{\, z^p \,\big\vert\; z \in \gerg
\,\big\} \big) \; $  as a  {\sl Poisson\/}  Hopf algebra, which
means that  $ H $  is a QrUEA at  $ \, \h := (\nu\!-\!1) \, $ 
(for  $ \u(\gerg) \, $)  and is a QFA at  $ \, \h := \nu \, $ (for 
$ F[\gerg^\star] \Big/ \big( \big\{\, z^p \,\big\vert\; z \in \gerg
\,\big\} \big) \, $).  Then we go and study Drinfeld's functors for 
$ H $  at  $ (\nu - 1) $  and at  $ (\nu) \, $.
                                          \par
   As for  $ \, H^{\vee_{\!(\nu)}} \, $,  it depends again
on the  $ p $--operation  of  $ \gerg \, $,  in short because
the  $ I $--filtration  of  $ \u_\nu(\gerg) $  depends on the
$ p $--filtration   of  $ \gerg \, $.  In the previous case
--- i.e.~when  $ \, \gerg = {\gerh}^{{[p\hskip0,7pt]}^\infty} \, $ 
for some Lie algebra  $ \gerh $  ---  the solution was plain, because
the  $ p $--filtration  of  $ \gerg $  is ``encoded'' in the standard
filtration of  $ U(\gerh) $;  but the general case will be more
complicated, and in consequence also the situation for  $ \,
{\big( H^{\vee_{\!(\nu)}} \big)}^{\prime_{(\nu)}} $,  since 
$ \, H^{\vee_{\!(\nu)}} \, $  will be different according to
the nature of  $ \gerg \, $.  Instead, proceeding exactly like
before one finds  $ \, H^{\prime_{(\nu)}} = \Big\langle \nu \,
P \big( u(\gerg_\nu) \big) \Big\rangle = \big\langle \nu \,
\gerg \big\rangle \, $,  whence, letting  $ \, \tilde{\gerg}
:= \nu \, \gerg \, $  and  $ \, \tilde{x} := \nu \, x \, $
   \hbox{for all  $ \, x \in \gerg \, $,  we have}   
  $$  H^{\prime_{(\nu)}} \; = \; T_{\Bbbk[\nu]} (\tilde{\gerg})
\bigg/ \Big( \Big\{\, \tilde{x} \, \tilde{y} - \tilde{y} \,
\tilde{x} - \nu^2 \, \widetilde{[x,y]} \, , \, \tilde{z}^p
- \nu^{p-1} \widetilde{z^{[p\hskip0,7pt]}} \;\Big\vert\;\,
x, y, z \in \gerg \,\Big\} \Big)  $$
so that  $ \; H^{\prime_{(\nu)}}
\;{\buildrel {\nu \rightarrow 0} \over
{\relbar\joinrel\llongrightarrow}}\;
T_\Bbbk (\tilde{\gerg}) \bigg/ \Big( \Big\{\, \tilde{x}
\, \tilde{y} - \tilde{y} \, \tilde{x} \, , \, \tilde{z}^p
\;\Big\vert\; \tilde{x}, \tilde{y}, \tilde{z} \in \tilde{\gerg}
\,\Big\} \Big) = S_\Bbbk (\gerg_{\text{\it ab}}) \Big/ \big(
\big\{\, z^p \,\big\vert\; z \in \gerg
                 \,\big\} \big) = $\break
      $ =F[\gerg^\star_{\delta-\text{\it ab}}] \Big/ \big(
\big\{\, z^p \,\big\vert\; z \in \gerg \,\big\} \big) \, $,
\, that is  $ \; H^{\prime_{(\nu)}}{\Big|}_{\nu=0} \!\! \cong
F[\gerg^\star_{\delta-\text{ab}}] \Big/ \big( \big\{\, z^p
\,\big\vert\; z \in \gerg \,\big\} \big) \; $  {\sl as Poisson
Hopf algebras\/}  (using notation as before).  Thus
$ H^{\prime_{(\nu)}} $  is a QFA (at  $ \h = \nu \, $)
for a non-reduced, zero-dimensional algebraic Poisson group
of height 1, whose cotangent Lie bialgebra is  $ \gerg $ 
with the trivial Lie bialgebra structure: so again we get
an analogue of part of  Theorem 2.2{\it (c)\/}.  Moreover,
iterating again one finds that all  $ \, \Big( \cdots
{\big( {(H)}^{\prime_{(\nu)}} \big)}^{\prime_{( \nu - 1 )}}
\cdots \Big)^{\prime_{(\nu - 1)}} \, $
are pairwise isomorphic, so
 \vskip-8pt
  $$  {{\Big( \! \cdots {\big( {(H)}^{\prime_{(\nu - 1)}}
\big)}^{\prime_{((\nu - 1)}} \! \cdots \Big)}^{\prime_{(\nu
- 1)}}}{\bigg|}_{(\nu - 1)=0}  \cong \;
S({\gerg}_{ab}) \Big/ \big( \big\{ z^p \,\big\vert\; z \! \in \!
\gerg \big\} \big) \; = \; F \big[ \gerg^\star_{\delta-{\text{\it ab}}}
\big] \Big/ \big( \big\{ z^p \,\big\vert\; z \! \in \! \gerg
\,\big\} \big) \; .  $$   
 \vskip-0pt
   \indent   As for Drinfeld's functors at  $ (\nu - 1) $,
the situation is more similar to the previous case of  $ \,
H = U_R(\gerg_\nu) \, $.  First  $ \; H^{\prime_{(\nu - 1)}}
= \Big\langle (\nu - 1) \cdot P \big( \u(\gerg_\nu) \big)
\Big\rangle = \big\langle (\nu - 1) \, \gerg_\nu \big\rangle
=: \big\langle {\gerg_\nu}^{\!\prime_{(\nu - 1)}} \big\rangle \, $, 
\, hence
 \vskip-10pt
  $$  H^{\prime_{(\nu - 1)}} = T_R \Big( {\gerg_\nu}^{\!\prime_{(
\nu - 1 )}} \Big) \bigg/ \! \Big( \Big\{\, x' \, y' - y' \, x' -
(\nu - 1) \, {\big[ x', y']}_* \, ,  \, {(x')}^p - {(\nu - 1)}^{p-1}
{\big( x^{[p\hskip0,7pt]} \big)}' \Big\}_{x, y \in \gerg_\nu \,}
\Big)  $$   
 \vskip1pt
\noindent
thus again  $ \; H^{\prime_{(\nu - 1)}}{\Big|}_{(\nu - 1) = 0} \!
\cong S_\Bbbk (\gerg) \Big/ \big( \big\{\, x^p \;\big|\; x \in
\gerg \,\big\} \big) = F[\gerg^\star] \Big/ \big( \big\{\, x^p
\;\big|\; x \in \gerg \,\big\} \big) \; $  as Poisson Hopf
algebras, that is  $ \; H^{\prime_{(\nu - 1)}} \,{\buildrel
{\, \nu \rightarrow 1 \,} \over \llongrightarrow}\, F[\gerg^\star]
\Big/ \big( \big\{\, x^p \;\big|\; x \in \gerg \,\big\} \big) \, $.
Iteration then shows that all  $ \, \Big( \cdots {\big(
{(H)}^{\prime_{(\nu)}} \big)}^{\prime_{(\nu - 1)}} \cdots
\Big)^{\prime_{(\nu)}} \, $  are pairwise isomorphic, so
that again   
 \vskip-9pt
  $$  {{\Big( \! \cdots {\big( {(H)}^{\prime_{(\nu - 1)}}
\big)}^{\prime_{((\nu - 1)}} \! \cdots \Big)}^{\prime_{(\nu
- 1)}}}{\bigg|}_{(\nu - 1)=0}  \cong \;  S({\gerg}_{ab}) \Big/
\big( \big\{ z^p \,\big\vert\; z \! \in \! \gerg \big\} \big) 
\; = \;  F \big[ \gerg^\star_{\delta-{\text{\it ab}}} \big]
\Big/ \big( \big\{ z^p \,\big\vert\; z \! \in \! \gerg
\,\big\} \big) \; .  $$
 \vskip-0pt
   \indent   Further, we have  $ \, {\big( H^{\prime_{(\nu - 1)}}
\big)}^{\vee_{(\nu - 1)}} = {\big\langle (\nu - 1) \, \gerg_\nu
\big\rangle}^{\vee_{(\nu - 1)}} = \big\langle \gerg_\nu \big\rangle
= \u_R(\gerg_\nu) =: H \, $,  \, which agrees at all with  Theorem
2.2{\it (b)}.  Finally,  $ H^{\vee_{(\nu - 1)}} $  again has the same
feature as in \S 3.7: in particular, the outcome strongly depends on
the properties of  {\sl both\/}  the lower central series  {\sl and\/} 
   \hbox{of the  $ p $--filtration  of  $ \gerg \, $.}

\vskip7pt

  {\bf 4.5 The hyperalgebra case.} \, Let  $ \Bbbk $  be again a field
with  $ \, \hbox{\it Char}\,(\Bbbk) = p > 0 \, $.  Like in \S 3.11, let
$ G $  be an algebraic group (finite-dimensional, for simplicity), and
let  $ \, \hyp(G) := {\big( {F[G]}^\circ \big)}_\epsilon = \big\{\, \phi
\in {F[G]}^\circ \,\big|\, \phi({\germ_e}^{\!n}) = 0 \, , \, \forall \;
n \gg 0 \,\big\} \, $  be the hyperalgebra associated to  $ G \, $  (see
\S 1.1).
                                               \par
   For each  $ \, \nu \in \Bbbk \, $,  \, let  $ \, \gerg_\nu := \big(
\gerg \, , {[\ ,\ ]}_\nu \big) \, $  be the Lie algebra given by  $ \,
\gerg \, $  endowed with the rescaled Lie bracket  $ \, {[\,\ ,\ ]}_\nu
:= \nu \, {[\,\ ,\ ]}_\gerg \; $.  By general theory, the algebraic group 
$ G $  is uniquely determined by a neighborhood of the identity together
with the formal group law uniquely determined by  $ {[\,\ ,\ ]}_\gerg
\; $.  Similarly, a neighborhood of the identity of  $ G $  together
with  $ {[\,\ ,\ ]}_\nu $  uniquely determines a new connected algebraic
group  $ G_\nu \, $,  whose hyperalgebra  $ \hyp(G_\nu) $  is an
algebraic deformation of  $ \hyp(G) \, $;  moreover,  $ G_\nu $  is
birationally equivalent to  $ G $,  and for  $ \, \nu \not= 0 \, $
they are also isomorphic as algebraic groups, via an isomorphism
induced by  $ \, \gerg \cong \gerg_\nu \, $,  $ \, x \mapsto \nu^{-1}
x \, $  (however, this may not be the case when  $ \, \nu = 0 \, $).
Note that  $ \hyp(G_0) $  is clearly commutative, because  $ G_0 $  is
Abelian: indeed,
%
%
%
%
we have
  $$  \hyp(G_0)  \; = \;  S_\Bbbk \Big( \gerg^{{(p\hskip0,7pt)^\infty}}
\Big) \hskip-2pt \bigg/ \hskip-2pt \Big( {\big\{\, x^p \,\big\}}_{x \in
\gerg^{{(p\hskip0,7pt)}^\infty}} \Big)  \, = \,  F \Big[ \Big( \gerg^{{(p
\hskip0,7pt)^\infty}} \Big)^{\!\star} \Big] \bigg/ \hskip-2pt \Big(
{\big\{\, y^p \,\big\}}_{y \in \gerg^{{(p\hskip0,7pt)}^\infty}}
\Big)  $$  
where  $ \; \gerg^{{(p\hskip0,7pt)}^\infty} \! := \text{\sl Span}\, \Big(
\Big\{\, x^{{(p^n\hskip0pt)}} \,\Big\vert\, x \in \gerg \, , n \in \N
\,\Big\} \Big) \, $;  \, here as usual  $ x^{(n)} $  denotes the
$ n $--th  divided power of  $ \, x \in \gerg \, $  (recall that
$ \hyp(G) $,  hence also  $ \hyp(G_\nu) $,  is generated as an
algebra by all the  $ x^{(n)} $'s,  some of which might be zero).
So  $ \, \hyp(G_0) = F[\varGamma] \, $  where  $ \varGamma $  is a
connected algebraic group of dimension zero and height 1: moreover,
$ \varGamma $  is a Poisson group, with cotangent Lie bialgebra
$ \gerg^{{(p\hskip0,7pt)^\infty}} $  and Poisson bracket induced
by the Lie bracket of  $ \gerg \, $.
                                  \par
%
%
%
%
%
   Now think at  $ \nu $  as a parameter in  $ \, R := \Bbbk[\nu] \, $
(as in \S 4.1), and set  $ \, H := \Bbbk[\nu] \otimes_\Bbbk \hyp(G_\nu)
\, $.  Then we find a situation much similar to that of \S 4.1,
which we shall shortly describe.
                                            \par
   Namely,  $ H $  is a free  $ \Bbbk[\nu] $--algebra,  thus  $ \,
H \in \HA \, $  and  $ \, H_F := \Bbbk(\nu) \otimes_{\Bbbk[\nu]} H
\in \HA_F \, $  (see \S 1.3); its specialization at  $ \, \nu = 1
\, $  is  $ \; H \Big/ (\nu\!-\!1) \, H \, = \, \hyp(G_1) \, = \,
\hyp(G) \, $,  \, and at  $ \, \nu = 0 \, $  is  $ \; H \Big/ \nu
\, H \, = \, \hyp(G_0) \, = \, F[\varGamma] \, $  (as a  {\sl Poisson}
Hopf algebra),  or  $ \; H \,{\buildrel {\, \nu \rightarrow 1 \,} \over
\llongrightarrow}\, \hyp(G) \, $  and  $ \; H \,{\buildrel {\, \nu
\rightarrow 0 \,} \over \llongrightarrow}\, F[\varGamma] \, $,
i.e.~$ H $  is a ``quantum hyperalgebra'' at  $ \, \h := (\nu-1) \, $
and a QFA at  $ \, \h := \nu \, $.  Now we study Drinfeld's functors
for  $ H $  at  $ \, \h = (\nu\!-\!1) \, $  and at  $ \h = \nu \, $.
                                        \par
   First, a straightforward analysis like in \S 4.2 yields  $ \,
H^{\vee_{\!(\nu)}} \cong \Bbbk[\nu] \otimes_\Bbbk \hyp(G) \, $
(induced by  $ \, \gerg \cong \gerg_\nu \, $,  $ \, x \mapsto \nu^{-1}
x \, $)  whence in particular  $ \, {\big( H^{\vee_{\!(\nu)}}
\big)}{\Big|}_{\nu=0} \cong \hyp(G) \, $,  \; that is  $ \;
H^{\vee_{\!(\nu)}} \,{\buildrel {\, \nu \rightarrow 0 \,}
\over \llongrightarrow}\, \hyp(G) \, $.  Second, one can also
see (essentially,  {\sl mutatis mutandis},  like in \S 4.2) that
$ \, {\big( H^{\vee_{\!(\nu)}} \big)}^{\prime_{(\nu)}} = H \, $,
\, whence  $ \; {\big( H^{\vee_{\!(\nu)}} \big)}^{\prime_{(\nu)}}
{\Big|}_{\nu=0} = \, H{\Big|}_{\nu=0} = \, \hyp(G_0) \, = \,
F[\varGamma] \; $  follows.
                                        \par
   At  $ \, \h = (\nu-1) \, $,  \, we can see by direct computation
that  $ \; H^{\prime_{(\nu-1)}} = \Big\langle \big( \gerg^{{(p
\hskip0,7pt)}^\infty} \big)^{\!\prime_{(\nu-1)}} \Big\rangle \; $
where  $ \, \big( \gerg^{{(p\hskip0,7pt)}^\infty} \big)^{\!\prime_{(\nu
- 1)}} := \text{\sl Span}\, \Big( \Big\{\, (\nu-1)^{p^n} x^{{(p^n)}}
\,\Big\vert\, x \in \gerg \, , n \in \N \,\Big\} \Big) \, $.  Indeed the
structure of  $ H^{\prime_{(\nu-1)}} $  depends only on the coproduct
of  $ H $,  in which  $ \nu $  plays no role; therefore we can do the same
analysis as in the trivial deformation case (see \S 3.11): the filtration
$ \underline{D} $  of  $ \hyp(G_\nu) $  is just the natural filtration
given by the order (of divided powers), and this yields the previous
description of  $ H^{\prime_{(\nu-1)}} \, $.  At  $ \, \nu = 1 \, $ 
we find   
  $$  H^{\prime_{(\nu-1)}} \Big/ (\nu-1) \, H^{\prime_{(\nu-1)}}
\; \cong \;  S_\Bbbk \Big( \gerg^{{(p\hskip0,7pt)^\infty}} \Big)
\hskip-2pt \bigg/ \hskip-2pt \Big( {\big\{\, x^p \,\big\}}_{x \in
\gerg^{{(p\hskip0,7pt)}^\infty}} \Big)  \; = \;  \hyp(G_0)  \; =
\;  F[\varGamma]  $$
as Poisson Hopf algebras: in short,  $ H^{\prime_{(\nu-1)}} $  is a
QFA, at  $ \, \h = \nu-1 \, $,  for the Poisson group  $ \varGamma $.
Similarly  $ \, H^{\prime_{(\nu)}} = \Big\langle \big( \gerg^{{(p
\hskip0,7pt)}^\infty} \big)^{\!\prime_{(\nu)}} \Big\rangle \, $  with 
$ \, \big( \gerg^{{(p\hskip0,7pt)}^\infty} \big)^{\!\prime_{(\nu)}} \!
:= \text{\sl Span}\, \Big( \Big\{\, \nu^{p^n} x^{{(p^n)}} \,\Big\vert\,
x \in \! \gerg \, , \, n \in \N \,\Big\} \Big) \, $;  \, therefore   
  $$  H^{\prime_{(\nu)}} \Big/ \nu \, H^{\prime_{(\nu)}}  \; \cong \;
S_\Bbbk \Big( \gerg_{\text{\it ab}}^{{(p\hskip0,7pt)^\infty}} \Big)
\hskip-2pt \bigg/ \hskip-2pt \Big( {\big\{\, x^p \,\big\}}_{x \in
\gerg^{{(p\hskip0,7pt)}^\infty}} \Big)  \; = \;  F \big[
\varGamma_{\text{\it ab\,}} \big]  $$   
where  $ \gerg_{\text{\it ab}} $  is simply  $ \gerg $  with trivialized
Lie bracket and  $ \varGamma_{\text{\it ab}} $  is the  {\sl same
algebraic group\/}  as  $ \varGamma $  but with  {\sl trivial\/}
Poisson bracket: this comes essentially like in \S 4.2, roughly
because  $ \, \big\{ \overline{\nu\,x} \, , \overline{\nu\,y} \big\}
:= \big( \nu^{-1} [\nu\,x \, , \nu\,y \,] \big){\Big|}_{\nu=0} =
\big( \nu^{-1} \cdot \nu^3 [x,y]_\gerg \big){\Big|}_{\nu=0} =
\big( \nu \cdot \nu \, [x,y]_\gerg \big){\Big|}_{\nu=0} = 0 \; $ 
(for all  $ \, x, y \in \gerg \, $).   
                                       \par
   Finally, we have  $ \, {\big( H^{\prime_{(\nu-1)}} \big)}^{\vee_{\!
(\nu-1)}} = \Big\langle \Big\{\, (\nu-1)^{p^n - 1\,} x^{{(p^n)}}
\,\Big\vert\, x \in \gerg \, , \, n \in \N \,\Big\} \Big\rangle \,
\varsubsetneqq H \, $  and  $ \, {\big( H^{\prime_{(\nu)}}
\big)}^{\vee_{\!(\nu)}} = \, \Big\langle \Big\{\, \nu^{p^n - 1\,}
x^{{(p^n)}} \,\Big\vert\, x \in \gerg \, , \, n \in \N \,\Big\}
\Big\rangle \, \varsubsetneqq H \, $,  \, by direct computation.
For  $ H^{\vee_{(\nu-1)}} $  we have the same features as in \S
3.7: the analysis therein can be repeated, with the final upshot
depending on the nature of  $ G $  (or of  $ \gerg \, $,  essentially,
in particular on its  $ p $--lower  central series).   

\vskip1,1truecm

\centerline {\bf \S \; 5 \  Second example:  $ {SL}_2 \, $,
$ {SL}_n $  and the semisimple case }

\vskip10pt

  {\bf 5.1 The classical setting.} \, Let  $ \Bbbk $  be any field of
characteristic  $ \, p \geq 0 \, $.  Let  $ \, G := {SL}_2(\Bbbk) \equiv
{SL}_2 \, $;  \, its tangent Lie algebra  $ \, \gerg = \gersl_2 \, $
is generated by  $ \, f $,  $ h $,  $ e \, $  (the  {\it Chevalley
generators\/})  with relations  $ \, [h,e] = 2 \, e $,  $ [h,f] =
-2 f $,  $ [e,f] = h \, $.  The formulas  $ \, \delta(f) = h \otimes f
- f \otimes h \, $,  $ \, \delta(h) = 0 \, $,  $ \, \delta(e) = h \otimes
e - e \otimes h \, $,  define a Lie cobracket on  $ \gerg $  which makes
it into a Lie bialgebra, corresponding to a structure of Poisson group
on  $ G $.  These formulas give also a presentation of the co-Poisson
Hopf algebra  $ U(\gerg) $  (with the standard Hopf structure).  If 
$ \, p > 0 \, $,  the  $ p $--operation  in  $ \gersl_2 $  is given
by  $ \, e^{[p\,]} = 0 \, $,  $ \, f^{[p\,]} = 0 \, $,  $ \, h^{[p\,]}
= h \, $.   
                                             \par
  On the other hand,  $ F[{SL}_2] $  is the unital associative
commutative  $ \Bbbk $--algebra  with generators  $ \, a $,
$ b $,  $ c $,  $ d \, $  and the relation  $ \, a d - b c
= 1 \, $,  and Poisson Hopf structure given by
  $$  \displaylines{
  \Delta(a) = a \otimes a + b \otimes c \, ,  \;\,
\Delta(b) = a \otimes b + b \otimes d \, ,  \;\,
\Delta(c) = c \otimes a + d \otimes c \, ,  \;\,
\Delta(d) = c \otimes b + d \otimes d  \cr
  \epsilon(a) = 1  \, ,  \hskip5pt  \epsilon(b) = 0 \, ,  \hskip5pt
\epsilon(c) = 0 \, ,  \hskip5pt  \epsilon(d) = 1 \, ,  \hskip13pt
S(a) = d \, ,  \hskip5pt  S(b) = -b \, ,  \hskip5pt  S(c) = - c \, ,
\hskip5pt  S(d) = a  \cr
  \{a,b\} = b \, a \, , \quad  \{a,c\} = c \, a \, ,  \quad  \{b,c\}
= 0 \, ,  \quad  \{d,b\} = - b \, d \, ,  \quad  \{d,c\} = - c \, d
\, ,  \quad  \{a,d\} = 2 \, b \, c \, .  \cr }  $$
                                             \par
  The dual Lie bialgebra  $ \, \gerg^* = {\gersl_2}^{\!*} \, $
is the Lie algebra with generators  $ \, \text{f} $,  $ \text{h} $,
$ \text{e} \, $,  and relations  $ \, [\text{h},\text{e}] = \text{e} $,
$ [\text{h},\text{f}\,] = \text{f} $,  $ [\text{e},\text{f}\,] = 0 \, $,
with Lie cobracket given by  $ \, \delta(\text{f}\,) = 2 (\text{f}
\otimes \text{h} - \text{h} \otimes \text{f}\,) $,  $ \, \delta(\text{h}) =
\text{e} \otimes \text{f} - \text{f} \otimes \text{e} $,  $ \,
\delta(\text{e}) = 2 (\text{h} \otimes \text{e} - \text{e} \otimes
\text{h}) \, $  (we choose as generators  $ \, \text{f} := f^* \, $,
$ \, \text{h} := h^* \, $,  $ \, \text{e} := e^* \, $,  where
$ \, \big\{ f^*, h^*, e^* \big\} \, $  is the basis of
$ {\gersl_2}^{\! *} $  which is the dual of the basis  $ \{ f, h, e \} $
of  $ \gersl_2 \, $).  This again yields also a presentation of  $ \,
U \left( {\gersl_2}^{\!*} \right) \, $.  If  $ \, p > 0 \, $,  the
$ p $--operation  in  $ {\gersl_2}^{\!*} $  is given by  $ \,
\text{e}^{[p\,]} = 0 \, $,  $ \, \text{f}^{\,[p\,]} = 0 \, $,  $ \,
\text{h}^{[p\,]} = \text{h} \, $.  The simply connected algebraic
Poisson group whose tangent Lie bialgebra is $ {\gersl_2}^{\!*} $
can be realized as the group of pairs of matrices (the left subscript
$ s $  meaning ``simply connected'')
  $$  {}_s{{SL}_2}^{\!*} = \Bigg\{\,
\bigg( \bigg( \begin{matrix}  z^{-1} & 0 \\  y & z  \end{matrix} \bigg) \, ,
\bigg( \begin{matrix}  z & x \\  0 & z^{-1}  \end{matrix} \bigg) \bigg)
\,\Bigg\vert\, x, y \in k, z \in \Bbbk \setminus \{0\} \,\Bigg\}
\; \leq \; {SL}_2 \times {SL}_2 \; .  $$
This group has centre  $ \, Z := \big\{ (I,I), (-I,-I) \big\} \, $,
so there is only one other (Poisson) group sharing the same Lie
(bi)algebra, namely the quotient  $ \, {}_a{{SL}_2}^{\!*} :=
{}_s{SL_2}^* \Big/ Z \, $  (the adjoint of  $ \, {}_s{{SL}_2}^{\!*}
\, $,  as the left subscript  $ a $  means).  Therefore  $ F\big[
{}_s{{SL}_2}^{\!*} \big] $  is the unital associative commutative
$ \Bbbk $--algebra  with generators  $ \, x $,  $ z^{\pm 1} $,
$ y $,  with Poisson Hopf structure given by
  $$  \displaylines{
   \Delta(x) = x \otimes z^{-1} + z \otimes x \, ,  \hskip21pt
\Delta\big(z^{\pm 1}\big) = z^{\pm 1} \otimes z^{\pm 1} \, ,
\hskip21pt  \Delta(y) = y \otimes z^{-1} + z \otimes y  \cr
  \epsilon(x) = 0  \, ,  \hskip10pt  \epsilon\big(z^{\pm 1}\big) = 1
\, ,  \hskip10pt  \epsilon(y) = 0 \, ,  \hskip31pt
  S(x) = -x \, ,  \hskip10pt  S\big(z^{\pm 1}\big) = z^{\mp 1} \, ,
\hskip10pt  S(y) = -y  \cr
  \{x,y\} = \big( z^2 - z^{-2} \big) \big/ 2 \, ,  \hskip27pt
\big\{z^{\pm 1},x\big\} = \pm \, x \, z^{\pm 1} \, ,  \hskip27pt
\big\{z^{\pm 1},y\big\} = \mp \, z^{\pm 1} y  \cr }  $$
(Remark: with respect to this presentation, we have  $ \, \text{f} =
{\partial_y}{\big\vert}_e \, $,  $ \, \text{h} = z \, {\partial_z}
{\big\vert}_e \, $,  $ \, \text{e} = {\partial_x}{\big\vert}_e \, $, 
\, where  $ e $  is the identity element of  $ {}_s{SL_2}^* \, $). 
Moreover,  $ F\big[{}_a{{SL}_2}^{\!*}\big] $ can be identified with
the Poisson Hopf subalgebra of  $ F\big[{}_s{SL_2}^*\big] $  spanned
by products of an even number of generators   --- i.e. monomials of
even degree: this is generated, as a unital subalgebra, by  $ \, x
z \, $,  $ \, z^{\pm 2} \, $,  \, and  $ \, z^{-1} y \, $.   
                                             \par   
   In general, we shall consider  $ \, \gerg = \gerg^\tau \, $  a
semisimple Lie algebra, endowed with the Lie cobracket   ---
depending on the parameter  $ \tau $  ---   given in [Ga1], \S 1.3;
in the following we shall also retain from  [{\it loc.~cit.}]  all
the notation we need: in particular, we denote by  $ Q $,
resp.~$ P $,  the root lattice, resp.~the weight lattice, of
$ \gerg \, $,  \, and by  $ r $  the rank of  $ \gerg \, $.

\vskip7pt

  {\bf 5.2 The\footnote
  {In \S\S 5--7 we should use
notation  $ \, U_{q-1}(\gerg) \, $  and  $ \, F_{q-1}[G] \, $,
after Remark 1.5 (for  $ \, \h = q-1 \, $);  instead, we write
$ \, U_q(\gerg) \, $  and  $ \, F_q[G] \, $  to be consistent
with the standard notation in use for these quantum algebras.}
           QrUEAs  $ \, U_q(\gerg) \, $.} \, We turn now to quantum
groups, starting with the  $ \gersl_2 $  case.  Let  $ R $   be any
1dD,  $ \, \h \in R \setminus \{0\} \, $  a prime such that  $ \, R
\big/ \h \, R = \Bbbk \, $;  \, moreover, letting  $ \, q := \h + 1
\, $  we assume that  $ q $  be invertible in  $ R $,  i.e.~there
exists  $ \, q^{-1} = {(\h + 1)}^{-1} \in R  \, $.  For instance,
one can pick  $ \, R := \Bbbk \! \left[ q, q^{-1} \right] \, $ 
for an indeterminate  $ q $  and  $ \, \h := q-1 \, $,  then 
$ \, F(R) = \Bbbk(q) \, $.   
                                            \par
   Let  $ \, \Bbb{U}_q(\gerg) = \Bbb{U}_q(\gersl_2) \, $  be
the associative unital  $ F(R) $--algebra  with (Chevalley-like)
generators  $ \, F $,  $ K^{\pm 1} $,  $ E $,  and relations
  $$  K K^{-1} = 1 = K^{-1} K \; ,  \;\;  K^{\pm 1} F = q^{\mp 2} F
K^{\pm 1} \; ,  \;\;  K^{\pm 1} E = q^{\pm 2} E K^{\pm 1} \; ,  \;\;
E F - F E = \frac{\, K - K^{-1} \,}{\, q - q^{-1} \,} \;\; .  $$
This is a Hopf algebra, with Hopf structure given by
  $$  \displaylines{
   \Delta(F) = F \otimes K^{-1} + 1 \otimes F \, ,  \hskip21pt
\Delta \big( K^{\pm 1} \big) = K^{\pm 1} \otimes K^{\pm 1} \, ,
\hskip21pt  \Delta(E) = E \otimes 1 + K \otimes E  \cr
   \epsilon(F) = 0  \, ,  \hskip3pt  \epsilon \big( K^{\pm 1} \big)
= 1 \, ,  \hskip3pt  \epsilon(E) = 0 \, ,  \hskip15pt  S(F) = - F K ,
\hskip3pt  S \big( K^{\pm 1} \big) = K^{\mp 1} ,  \hskip3pt
S(E) = - K^{-1} E \, .  \cr }  $$
Then let  $ \, U_q(\gerg) \, $  be the  $ R $--subalgebra  of
$ \Bbb{U}_q(\gerg) $  generated by  $ \, F \, $,  $ \, H :=
\displaystyle{\, K - 1 \, \over \, q - 1 \,} \, $,  $ \,
\varGamma := \displaystyle{\, K - K^{-1} \, \over \, q -
q^{-1} \,} \, $,  $ K^{\pm 1} \, $,  $ E \, $.  From the definition
of  $ \Bbb{U}_q(\gerg) $  one gets a presentation of  $ U_q(\gerg) $
as the associative unital algebra with generators  $ \, F $,  $ H $,
$ \varGamma $,  $ K^{\pm 1} $,  $ E $  and relations
  $$  \displaylines{
   K K^{-1} = 1 = K^{-1} K \, ,  \;\quad  K^{\pm 1} H = H K^{\pm 1} \, ,
\;\quad  K^{\pm 1} \varGamma = \varGamma K^{\pm 1} \, ,  \;\quad  H
\varGamma = \varGamma H  \cr
     ( q - 1 ) H = K - 1 \, ,  \quad  \big( q - q^{-1} \big) \varGamma = K -
K^{-1} \, ,  \quad  H \big( 1 + K^{-1} \big) = \big( 1 + q^{-1} \big)
\varGamma \, ,  \quad  E F - F E = \varGamma  \cr
 }  $$   
  $$  \displaylines{
   K^{\pm 1} F = q^{\mp 2} F K^{\pm 1} \, ,  \,\quad  H F = q^{-2} F H -
(q+1) F \, ,  \,\quad  \varGamma F = q^{-2} F \varGamma - \big( q + q^{-1}
\big) F  \cr
   K^{\pm 1} E = q^{\pm 2} E K^{\pm 1} \, ,  \,\quad  H E = q^{+2} E H +
(q+1) E \, ,  \,\quad  \varGamma E = q^{+2} E \varGamma + \big( q + q^{-1}
\big) E  \cr }  $$   
and with a Hopf structure given by the same formulas as above
for  $ \, F $,  $ K^{\pm 1} $,  and  $ E \, $  plus
  $$  \begin{matrix}
      \Delta(\varGamma) = \varGamma \otimes K + K^{-1} \otimes \varGamma
\, , &  \hskip10pt  \epsilon(\varGamma) = 0  \, ,  &  \hskip10pt
S(\varGamma) = - \varGamma  \\
      \Delta(H) = H \otimes 1 + K \otimes H \, ,  &  \hskip10pt
\epsilon(H) = 0 \, ,  &  \hskip10pt  S(H) = - K^{-1} H \, .  \\
     \end{matrix}  $$
Note also that  $ \, K = 1 + (q-1) H \, $  and  $ \, K^{-1} = K -
\big( q - q^{-1} \big) \varGamma = 1 + (q-1) H - \big( q - q^{-1}
\big) \varGamma \, $,  hence  $ \, U_q(\gerg) \, $  is generated
even by  $ \, F $,  $ H $, $ \varGamma $  and $ E $  alone.  Further,
notice also that
  $$  \begin{array}{rcl}
   \hskip2,5cm  \Bbb{U}_q(\gerg) \!\!  &  \! = \; \hbox{free 
$ F(R) $--module  over}  \;\;  \Big\{\, F^a K^z E^d \,\Big\vert\,
a, d \in \N, z \in \Z \,\Big\}  &   \hskip1,6cm \hskip-0.7pt   (5.1)  \\
   \hskip2,7cm  U_q(\gerg) \hskip-4pt  &  \, = \; \hbox{$ R $--span of} 
\;\;  \Big\{\, F^a H^b \varGamma^c E^d \,\Big\vert\, a, b, c, d \in \N
\,\Big\} \;\; \hbox{inside} \;\; \Bbb{U}_q(\gerg)  &   \hskip1,6cm
\hskip-0.7pt   (5.2)  \end{array}  $$
which implies that  $ \, F(R) \otimes_R U_q(\gerg) = \Bbb{U}_q(\gerg)
\, $.  Moreover, definitions imply at once that  $ \, U_q(\gerg) \, $
is torsion-free, and also that it is a Hopf  $ R $--subalgebra  of
$ \Bbb{U}_q(\gerg) \, $.  Therefore  $ \, U_q(\gerg) \in \HA \, $,  \,
and in fact  $ \, U_q(\gerg) \, $  is even a QrUEA, whose semiclassical
limit is  $ \, U(\gerg) = U(\gersl_2) \, $,  \, with the generators
$ \, F $,  $ K^{\pm 1} $,  $ H $,  $ \varGamma $,  $ E \, $
respectively mapping to  $ \, f $,  $ 1 $,  $ h $,  $ h $,
$ e \in U(\gersl_2) \, $.
                                       \par
   It is also possible to define a ``simply connected" version
of  $ \Bbb{U}_q(\gerg) $  and  $ U_q(\gerg) $,  obtained from
the previous ones   --- called ``adjoint'' ---   as follows.  For 
$ \Bbb{U}_q(\gerg) $,  one adds a square root of  $ K^{\pm 1} $, 
call it  $ L^{\pm 1} $,  as new generator; for  $ U_q(\gerg) $ 
one adds the new generators  $ L^{\pm 1} $  and also  $ \, D :=
\displaystyle{ {\, L - 1 \,} \over {\, q - 1 \,} } \, $.  Then
the same analysis as before shows that  $ U_q(\gerg) $  is another
quantization (containing the ``adjoint'' one) of  $ U(\gerg) \, $.
                                               \par
   In the general case of semisimple  $ \gerg \, $,  let  $ \Bbb{U}_q
(\gerg) $  be the Lusztig-like quantum group   --- over  $ R $  ---  
associated to  $ \, \gerg = \gerg^\tau \, $  as in [Ga1], namely  $ \,
\Bbb{U}_q(\gerg) := U_{q,\varphi}^{\scriptscriptstyle M}(\gerg) \, $ 
with respect to the notation in  [{\it loc.~cit.}],  where  $ M $  is
any intermediate lattice such that  $ \, Q \leq M \leq P \, $  (this
is just a matter of choice, of the type mentioned in the statement of 
Theorem 2.2{\it (c)}):  this is a Hopf algebra over  $ F(R) $,  generated
by elements  $ \, F_i \, $,  $ M_i $,  $ E_i $  for  $ \, i=1, \dots, r
=: \hbox{\it rank}\,(\gerg) \, $.  Then let  $ \, U_q(\gerg) \, $  be the
unital  $ R $--subalgebra  of  $ \Bbb{U}_q(\gerg) $  generated by the
elements  $ \, F_i \, $,  $ \, H_i := \displaystyle{\, M_i - 1 \, \over
\, q - 1 \,} \, $,  $ \, \varGamma_i := \displaystyle{\, K_i - K_i^{-1}
\, \over \, q - q^{-1} \,} \, $,  $ \, M_i^{\pm 1} \, $,  $ \, E_i
\, $,  where the  $ \, K_i = M_{\alpha_i} \, $  are suitable product
of  $ M_j $'s,  defined as in [Ga1], \S 2.2 (whence  $ \, K_i \, $, 
$ \, K_i^{-1} \in U_q(\gerg) \, $).  From [Ga1], \S\S 2.5, 3.3, we
have that  $ \Bbb{U}_q(\gerg) $  is the free  $ F(R) $--module 
with basis the set of monomials
  $$  \bigg\{\, \textstyle{\prod\limits_{\alpha \in \Phi^+}}
F_\alpha^{f_\alpha} \cdot \textstyle{\prod\limits_{i=1}^n}
K_i^{z_i} \cdot \textstyle{\prod\limits_{\alpha \in \Phi^+}}
E_\alpha^{e_\alpha} \,\bigg\vert\, f_\alpha, e_\alpha \in \N \, ,
\, z_i \in \Z \, , \;\; \forall \; \alpha \in \Phi^+, \, i = 1,
\dots, n \,\bigg\}  $$
while  $ U_q(\gerg) $  is the  $ R $--span  inside  $ \Bbb{U}_q(\gerg) $
of the set of monomials
  $$  \bigg\{\, \textstyle{\prod\limits_{\alpha \in \Phi^+}}
F_\alpha^{f_\alpha} \cdot \textstyle{\prod\limits_{i=1}^n}
H_i^{t_i} \cdot \textstyle{\prod\limits_{j=1}^n} \varGamma_j^{c_j}
\cdot \textstyle{\prod\limits_{\alpha \in \Phi^+}} E_\alpha^{e_\alpha}
\,\bigg\vert\, f_\alpha, t_i, c_j, e_\alpha \in \N \;\; \forall \;
\alpha \in \Phi^+, \, i, j= 1, \dots, n \,\bigg\}  $$
(hereafter,  $ \Phi^+ $  is the set of positive roots of  $ \gerg \, $,
each  $ \, E_\alpha \, $,  resp.~$ \, F_\alpha \, $,  is a root vector
attached to  $ \, \alpha \in \Phi^+ $,  resp.~to  $ \, -\alpha \in
(-\Phi^+) $,  and the products of factors indexed by  $ \Phi^+ $
are ordered with respect to a fixed convex order of  $ \Phi^+ $,
see [Ga1]), whence (as for  $ \, n = 2 \, $)  $ U_q(\gerg) $  is
a free  $ R $--module.  In this case again  $ \, U_q(\gerg) \, $
is a QrUEA, with semiclassical limit  $ \, U(\gerg) \, $.

 \vskip7pt

  {\bf 5.3 Computation of  $ \, {U_q(\gerg)}' \, $  and
specialization  $ \, {U_q(\gerg)}' \,{\buildrel {q \rightarrow 1}
\over \llongrightarrow}\, F[G^\star] \, $.} \, We begin with the
simplest case  $ \, \gerg = \gersl_2 \, $.  From the definition of 
$ \, U_q(\gerg) = U_q(\gersl_2) \, $  we have  $ \; \delta_n(E) =
{(\text{id} - \epsilon)}^{\otimes n} \big( \Delta^n(E) \big) =
{(\text{id} - \epsilon)}^{\otimes n} \left( \sum\limits_{s=1}^n
K^{\otimes (s-1)} \otimes E \otimes 1^{\otimes (n-s)} \! \right)
%
%
    \hbox{$ = {(q-1)}^{n-1} H^{\otimes (n-1)} \otimes E \; $} $
%
%
 from which  $ \; \delta_n\big( (q-1) E \big) \in {(q-1)}^n
U_q(\gerg) \setminus {(q-1)}^{n+1} U_q(\gerg) \, $   (for all
$ \, n \in \N \, $),  whence  $ \, (q-1) E \in {U_q(\gerg)}' $,
whereas  $ \, E \notin {U_q(\gerg)}' $.  Similarly,  $ \, (q-1)
F \in {U_q(\gerg)}' $,  whilst  $ \, F \notin {U_q(\gerg)}' $.  As
for generators  $ H $,  $ \varGamma $,  $ K^{\pm 1} $,  we have
$ \; \Delta^n(H) = \sum_{s=1}^n K^{\otimes (s-1)} \otimes H \otimes
1^{\otimes (n-s)} $,  $ \, \Delta^n \big( K^{\pm 1} \big) = {\big(
K^{\pm 1} \big)}^{\otimes n} $,  $ \, \Delta^n(\varGamma) =
\sum_{s=1}^n K^{\otimes (s-1)} \otimes \varGamma \otimes {\big(
K^{-1} \big)}^{\otimes (n-s)} $,  \; hence for  $ \, \delta_n
= {(\text{id} - \epsilon)}^{\otimes n} \circ \Delta^n \, $  we
have  $ \; \delta_n(H) = {(q-1})^{n-1} \cdot H^{\otimes n} \, $,
$ \, \delta^n \big( K^{-1} \big) = {(q-1})^n \cdot {(- K^{-1}
H)}^{\otimes n} \, $,  $ \, \delta^n(K) = {(q-1})^n \cdot
H^{\otimes n} \, $,  $ \, \delta^n(\varGamma) = {(q-1})^{n-1}
\cdot \sum\limits_{s=1}^n {(-1)}^{n-s} H^{\otimes (s-1)} \otimes
\varGamma \otimes {\big( H K^{-1} \big)}^{\otimes (n-s)} \, $
%
%
for all  $ \, n \in \N \, $,  so that  $ \, (q-1) H $,  $ (q-1)
\varGamma $,  $ K^{\pm 1} \in {U_q(\gerg)}' \setminus (q-1)
{U_q(\gerg)}' \, $.  Therefore  $ {U_q(\gerg)}' $  contains
the subalgebra  $ U' $  generated by  $ \, (q-1) F $,  $ K $,
$ K^{-1} $,  $ (q-1) H $,  $ (q-1) \varGamma $,  $ (q-1) E \, $.  On
the other hand, using (5.2) a thorough   --- but straightforward ---
computation along the same lines as above shows that any element in
$ {U_q(\gerg)}' $  does necessarily lie in  $ U' $  (details are
left to the reader: everything follows from definitions and the
formulas above for  $ \, \Delta^n \, $).  Thus  $ {U_q(\gerg)}' $
is nothing but the subalgebra of  $ U_q(\gerg) $  generated by
$ \, \dot{F} := (q-1) F $,  $ K $,  $ K^{-1} $,  $ \dot{H} :=
(q-1) H $,  $ \dot{\varGamma} := (q-1) \varGamma $,  $ \dot{E} :=
(q-1) E \, $;  notice also that the generator  $ \, \dot{H} $  is
unnecessary, for  $ \, \dot{H} = K - 1 \, $.  Then  $ {U_q(\gerg)}' $ 
can be presented as the unital associative  $ R $--algebra  with
generators  $ \, \dot{F} \, $,  $ \dot{\varGamma} $,  $ K^{\pm 1} $, 
$ \dot{E} \, $  and relations
  $$  \displaylines{
   K K^{-1} = 1 = K^{-1} K ,  \,\;  K^{\pm 1} \dot{\varGamma} =
\dot{\varGamma} K^{\pm 1} ,  \,\;  \big( 1 + q^{-1} \big)
\dot{\varGamma} = K - K^{-1} ,  \,\;  \dot{E} \dot{F} - \dot{F}
\dot{E} = (q-1) \dot{\varGamma}  \cr
   K - K^{-1} = \big( 1 + q^{-1}\big) \dot{\varGamma} \, ,
\quad  K^{\pm 1} \dot{F} = q^{\mp 2} \dot{F} K^{\pm 1} \, ,
\quad  K^{\pm 1} \dot{E} = q^{\pm 2} \dot{E} K^{\pm 1}  \cr
   \dot{\varGamma} \dot{F} = q^{-2} \dot{F} \dot{\varGamma} -
(q-1) \big( q + q^{-1} \big) \dot{F} \, ,  \qquad  \dot{\varGamma}
\dot{E} = q^{+2} \dot{E} \dot{\varGamma} + (q-1) \big( q + q^{-1} \big)
\dot{E}  \cr }  $$  
with Hopf structure given by
  $$  \begin{matrix}
      \Delta\big(\dot{F}\big) = \dot{F} \otimes K^{-1} + 1 \otimes
\dot{F} \, ,  &  \hskip25pt  \epsilon\big(\dot{F}\big) = 0 \, ,  &
\hskip25pt  S\big(\dot{F}\big) = - \dot{F} K \, \phantom{.}  \\
      \Delta\big(\dot{\varGamma}\big) = \dot{\varGamma} \otimes K +
K^{-1} \otimes \dot{\varGamma} \, ,  &  \hskip25pt  \epsilon\big(
\dot{\varGamma}\big) = 0 \, ,  &  \hskip25pt  S\big(\dot{\varGamma}\big)
= - \dot{\varGamma}  \\
      \Delta \big( K^{\pm 1} \big) = K^{\pm 1} \otimes K^{\pm 1} \, ,
&  \hskip25pt  \epsilon \big( K^{\pm 1} \big) = 1 \, ,  &  \hskip25pt
S \big( K^{\pm 1} \big) = K^{\mp 1} \, \phantom{.}  \\
      \Delta\big(\dot{E}\big) = \dot{E} \otimes 1 + K \otimes \dot{E}
\, ,  &  \hskip25pt  \epsilon\big(\dot{E}\big) = 0 \, ,  &
\hskip25pt  S\big(\dot{E}\big) = - K^{-1} \dot{E} \, .  \\
     \end{matrix}  $$
   \indent   When  $ \, q \rightarrow 1 \, $,  a direct
computation shows that this gives a presentation of 
$ \, F \big[ {}_a{{SL}_2}^{\!*} \big] $,  and the Poisson
structure that $ \, F \big[ {}_a{{SL}_2}^{\!*} \big] \, $  inherits
from this quantization process is exactly the one coming from the
Poisson structure on  $ \, {}_a{{SL}_2}^{\!*} \, $:  in fact,
there is a Poisson Hopf algebra isomorphism
 $$  {U_q(\gerg)}' \Big/ (q-1) \, {U_q(\gerg)}' \,{\buildrel \cong
\over \llongrightarrow}\, F \big[ {}_a{{SL}_2}^{\!*} \big] \quad
\Big( \, \subseteq F \big[ {}_s{{SL}_2}^{\!*} \big] \, \Big)  $$
given by:  $ \;\, \dot{E} \mod (q-1) \mapsto x \, z \, $,
$ \, K^{\pm 1} \mod (q-1) \mapsto z^{\pm 2} \, $,
$ \, \dot{H} \mod (q-1)
                \mapsto z^2 - 1 \, $,\break
$ \, \dot{\varGamma} \mod (q-1) \mapsto \big( z^2 - z^{-2}
\big) \Big/ 2 \, $,  $ \, \dot{F} \mod (q-1) \mapsto z^{-1} y \, $.
In other words, $ \, {U_q(\gerg)}' \, $ specializes to  $ \, F \big[
{}_a{{SL}_2}^{\!*} \big] \, $  as a Poisson Hopf algebra.  {\sl Note
that this was predicted by Theorem 2.2{\it (c)\/}  when  $ \,
\Char(\Bbbk) = 0 \, $,  \, but our analysis now proves it also
for  $ \, \Char(\Bbbk) > 0 \, $.}
                                            \par
  Note that we got the  {\sl adjoint}  Poisson group dual of
$ \, G = {SL}_2 \, $,  that is  $ {}_a{{SL}_2}^{\!*} \, $; a
different choice of the initial QrUEA leads us to the {\sl simply
connected}  one, i.e.~$ \, {}_s{{SL}_2}^{\!*} $.  Indeed, if we
start from the ``simply connected" version of  $ U_q(\gerg) $  (see
\S 5.2) the same analysis shows that  $ {U_q(\gerg)}' $  is like
above but for containing also the new generators  $ L^{\pm 1} $,
and similarly when specializing  $ q $  at  $ 1 $:  thus we get the
function algebra of a Poisson group which is a double covering of
$ {}_a{{SL}_2}^{\!*} $,  namely  $ {}_s{{SL}_2}^{\!*} $.  So changing
the QrUEA quantizing  $ \gerg $  we get two different QFAs, one for each
of the two connected Poisson algebraic groups dual of  $ {SL}_2 $,
i.e.~with tangent Lie bialgebra  $ \, {\gersl_2}^{\!*} \, $;  this
shows the dependence of  $ G^\star $  (here denoted  $ G^* $  since 
$ \, \gerg^\times = \gerg^* \, $)  in  Theorem 2.2{\it (c)}  on the
choice of the QrUEA  $ U_q(\gerg) \, $,  for fixed  $ \gerg \, $.
                                                 \par
   With a bit more careful study, exploiting the analysis in [Ga1],
one can treat the general case too: we sketch briefly our arguments
--- restricting to the simply laced case, to simplify the exposition
---   leaving to the reader the straightforward task of filling in
details.   
                                                 \par
   So now let  $ \, \gerg = \gerg^\tau $  be a semisimple Lie algebra,
as in \S 5.1, and let  $ \, U_q(\gerg) \, $  be the QrUEA introduced in
\S 5.2: our aim again is to compute the QFA  $ \, {U_q(\gerg)}' \, $.
                                               \par
   The same computations as for  $ \, \gerg = {\frak s}{\frak l}(2)
\, $  show that  $ \; \delta_n(H_i) = {(q-1})^{n-1} \cdot
H_i^{\otimes n} \, $  and  $ \, \delta^n(\varGamma_i) = {(q-1})^{n-1}
\cdot \sum_{s=1}^n {(-1)}^{n-s} H_i^{\otimes (s-1)} \otimes
\varGamma_i \otimes {\big( H_i K_i^{-1} \big)}^{\otimes (n-s)} $,
which gives
  $$  \dot{H}_i := (q-1) H_i \in {U_q(\gerg)}' \setminus (q-1) \,
{U_q(\gerg)}'  \quad   \hbox{and}  \quad  \dot{\varGamma}_i :=
(q-1) \, \varGamma_i \in {U_q(\gerg)}' \setminus (q-1) \,
{U_q(\gerg)}' \, .  $$
                                                 \par
   As for root vectors, let  $ \, \dot{E}_\gamma := (q-1) E_\gamma
\, $  and  $ \, \dot{F}_\gamma := (q-1) F_\gamma \, $  for all
$ \, \gamma \in \Phi^+ \, $:  using the same type of arguments
          as in [Ga1]\footnote
          {In [Ga1] one assumes
$ \, \Char(\Bbbk) = 0 \, $:  \, however,  {\sl this is not
necesary\/}  for the present analysis.}
%
%
%
   \S 5.16, we can prove that  $ \, E_\alpha \not\in {U_q(\gerg)}'
\, $  but  $ \, \dot{E}_\alpha \in {U_q(\gerg)}' \setminus (q-1) \,
{U_q(\gerg)}' \, $.  In fact, let  $ \, \Bbb{U}_q({\frak b}_+) \, $  and
$ \, \Bbb{U}_q({\frak b}_-) \, $  be quantum Borel subalgebras, and
$ \, {\frak U}_{\varphi,\geq}^{\scriptscriptstyle M} \, $,
$ \, {\mathcal U}_{\varphi,\geq}^{\scriptscriptstyle M} \, $,
$ \, {\frak U}_{\varphi,\leq}^{\scriptscriptstyle M} \, $,
$ \, {\mathcal U}_{\varphi,\leq}^{\scriptscriptstyle M} \, $
their  $ R $--subalgebras  defined in [Ga1], \S 2: then both  $ \,
\Bbb{U}_q({\frak b}_+) \, $  and  $ \, \Bbb{U}_q({\frak b}_-) \, $  are
Hopf subalgebras  of  $ \Bbb{U}_q(\gerg) \, $.  In addition, letting
$ M' $  be the lattice between  $ Q $  and  $ P $  dual of  $ M $
(in the sense of [Ga1], \S 1.1, there exists an  $ F(R) $--valued
perfect Hopf pairing between  $ \Bbb{U}_q({\frak b}_\pm) $  and
$ \Bbb{U}_q({\frak b}_\mp) $   --- one built up on  $ M $  and
the other on  $ M' $  ---   such that  $ \,
{\frak U}_{\varphi,\geq}^{\scriptscriptstyle M} = {\Big(
{\mathcal U}_{\varphi,\leq}^{\scriptscriptstyle M'} \Big)}^{\!\bullet}
\, $,  $ \, {\frak U}_{\varphi,\leq}^{\scriptscriptstyle M} = {\Big(
{\mathcal U}_{\varphi,\geq}^{\scriptscriptstyle M'} \Big)}^{\!\bullet}
\, $,  $ \, {\mathcal U}_{\varphi,\geq}^{\scriptscriptstyle M} = {\Big(
{\frak U}_{\varphi,\leq}^{\scriptscriptstyle M'} \Big)}^{\!\bullet}
\, $,  and  $ \, {\mathcal U}_{\varphi,\leq}^{\scriptscriptstyle M} =
{\Big( {\frak U}_{\varphi,\geq}^{\scriptscriptstyle M'}
\Big)}^{\!\bullet} \, $.  Now,  $ \, \big( q - q^{-1} \big) E_\alpha
\in {\mathcal U}_{\varphi,\geq}^{\scriptscriptstyle M} = {\Big(
{\frak U}_{\varphi,\leq}^{\scriptscriptstyle M'}
\Big)}^{\!\bullet} \, $,  hence   --- since
$ {\frak U}_{\varphi,\leq}^{\scriptscriptstyle M'} $
is an algebra ---   we have  $ \, \Delta \Big( \big(
q - q^{-1} \big) E_\alpha \Big) \in {\Big(
{\frak U}_{\varphi,\leq}^{\scriptscriptstyle M'} \otimes
{\frak U}_{\varphi,\leq}^{\scriptscriptstyle M'} \Big)}^{\!\bullet}
= {\Big( {\frak U}_{\varphi,\leq}^{\scriptscriptstyle M'}
\Big)}^{\!\bullet} \otimes
{\Big( {\frak U}_{\varphi,\leq}^{\scriptscriptstyle M'}
\Big)}^{\!\bullet} = {\mathcal U}_{\varphi,\geq}^{\scriptscriptstyle M}
\otimes {\mathcal U}_{\varphi,\geq}^{\scriptscriptstyle M} \, $.
Therefore, by definition of
$ {\mathcal U}_{\varphi,\geq}^{\scriptscriptstyle M} $
and by the PBW theorem for it and for
$ {\frak U}_{\varphi,\leq}^{\scriptscriptstyle M'} $
(cf.~[Ga1], \S 2.5) we have that  $ \, \Delta \Big( \big(
q - q^{-1} \big) E_\alpha \Big) \, $  is an  $ R $--linear
combination like  $ \; \Delta \Big( \big( q - q^{-1} \big)
E_\alpha \Big) = \sum_r \, A^{(1)}_r \otimes A^{(2)}_r \; $
in which the  $ A^{(j)}_r $'s  are monomials in the  $ M_j $'s
and in the  $ \, \overline{E}_\gamma $'s,  where  $ \,
\overline{E}_\gamma := \big( q - q^{-1} \big) E_\gamma \, $
for all  $ \, \gamma \in \Phi^+ \, $:  iterating, we find that
$ \, \Delta^\ell \Big( \big( q - q^{-1} \big) E_\alpha \Big) \, $
is an  $ R $--linear  combination
  $$  \Delta^\ell \Big( \big( q - q^{-1} \big) E_\alpha \Big) \; =
\; {\textstyle \sum_r} \, A^{(1)}_r \otimes A^{(2)}_r \otimes \cdots
\otimes A^{(\ell)}_r   \eqno (5.3)  $$   
in which the  $ A^{(j)}_r $'s  are again monomials in the  $ M_j $'s
and in the  $ \, \overline{E}_\gamma $'s.  Now, we distinguish
two cases: either  $ A^{(j)}_r $  does contain some  $ \,
\overline{E}_\gamma \, (\, \in \big( q - q^{-1} \big) \, U_q(\gerg)
\big) \, $,  thus  $ \, \epsilon \Big( A^{(j)}_r \Big) = A^{(j)}_r
\in (q-1) \, U_q(\gerg) \, $  whence  $ \, (\text{id} - \epsilon)
\Big( A^{(j)}_r \Big) = 0 \, $;  or  $ A^{(j)}_r $  does not contain
any  $ \, \overline{E}_\gamma \, $  and is only a monomial in the
$ M_t $'s,  say  $ \, A^{(j)}_r = \prod_{t=1}^n M_t^{m_t} \, $:  then
$ \, (\text{id} - \epsilon) \Big( A^{(j)}_r \Big) = \prod_{t=1}^n
M_t^{m_t} - 1 = \prod_{t=1}^n {\big( (q-1) \, H_t + 1 \big)}^{m_t} -
1 \in (q-1) \, U_q(\gerg) \, $.  In addition, for some  ``$ Q $--grading
reasons" (as in [Ga1], \S 3.16), in each one of the summands in (5.3)
the sum of all the  $ \gamma $'s  such that the (rescaled) root
vectors  $ \overline{E}_\gamma $  occur in any of the factors
$ A^{(1)}_r $,  $ A^{(2)}_r $,  $ \dots $,  $ A^{(n)}_r $  must
be equal to  $ \alpha $:  therefore, in each of these summands at
least one factor  $ \overline{E}_\gamma $  does occur.  The conclusion
is that  $ \; \delta_\ell \big( \overline{E}_\alpha \big) \in \big(
1 + q^{-1} \big) (q-1)^\ell \, {U_q(\gerg)}^{\otimes \ell} \; $
(the factor  $ \, \big( 1 + q^{-1} \big) \, $  being there because at
least one rescaled root vector  $ \overline{E}_\gamma $  occurs in each
summand of  $ \, \delta_\ell \big( \overline{E}_\alpha \big) \, $,
thus providing a coefficient  $ \, \big( q - q^{-1} \big) \, $  the
term  $ \, \big( 1 + q^{-1} \big) \, $  is factored out of), whence
$ \; \delta_\ell \big( \dot{E}_\alpha \big) \in (q-1)^\ell \,
{U_q(\gerg)}^{\otimes \ell} \, $.  More precisely, we have also
$ \; \delta_\ell \big( \dot{E}_\alpha \big) \not\in (q-1)^{\ell+1} \,
{U_q(\gerg)}^{\otimes \ell} \, $,  \; for we can easily check that
$ \, \Delta^\ell \big( \dot{E}_\alpha \big) \, $  is the sum of  $ \,
M_\alpha \otimes M_\alpha \otimes \cdots \otimes M_\alpha \otimes
\dot{E}_\alpha \, $  plus other summands which are  $ R $--linearly
independent of this first term: but then  $ \, \delta_\ell \big(
\dot{E}_\alpha \big) \, $  is the sum of  $ \, {(q-1)}^{\ell-1}
H_\alpha \otimes H_\alpha \otimes \cdots \otimes H_\alpha \otimes
\dot{E}_\alpha \, $  (where  $ \, H_\alpha := {\, M_\alpha - 1 \,
\over \, q - 1 \,} \, $  is equal to an  $ R $--linear  combination
of products of  $ M_j $'s  and  $ H_t $'s)  plus other summands which
are  $ R $--linearly  independent of the first one, and since
$ \, H_\alpha \otimes H_\alpha \otimes \cdots \otimes H_\alpha \otimes
\dot{E}_\alpha \not\in (q-1)^2 \, {U_q(\gerg)}^{\otimes \ell} \, $
we can conclude as claimed.  Therefore  $ \; \delta_\ell \big(
\dot{E}_\alpha \big) \in (q-1)^\ell \, {U_q(\gerg)}^{\otimes \ell}
\setminus (q-1)^{\ell+1} \, {U_q(\gerg)}^{\otimes \ell} \, $,
\; whence we get  $ \; \displaystyle{ \dot{E}_\alpha := (q-1)
E_\alpha \in {U_q(\gerg)}' \setminus (q-1) \, {U_q(\gerg)}' \;\;
\forall \; \alpha \in \Phi^+ } \, $.  \; An entirely similar analysis
yields also  $ \; \displaystyle{ \dot{F}_\alpha := (q-1) F_\alpha
\in {U_q(\gerg)}' \setminus (q-1) \, {U_q(\gerg)}' \;\; \forall \;
\alpha \in \Phi^+ } \, $.
                                                 \par
   Summing up, we have found that  $ \, {U_q(\gerg)}' \, $  contains
the subalgebra  $ U' $  generated by  $ \, \dot{F}_\alpha \, $,  $ \,
\dot{H}_i \, $,  $ \, \dot{\varGamma}_i \, $,  $ \, \dot{E}_\alpha \, $ 
for all  $ \, \alpha \in \Phi^+ \, $  and all  $ \, i = 1, \dots, n \, $. 
On the other hand, using (5.2) a thorough   --- but straightforward ---  
computation along the same lines as above shows that any element in 
$ {U_q(\gerg)}' $  must lie in  $ U' $  (details are left to the reader). 
Thus finally  $ \, {U_q(\gerg)}' = U' \, $,  so we have a concrete
description of  $ {U_q(\gerg)}' $.
                                                \par
   Now compare  $ \, U' = {U_q(\gerg)}' \, $  with the algebra
$ \, {\mathcal U}_\varphi^{\scriptscriptstyle M}(\gerg) \, $  in [Ga1],
\S 3.4 (for  $ \, \varphi = 0 \, $),  the latter being just the
$ R $--subalgebra  of  $ \Bbb{U}_q(\gerg) $  generated by the set
$ \, \big\{\, \overline{F}_\alpha, M_i, \overline{E}_\alpha
\,\big\vert\, \alpha \in \Phi^+, i=1, \dots, n \,\big\} \, $.
First of all, by definition, we have  $ \,
{\mathcal U}_\varphi^{\scriptscriptstyle M}(\gerg)
\subseteq U' = {U_q(\gerg)}' \, $;  moreover,
$ \; \dot{F}_\alpha \equiv {\, 1 \, \over \, 2 \,}
\overline{F}_\alpha \, $,  $ \, \dot{E}_\alpha \equiv
{\, 1 \, \over \, 2 \,} \overline{E}_\alpha \, $,
$ \, \dot{\varGamma}_i \equiv {\, 1 \, \over \, 2 \,}
\big( K_i - K_i^{-1} \big) \, \mod\, (q-1) \, {\mathcal
U}_\varphi^{\scriptscriptstyle M}(\gerg) \; $  for all
$ \, \alpha $,  $ i \, $.  Then
  $$  {\big(U_q(\gerg)'\big)}_1 \; := \; {U_q(\gerg)}' \Big/ (q-1) \,
{U_q(\gerg)}' \; = \; {\mathcal U}_\varphi^{\scriptscriptstyle M}
(\gerg) \Big/ (q-1) \, {\mathcal U}_\varphi^{\scriptscriptstyle M}
(\gerg) \; \cong \; F \big[ G^*_{\scriptscriptstyle M} \big]  $$
where  $ \, G^*_{\scriptscriptstyle M} \, $  is the Poisson
group dual of  $ \, G = G^\tau \, $  with centre  $ \,
Z(G^*_{\scriptscriptstyle M}) \cong M \big/ Q \, $  and
fundamental group  $ \, \pi_1(G^*_{\scriptscriptstyle M})
\cong P \big/ M \, $,  and the isomorphism (of Poisson Hopf
algebras) on the right is given by [Ga1], Theorem 7.4 (see
also references therein for the original statement and proof). 
In other words,  $ \, {U_q(\gerg)}' \, $  specializes to  $ \,
F \big[ G^*_{\scriptscriptstyle M} \big] \, $  {\sl as a Poisson
Hopf algebra},  as prescribed by Theorem 2.2.  By the way, notice
that in the present case the dependence of the dual group  $ \,
G^\star = G^*_{\scriptscriptstyle M} \, $  on the choice of the
initial QrUEA (for fixed  $ \gerg $)   --- mentioned in the last
part of the statement of Theorem 2.2{\it (c)}  ---   is evident.  
                                                 \par
   By the way, the previous discussion applies as well to the case
of  $ \gerg $  {\it an untwisted affine Kac-Moody algebra},  just
replacing quotations from [Ga1]   --- referring to results about 
{\sl finite}  Kac-Moody algebras ---   with similar quotations
from [Ga3]   --- referring to untwisted  {\sl affine}  Kac-Moody
algebras. 

\vskip7pt

  {\bf 5.4 The identity  $ \, {\big({U_q(\gerg)}'\big)}^\vee
= U_q(\gerg) \, $.} \,  In this section we check the part of
Theorem 2.2{\it (b)}  claiming that, when  $ \, p = 0 \, $,  \, one
has  $ \;  H \in \QrUEA \,\Longrightarrow\, {\big( H' \big)}^\vee =
H \; $  for  $ \, H = U_q(\gerg) \, $  as above.  In addition, our
proof now will work for  $ \, p > 0 \, $  as well.  Of course, we
start once again from  $ \, \gerg = \gersl_2 \, $.
                                                     \par
   Since  $ \, \epsilon\big(\dot{F}\big) = \epsilon\big(\dot{H}\big)
= \epsilon\big(\dot{\varGamma}\big) = \epsilon\big(\dot{E}\big) = 0
\, $,  the ideal  $ \, J := \text{\sl Ker}\,\big( \epsilon \, \colon
\, {U_q(\gerg)}' \! \loongrightarrow R \,\big) \, $  is generated by
$ \, \dot{F} $,  $ \dot{H} $,  $ \dot{\varGamma} $,  and  $ \dot{E}
\, $.  This implies that  $ J $  is the  $ R $--span  of  $ \,
\Big\{\, {\dot{F}}^\varphi {\dot{H}}^\kappa {\dot{\varGamma}}^\gamma
{\dot{E}}^\eta \,\Big\vert\, (\varphi, \kappa, \gamma, \eta) \in \N^4
\setminus \{(0,0,0,0)\} \! \Big\} \, $.  
%
%
%
    Therefore  $ \; {\big({U_q(\gerg)}'\big)}^\vee := \sum_{n \geq 0}
  {\Big( {(q-1)}^{-1} J \Big)}^n \; $  
is generated, as a unital  $ R $-subalgebra 
of  $ \Bbb{U}_q(\gerg) \, $,  by the elements  $ \, {(q-1)}^{-1} \dot{F}
= F \, $,  $ \, {(q-1)}^{-1} \dot{H} = H \, $,  $ \, {(q-1)}^{-1}
\dot{\varGamma} = \varGamma \, $,  $ \, {(q-1)}^{-1} \dot{E} = E \, $, 
\, hence it coincides with  $ U_q(\gerg) \, $,  q.e.d.  A similar
analysis works in the ``adjoint" case as well, and also for the
general semisimple or affine Kac-Moody case.

\vskip7pt

  {\bf 5.5 The quantum hyperalgebra  $ \hyp_q(\gerg) $.} \, Let  $ G $
be a semisimple (affine) algebraic group, with Lie algebra  $ \gerg \, $,
and let  $ \Bbb{U}_q(\gerg) $  be the quantum group considered in the
previous sections.  Lusztig introduced (cf.~[Lu1-2]) a ``quantum
hyperalgebra'', i.e.~a Hopf subalgebra of  $ \Bbb{U}_q(\gerg) $
over  $ \Z \big[ q, q^{-1} \big] $  whose specialization at  $ \,
q = 1 \, $  is exactly the Kostant's  $ \Z $--integer  form
$ U_\Z(\gerg) $  of  $ U(\gerg) $  from which one gets the
hyperalgebra  $ \hyp(\gerg) $  over any field  $ \Bbbk $  of
characteristic  $ \, p > 0 \, $  by scalar extension, namely
$ \, \hyp(\gerg) = \Bbbk \otimes_\Z U_\Z(\gerg) \, $.  In fact,
to be precise one needs a suitable enlargement of the algebra
given by Lusztig, which is provided in [DL], \S 3.4, and denoted
by  $ \, \varGamma(\gerg) $.  Now we study Drinfeld's functors
(at  $ \, \h = q - 1 \, $)  on  $ \, \hyp_q(\gerg) := R
\otimes_{\Z[q,q^{-1}]} \varGamma(\gerg) \, $  (with  $ R $
like in \S 5.2), taking as sample the case  $ \, \gerg
= \gersl_2 \, $.
                                             \par
   Let  $ \, \gerg = \gersl_2 \, $.  Let  $ \hyp^\Z_q(\gerg) $  be
the unital  $ \Z\big[q,q^{-1}\big] $--subalgebra  of  $ \Bbb{U}_q(\gerg) $
(say the one of ``adjoint type'' defined like above  {\sl but over\/}
$ \Z\big[q,q^{-1}\big] $)  generated the ``quantum divided powers''
 \vskip-11pt  
  $$  F^{(n)} := F^n \! \Big/ {[n]}_q! \quad ,  \qquad 
\left( {{K \, ; \, c} \atop {n}} \right) \! := \prod_{s=1}^n
{{\; q^{c+1-s} K - 1 \,} \over {\; q^s - 1 \;}} \quad , 
\qquad  E^{(n)} := E^n \! \Big/ {[n]}_q!  $$  
 \vskip-5pt  
(for all  $ \, n \in \N \, $,  $ \, c \in \Z \, $)  and by  $ K^{-1} \, $, 
where  $ \, {[n]}_q! := \prod_{s=1}^n {[s]}_q \, $  and  $ \, {[s]}_q =
\big( q^s - q^{-s} \big) \Big/ \big( q - q^{-1} \big) \, $  for all 
$ \, n $,  $ s \in \N \, $.  Then (cf.~[DL]) this is a Hopf subalgebra
of  $ \Bbb{U}_q(\gerg) $,  and  $ \, \hyp^\Z_q(\gerg){\Big|}_{q=1} \cong
\, U_\Z(\gerg) \, $;  \, therefore  $ \, \hyp_q(\gerg) := R \otimes_{\Z[q,
q^{-1}]} \hyp^\Z_q(\gerg) \, $  (for any  $ R $  like in \S 5.2, with 
$ \, \Bbbk := R \big/ \h \, R \, $  and  $ \, p := \Char(\Bbbk) \, $)
specializes at  $ \, q = 1 \, $  to the  $ \Bbbk $--hyperalgebra
$ \hyp(\gerg) $.  Moreover, among all the  $ \left( {{K ; \, c}
\atop {n}} \right) $'s  it is enough to take only those with
$ \, c = 0 \, $.  {\sl From now on we assume  $ \, p > 0 \, $.}
                                             \par
   Using formulas for the iterated coproduct in [DL], Corollary 3.3
(which uses the opposite coproduct than ours, but this doesn't matter),
and exploiting the PBW-like theorem for  $ \hyp_q(\gerg) $  (see [DL]
again) we see by direct inspection that  $ \, {\hyp_q(\gerg)}' \, $
is the unital  $ R $--subalgebra  of  $ \hyp_q(\gerg) $  generated by
$ K^{-1} $  and the ``rescaled quantum divided powers''  $ \, {(q-1)}^n
F^{(n)} \, $,  $ \, {(q-1)}^n \left( {{K ; \, 0} \atop {n}} \right)
\, $  and  $ \, {(q-1)}^n E^{(n)} \, $  for all  $ \, n \in \N \, $.
Since  $ \, {[n]}_q!{\Big|}_{q=1} \! = n! = 0 \, $  iff  $ \, p \,\Big|
n \, $,  we argue that  $ \, {\hyp_q(\gerg)}'{\Big|}_{q=1} \, $  is
generated by the corresponding specializations of  $ \, {(q-1)}^{p^s}
F^{(p^s)} \, $,  $ \, {(q-1)}^{p^s} \left( {{K ; \, 0} \atop {p^s}}
\right) \, $  and  $ \, {(q-1)}^{p^s} E^{(p^s)} \, $  for all  $ \, s
\! \in \! \N \, $:  in particular this shows 
   \hbox{that the spectrum of  $ \,
{\hyp_q(\gerg)}'{\Big|}_{q=1} \, $  has dimension 0 and height 1, and}
its cotangent Lie algebra  $ \, J \Big/ J^{\,2} \, $   --- where  $ J $
is the augmentation ideal of  $ {\hyp_q(\gerg)}' {\Big|}_{q=1} $  ---
has basis  $ \, \Big\{ {(q\!-\!1)}^{p^s} F^{(p^s)}, \, {(q\!-\!1)}^{p^s}
\! \left( {{K ; \, 0} \atop {p^s}} \right) , \, {(q\!-\!1)}^{p^s} E^{(p^s)}
\, \mod (q\!-\!1) \, {\hyp_q(\gerg)}'\, , \mod J^{\,2} \;\Big|\; s \in \N
\,\Big\} \, $.  Furthermore,  $ \, \big( {\hyp_q(\gerg)}' \big)^{\!\vee}
\, $  is generated by the elements  $ \, {(q-1)}^{p^s-1} F^{(p^s)} \, $, 
$ \, {(q-1)}^{p^s-1} \left( {{K ; \, 0} \atop {p^s}} \right) \, $,  $ \,
K^{-1} \, $  and  $ \, {(q-1)}^{p^s-1} E^{(p^s)} \, $  for all  $ \,
s \in \N \, $:  \, in particular we have that  $ \, \big( {\hyp_q(\gerg)}'
\big)^{\!\vee} \subsetneqq \hyp_q(\gerg) \, $,  \, and  $ \, \big(
{\hyp_q(\gerg)}' \big)^{\!\vee}{\Big|}_{q=1} \, $  is generated by the
cosets modulo  $ (q-1) $  of the previous elements, which do form a
basis of the restricted Lie bialgebra  $ \gerk $  such that  $ \, \big(
{\hyp_q(\gerg)}' \big)^{\!\vee}{\Big|}_{q=1} = \, \u(\gerk) \; $.
                                                \par
   We performed the previous study using the ``adjoint'' version of
$ U_q(\gerg) $  as starting point: instead, we can use as well its
``simply connected'' version, thus obtaining a ``simply connected
version of  $ \hyp_q(\gerg) $''  which is defined like before but
for using  $ L^{\pm 1} $  instead of  $ K^{\pm 1} $;  up to these
changes, the analysis and its outcome will be exactly
the same.  Note that all quantum objects involved   --- namely,
$ \hyp_q(\gerg) $,  $ {\hyp_q(\gerg)}' $  and  $ \big( {\hyp_q(\gerg)}'
\big)^{\!\vee} $  ---   will strictly contain the corresponding
``adjoint'' quantum objects; on the other hand, the semiclassical
limit is the same in the case of  $ \hyp_q(\gerg) $  (giving
$ \hyp(\gerg) $,  in both cases)  and in the case of  $ \big(
{\hyp_q(\gerg)}' \big)^{\!\vee} $  (giving  $ \u(\gerk) $,  in both
cases), whereas the semiclassical limit of  $ {\hyp_q(\gerg)}' $ 
in the ``simply connected'' case is a (countable) covering of the
``adjoint'' one.
                                             \par
   The general case of semisimple or affine Kac-Moody  $ \gerg $
can be dealt with similarly, with analogous outcome.  Indeed,
$ \hyp^\Z_q(\gerg) $  is defined as the unital  $ \Z \big[ q,
q^{-1} \big] $--subalgebra  of  $ \Bbb{U}_q(\gerg) $  (defined like
before  {\sl but over\/}  $ \Z\big[q,q^{-1}\big] $) generated by
$ K_i^{-1} $  and the ``quantum divided powers'' (in the above
sense)  $ \, F_i^{(n)} \, $,  $ \, \left( {{K_i ; \, c} \atop {n}}
\right) \, $,  $ \, E_i^{(n)} \, $  for all  $ \, n \in \N \, $,
$ \, c \in \Z \, $  and  $ \, i =1, \dots, \text{\it rank}\,(\gerg)
\, $  (notation of \S 5.2, but now each divided power relative to
$ i $  is built upon  $ q_i $,  see [Ga1]).  Then (cf.~[DL]) this
is a Hopf subalgebra of  $ \Bbb{U}_q(\gerg) $  with  $ \, \hyp^\Z_q(\gerg)
{\Big|}_{q=1} \cong\, U_\Z(\gerg) \, $,  \, so  $ \, \hyp_q(\gerg)
:= R \otimes_{\Z[q,q^{-1}]} \hyp^\Z_q(\gerg) \, $  (for any
$ R $  like before)  specializes at  $ \, q = 1 \, $  to the
$ \Bbbk $--hyperalgebra  $ \hyp(\gerg) $;  \, and among the
$ \left( {{K_i ; \, c} \atop {n}} \right) $'s  it is enough
to take those with  $ \, c = 0 \, $.
                                             \par
   Again a PBW-like theorem holds for  $ \hyp_q(\gerg) $  (see [DL]),
where powers of root vectors are replaced by quantum divided powers
like  $ \, F_\alpha^{(n)} \, $,  $ \, \left( {{K_i ; \, c} \atop {n}}
\right) \cdot K_i^{-\text{\it Ent}(n/2)} \, $  and  $ \, E_\alpha^{(n)}
\, $,  \, for all positive roots  $ \alpha $  of  $ \gerg $  (each
divided power being relative to  $ q_\alpha $,  see [Ga1])  both in
the finite and in the affine case.  Using this and the same type of
arguments as in \S 5.3   --- i.e.~the perfect graded Hopf pairing
between quantum Borel subalgebras ---   we see by direct inspection
that  $ \, {\hyp_q(\gerg)}' \, $  is the unital  $ R $--subalgebra
of  $ \hyp_q(\gerg) $  generated by the  $ K_i^{-1} $'s  and the
``rescaled quantum divided powers''  $ \, {(q_\alpha - 1)}^n
F_\alpha^{(n)} \, $,  $ \, {(q_i - 1)}^n \left( {{K_i ; \, 0} \atop
{n}} \right) \, $  and  $ \, {(q_\alpha - 1)}^n E_\alpha^{(n)} \, $
for all  $ \, n \in \N \, $.  Since  $ \, {[n]}_{q_\alpha\!}!
{\,\Big|}_{q=1} \! = n! = 0 \, $  iff  $ \, p \,\Big| n \, $, 
one argues like before that  $ \, {\hyp_q(\gerg)}'{\Big|}_{q=1}
\, $  is generated by the corresponding specializations of  $ \,
{(q_\alpha - 1)}^{p^s} F_\alpha^{(p^s)} \, $,  $ \, {(q_i-1)}^{p^s}
\left( {{K_i ; \, 0} \atop {p^s}} \right) \, $  and  $ \, {(q_\alpha
- 1)}^{p^s} E_\alpha^{(p^s)} \, $  for all  $ \, s \in \N \, $  and all
positive roots  $ \alpha \, $.  Again, this shows that the spectrum of 
$ \, {\hyp_q(\gerg)}'{\Big|}_{q=1} \, $  has (dimension 0 and) height 1,
and its cotangent Lie algebra  $ \, J \Big/ J^{\,2} \, $  (where  $ J $
is the augmentation ideal of  $ {\hyp_q(\gerg)}' {\Big|}_{q=1} $)  has
basis  $ \, \Big\{ {(q_\alpha \!\! - \! 1)}^{p^s} \! F_\alpha^{(p^s)},
\, {(q_i \! - \! 1)}^{p^s} \!\! \left( {{K_i ; \, 0} \atop {p^s}}
\right) , \, {(q_\alpha \!\! - \! 1)}^{p^s} E^{(p^s)} \mod (q \! - \! 1)
{\hyp_q(\gerg)}' \mod \! J^{\,2\!} \;\Big|\, s \in \N \Big\} \, $. 
Moreover,  $ \, \big( {\hyp_q(\gerg)}' \big)^{\!\vee} \, $  is
generated by  $ {(q_\alpha - 1)}^{p^s-1} F_\alpha^{(p^s)} $, 
$ {(q_i-1)}^{p^s-1} \! \left( {{K_i ; \, 0} \atop {p^s}} \right) $, 
$ K_i^{-1} $  and  $ {(q_\alpha - 1)}^{p^s-1} E_\alpha^{(p^s)} $  for
all  $ s \, $,  $ i \, $  and  $ \alpha \, $:  \, in particular  $ \,
\big( {\hyp_q(\gerg)}' \big)^{\!\vee} \subsetneqq \hyp_q(\gerg) \, $, 
\, and  $ \, \big( {\hyp_q(\gerg)}' \big)^{\!\vee}{\Big|}_{q=1} \, $ 
is generated by the cosets modulo  $ (q-1) $  of the previous elements,
which in fact do form a basis of the restricted Lie bialgebra  $ \gerk $ 
such that  $ \, \big( {\hyp_q(\gerg)}' \big)^{\!\vee}{\Big|}_{q=1} = \,
\u(\gerk) \; $.

\vskip7pt

  {\bf 5.6 The QFA  $ F_q[G] \, $.} \, In this and the following
sections we pass to look at Theorem 2.2 the other way round: namely,
we start from QFAs and produce QrUEAs.
                                             \par
   We begin with  $ \, G = {SL}_n \, $,  with the standard Poisson
structure,  for which an especially explicit description of the QFA
is available.  Namely, let  $ \, F_q[{SL}_n] \, $  be the unital
associative  \hbox{$R$--alge}bra  generated by \,  $ \{\,
\rho_{ij} \mid i, j = 1, \ldots, n \,\} $  \, with relations
  $$  \begin{array}{lll}
   \rho_{ij} \rho_{ik} = q \, \rho_{ik} \rho_{ij} \; , \quad \quad
\rho_{ik} \rho_{hk} = q \, \rho_{hk} \rho_{ik}  & &  \forall\, j<k,
\, i<h  \quad  \\
   \rho_{il} \rho_{jk} = \rho_{jk} \rho_{il} \; , \quad \quad
\rho_{ik} \rho_{jl} - \rho_{jl} \rho_{ik} = \left( q - q^{-1}
\right) \, \rho_{il} \rho_{jk} & & \forall\, i<j, \, k<l \\
   {det}_q (\rho_{ij}) := {\textstyle \sum_{\sigma \in S_n}}
{(-q)}^{l(\sigma)} \rho_{1,\sigma(1)} \rho_{2,\sigma(2)} \cdots
\rho_{n,\sigma(n)} = 1 \, .  & &  \end{array} $$
  \indent   This is a Hopf algebra, with comultiplication,
counit and antipode given by
  $$  \Delta (\rho_{ij}) = {\textstyle \sum_{k=1}^n} \rho_{ik}
\otimes \rho_{kj} \, ,  \quad  \epsilon(\rho_{ij}) = \delta_{ij} \, ,
\quad  S(\rho_{ij}) = {(-q)}^{i-j} \, {det}_q \! \left(
{(\rho_{hk})}_{h \neq j}^{k \neq i} \right)  $$
for all  $ \, i, j = 1, \dots, n \, $.  Let  $ \, \F_q[{SL}_n] :=
F(R) \otimes_R F_q[{SL}_n] \, $.  The set of ordered monomials
  $$  M \; := \; \bigg\{\; {\textstyle \prod\limits_{i>j}}
\rho_{ij}^{N_{ij}} {\textstyle \prod\limits_{h=k}} \rho_{hk}^{N_{hk}}
{\textstyle \prod\limits_{l<m}} \rho_{lm}^{N_{lm}} \;\bigg\vert\;
N_{st} \in \N  \;\; \forall \; s, t \; ; \; \min \big\{ N_{1,1},
\dots, N_{n,n} \big\} = 0 \;\bigg\}   \eqno (5.4)  $$
is an  $ R $--basis  of  $ F_q[{SL}_n] $  and an  $ F(R) $--basis
of  $ \F_q[{SL}_n] $  (cf.~[Ga2], Theorem 7.4, suitably adapted to 
$ F_q[{SL}_n] \, $).  Moreover,  $ \, F_q[{SL}_n] \, $  is a QFA (at 
$ \, \h = q - 1 \, $),  with  $ \; F_q[{SL}_n] \, {\buildrel \, q
\rightarrow 1 \, \over \llongrightarrow} \, F[{SL}_n] \; $.

\vskip7pt

  {\bf 5.7 Computation of  $ \, {F_q[G]}^\vee \, $  and specialization
$ \, {F_q[G]}^\vee \,{\buildrel {q \rightarrow 1} \over \llongrightarrow}\,
U(\gerg^\times) \, $.} \, In this section we compute  $ \, {F_q[G]}^\vee
\, $  and its semiclassical limit (= specialization at  $ \, q = 1 \, $). 
Note that   
  $$  M' \, := \, \bigg\{\; {\textstyle \prod\limits_{i>j}}
\rho_{ij}^{N_{ij}} {\textstyle \prod\limits_{h=k}} {(\rho_{hk} -
1)}^{N_{hk}} {\textstyle \prod\limits_{l<m}} \rho_{lm}^{N_{lm}}
\;\bigg\vert\; N_{st} \in \N \;\, \forall \; s, t \; ; \; \min
\big\{ N_{1,1}, \dots, N_{n,n} \big\} = 0 \,\bigg\}  $$  
is an  $ R $--basis  of  $ F_q[{SL}_n] $  and an  $ F(R) $--basis
of  $ \F_q[{SL}_n] $;  then, from the definition of the counit, it
follows that  $ \, M' \setminus \{1\} \, $  is an  $ R $--basis  of
$ \, \hbox{\sl Ker}\,\big( \epsilon: F_q[{SL}_n] \longrightarrow R
\,\big) \, $.  Now, by definition  $ \, I := \hbox{\sl Ker} \left(
F_q[{SL}_n] {\buildrel \epsilon \over \llongtwoheadrightarrow} R
{\buildrel {q \mapsto 1} \over \llongtwoheadrightarrow} \Bbbk \right)
\, $,  \, whence  $ \, I = \hbox{\sl Ker}\,(\epsilon) + (q-1) \cdot
F_q[{SL}_n] \, $;  \, therefore  $ \, \big( M' \setminus \{1\}
\big) \cup \big\{ (q-1) \cdot 1 \big\} \, $  is an  $ R $--basis
of  $ I $,  hence  $ \; {(q-1)}^{-1} I \; $  has  $ R $--basis
$ \; {(q-1)}^{-1} \cdot \big( M' \setminus \{1\} \big) \cup
\{1\} \, $.  The outcome is that  $ \, {F_q[{SL}_n]}^\vee :=
\sum_{n \geq 0} {\Big( {(q-1)}^{-1} I \Big)}^n \, $  is just
the unital  $ R $--subalgebra of  $ \F_q[{SL}_n] $  generated
by  $ \; \left\{\, \displaystyle{r_{i{}j} := {\, \rho_{i{}j} -
\delta_{i{}j} \, \over \, q-1 \,}} \;\bigg\vert\; i, j= 1, \dots,
n \,\right\} \, $.  Then one can directly show that this is a Hopf
algebra, and that  $ \, {F_q[{SL}_n]}^\vee \! {\buildrel \, q
\rightarrow 1 \, \over \llongrightarrow} \, U({\gersl_n}^{\!*}) $
as predicted by Theorem 2.2.  Details can be found in [Ga2], \S\S \,
2, 4, looking at the algebra  $ \widetilde{F}_q[{SL}_n] $  considered
therein, up to the following changes.  The algebra which is considered
in  [{\it loc.~cit.}]  has generators  $ \, \displaystyle{ \big( 1 +
q^{-1} \big)^{\delta_{i{}j}} {\, \rho_{i{}j} - \delta_{i{}j} \, \over
\, q - q^{-1} \,}} \, $  ($ \, i, j = 1, \dots, n \, $)  instead of our
$ \, r_{i{}j} \, $'s  (they coincide iff  $ \, i = j \, $)  and also
generators  $ \, \rho_{i{}i} = 1 + (q-1) \, r_{i{}i} \, $  ($ \, i =
1, \dots, n \, $);  then the presentation in \S 2.8 of  [{\it loc.~cit.}]
must be changed accordingly; computing the specialization then goes
exactly the same, and gives the same result   --- specialized
generators are rescaled, though, compared with the standard
ones given in  [{\it loc.~cit.}],  \S 1.
                                             \par
  We sketch the case of  $ \, n = 2 \, $  (see also [FG]).   
                                             \par
   Using notation $ \, \text{a} := \rho_{1,1} \, $,  $ \, \text{b} :=
\rho_{1,2} \, $,  $ \, \text{c} := \rho_{2,1} \, $,  $ \, \text{d} :=
\rho_{2,2} \, $,  we have the relations
 $$  \displaylines{
  \text{a} \, \text{b} = q \, \text{b} \, \text{a} \, ,  \;\hskip27pt
\text{a} \, \text{c} = q \, \text{c} \, \text{a} \, , \;\hskip27pt
\text{b} \, \text{d} = q \, \text{d} \, \text{b} \, ,  \;\hskip27pt
\text{c} \, \text{d} = q \, \text{d} \, \text{c} \, , \cr
  {}  \hskip9pt  \text{b} \, \text{c} = \text{c} \, \text{b} \, ,
\;\hskip41pt  \text{a} \, \text{d} - \text{d} \, \text{a} = \big( q -
q^{-1} \big) \text{b} \, \text{c} \, ,  \;\hskip31pt  \text{a} \, \text{d}
- q \, \text{b} \, \text{c} = 1 \cr }  $$   
holding in  $ F_q[{SL}_2] $  and in  $ \F_q[{SL}_2] $,  with
  $$  \displaylines{
   \Delta(\text{a}) = \text{a} \otimes \text{a} + \text{b} \otimes \text{c}
\, , \; \Delta(\text{b}) = \text{a} \otimes \text{b} + \text{b} \otimes
\text{d} \, , \;  \Delta(\text{c}) = \text{c} \otimes \text{a} + \text{d}
\otimes \text{c} \, , \;  \Delta(\text{d}) = \text{c} \otimes \text{b} +
\text{d} \otimes \text{d}  \cr
   \epsilon(\text{a}) = 1 \, ,  \hskip3pt  \epsilon(\text{b}) = 0 \, ,
\hskip3pt  \epsilon(\text{c}) = 0 \, ,  \hskip3pt  \epsilon(\text{d})
= 1 \, ,  \hskip9,7pt  S(\text{a}) = \text{d} \, ,  \hskip3pt  S(\text{b})
= - q^{-1} \text{b} \, ,  \hskip3pt  S(\text{c}) = - q^{+1} \text{c} \, ,
\hskip1pt  S(\text{d}) = \text{a} \, .  \cr }  $$
Then the elements  $ \, H_+ := r_{1,1} = \displaystyle{\, \text{a} - 1 \,
\over \, q - 1\,} \, $,  $ \, E := r_{1,2} = \displaystyle{\, \text{b} \,
\over \, q - 1\,} \, $, $ \, F := r_{2,1} = \displaystyle{\, \text{c} \,
\over \, q - 1\,} \, $ and $ \, H_- := r_{2,2} = \displaystyle{\, \text{d}
- 1 \, \over \, q - 1\,} \, $  generate  $ {F_q[{SL}_2]}^\vee \, $. 
Moreover, these generators have relations
  $$  \displaylines{
   H_+ E = q \, E H_+ + E \, ,  \;\;  H_+ F = q \, F H_+ + F \, ,  \;\;
E H_- = q \, H_- E + E \, ,  \;\; F H_- = q \, H_- F + F \, ,  \cr
   E F = F E \, ,  \;\;\;  H_+ H_- - H_- H_+ = \big( q - q^{-1} \big) E F
\, , \;\;\;  H_- + H_+ = (q-1) \big( q \, E F - H_+ H_- \big)  \cr }  $$
and Hopf operations given by
  $$  \displaylines{
   \Delta(H_+) = H_+ \otimes 1 + 1 \otimes H_+ + (q-1) \big( H_+ \otimes
H_+ + E \otimes F \big) \, ,  \quad  \epsilon(H_+) = 0 \, ,  \quad  S(H_+)
= H_-  \cr
   {} \, \Delta(E) = E \otimes 1 + 1 \otimes E + (q-1) \big( H_+
\otimes E + E \otimes H_- \big) \, ,  \,\;\;\quad  \epsilon(E) = 0 \, ,
\;\;\quad  S(E) = - q^{-1} E  \cr
   {} \, \Delta(F) = F \otimes 1 + 1 \otimes F + (q-1) \big( F \otimes
H_+ + H_- \otimes F \big) \, ,  \,\;\;\quad  \epsilon(F) = 0 \, ,
\;\;\quad  S(F) = - q^{+1} F  \cr
   \Delta(H_-) = H_- \otimes 1 + 1 \otimes H_- + (q-1) \big( H_- \otimes
H_- + F \otimes E \big) \, ,  \quad  \epsilon(H_-) = 0 \, ,  \quad
S(H_-) = H_+  \cr }  $$
from which one easily checks that  $ \, {F_q[{SL}_2]}^\vee
\,{\buildrel \, q \rightarrow 1 \, \over \llongrightarrow}\,
U({\gersl_2}^{\! *}) \, $  as co-Poisson Hopf algebras, for a
co-Poisson Hopf algebra isomorphism
  $$  {F_q[{SL}_2]}^\vee \Big/ (q-1) \, {F_q[{SL}_2]}^\vee
\;{\buildrel \cong \over \llongrightarrow}\; U({\gersl_2}^{\! *})  $$
exists, given by:  $ \;\, H_\pm \mod (q-1) \mapsto \pm \text{h} \, $,
$ \, E \mod (q-1) \mapsto \text{e} \, $,  $ \, F \mod (q-1) \mapsto
\text{f} \, $;  that is, $ \, {F_q[{SL}_2]}^\vee \, $  specializes
to  $ \, U({\gersl_2}^{\! *}) \, $  {\sl as a co-Poisson Hopf
algebra},  q.e.d.
                                             \par
   Finally, the general case of any semisimple group  $ \, G = G^\tau
\, $,  with the Poisson structure induced from the Lie bialgebra
structure of  $ \, \gerg = \gerg^\tau \, $,  can be treated in a
different way.  Following [Ga1], \S\S 5--6,  $ \F_q[G] $  can be
embedded into a (topological) Hopf algebra  $ \, \Bbb{U}_q(\gerg^*) =
{\Bbb U}_{q,\varphi}^{\scriptscriptstyle M}(\gerg^*) \, $,  so
that the image of the integer form  $ F_q[G] $
lies into a suitable (topological) integer form  $ \,
{\mathcal U}_{q,\varphi}^{\scriptscriptstyle M}(\gerg^*) \, $
of  $ \, \Bbb{U}_q(\gerg^*) \, $.  Now, the analysis given
in  [{\it loc.~cit.\/}],  when carefully read, shows that 
$ \, {F_q[G]}^\vee = \F_q[G] \cap {{\mathcal U}_{q,
\varphi}^{\scriptscriptstyle M}(\gerg^*)}^\vee \, $; 
moreover, the latter (intersection) algebra ``almost''
coincides   --- it is its closure in a suitable topology ---
with the integer form  $ \, {\mathcal F}_q[G] \, $  considered in 
[{\it loc.~cit.\/}]:  in particular, they have the same specialization
at  $ \, q = 1 \, $.  Since in addition  $ \, {\mathcal F}_q[G] \, $
does specialize to  $ U(\gerg^*) $,  the same is true for
$ {F_q[G]}^\vee $,  q.e.d.
                                             \par
   The last point to stress is that, once more, the whole analysis
above is valid for  $ \, p := \Char(\Bbbk) \geq 0 \, $, i.e.~also
for  $ \, p > 0 \, $,  \, which was not granted by Theorem 2.2.

\vskip7pt

  {\bf 5.8 The identity  $ \; {\Big({F_q[G]}^\vee\Big)}' = F_q[G]
\, $.} \, In this section we verify the validity of that part
of Theorem 2.2{\it (b)}  claiming that  $ \; H \in \QFA
\,\Longrightarrow\, {\big(H^\vee\big)}' = H \; $  for
$ \, H = F_q[G] \, $  as above; moreover we show that
this holds for  $ \, p > 0 \, $  too.  We begin with
$ \, G = {SL}_n \, $.
                                                     \par
   From  $ \; \Delta(\rho_{i{}j}) = \sum\limits_{k=1}^n \rho_{i,k}
\otimes \rho_{k,j} \, $,  \; we get  $ \; \Delta^{N} (\rho_{i{}j})
= \!\! \sum\limits_{k_1, \ldots, k_{N-1} = 1}^n \hskip-15pt \rho_{i,k_1}
\otimes \rho_{k_1,k_2} \otimes \cdots \otimes \rho_{k_{N-1},j} \, $,
\, by repeated iteration, whence a simple computation yields
  $$  \delta_N(r_{i{}j}) = {\textstyle \sum\limits_{k_1, \ldots,
k_{N-1} = 1}^n} {(q-1)}^{-1} \cdot \big( (q-1) \, r_{i,k_1} \otimes
(q-1) \, r_{k_1,k_2} \otimes \cdots \otimes (q-1) \, r_{k_{N-1},j}
\big) \qquad \;  \forall \hskip5pt  i, j  $$
%
%
%
%
so that
  $$  \delta_N \big( (q-1) r_{i{}j} \big) \in {(q-1)}^N {F_q[{SL}_n]}^\vee
\setminus {(q-1)}^{N+1} {F_q[{SL}_n]}^\vee  \qquad \hskip11pt \forall
\hskip5pt i, j \, .   \eqno (5.5)  $$
   \indent   Now consider  $ \; M' := \bigg\{\, \prod\limits_{i>j}
\rho_{ij}^{N_{ij}} \prod\limits_{h=k} {(\rho_{hk} -1)}^{N_{hk}}
\prod\limits_{l<m} \rho_{lm}^{N_{lm}} \,\bigg\vert\, N_{st} \in \N
\; \forall \, s, t \, ; \, \min_i \big\{ N_{i,i} \big\} = 0 \,\bigg\}
\, $:  \; since this is an  $ R $--basis  of $ F_q[{SL}_n] \, $,  we
have also that
 $$  M'' \, := \; \left\{\; {\textstyle \prod\limits_{i>j}}
r_{ij}^{N_{ij}} {\textstyle \prod\limits_{h=k}} r_{hk}^{N_{hk}}
{\textstyle \prod\limits_{l<m}} r_{lm}^{N_{lm}} \;\bigg\vert\;
N_{st} \in \N  \;\; \forall \; s, t \, ; \; \min \big\{\, N_{1,1},
\dots, N_{n,n} \,\big\} = 0 \;\right\}  $$
is an  $ R $--basis  of  $ {F_q[{SL}_n]}^\vee $.  This and  $ (5.5) $
above imply that  $ \, {\big( {F_q[{SL}_n]}^\vee \big)}' $  is the
unital  $ R $--subalgebra of  $ \, \F_q[{SL}_n] \, $  generated by
the set  $ \, \big\{\, (q-1) r_{i{}j} \,\big\vert\, i, j= 1, \dots,
n \,\big\} \, $;  since  $ \, (q-1) \, r_{i{}j} = \rho_{i{}j} -
\delta_{i{}j} \, $,  the latter algebra does coincide with  $ \,
F_q[{SL}_n] \, $,  as expected. 
                                           \par
   For the general case of any semisimple group  $ \, G = G^\tau
\, $,  the result can be obtained again by looking at the immersions
$ \, \F_q[G] \subseteq \Bbb{U}_q(\gerg^*) \, $  and  $ \, F_q[G]
\subseteq {\mathcal U}_{q,\varphi}^{\scriptscriptstyle M}(\gerg^*) \, $,
and at the identity  $ \, {F_q[G]}^\vee = \F_q[G] \cap {{\mathcal U}_{q,
\varphi}^{\scriptscriptstyle M}(\gerg^*)}^\vee \, $  (cf.~\S 5.6).  If we
go and compute  $ \, {\Big( {{\mathcal U}_{q,\varphi}^{\scriptscriptstyle
M}(\gerg^*)}^\vee \Big)}' \, $  (noting that  $ \, {\big( {\mathcal U}_{q,
\varphi}^{\scriptscriptstyle M}(\gerg^*) \big)}^\vee \, $  is a QrUEA),
we have to apply the like methods as for  $ \, {U_q(\gerg)}' \, $,
\, thus finding a similar result; this and the identity  $ \;
{F_q[G]}^\vee = \F_q[G] \cap {{\mathcal U}_{q,\varphi}^{\scriptscriptstyle
M}(\gerg^*)}^\vee \; $  eventually yield  $ \; {\Big({F_q[G]}^\vee\Big)}'
= \, F_q[G] \; $.   
                                             \par
   Is is worth pointing out once more that the previous analysis is valid
for  $ \, p := \Char(\Bbbk) \geq 0 \, $,  i.e.~also for  $ \, p > 0 \, $,
\, so the outcome is stronger than what ensured by Theorem 2.2.

\vskip7pt

   {\sl $ \underline{\hbox{\it Remark}} $:}  \; Formula (5.4)
gives an explicit  $ R $--basis  $ M $  of  $ F_q[{SL}_2] $.
By direct computation one sees that  $ \, \delta_n(\mu)
\in {F_q[{SL}_2]}^{\otimes n} \setminus (q-1) \,
{F_q[{SL}_2]}^{\otimes n} $  for all  $ \, \mu \in M \setminus
\{1\} \, $  and  $ \, n \in \N \, $,  whence  $ \, {F_q[{SL}_2]}'
= R \cdot 1 \, $,  which implies  $ \, {\big( {F_q[{SL}_2]}' \big)}_F
= F(R) \cdot 1 \subsetneqq \F_q[{SL}_2] \, $  and also  $ \, {\big(
{F_q[{SL}_2]}' \big)}^\vee = \, R \cdot 1 \, \subsetneqq \, F_q[{SL}_2]
\, $.  This yields a counterexample to part of  Theorem 2.2{\it (b)}.  

\vskip7pt

  {\bf 5.9 Drinfeld's functors and  $ L $--operators  in
$ U_q(\gerg) $  for classical  $ \gerg \, $.} \, Let now
$ \Bbbk $  have zero characteristic,  and let  $ \gerg $  be
a finite dimensional semisimple Lie algebra over  $ \Bbbk $
whose simple Lie subalgebra are all of classical type.  It is
known from [FRT2] that in this case  $ \Bbb{U}^P_q(\gerg) $
(where the subscript  $ P $  means that we are taking a
``simply-connected'' quantum group) admits an alternative
presentation, in which the generators are the so-called
$ L $--operators,  denoted  $ l_{i,j}^{(\varepsilon)} $  with
$ \, \varepsilon = \pm 1 \, $  and  $ i $,  $ j $  ranging
in a suitable set of indices (see [FRT2], \S 2).  Now, if we
consider instead the subalgebra of  $ \Bbb{U}^P_q(\gerg) \, $, 
call it  $ H \, $,  generated by the  $ L $--operators  {\it over 
$ R \, $},  we get at once from the very description of the relations
between the  $ l^{(\varepsilon)}_{i,j} $'s  given in [FRT2] that 
$ H $  is a  {\sl Hopf\/}  $ R $--subalgebra  of  $ \Bbb{U}^P_q
(\gerg) \, $,  and more precisely it is a QFA for the connected
simply-connected dual Poisson group  $ G^\star \, $.
                                          \par
   When computing  $ H^\vee $,  it is generated by the elements
$ \, {(q-1)}^{-1} l^{(\varepsilon)}_{i,j} \, $;  \, even more,
the elements  $ \, {(q-1)}^{-1} l^{(+)}_{i,i+1} \, $  and  $ \,
{(q-1)}^{-1} l^{(-)}_{i+1,i} \, $  are enough to generate.  Now,
Theorem 12 in [FRT2] shows that these latter generators are simply
multiples of the Chevalley generators of  $ U^P_q(\gerg) $  (in
the sense of Jimbo, Drinfeld, etc.), by a coefficient of type
$ \, \pm q^s \big( 1 + q^{-1} \big) \, $  for some  $ \, s \in
\Z \, $:  \, this proves directly that  $ H^\vee $  is a QrUEA
associated to  $ \gerg \, $,  that is the dual Lie bialgebra of
$ G^\star $,  as prescribed by Theorem 2.2.  Conversely, if we
start from  $ U_q^P(\gerg) \, $,  again Theorem 12 of [FRT2]
shows that the  $ {\big( q - q^{-1} \big)}^{-1}
l^{(\varepsilon)}_{i,j} $'s  are quantum root vectors
in  $ U_q^P(\gerg) \, $.  Then when computing  $ {U_q^P(\gerg)}' $
we can shorten a lot the analysis in \S 5.3, because the explicit
expression of the coproduct on the  $ L $--operators  given in
[FRT2]   --- roughly,  $ \Delta $  is given on them by a standard
``matrix coproduct'' ---   tells us directly that all the
$ {\big( 1 + q^{-1} \big)}^{-1} l^{(\varepsilon)}_{i,j} $'s
do belong to  $ {U_q^P(\gerg)}' $,  and again by a PBW argument
we conclude that  $ {U_q^P(\gerg)}' $  is generated by these
rescaled  $ L $--operators,  i.e.~the  $ {\big( 1 + q^{-1}
\big)}^{-1} l^{(\varepsilon)}_{i,j} \; $.
                                          \par
   Therefore, we can say in short that shifting from  $ H $  to
$ H^\vee $  or from  $ U_q^P(\gerg) $  to  $ {U_q^P(\gerg)}' $
essentially amounts   --- up to rescaling by irrelevant factors
(in that they do not vanish at  $ \, q = 1 \, $)  ---   to
switching from the presentation of  $ \Bbb{U}^P_q(\gerg) $
via  $ L $--operators  (after [FRT2]) to the presentation of
Serre-Chevalley type (after Drinfeld and Jimbo), and conversely. 
See also [Ga8] for the cases  $ \, \gerg = \frak{gl}_n \, $  and 
$ \, \gerg = \frak{sl}_n \; $.

\vskip7pt

  {\bf 5.10 The cases  $ U_q(\frak{gl}_n) \, $,  $ F_q[{GL}_n] \, $
and  $ F_q[M_n] \, $.}  In [Ga2], \S 3.2, a certain algebra
$ U_q(\frak{gl}_n) $  is considered as a quantization of
$ \frak{gl}_n \, $;  \, due to their strict relationship, from
the analysis for  $ \gersl_n $  one argues a description of 
$ {U_q(\frak{gl}_n)}' $  and its specialization at  $ \, q = 1 \, $, 
and also verifies that  $ \, \big({U_q(\frak{gl}_n)}'\big)^{\!\vee}
\! = \, U_q(\frak{gl}_n) \, $.   
                                             \par
   Similarly, we can consider the unital associative $R$--algebra  $ \,
F_q[M_n] \, $  with generators  $ \, \rho_{ij} \, $  ($ i $,  $ j = 1,
\ldots, n \, $)  and relations  $ \; \rho_{ij} \rho_{ik} = q \, \rho_{ik}
\rho_{ij} \, $,  $ \, \rho_{ik} \rho_{hk} = q \, \rho_{hk} \rho_{ik} \, $
(for all  $ \, j < k \, $,  $ \, i < h \, $),  $ \; \rho_{il} \rho_{jk} =
\rho_{jk} \rho_{il} \, $,  $ \, \rho_{ik} \rho_{jl} - \rho_{jl} \rho_{ik}
= \left( q - q^{-1} \right) \, \rho_{il} \rho_{jk} \, $  (for all  $ \,
i < j \, $,  $ \, k < l \, $)   --- i.e.~like for  $ {SL}_n \, $,  \,
but for skipping the last relation.  This is the celebrated standard
quantization of  $ F[M_n] $,  the function algebra of the variety
$ M_n $  of  ($ n \times n $)--matrices  over  $ \Bbbk \, $:  \, it
is a  $ \Bbbk $--bialgebra,  whose structure is given by formulas
$ \; \Delta (\rho_{ij}) = \sum_{k=1}^n \rho_{ik} \otimes \rho_{kj} \, $,
$ \, \epsilon(\rho_{ij}) = \delta_{ij} \, $  (for all  $ \, i $,  $ j
= 1, \dots, n \, $)  again, but it is  {\sl not a Hopf algebra}.  The
quantum determinant  $ \; {det}_q (\rho_{ij}) := \sum_{\sigma \in S_n}
{(-q)}^{l(\sigma)} \rho_{1,\sigma(1)} \, \rho_{2,\sigma(2)} \cdots
\rho_{n,\sigma(n)} \; $  is central in  $ F_q[M_n] $,  so by standard
theory we can extend  $ F_q[M_n] $  by adding a formal inverse to
$ \, {det}_q (\rho_{ij}) \, $,  \, thus getting a larger algebra
$ \, F_q[{GL}_n] := F_q[M_n] \big[ {{det}_q(\rho_{ij})}^{-1} \big]
\, $:  \, this is now a Hopf algebra, with antipode  $ \; S(\rho_{ij})
= {(-q)}^{i-j} \, {det}_q \! \left( {(\rho_{hk})}_{h \neq j}^{k \neq i}
\right) \; $  (for all  $ \, i $,  $ j = 1, \dots, n \, $),  the
well-known standard quantization of  $ F[{GL}_n] $,  due to Manin
(see [Ma]).
                                             \par
   Applying Drinfeld's functor  $ (\ )^\vee $  w.r.t.~$ \, \h := (q - 1)
\, $  at  $ F_q[{GL}_n] $  we can repeat step by step the analysis made
for  $ F_q[{SL}_n] $:  \, then  $ {F_q[{GL}_n]}^\vee $  is generated by
the  $ r_{i{}j} $'s  and  $ \, {(q-1)}^{-1} \big( {det}_q(\rho_{ij}) -
1 \big) \, $,  the sole real difference being the lack of the relation 
$ \, {det}_q(\rho_{ij}) = 1 \, $,  \, which implies one relation less
among the  $ r_{i{}j} $'s  inside  $ {F_q[{GL}_n]}^\vee $,  hence also
one relation less among their cosets modulo  $ (q-1) $.  The outcome is
pretty similar, in particular  $ \, {F_q[{GL}_n]}^\vee{\Big|}_{q=1} \!
= \, U({\frak{gl}_n}^{\!\! *}) \, $  (cf.~[Ga2], \S 6.2).  Even more,
we can do the same with  $ F_q[M_n] \, $:  \, things are even easier,
because we have only the  $ r_{i{}j} $'s  alone which generate 
$ {F_q[M_n]}^\vee $,  with no relation coming from the relation  $ \,
{det}_q(\rho_{ij}) = 1 \, $;  \, nevertheless at  $ \, q = 1 \, $  the
relations among the cosets of the  $ r_{i{}j} $'s  are exactly the same
as in the case of  $ {F_q[{GL}_n]}^\vee{\Big|}_{q=1} $,  whence we get 
$ \, {F_q[M_n]}^\vee{\Big|}_{q=1} \! = \, U({\frak{gl}_n}^{\!\! *}) \, $.
{\sl In particular,  $ {F_q[M_n]}^\vee{\Big|}_{q=1} $  is a Hopf algebra,
although both  $ F_q[M_n] $  and  $ {F_q[M_n]}^\vee $  are only
bialgebras,  {\sl not\/}  Hopf algebras: this gives a non-trivial
explicit example of how Theorem 2.2 may be improved.}  The general
result in this sense is Theorem 4.9 in [Ga5].   
                                             \par
   Finally, an analysis of the relationship between Drinfeld functors
and  $ L $--operators  about  $ \Bbb{U}^P_q(\frak{gl}_n) $  can be
done again, exactly like in \S 7.9, leading to entirely similar
results.

\vskip1,1truecm

\centerline {\bf \S \; 6 \  Third example: the three-dimensional
Euclidean group  $ \, E_2 \, $ }

\vskip10pt

  {\bf 6.1  The classical setting.} \, Let  $ \Bbbk $  be any field of
characteristic  $ \, p \geq 0 \, $.  Let  $ \, G := {E}_2(\Bbbk) \equiv
{E}_2 \, $,  the three-dimensional Euclidean group; its tangent Lie
algebra  $ \, \gerg = \gere_2 \, $  is generated by  $ \, f \, $,
$ h \, $,  $ e \, $  with relations  $ \, [h,e] = 2 e \, $,  $ \,
[h,f] = -2 f \, $,  $ \, [e,f] = 0 \, $.  The formulas  $ \, \delta(f)
= h \otimes f - f \otimes h \, $,  $ \, \delta(h) = 0 \, $,  $ \,
\delta(e) = h \otimes e - e \otimes h \, $,  make  $ {\gere}_2 $  into
a Lie bialgebra, hence  $ {E}_2 $  into a Poisson group.  These also give
a presentation of the co-Poisson Hopf algebra  $ U({\gere}_2) $  (with
standard Hopf structure).  If  $ \, p > 0 \, $,  \, we consider on 
$ \gere_2 $  the  $ p $--operation  given by  $ \, e^{[p\,]} = 0 \, $, 
$ \, f^{[p\,]} = 0 \, $,  $ \, h^{[p\,]} = h \, $.
                                             \par
  On the other hand, the function algebra  $ F[{E}_2] $  is the unital
associative commutative  $ \Bbbk $--algebra  with generators  $ \, b $,
$ a^{\pm 1} $,  $ c $,  with Poisson Hopf algebra structure given by
  $$  \displaylines{
   \Delta(b) = b \otimes a^{-1} + a \otimes b \; ,  \hskip29pt
\Delta\big(a^{\pm 1}\big) = a^{\pm 1} \otimes a^{\pm 1} \; ,
\hskip29pt  \Delta(c) = c \otimes a + a^{-1} \otimes c  \cr
   \epsilon(b) = 0  \; ,  \hskip7pt  \epsilon\big(a^{\pm 1}\big) = 1
\; ,  \hskip7pt  \epsilon(c) = 0 \; ,  \hskip39pt  S(b) = -b \; ,
\hskip7pt  S\big(a^{\pm 1}\big) = a^{\mp 1} \; ,  \hskip7pt
S(c) = -c  \cr
   \big\{a^{\pm 1},b\big\} \, = \, \pm \, a^{\pm 1} b \; ,  \hskip33pt
\big\{a^{\pm 1},c\big\} \, = \, \pm \, a^{\pm 1} c \; ,  \hskip33pt
\{b,c\} \, = \, 0  \cr }  $$   
  \indent   We can realize  $ E_2 $  as  $ \, E_2 = \big\{ (b,a,c)
\,\big\vert\, b \, , c \in k \, , a \in \Bbbk \setminus \{0\} \big\}
\, $,  \, with group operation    
  $$  (b_1,a_1,c_1) \cdot (b_2,a_2,c_2) \; = \; \big( b_1 a_2^{-1} +
a_1 b_2 \, , \ a_1 a_2 \, , \, c_1 a_2 + a_1^{-1} c_2 \big) \; ;  $$  
in particular the centre of  $ E_2 $  is simply  $ \, Z := \big\{
(0,1,0), (0,-1,0) \big\} \, $,  so there is only one other connected
Poisson group having  $ \gere_2 $  as Lie bialgebra, namely the
adjoint group  $ \, {}_aE_2 := E_2 \Big/ Z \, $  (the left subscript
$ a $  stands for ``adjoint").  Then  $ F[{}_aE_2] $  coincides
with the Poisson Hopf subalgebra of  $ F[{}_aE_2] $  spanned by
products of an even number of generators, i.e. monomials of even
degree: as a unital subalgebra, this is generated by  $ \, b \,
a \, $,  $ \, a^{\pm 2} \, $,  and  $ \, a^{-1} c \, $.
                                                \par
  The dual Lie bialgebra  $ \, \gerg^* = {\gere_2}^{\!*} \, $  is
the Lie algebra with generators  $ \, \text{f} $,  $ \text{h} $,
$ \text{e} \, $,  and relations  $ \, [\text{h},\text{e}] = 2 \text{e}
\,  $,  $ \, [\text{h},\text{f}\,] = 2 \text{f} \, $,  $ [\text{e},
\text{f}\,] = 0 \, $,  \, with Lie cobracket given by  $ \,
\delta(\text{f}\,) = \text{f} \otimes \text{h} - \text{h} \otimes
\text{f} $,  $ \, \delta(\text{h}) = 0 \, $,  $ \, \delta(\text{e})
= \text{h} \otimes \text{e} - \text{e} \otimes \text{h} \, $  (we choose
as generators  $ \, \text{f} := f^* \, $,  $ \, \text{h} := 2 h^* \, $, 
$ \, \text{e} := e^* \, $,  where  $ \, \big\{ f^*, h^*, e^* \big\} \, $ 
is the basis of  $ {\gere_2}^{\! *} $  which is dual to the basis 
$ \{ f, h, e \} $  of  $ \gere_2 \, $).  If  $ \, p > 0 \, $,  the 
$ p $--operation  of  $ {\gere_2}^{\!*} $  reads  $ \, \text{e}^{[p\,]}
= 0 \, $,  $ \, \text{f}^{\,[p\,]} = 0 \, $,  $ \, \text{h}^{[p\,]} =
\text{h} \, $.  This again gives a presentation of  $ \, U \left(
{{\gere}_2}^{\!*} \right) \, $  too.  The simply connected algebraic
Poisson group with tangent Lie bialgebra  $ {{\gere}_2}^{\!*} $  can
be realized as the group of pairs of matrices   
 $$  {{}_s{E}_2}^{\! \star} \, \equiv \, {{}_s{E}_2}^{\!*} \, := \Bigg\{\,
\bigg( \bigg( \begin{matrix}  z^{-1} & 0 \\  y & z  \end{matrix} \bigg) \,
, \bigg( \begin{matrix}  z & x \\  0 & z^{-1}  \end{matrix} \bigg) \bigg)
\,\Bigg\vert\, x, y \in k, z \in \Bbbk \setminus \{0\} \,\Bigg\} \; ;  $$ 
this group has centre  $ \, Z := \big\{ (I,I), (-I,-I) \big\} \, $,
so there is only one other (Poisson) group with Lie (bi)algebra
$ {{\gere}_2}^{\! *} \, $,  namely the adjoint group  $ \,
{{}_a{E}_2}^{\! *} := {{}_s{E}_2}^{\! *} \Big/ Z \; $.
                                                  \par
  Therefore  $ F\big[{{}_s{E}_2}^{\! *}\big] $  is the unital
associative commutative  $ \Bbbk $--algebra  with generators  $ \, x $,
$ z^{\pm 1} $,  $ y $,  with Poisson Hopf structure given by
  $$  \displaylines{
   \Delta(x) = x \otimes z^{-1} + z \otimes x \; ,  \hskip29pt
\Delta\big(z^{\pm 1}\big) = z^{\pm 1} \otimes z^{\pm 1} \; ,
\hskip29pt  \Delta(y) = y \otimes z^{-1} + z \otimes y  \cr
   \epsilon(x) = 0  \; ,  \hskip7pt  \epsilon\big(z^{\pm 1}\big) = 1
\; ,  \hskip7pt  \epsilon(y) = 0 \; ,  \hskip39pt  S(x) = -x \; ,
\hskip7pt  S\big(z^{\pm 1}\big) = z^{\mp 1} \, ,  \hskip7pt
S(y) = -y  \cr
   \{x,y\} \, = \, 0 \; ,  \hskip33pt  \big\{z^{\pm 1},x\big\} \, =
\, \pm \, z^{\pm 1} x \; ,  \hskip33pt  \big\{z^{\pm 1},y\big\} \, =
\, \mp \, z^{\pm 1} y  \cr }  $$   
(Remark: with respect to this presentation, we have  $ \, \text{f}
= {\partial_y}{\big\vert}_e \, $,  $ \, \text{h} = z \,
{\partial_z}{\big\vert}_e \, $,  $ \, \text{e} =
{\partial_x}{\big\vert}_e \, $,  where  $ e $  is the
identity element of  $ \, {{}_s{E}_2}^{\! *} \, $).  Moreover,
$ F\big[{{}_a{E}_2}^{\! *}\big] $  can be identified with the
Poisson Hopf subalgebra of  $ F\big[{{}_s{E}_2}^{\! *}\big] $
spanned by products of an even number of generators, i.e.{}
monomials of even degree: this is generated, as a unital
subalgebra, by  $ \, x \, z \, $,  $ \, z^{\pm 2} \, $ 
and  $ \, z^{-1} y \, $.

\vskip7pt

  {\bf 6.2 The QrUEAs  $ \, U_q^s(\gere_2) \, $  and  $ \,
U_q^a(\gere_2) \, $.} \, We turn now to quantizations: the
situation is much similar to the case of  $ \gersl_2 \, $, 
so we follow the same pattern; nevertheless, now we stress
a bit more the occurrence of different groups sharing the
same tangent Lie bialgebra.   
                                                  \par
  Let  $ R $  be a 1dD, and let  $ \, \h \in R \setminus \{0\} \, $
and  $ \, q := \h + 1 \in R \, $  be like in \S 5.2.
                                                  \par
   Let $ \, \Bbb{U}_q(\gerg) = \Bbb{U}_q^s({\gere}_2) \, $  (where the
superscript  $ s $  stands for ``simply connected") be the
associative unital  $ F(R) $--algebra  with generators
$ \, F $,  $ L^{\!\pm 1} $,  $ E $,  and relations
  $$  L L^{-1} = 1 = L^{-1} L \; ,  \;\quad  L^{\!\pm 1} F =
q^{\mp 1} F L^{\pm 1} \; ,  \;\quad  L^{\pm 1} E = q^{\pm 1}
E L^{\pm 1} \; ,  \;\quad  E F = F E \; .  $$   
This is a Hopf algebra, with Hopf structure given by
  $$  \displaylines{
     \Delta(F) = F \otimes L^{-2} + 1 \otimes F \; ,  \hskip19pt
\Delta \big( L^{\pm 1} \big) = L^{\pm 1} \otimes L^{\pm 1} \; ,
\hskip19pt  \Delta(E) = E \otimes 1 + L^2 \otimes E  \cr
\epsilon(F) = 0 \; ,  \hskip5pt  \epsilon \big( L^{\pm 1} \big) = 1
\; ,  \hskip5pt  \epsilon(E) = 0 \; ,  \hskip15pt  S(F) = - F L^2
\; ,  \hskip5pt  S \big( L^{\pm 1} \big) = L^{\mp 1} \; , 
\hskip5pt  S(E) = - L^{-2} E \; .  \cr }  $$
Then let  $ \, U_q^s(\gere_2) \, $  be the  $ R $--subalgebra of 
$ \Bbb{U}_q^s(\gere_2) $  generated by  $ \, F $,  $ D_\pm :=
\displaystyle{\, L^{\pm 1} - 1 \, \over \, q - 1 \,} \, $,  $ \,
E \, $.  From the definition of  $ \Bbb{U}_q^s(\gere_2) $  one gets
a presentation of  $ U_q^s(\gere_2) $  as the associative unital
algebra with generators  $ \, F \, $,  $ \, D_\pm \, $,  $ \,
E \, $  and relations
  $$  \displaylines{
   D_+ E = q E D_+ + E \; ,  \quad  F D_+ = q D_+ F + F \; , \quad
E D_- = q D_- E + E \; ,  \quad  D_- F = q F D_- + F  \cr
   E F = F E \; ,  \;\;\qquad  D_+ D_- = D_- D_+ \; ,  \;\;\qquad
D_+ + D_- + (q-1) D_+ D_- = 0  \cr }  $$
with a Hopf structure given by
  $$  \displaylines{
   \Delta(E) = E \otimes 1 + 1 \otimes E + 2(q-1) D_+ \otimes E +
{(q-1)}^2 \cdot D_+^2 \otimes E  \cr
  \Delta(D_\pm) = D_\pm \otimes 1 + 1 \otimes D_\pm + (q-1) \cdot D_\pm
\otimes D_\pm  \cr
  \Delta(F) = F \otimes 1 + 1 \otimes F + 2(q-1) F \otimes D_- +
{(q-1)}^2 \cdot F \otimes D_-^2  \cr  }  $$
 \vskip-25pt
  $$  \begin{matrix}
   \epsilon(E) = 0 \; ,  &  \quad \qquad  S(E) = -E - 2(q-1) D_- E
-{(q-1)}^2 D_-^2 E \; \phantom{.}  \\
   \epsilon(D_\pm) = 0 \; ,  &  \quad \qquad  S(D_\pm) = D_\mp \;
\phantom{.} \\
   \epsilon(F) = 0 \; ,  &  \quad \qquad  S(F) = -F - 2(q-1) F D_+ -
{(q-1)}^2 F D_+^2 \; .  \\
      \end{matrix}  $$
  \indent   The ``adjoint version'' of  $ \Bbb{U}_q^s(\gere_2) $
is the unital subalgebra  $ \Bbb{U}_q^a(\gere_2) $  generated by
$ \, F \, $,  $ \, K^{\pm 1} := L^{\pm 2} \, $,  $ \, E \, $,  which
is clearly a Hopf subalgebra.  It also has an  $ R $--integer  form
$ U_q^a(\gere_2) \, $,  the  unital  $ R $--subalgebra  generated
by  $ \, F \, $,  $ \, H_\pm := \displaystyle{\, K^{\pm 1} - 1 \,
\over \, q - 1 \,} \, $,  $ \, E \, $:  this has relations
  $$  \displaylines{
   E F = F E \; ,  \,\;  H_+ E = q^2 E H_+ + (q+1) E \; ,  \,\;
F H_+ = q^2 H_+ F + (q+1) F \; ,  \,\; H_+ H_- = H_- H_+  \cr
   E H_- = q^2 H_- E + (q+1) E \; ,  \; H_- F = q^2 F H_- +
(q+1) F \; ,  \; H_+ + H_- + (q-1) H_+ H_- = 0  \cr }  $$
and it is a Hopf subalgebra, with Hopf operations given by
  $$  \displaylines{
   \Delta(E) = E \otimes 1 + 1 \otimes E + (q-1) \cdot H_+
\otimes E \; ,  \hskip19pt  \epsilon(E) = 0 \; ,  \hskip19pt 
S(E) = - E - (q-1) H_- E \; \phantom{.}  \cr
   \Delta(H_\pm) = H_\pm \otimes 1 + 1 \otimes H_\pm + (q-1) \cdot H_\pm
\otimes H_\pm \; ,  \hskip17pt  \epsilon(H_\pm) = 0 \; ,  \hskip17pt
S(H_\pm) = H_\mp \;  \hskip25pt \phantom{.}  \cr  
   \Delta(F) = F \otimes 1 + 1 \otimes F + (q-1) \cdot F \otimes H_-
\; ,  \hskip19pt  \epsilon(F) = 0 \; ,  \hskip19pt  S(F) = - F - (q-1)
F H_+ \; .  \cr }  $$   
  \indent   It is easy to check that  $ \, U_q^s(\gere_2) \, $  is a
QrUEA, whose semiclassical limit is $ \, U(\gere_2) \, $:  in fact,
mapping the generators  $ \, F \mod (q-1) $,  $ \, D_\pm \mod (q-1) $,
$ \, E \mod (q-1) \, $  respectively to $ \, f \, $, $ \, \pm h\big/2
\,  $,  $ \, e \in U(\gere_2) $  gives an isomorphism  $ \, U_q^s(\gere_2)
\Big/ (q-1) \, U_q^s(\gere_2) \,{\buildrel \cong \over \longrightarrow}\,
U(\gere_2) \, $  of co-Poisson Hopf algebras.  Similarly, $ \,
U_q^a(\gere_2) \, $  is a QrUEA too, with semiclassical limit
$ \, U(\gere_2) \, $  again: here a co-Poisson Hopf algebra
isomorphism  $ \, U_q^a(\gere_2) \Big/ (q\!-\!1) \, U_q^a(\gere_2)
\, \cong \, U(\gere_2) \, $  is given mapping  $ \, F \! \mod
(q\!-\!1) \, $,  $ \, H_\pm \! \mod (q\!-\!1) \, $,  $ \, E \!
\mod (q\!-\!1) \, $  respectively to  $ \, f \, $,  $ \, \pm h
\, $,  $ \, e \in U(\gere_2) \, $.

\vskip7pt

  {\bf 6.3 Computation of  $ \, {U_q(\gere_2)}' \, $  and
specialization  $ \, {U_q(\gere_2)}' \,{\buildrel {q \rightarrow 1}
\over \llongrightarrow}\, F\big[{E_2}^{\!\star}\big] \, $.} \, This
section is devoted to compute  $ {U_q^s(\gere_2)}' $  and
$ {U_q^a(\gere_2)}' \, $,  and their specialization at  $ \, q =
1 \, $:  everything goes on as in \S 5.3, so we can be more sketchy.
From definitions we have, for any  $ \, n \in \N \, $,  $ \, \Delta^n(E)
= \sum_{s=1}^n K^{\otimes (s-1)} \otimes E \otimes 1^{\otimes (n-s)} $,
\, so  $ \; \delta_n(E) = {(K-1)}^{\otimes (n-1)} \otimes E =
{(q-1)}^{n-1} \cdot H_+^{\otimes (n-1)} \otimes E \, $,  \; whence
$ \; \delta_n\big( (q-1) E \big) \in {(q-1)}^n \, U_q^a(\gere_2)
\setminus {(q-1)}^{n+1} \, U_q^a(\gere_2) \; $  thus  $ \,
(q-1) E \in {U_q^a(\gere_2)}' \, $,  whereas  $ \, E \notin
{U_q^a(\gere_2)}' \, $.  Similarly, we have  $ \, (q-1) F \, $,
$ \, (q-1) H_\pm \in {U_q^a(\gere_2)}' \setminus (q-1) \,
{U_q^a(\gere_2)}' \, $.  Therefore  $ {U_q^a(\gere_2)} $
contains the subalgebra  $ U' $  generated by  $ \, \dot{F}
:= (q-1) F \, $,  $ \, \dot{H}_\pm := (q-1) H_\pm \, $,  $ \,
\dot{E} := (q-1) E \, $.  On the other hand,  $ {U_q^a(\gere_2)}' $ 
is clearly the  $ R $--span  of the set  $ \; \Big\{\, F^a H_+^b H_-^c
E^d \,\Big\vert\, a, b, c, d \in \N \,\Big\} \, $:  \; to be precise,
the set   
  $$  \Big\{\, F^a H_+^b K^{-[b/2]} E^d \,\Big\vert\, a, b, d \in \N
\,\Big\} \; = \; \Big\{\, F^a H_+^b {\big( 1 + (q-1) H_- \big)}^{[b/2]}
E^d \,\Big\vert\, a, b, d \in \N \,\Big\}  $$   
is an  $ R $--basis  of  $ {U_q^a(\gere_2)}' $;  therefore,
a straightforward computation shows that any element in
$ {U_q^a(\gere_2)}' $  does necessarily lie in  $ U' \, $, 
thus  $ {U_q^a(\gere_2)}' $  coincides with  $ U ' \, $. 
Moreover, since  $ \, \dot{H}_\pm = K^{\pm 1} - 1 \, $, 
the unital algebra  $ {U_q^a(\gere_2)}' $  is generated by 
$ \dot{F} \, $,  $ K^{\pm 1} $  and  $ \dot{E} $  as well.
                                       \par
  The previous analysis   --- {\it mutatis mutandis} ---   ensures
also that $ \, {U_q^s(\gere_2)}' \, $  coincides with the unital
$ R $--subalgebra  $ U'' $  of  $ \Bbb{U}_q^s(\gere_2) $  generated by
$ \, \dot{F} := (q-1) F \, $,  $ \, \dot{D}_\pm := (q-1) D_\pm \, $,
$ \, \dot{E} := (q-1) E \, $;  in particular,  $ \, {U_q^s(\gere_2)}'
\supset {U_q^a(\gere_2)}' \, $.  Moreover, as  $ \, \dot{D}_\pm =
L^{\pm 1} - 1 \, $,  the unital algebra  $ {U_q^s(\gere_2)}' $  is
generated by  $ \dot{F} \, $,  $ L^{\pm 1} $  and  $ \dot{E} $  as
well.  Thus  $ {U_q^s(\gere_2)}' $  is the unital associative
$ R $--algebra  with generators  $ \, {\mathcal F} := L \dot{F} \, $,
$ \, {\mathcal L}^{\pm 1} := L^{\pm 1} \, $,  $ \, {\mathcal E} :=
\dot{E} L^{-1} \, $  and relations
  $$  \displaylines{
 {\mathcal L} {\mathcal L}^{-1} = 1 = {\mathcal L}^{-1} {\mathcal L} \; , 
\qquad  {\mathcal E} \, {\mathcal F} = {\mathcal F} \, {\mathcal E} \; , 
\qquad  {\mathcal L}^{\pm 1} {\mathcal F} = q^{\mp 1} {\mathcal F}
{\mathcal L}^{\pm 1} \; ,  \qquad  {\mathcal L}^{\pm 1} {\mathcal E}
= q^{\pm 1} {\mathcal E} {\mathcal L}^{\pm 1}  \cr }  $$   
with Hopf structure given by
  $$  \displaylines{
   \Delta({\mathcal F}) = {\mathcal F} \otimes {\mathcal L}^{-1}
+ {\mathcal L} \otimes {\mathcal F} \; ,  \hskip19pt 
 \Delta \big( {\mathcal L}^{\pm 1} \big) = {\mathcal L}^{\pm 1}
\otimes {\mathcal L}^{\pm 1} \; ,  \hskip19pt 
 \Delta({\mathcal E}) = {\mathcal E} \otimes {\mathcal L}^{-1}
+ {\mathcal L} \otimes {\mathcal E}  \cr
   \epsilon({\mathcal F}) = 0 \; ,  \hskip7pt  \epsilon \big(
{\mathcal L}^{\pm 1} \big) = 1 \; ,  \hskip7pt  \epsilon({\mathcal E}) = 0
\; ,  \hskip29pt  S({\mathcal F}) = - {\mathcal F} \; ,  \hskip7pt
S \big( {\mathcal L}^{\pm 1} \big) = {\mathcal L}^{\mp 1} \; , 
\hskip7pt S({\mathcal E}) = - {\mathcal E} \, .  \cr }  $$
   As  $ \, q \rightarrow 1 \, $,  this yields a presentation of the
function algebra  $ \, F \big[ {}_s{E_2}^{\!*} \big] $,  and the
Poisson bracket that  $ \, F \big[ {}_s{E_2}^{\!*} \big] \, $
earns from this quantization process coincides with the one
coming from the Poisson structure on  $ \, {}_s{E_2}^{\!*} \, $:
namely, there is a Poisson Hopf algebra isomorphism
 $$  {U_q^s(\gere_2)}' \Big/ (q-1) \, {U_q^s(\gere_2)}'
\,{\buildrel \cong \over \llongrightarrow}\, F \big[
{}_s{E_2}^{\!*} \big]  $$  
given by  $ \;\, {\mathcal E} \mod (q-1) \mapsto x \, $,
$ \, {\mathcal L}^{\pm 1} \mod (q-1) \mapsto z^{\pm 1} \, $,
$ {\mathcal F} \mod (q-1) \mapsto y \, $.  That is, $ \,
{U_q^s(\gere_2)}' \, $  specializes to  $ \, F \big[
{}_s{E_2}^{\!*} \big] \, $  {\sl as a Poisson Hopf
algebra},  as predicted by Theorem 2.2.
                                            \par
  In the ``adjoint case", from the definition of  $ U' $  and
from  $ \, {U_q^a(\gere_2)}' = U' \, $  we find that
$ {U_q^a(\gere_2)}' $  is the unital associative
$ R $--algebra  with generators  $ \, \dot{F} $,
$ K^{\pm 1} $,  $ \dot{E} \, $  and relations
  $$  K K^{-1} = 1 = K^{-1} K \; ,  \;\quad  \dot{E} \dot{F} =
\dot{F} \dot{E} \; ,  \;\quad  K^{\pm 1} \dot{F} = q^{\mp 2} \dot{F}
K^{\pm 1} \; ,  \;\quad  K^{\pm 1} \dot{E} = q^{\pm 2} \dot{E}
K^{\pm 1}  $$   
with Hopf structure given by
  $$  \displaylines{
   \Delta\big(\dot{F}\big) = \dot{F} \otimes K^{-1} + 1 \otimes
\dot{F} \; ,  \hskip19pt  \Delta \big( K^{\pm 1} \big) = K^{\pm 1}
\otimes K^{\pm 1} \; ,  \hskip19pt  \Delta\big(\dot{E}\big) = \dot{E}
\otimes 1 + K \otimes \dot{E}  \cr
   \epsilon\big(\dot{F}\big) = 0  \; ,  \hskip3pt
\epsilon \big( K^{\pm 1} \big) = 1 \; ,  \hskip3pt
\epsilon\big(\dot{E}\big) = 0 \; ,  \hskip11pt
S\big(\dot{F}\big) = - \dot{F} K \; ,  \hskip3pt
S \big( K^{\pm 1} \big) = K^{\mp 1} \; ,  \hskip3pt
S\big(\dot{E}\big) = - K^{-1} \dot{E} \; .  \cr }  $$
   \indent   The outcome is that there is a Poisson Hopf algebra
isomorphism   
  $$  {U_q^a(\gere_2)}' \Big/ (q-1) \, {U_q^a(\gere_2)}' \,{\buildrel
\cong \over \llongrightarrow}\, F \big[ {}_a{E_2}^{\!*} \big]  \quad
\Big( \subset F \big[ {}_s{E_2}^{\!*} \big] \Big) $$
given by  $ \;\, \dot{E} \mod (q-1) \mapsto x \, z \, $,  $ \, K^{\pm 1}
\mod (q-1) \mapsto z^{\pm 2} \, $,  $ \dot{F} \mod (q-1) \mapsto z^{-1}
y \, $,  which means  $ \, {U_q^a(\gere_2)}' \, $  specializes to  $ \,
F \big[ {}_a{E_2}^{\!*} \big] \, $  {\sl as a Poisson Hopf algebra}, 
according to Theorem 2.2.
                                            \par
   To finish with, note that  {\sl all this analysis (and its
outcome) is entirely characteristic-free}.

\vskip7pt

  {\bf 6.4 The identity  $ \, {\big({U_q(\gere_2)}'\big)}^{\!\vee}
= U_q(\gere_2) \, $.} \,  This section goal is to check the part
of Theorem 2.2{\it (b)}  claiming that  $ \; H \in \QrUEA
\,\Longrightarrow\, {\big(H'\big)}^{\!\vee} = H \; $  both for
$ \, H = U_q^s(\gere_2) \, $  and  $ \, H = U_q^a(\gere_2) \, $.
{\sl In addition, our analysis work for all  $ \, p:= \Char(\Bbbk)
\, $,  thus giving a stronger result than Theorem 2.2{\it (b)}}.
                                            \par   
   First,  $ \, {U_q^s(\gere_2)}' $  is clearly a free
$ R $--module,  with basis  $ \, \Big\{\, {\mathcal F}^a
{\mathcal L}^d {\mathcal E}^c \,\Big\vert\, a, c \in \N, d \in \Z
\,\Big\} \, $,  hence  $ \, {\Bbb B} := \Big\{\, {\mathcal F}^a
{({\mathcal L}^{\pm 1} - 1)}^b {\mathcal E}^c \,\Big\vert\, a, b,
c \in \N \,\Big\} \, $,  is an  $ R $--basis  as well.  Second,
since  $ \, \epsilon({\mathcal F}) = \epsilon\big({\mathcal L}^{\pm 1}
- 1\big) = \epsilon({\mathcal E}) = 0 \, $,  the ideal  $ \, J :=
\hbox{\sl Ker} \Big( \epsilon \, \colon \, {U_q^s(\gere_2)}'
\longrightarrow R \Big) \, $  is the span of  $ \, {\Bbb B}
\setminus \{1\} $.   
%
%
    Therefore  $ \; {\big({U_q^s(\gere_2)}'\big)}^{\!\vee} =
  \sum_{n \geq 0} {\Big( {(q-1)}^{-1} J \Big)}^n \; $  is generated  
%
%
by  $ \,
{(q-1)}^{-1} {\mathcal F} = L F \, $,  $ \, {(q-1)}^{-1} ({\mathcal L}
- 1) = D_+ \, $,  $ \, {(q-1)}^{-1} \big( {\mathcal L}^{-1} - 1 \big)
= D_- \, $,  $ \, {(q-1)}^{-1} {\mathcal E} = E L^{-1} \, $,  hence
by  $ \, F \, $,  $ \, D_\pm \, $,  $ \, E \, $,  so it coincides
with  $ U_q^s(\gere_2) \, $.   
                                                 \par
  The situation is entirely similar for the adjoint case: one simply
has to change  $ \, {\mathcal F} \, $,  $ \, {\mathcal L}^{\pm 1} \, $, 
$ \, {\mathcal E} \, $  respectively with  $ \, \dot{F} \, $,  $ \,
K^{\pm 1} \, $,  $ \, \dot{E} \, $,  and  $ \, D_\pm \, $  with 
$ \, H_\pm \, $,  then everything goes through as above.

\vskip7pt

  {\bf 6.5 The quantum hyperalgebra  $ \hyp_q(\gere_2) $.} \, Like
for semisimple groups, we can define ``quantum hyperalgebras'' for 
$ \gere_2 $  mimicking what done in \S 5.5.  Namely, we can first
define a Hopf  $ \Z \big[ q, q^{-1} \big] $--subalgebra  of 
$ \Bbb{U}_q^s(\gere_2) $  whose specialization at  $ \, q = 1 \, $ 
is the Kostant-like  $ \Z $--integer  form  $ U_\Z(\gere_2) $  of 
$ U(\gere_2) $  (generated by divided powers, and giving the
hyperalgebra  $ \hyp(\gere_2) $  over any field  $ \Bbbk $  by
scalar extension, namely  $ \, \hyp(\gere_2) = \Bbbk \otimes_\Z
U_\Z(\gere_2) \, $),  \, and then take its scalar extension over
$ R \, $.
                                             \par
   To be precise, let  $ \hyp^{s,\Z}_q(\gere_2) $  be the unital
$ \Z\big[q,q^{-1}\big] $--subalgebra  of  $ \Bbb{U}^s_q(\gere_2) $
(defined like above  {\sl but over\/}  $ \Z\big[q,q^{-1}\big] $)
generated by the ``quantum divided powers''  
  $$  F^{(n)} \! := F^n \! \Big/ {[n]}_q! \; ,  \quad 
\left( {{L \, ; \, c} \atop {n}} \right) \! := \prod_{r=1}^n {{\;
q^{c+1-r} L - 1 \,} \over {\; q^r - 1 \;}} \; ,  \quad 
E^{(n)} \! := E^n \! \Big/ {[n]}_q!  $$  
(for all  $ \, n \in \N \, $  and  $ \, c \in \Z \, $,  with notation
of \S 5.5)  and by  $ L^{-1} \, $.  Comparing with the case of 
$ \gersl_2 $  one easily sees that this is a Hopf subalgebra of 
$ \Bbb{U}_q^s(\gere_2) $,  and  $ \, \hyp^{s,\Z}_q(\gere_2)
{\Big|}_{q=1} \cong\, U_\Z(\gere_2) \, $;  \, thus  $ \,
\hyp^s_q(\gere_2) := R \otimes_{\Z[q,q^{-1}]} \hyp^{s,\Z}_q(\gere_2)
\, $  (for any  $ R $  like in \S 6.2, with  $ \, \Bbbk := R \big/
\h \, R \, $  and  $ \, p := \Char(\Bbbk) \, $)  specializes at  $ \,
q = 1 \, $  to the  $ \Bbbk $--hyperalgebra  $ \hyp(\gere_2) $.  In
addition, among all the  $ \left( {{L \, ; \, c} \atop {n}} \right) $'s
it is enough to take only those with  $ \, c = 0 \, $.  {\sl From now
on we assume  $ \, p > 0 \, $.}
                                             \par
   Again a strict comparison with the  $ \gersl_2 $  case shows us
that  $ \, {\hyp^s_q(\gere_2)}' \, $  is the unital  $ R $--subalgebra
of  $ \hyp^s_q(\gere_2) $  generated by  $ L^{-1} $  and the ``rescaled
quantum divided powers''  $ \, {(q \! - \! 1)}^n F^{(n)} \, $,  $ \,
{(q\!-\!1)}^n \! \left( {{L \, ; \, 0} \atop {n}} \right) \, $  and
$ \, {(q\!-\!1)}^n E^{(n)} \, $  for all  $ \, n \in \N \, $.  It
follows that  $ \, {\hyp^s_q(\gere_2)}'{\Big|}_{q=1} \, $  is
generated by the corresponding specializations of  $ \, {(q-1)}^{p^r}
F^{(p^r)} \, $,  $ \, {(q-1)}^{p^r} \! \left( {{L \, ; \, 0} \atop
{p^r}} \right) \, $  and  $ \, {(q-1)}^{p^r} E^{(p^r)} \, $  for
all  $ \, r \in \N \, $:  \, this proves that the spectrum of  $ \,
{\hyp^s_q(\gere_2)}'{\Big|}_{q=1} \, $  has dimension 0 and height 1,
and its cotangent Lie algebra has basis  $ \, \Big\{\, {(q\!-\!1)}^{p^r}
F^{(p^r)}, \, {(q\!-\!1)}^{p^r} \! \left( {{L \, ; \, 0} \atop {p^r}}
\right) \! , \, {(q\!-\!1)}^{p^r} E^{(p^r)} \, \mod (q\!-\!1) \,
{\hyp^s_q(\gerg)}' \, \mod J^{\,2} \;\Big|\; r \!\in\! \N \,\Big\} \, $ 
(where  $ J $  is the augmentation ideal of  $ {\hyp^s_q(\gere_2)}'
{\Big|}_{q=1} \, $,  so that  $ \, J \Big/ J^{\,2} \, $  is the
aforementioned cotangent Lie bialgebra).  Moreover,  $ \, \big(
{\hyp^s_q(\gere_2)}' \big)^{\!\vee} \, $  is generated by  $ \,
{(q-1)}^{p^r-1} F^{(p^r)} \, $,  $ \, {(q-1)}^{p^r-1} \left(
{{L \, ; \, 0} \atop {p^r}} \right) \, $,  $ L^{-1} $  and  $ \,
{(q-1)}^{p^r-1} E^{(p^r)} \, $  (for all  $ \, r \in \N \, $): 
\, in particular  $ \, \big( {\hyp^s_q(\gere_2)}' \big)^{\!\vee}
\subsetneqq \hyp^s_q(\gere_2) \, $,  \, and finally  $ \, \big(
{\hyp^s_q(\gere_2)}' \big)^{\!\vee}{\Big|}_{q=1} \, $  is generated
by the cosets modulo  $ (q-1) $  of the elements above, which in fact
form a basis of the restricted Lie bialgebra  $ \gerk $  such that 
$ \, \big( {\hyp^s_q(\gere_2)}' \big)^{\!\vee} {\Big|}_{q=1} =
\, \u(\gerk) \; $.   
                                             \par
   All this analysis was made starting from  $ \Bbb{U}_q^s(\gere_2) \, $, 
which gave ``simply connected quantum objects''.  If we start instead
from  $ \Bbb{U}_q^a(\gere_2) \, $,  we get ``adjoint quantum objects''
following the same pattern but for replacing everywhere  $ L^{\pm 1} $ 
by  $ K^{\pm 1} \, $:  apart from these changes, the analysis and its
outcome will be exactly the same.  Like for  $ \, \gersl_2 $  (cf.~\S
5.5), all the adjoint quantum objects   ---  i.e.~$ \hyp^a_q(\gere_2)
\, $,  $ {\hyp^a_q(\gere_2)}' $  and  $ \big( {\hyp^a_q(\gere_2)}'
\big)^{\!\vee} $  ---   will be strictly contained in the corresponding
simply connected quantum objects; nevertheless, the semiclassical limits
will be the same in the case of  $ \hyp_q(\gere_2) $  (always yielding
$ \hyp(\gere_2) \, $)  and in the case of  $ \big( {\hyp_q(\gere_2)}'
\big)^{\!\vee} $  (giving  $ \u(\gerk) $,  in both cases), while the
semiclassical limit of  $ {\hyp_q(\gere_2)}' $  in the simply connected
case will be a (countable) covering of that in the adjoint case.

\vskip7pt

  {\bf 6.6 The QFAs  $ \, F_q[E_2] \, $  and  $ \, F_q[{}_aE_2] \, $.}
\, In this and the following sections we look at Theorem 2.2 starting
from QFAs, to get QrUEAs out of them.
                                            \par
   We begin by introducing a QFA for the Euclidean groups  $ E_2 $ 
and  $ {}_aE_2 \, $.  Let  $ \, F_q[E_2] \, $  be the unital associative
\hbox{$R$--alge}bra  with generators  $ \, \text{a}^{\pm 1} $, 
$ \text{b} $,  $ \text{c} \, $  and relations  
 \vskip-11pt  
  $$  \text{a} \, \text{b} = q \, \text{b} \, \text{a} \, ,  \qquad\qquad
\text{a} \, \text{c} = q \, \text{c} \, \text{a} \, ,  \qquad\qquad
\text{b} \, \text{c} = \text{c} \, \text{b}  $$  
endowed with the Hopf algebra structure given by
  $$  \displaylines{
   \Delta\big(\text{a}^{\pm 1}\big) = \text{a}^{\pm 1} \otimes
\text{a}^{\pm 1} \, ,  \hskip19pt    \Delta(\text{b}) = \text{b}
\otimes \text{a}^{-1} + \text{a} \otimes \text{b} \, ,  \hskip19pt 
\Delta(\text{c}) = \text{c} \otimes \text{a} + \text{a}^{-1}
\otimes \text{c}  \cr
   \epsilon\big(\text{a}^{\pm 1}\big) = 1 \, ,  \hskip5pt 
\epsilon(\text{b}) = 0 \, ,  \hskip5pt  \epsilon(\text{c}) = 0 \, , 
\hskip21pt  S\big(\text{a}^{\pm 1}\big) = \text{a}^{\mp 1} \, ,
\hskip5pt  S(\text{b}) = - q^{-1} \, \text{b} \, ,  \hskip5pt 
S(\text{c}) = - q^{+1} \, \text{c} \, .  \cr }  $$  
   \indent   Define  $ \, F_q[{}_aE_2] \, $  as the
$ R $--submodule  of  $ F_q[E_2] $  spanned by the products
of an even number of generators, i.e.~monomials of even degree
in  $ \text{a}^{\pm 1} \, $,  $ \text{b} \, $,  $ \text{c} \, $: 
this is a unital subalgebra of  $ F_q[E_2] \, $,  generated by 
$ \, \beta := \text{b} \, \text{a} \, $,  $ \, \alpha^{\pm 1} :=
\text{a}^{\pm 2} \, $,  and  $ \, \gamma := \text{a}^{-1} \text{c}
\, $.  Let also  $ \, \F_q[E_2] := {\big(F_q[E_2]\big)}_F \, $  and 
$ \, \F_q[{}_aE_2] := {\big(F_q[{}_aE_2]\big)}_F \, $,  \, having the
same presentation than  $ F_q[E_2] $  and  $ F_q[{}_aE_2] $  but over 
$ F(R) \, $.  By construction  $ \, F_q[E_2] \, $  and  $ \, F_q[{}_aE_2]
\, $  are QFAs (at  $ \, \h = q - 1 \, $),  with semiclassical limit 
$ F[E_2] $  and  $ F[{}_aE_2] $  respectively.  

\vskip7pt

  {\bf 6.7 Computation of  $ {F_q[E_2]}^\vee $  and
$ {F_q[{}_aE_2]}^\vee $  and specializations  ${F_q[E_2]}^\vee 
\!\!{\buildrel {q \rightarrow 1} \over \llongrightarrow}\,
U(\gerg^\times) $  and  \break $ \, {F_q[{}_aE_2]}^\vee \!\!{\buildrel
{q \rightarrow 1} \over \llongrightarrow}\, U(\gerg^\times) \, $.}
\,  In this section we go and compute  $ \, {F_q[G]}^\vee \, $  and
its semiclassical limit (i.e.~its specialization at  $ \, q = 1
\, $),  both for  $ \, G = E_2 \, $  and  $ \, G = {}_aE_2 \, $.
                                             \par
   First,  $ F_q[E_2] $  is free over  $ R \, $,  with basis  $ \,
\Big\{\, \text{b}^b \text{a}^a \text{c}^c \,\Big\vert\, a \in \Z,
b, c \in \N \,\Big\} \, $,  and so also the set  $ \, {\Bbb B}_s
:= \Big\{\, \text{b}^b {(\text{a}^{\pm 1} - 1)}^a \text{c}^c
\,\Big\vert\, a, b, c \in \N \,\Big\} \, $  is an  $ R $--basis. 
Second, as  $ \, \epsilon(\text{b}) =   \hbox{$ \epsilon \big(
\text{a}^{\pm 1} - 1 \big) $}   = \epsilon(\text{c}) = 0 \, $, 
the ideal  $ \, J := \text{\sl Ker} \, \Big( \epsilon \, \colon
\, F_q[E_2] \loongrightarrow R \Big) \, $  is the span of  $ \,
{\Bbb B}_s \setminus \{1\} \, $.    
%
%
    Then  $ \; {F_q[E_2]}^\vee = \sum_{n \geq 0} {\Big( \!
{(q\!-\!1)}^{-1} J \Big)}^n \, $  
%
%
 is
the unital  $ R $--algebra
%
%
with generators  $ \, D_\pm :=
\displaystyle{\, \text{a}^{\pm 1} - 1 \, \over \, q - 1\,} \, $,
$ \, E := \displaystyle{\, \text{b} \, \over \, q - 1\,} \, $,
and  $ \, F := \displaystyle{\, \text{c} \, \over \, q - 1\,} \, $
and relations
  $$  \displaylines{
   D_+ E = q E D_+ + E \, ,  \quad  D_+ F = q F D_+ + F \, , \quad
E D_- = q D_- E + E \, ,  \quad  F D_- = q D_- F + F  \cr
   E F = F E \, ,  \;\qquad  D_+ D_- = D_- D_+ \, ,  \;\qquad
D_+ + D_- + (q-1) D_+ D_- = 0  \cr }  $$
with a Hopf structure given by
  $$  \Delta(E) = E \otimes 1 + 1 \otimes E + (q-1) \big( E \otimes D_-
+ D_+ \otimes E \big) \; ,  \quad  \epsilon(E) = 0 \; ,  \quad  S(E) =
- q^{-1} E  $$
  $$   \begin{matrix}
%
%
   \Delta(D_\pm) = D_\pm \otimes 1 + 1 \otimes D_\pm +
(q-1) \cdot D_\pm \otimes D_\pm \; ,  &  \;  \epsilon(D_\pm) = 0 \; ,
&  \;  S(D_\pm) = D_\mp \; \phantom{.}  \\
   \Delta(F) = F \otimes 1 + 1 \otimes F + (q-1) \big( F \otimes D_+
+ D_- \otimes F \big) \; ,  &  \;  \epsilon(F) = 0 \; ,  &  \;
S(F) = - q^{+1} F \; .  \\
     \end{matrix}  $$
This implies that  $ \, {F_q[E_2]}^\vee \,{\buildrel \, q
\rightarrow 1 \, \over \llongrightarrow}\, U({\gere_2}^{\! *}) \, $ 
as co-Poisson Hopf algebras, for an isomorphism   
  $$  {F_q[E_2]}^\vee \Big/ (q-1) \, {F_q[E_2]}^\vee
\,{\buildrel \cong \over \llongrightarrow}\; U({\gere_2}^{\! *})  $$
of co-Poisson Hopf algebra exists, given by  $ \;\, D_\pm \mod
(q-1) \mapsto \pm \, \text{h} \big/ 2 \, $,  $ \, E \mod (q-1) \mapsto
\text{e} \, $,  \, and  $ \, F \mod (q-1) \mapsto \text{f} \, $;  so 
$ \, {F_q[E_2]}^\vee \, $  does specialize to  $ \, U({\gere_2}^{\! *})
\, $ {\sl as a co-Poisson Hopf algebra},  q.e.d.
                                             \par
  Similarly, if we consider  $ F_q[{}_aE_2] $  the same
analysis works again.  In fact,  $ F_q[{}_aE_2] $  is free
over $ R $,  with basis  $ \, {\Bbb B}_a :=  \Big\{\, \beta^b
{(\alpha^{\pm 1} - 1)}^a \gamma^c \,\Big\vert\, a, b, c \in \N
\,\Big\} \, $;  then as above $ \, J := \text{\sl Ker} \, \Big(
\epsilon \, \colon \, F_q[{}_aE_2] \rightarrow R \Big) \, $  is the
span of  $ \, {\Bbb B}_a \setminus \{1\} \, $.  
%
%
    $ \; {F_q[{}_aE_2]}^\vee = \sum_{n \geq 0} \! {\Big( \!
  {(q\!-\!1)}^{-1} J \Big)}^n $  \; is nothing but the unital
$ R $--algebra  (inside  $ \F_q[{}_aE_2] \, $)  with
generators  $ \, H_\pm := \displaystyle{\, \alpha^{\pm 1} -
1 \, \over \, q - 1\,} \, $,  $ \, E' := \displaystyle{\,
\beta \, \over \, q - 1\,} \, $,  and  $ \, F' :=
\displaystyle{\, \gamma \, \over \, q - 1\,} \, $
and relations
  $$  \displaylines{
   E' F' \! = \! q^{-2} F' E' \! ,  \, H_+ E' \! = \! q^2 E' H_+
+ (q+1) E' \! ,  \, H_+ F' \! = \! q^2 F' H_+ + (q+1) F' \! ,  \,
H_+ H_- \! = \! H_- H_+  \cr
   E' H_- = q^2 H_- E' + (q+1) E' \! ,  \, F' H_- = q^2 H_- F'
+ (q+1) F' \! ,  \, H_+ + H_- + (q-1) H_+ H_- = 0  \cr }  $$
with a Hopf structure given by
  $$  \displaylines{
   \Delta(E') = E' \otimes 1 + 1 \otimes E' + (q-1) \cdot H_+
\otimes E' \, ,  \hskip19pt  \epsilon(E') = 0 \, ,  \hskip19pt
S(E') = - E' - (q-1) H_- E' \, \phantom {.}  \cr
   \;\;  \Delta(H_\pm) = H_\pm \otimes 1 + 1 \otimes H_\pm +
(q-1) \cdot H_\pm \otimes H_\pm \, ,  \hskip17pt
\epsilon(H_\pm) = 0 \, ,  \hskip17pt S(H_\pm) = H_\mp \,
\hskip25pt  \phantom{.}  \cr
   \Delta(F') = F' \otimes 1 + 1 \otimes F' + (q-1) \cdot H_- \otimes F'
\, ,  \hskip19pt  \epsilon(F') = 0 \, ,  \hskip19pt  S(F') = - F' - (q-1)
H_+ F' \, .  \cr }  $$
This implies that  $ \, {F_q[{}_aE_2]}^\vee \,{\buildrel \, q
\rightarrow 1 \, \over \llongrightarrow}\, U({\gere_2}^{\! *}) \, $  as
co-Poisson Hopf algebras, for an isomorphism
  $$  {F_q[{}_aE_2]}^\vee \Big/ (q-1) \, {F_q[{}_aE_2]}^\vee
\,{\buildrel \cong \over \llongrightarrow}\, U({\gere_2}^{\! *})  $$
of co-Poisson Hopf algebras is given by  $ \;\, H_\pm \mod (q-1) \mapsto
\pm \text{h} \, $,  $ \, E' \mod (q-1) \mapsto \text{e} \, $,  \, and 
$ \, F' \mod (q-1) \mapsto \text{f} \, $;  so  $ \, {F_q[{}_aE_2]}^\vee
\, $  too specializes to  $ \, U({\gere_2}^{\! *}) \, $  {\sl as a
co-Poisson Hopf algebra},  as expected.
                                            \par
   We finish noting that, once more,  {\sl this analysis (and its
outcome) is characteristic-free}.

\vskip7pt

  {\bf 6.8 The identities  $ \, {\big({F_q[E_2]}^\vee\big)}' =
F_q[E_2] \, $  and  $ \, {\big({F_q[{}_aE_2]}^\vee\big)}' =
F_q[{}_aE_2] \, $.} \, In this section we verify   for the QFAs
$ \, H = F_q[E_2] \, $  and  $ \, H = F_q[{}_aE_2] \, $  the
validity of the part of Theorem 2.2{\it (b)}  claiming that
$ \; H \in \QFA \,\Longrightarrow\, {\big(H^\vee\big)}' = H \, $.
Once more, our arguments will prove this result for  $ \, \Char(\Bbbk)
\geq 0 \, $,  \, thus going beyond what forecasted by Theorem 2.2.
                                           \par
   Formulas  $ \; \Delta^n(E) = \sum\limits_{r+s+1=n} \text{a}^{\otimes r}
\otimes E \otimes {\big( \text{a}^{-1} \big)}^{\otimes s} \, $,  $ \;
\Delta^n(D_\pm) = \sum\limits_{r+s+1=n} {\big( \text{a}^{\pm 1}
\big)}^{\otimes r} \otimes D_\pm \otimes 1^{\otimes s} \; $  and 
$ \; \Delta^n(F) = \sum\limits_{r+s+1=n} {\big( \text{a}^{-1}
\big)}^{\otimes r} \otimes E \otimes \text{a}^{\otimes s} \; $ 
are found by induction.  
%
%
These identities imply the following  
  $$  \begin{array}{rcl}
   \delta_n(E)  &  = & {\textstyle \sum\limits_{r+s+1=n}} \!
{(\text{a} - 1)}^{\otimes r} \otimes E \otimes {\big( \text{a}^{-1}
- 1 \big)}^{\otimes s} = {(q-1)}^{n-1} \! {\textstyle
\sum\limits_{r+s+1=n}} \! {D_+}^{\! \otimes r}
\otimes E \otimes {D_-}^{\! \otimes s}  \\
  \delta_n (D_\pm) & =  & {\big( \text{a}^{\pm 1} - 1 \big)}^{\otimes
(n-1)} \otimes D_\pm = {(q-1)}^{n-1} {D_\pm}^{\!\otimes n}  \\
  \delta_n(F)  &  = & {\textstyle \sum\limits_{r+s+1=n}} \!
{\big( \text{a}^{-1} - 1 \big)}^{\otimes r} \otimes E \otimes
{(\text{a} - 1)}^{\otimes s} = {(q-1)}^{n-1} \! {\textstyle
\sum\limits_{r+s+1=n}} \! {D_-}^{\! \otimes r} \otimes
E \otimes {D_+}^{\! \otimes s}  \end{array}  $$
which give   $ \; \dot{E} := (q-1) E $,  $ \, \dot{D}_\pm
:= (q-1) D_\pm $,  $ \, \dot{F} := (q-1) F \in {\big( {F_q[E_2]}^\vee
\big)}' \setminus (q-1) \cdot {\big({F_q[E_2]}^\vee\big)}' \, $. 
\; So  $ {\big( {F_q[E_2]}^\vee \big)}' $  contains the unital 
$ R $--subalgebra  $ A' $  generated (inside  $ \, \F_q[E_2] \, $) 
by  $ \dot{E} \, $, $ \dot{D}_\pm $  and  $ \dot{F} \, $;  but  $ \,
\dot{E} = \text{b} \, $,  $ \, \dot{D}_\pm = \text{a}^{\pm 1} - 1 \, $, 
\, and  $ \, \dot{F} = \text{c} \, $,  thus  $ A' $  is just  $ F_q[E_2]
\, $.  Since  $ {F_q[E_2]}^\vee $  is the  $ R $--span  of  $ \, \Big\{\,
E^e D_+^{d_+} D_-^{d_-} F^f \,\Big\vert\, e, d_+, d_-, f \in \N \,\Big\}
\, $,  one easily sees   --- using the previous formulas for  $ \Delta^n $ 
---   that in fact  $ \, {\big({F_q[E_2]}^\vee\big)}' = A' = F_q[E_2]
\, $,  \, q.e.d.   
                                            \par
   When dealing with the adjoint case, the previous arguments
go through again: in fact,  $ {\big({F_q[{}_aE_2]}^\vee\big)}' $
turns out to coincide with the unital  $ R $--subalgebra  $ A'' $
generated (inside  $ \, \F_q[{}_aE_2] \, $)  by  $ \, \dot{E}' :=
(q-1) E' = \beta \, $,  $ \, \dot{H}_\pm := (q-1) H_\pm = \alpha^{\pm 1}
- 1 \, $,  and  $ \, \dot{F}' := (q-1) F' = \gamma \, $;  but this is
also generated by  $ \beta \, $,  $ \alpha^{\pm 1} $  and  $ \gamma
\, $,  \, thus it coincides with  $ F_q[{}_aE_2] \, $,  \, q.e.d.

 \vskip1,1truecm

\centerline {\bf \S \; 7 \  Fourth example: the Heisenberg
group  $ \, H_n \, $ }

\vskip10pt

  {\bf 7.1  The classical setting.} \, Let  $ \Bbbk $  be any field
of characteristic  $ \, p \geq 0 \, $.  Let  $ \, G := H_n(\Bbbk) =
H_n \, $,  the  $ (2 \, n + 1) $--dimensional  Heisenberg group; its
tangent Lie algebra  $ \, \gerg = \gerh_n \, $  is generated by
$ \, \{\, f_i, h, e_i \,\vert\, i = 1, \dots, n \,\} \, $  with
relations  $ \, [e_i,f_j] = \delta_{i{}j} h \, $,  $ \, [e_i,e_j]
= [f_i,f_j] = [h,e_i] = [h,f_j] = 0 \, $  ($ \forall \, i, j = 1,
\dots n \, $).  The formulas  $ \, \delta(f_i) = h \otimes f_i -
f_i \otimes h \, $,  $ \, \delta(h) = 0 \, $,  $ \, \delta(e_i) =
h \otimes e_i - e_i \otimes h \, $  ($ \forall \, i = 1, \dots n
\, $)  make  $ \gerh_n $  into a Lie bialgebra, which provides  $ H_n $
with a structure of Poisson group; these same formulas give also a
presentation of the co-Poisson Hopf algebra  $ U({\gerh}_n) $  (with
the standard Hopf structure).  When  $ \, p > 0 \, $  we consider on
$ \gerh_n $  the  $ p $--operation  uniquely defined by  $ \, e_i^{\,
[p\,]} = 0 \, $,  $ \, f_i^{\,[p\,]} = 0 \, $,  $ \, h^{[p\,]} = h
\, $  (for all  $ \, i = 1, \dots, n \, $),  which makes it into a
restricted Lie bialgebra.  The group  $ H_n $  is usually realized
as the group of all square matrices  $ \, {\big( a_{i{}j} \big)}_{i,j
= 1, \dots, n+2;} \, $  such that  $ \, a_{i{}i} = 1 \; \forall \, i
\, $  and  $ \, a_{i{}j} = 0 \; \forall\, i, j \, $  such that either
$ \, i > j \, $  or  $ \, 1 \not= i < j \, $  or  $ \, i < j \not= n+2
\, $;  it can also be realized as  $ \, H_n = \Bbbk^n \times \Bbbk
\times \Bbbk^n \, $  with group operation given by  $ \; \big(
\underline{a}', c', \underline{b}' \big) \cdot \big( \underline{a}'',
c'', \underline{b}'' \big) = \big( \underline{a}' + \underline{a}'',
c' + c'' + \underline{a}' \ast \underline{b}'', \underline{b}' +
\underline{b}'' \big) \, $,  \;  where we use vector notation  $ \,
\underline{v} = (v_1, \dots, v_n) \in k^n \, $  and  $ \, \underline{a}'
\ast \underline{b}'' := \sum_{i=1}^n a'_i b''_i \, $  is the standard
scalar product in  $ \, k^n \, $;  in particular the identity of
$ H_n $  is  $ \, e = (\underline{0}, 0, \underline{0}) \, $  and the
inverse of a generic element is given by  $ \, {\big( \underline{a},
c, \underline{b} \big)}^{-1} = \big( \! - \! \underline{a} \, , - c
+ \underline{a} \ast \underline{b} \, , \! - \underline{b} \, \big)
\, $.  Therefore  $ F[{H}_n] $  is the unital associative commutative 
$ \Bbbk $--algebra  with generators  $ \, a_1 \, $,  $ \dots $, 
$ a_n \, $,  $ c \, $,  $ b_1 \, $,  $ \dots $,  $ b_n \, $,  and
with Poisson Hopf structure given by
  $$  \displaylines{
   \Delta(a_i) = a_i \otimes 1 + 1 \otimes a_i \; ,  \hskip13pt
\Delta(c) = c \otimes 1 + 1 \otimes c + {\textstyle \sum_{\ell=1}^{n}}
a_\ell \otimes b_\ell \; ,  \hskip13pt  \Delta(b_i) =
b_i \otimes 1 + 1 \otimes b_i  \cr   
  \epsilon(a_i) = 0 \; ,  \hskip7pt  \epsilon(c) = 0 \; ,  \hskip7pt
\epsilon(b_i) = 0 \; ,  \hskip21pt  S(a_i) = - a_i \; ,  \hskip7pt
S(c) = - c + {\textstyle \sum_{\ell=1}^{n}} a_\ell b_\ell \; ,
\hskip7pt  S(b_i) = - b_i  \cr }  $$
  $$  \{a_i,a_j\} = 0 \; ,  \hskip17pt  \{a_i,b_j\} = 0 \; ,  \hskip17pt
\{b_i,b_j\} = 0 \; , \hskip17pt  \{c \, , a_i\} = a_i \; ,  \hskip17pt
\{c \, , b_i\} = b_i  $$
for all  $ \; i, j= 1, \dots, n \, $.  (Remark: with respect to this
presentation, we have  $ \, f_i = {\partial_{b_i}}{\big\vert}_e \, $,
$ \, h = {\partial_c}{\big\vert}_e \, $,  $ \, e_i = {\partial_{a_i}}
{\big\vert}_e \, $,  where  $ e $  is the identity element of  $ H_n
\, $).  The dual Lie bialgebra  $ \, \gerg^* = {\gerh_n}^{\!*} \, $
is the Lie algebra with generators  $ \, \text{f}_i \, $,  $ \text{h}
\, $,  $ \text{e}_i \, $,  and relations  $ \, [\text{h},\text{e}_i]
= \text{e}_i \, $,  $ [\text{h},\text{f}_i] = \text{f}_i \, $, 
$ [\text{e}_i,\text{e}_j] = [\text{e}_i,\text{f}_j] = [\text{f}_i,
\text{f}_j] = 0 \, $,  with Lie cobracket given by  $ \, \delta
(\text{f}_i) = 0 \, $, $ \, \delta(\text{h}) = \sum_{j=1}^{n}
(\text{e}_j \otimes \text{f}_j - \text{f}_j \otimes \text{e}_j) \, $, 
$ \, \delta(\text{e}_i) = 0 \, $  for all  $ \, i= 1, \dots, n \, $ 
(we take  $ \, \text{f}_i := f_i^* \, $,  $ \, \text{h} := h^* \, $, 
$ \, \text{e}_i := e_i^* \, $,  where  $ \, \big\{\, f_i^*, h^*,
e_i^* \,\vert\, i= 1, \dots, n \,\big\} \, $  is the basis of 
$ {\gerh_n}^{\! *} $  which is the dual of the basis  $ \{\, f_i,
h, e_i \,\vert\, i= 1, \dots, n \,\} $  of  $ \gerh_n \, $).  This
again gives a presentation of  $ \, U({\gerh_n}^{\!*}) \, $  too. 
If  $ \, p > 0 \, $  then  $ {\gerh_n}^{\!*} $  is a restricted Lie
bialgebra with respect to the  $ p $--operation  given by  $ \,
\text{e}_i^{\,[p\,]} = 0 \, $,  $ \, \text{f}_i^{\,\,[p\,]} = 0 \, $, 
$ \, \text{h}^{[p\,]} = \text{h} \, $  (for all  $ \, i = 1, \dots,
n \, $).  The simply connected algebraic Poisson group with tangent
Lie bialgebra  $ {{\gerh}_n}^{\!*} $  can be realized (with  $ \,
\Bbbk^\star := \Bbbk \setminus \{0\} \, $)  as  $ \; {{}_s{H}_n}^{\! *}
= \Bbbk^n \times \Bbbk^\star \times \Bbbk^n \, $,  with group operation 
$ \; \big( \underline{\dot{\alpha}}, \underline{\dot{\gamma}},
\underline{\dot{\beta}} \,\big) \cdot \big( \underline{\check{\alpha}},
\underline{\check{\gamma}}, \underline{\check{\beta}} \,\big) = \big(
\check{\gamma} \underline{\dot{\alpha}} + {\dot{\gamma}}^{-1}
\underline{\check{\alpha}}, \dot{\gamma} \check{\gamma},
\check{\gamma} \underline{\dot{\beta}} + {\dot{\gamma}}^{-1}
\underline{\check{\beta}} \,\big) \, $;  so the identity of
$ {{}_s{H}_n}^{\! *} $  is  $ \, e = (\underline{0}, 1, \underline{0})
\, $  and the inverse is given by  $ \, {\big( \underline{\alpha},
\gamma, \underline{\beta} \,\big)}^{-1} = \big( \! - \!
\underline{\alpha}, \gamma^{-1}, - \underline{\beta} \,\big) $.
Its centre is  $ \, Z \big( {{}_s{H}_n}^{\! *} \big) = \big\{
(\underline{0}, 1, \underline{0}), (\underline{0}, -1, \underline{0})
\big\} =: Z \, $,  so there is only one other (Poisson) group with
tangent Lie bialgebra  $ {{\gerh}_n}^{\! *} \, $,  that is the
adjoint group  $ \, {{}_a{H}_n}^{\! *} := {{}_s{H}_n}^{\! *}
\Big/ Z \; $.
                                                  \par
   It is clear that  $ F\big[{{}_s{H}_n}^{\! *}] $  is the unital
associative commutative  $ \Bbbk $--algebra  with generators
$ \, \alpha_1 \, $,  $ \dots $,  $ \alpha_n \, $,  $ \gamma^{\pm 1}
\, $,  $ \beta_1 \, $,  $ \dots $,  $ \beta_n \, $,  \, and with
Poisson Hopf algebra structure given by
  $$  \displaylines{
   \Delta(\alpha_i) = \alpha_i \otimes \gamma + \gamma^{-1}
\otimes \alpha_i \; ,  \hskip17pt  \Delta\big(\gamma^{\pm 1}\big) =
\gamma^{\pm 1} \otimes \gamma^{\pm 1} \; ,  \hskip17pt  \Delta(\beta_i)
= \beta_i \otimes \gamma + \gamma^{-1} \otimes \beta_i  \cr
   \epsilon(\alpha_i) = 0 \; ,  \hskip5pt  \epsilon \big(
\gamma^{\pm 1} \big) = 1 \; ,  \hskip5pt  \epsilon(\beta_i) = 0 \; ,
\hskip29pt  S(\alpha_i) =  - \alpha_i \; ,  \hskip5pt  S \big(
\gamma^{\pm 1} \big) = \gamma^{\mp 1} \; ,  \hskip5pt  S(\beta_i)
= - \beta_i  \cr
   \{\alpha_i,\alpha_j\} = \{\alpha_i,\beta_j\} = \{\beta_i,
\beta_j\} = \{\alpha_i,\gamma\} = \{\beta_i,\gamma\} = 0 \; ,
\hskip13pt  \{\alpha_i,\beta_j\} = \delta_{i{}j} \big(
\gamma^2 - \gamma^{-2} \big) \big/ 2  \cr }  $$
for all  $ \; i, j= 1, \dots, n \, $  (Remark: with respect to this
presentation, we have  $ \, \text{f}_i = {\partial_{\beta_i}}
{\big\vert}_e \, $,  $ \, \text{h} = {\,1\, \over \,2\,} \, \gamma \,
{\partial_\gamma}{\big\vert}_e \, $,  $ \, \text{e}_i =
{\partial_{\alpha_i}}{\big\vert}_e \, $,  where  $ e $
is the identity element of  $ {{}_s{H}_n}^{\! *} \, $),  and
$ F\big[{{}_a{H}_n}^{\! *}\big] $  can be identified   --- as in the
case of the Euclidean group ---   with the Poisson Hopf subalgebra of
$ F\big[{{}_a{H}_n}^{\! *}\big] $  which is spanned by products of an
even number of generators: this is generated
%
%
by  $ \, \alpha_i \, \gamma \, $,  $ \, \gamma^{\pm 2} \, $,  and 
$ \, \gamma^{-1} \, \beta_i \; $  ($ \, i= 1, \dots, n \, $).

\vskip7pt

  {\bf 7.2 The QrUEAs  $ \, U_q^s({\gerh}_n) \, $  and  $ \, U_q^a
({\gerh}_n) \, $.} \, We switch now to quantizations.  Once again,
let  $ R $  be a 1dD and let  $ \, \h \in R \setminus \{0\} \, $  and
assume  $ \, q := 1 + \h \in R \, $  be invertible, like in \S 5.2.
                                                  \par
   Let $ \, \Bbb{U}_q(\gerg) = \Bbb{U}_q^s ({\gerh}_n) \, $  be the
unital associative  $ F(R) $--algebra  with generators  $ \, F_i \, $,
$ L^{\!\pm 1} \, $,  $ E_i \, $  (for  $ \, i = 1, \dots, n \, $)  and
relations   
  $$  L L^{-1} = 1 = L^{-1} L \; ,  \,\quad  L^{\!\pm 1} F = F L^{\pm 1}
\; ,  \,\quad  L^{\pm 1} E = E L^{\pm 1} \; ,  \,\quad  E_i F_j - F_j E_i
= \delta_{i{}j} {{\, L^2 - L^{-2} \, \over \, q - q^{-1} \,}}   $$
for all  $ \, i, j= 1, \dots, n \, $;  we give it a structure of Hopf
algebra, by setting  ($ \forall \; i, j = 1, \dots, n \, $)
  $$  \displaylines{
   \Delta(E_i) = E_i \otimes 1 + L^2 \otimes E_i \; ,  \hskip15pt
\Delta\big(L^{\pm 1}\big) = L^{\pm 1} \otimes L^{\pm 1} \; , 
\hskip15pt  \Delta(F_i) = F_i \otimes L^{-2} + 1 \otimes F_i  \cr
   \epsilon(E_i) = 0 \; ,  \hskip9pt  \epsilon\big(L^{\pm 1}\big) = 1 \; ,
\hskip9pt  \epsilon(F_i) = 0 \; ,  \hskip15pt  S(E_i) = - L^{-2} E_i \; ,
\hskip9pt  S\big(L^{\pm 1}\big) =  L^{\mp 1} \; ,  \hskip9pt
S(F_i) = - F_i L^2  \cr }  $$
Note then that  $ \; \Big\{\, \prod_{i=1}^n F_i^{a_i} \! \cdot \! L^z
\! \cdot \! \prod_{i=1}^n E_i^{d_i} \,\Big\vert\, z \in \Z, \, a_i,
d_i \in \N, \, \forall\, i \,\Big\} \; $  is an  $ F(R) $--basis
of  $ \Bbb{U}_q^s(\gerh_n) \, $.   
                                             \par
   Now, let  $ \, U_q^s(\gerh_n) \, $  be the unital
$ R $--subalgebra  of  $ \Bbb{U}_q^s(\gerh_n) $  generated by
the elements  $ \, F_1 \, $,  $ \dots $,  $ F_n \, $,  $ D :=
\displaystyle{\, L - 1 \, \over \, q - 1 \,} \, $,  $ \varGamma
:= \displaystyle{{\, L - L^{-2} \, \over \, q - q^{-1} \,}} \, $, 
$ E_1 \, $,  $ \dots $,  $ E_n \, $.  Then  $ U_q^s(\gerh_n) $  can
be presented as the associative unital algebra with generators  $ \,
F_1 \, $,   $ \dots $,  $ F_n \, $,  $ L^{\pm 1} \, $,  $ D \, $, 
$ \varGamma \, $,  $ E_1 \, $,  $ \dots $,  $ E_n \, $ and relations
  $$  \displaylines{
   \hskip7pt  D X = X D \; ,  \hskip25pt  L^{\pm 1} X = X L^{\pm 1} \; ,
\hskip25pt  \varGamma X = X \varGamma \; ,  \hskip25pt  E_i F_j - F_j E_i
= \delta_{i{}j} \varGamma  \cr
  L = 1 + (q-1) D \; ,  \hskip15pt  L^2 - L^{-2} = \big( q - q^{-1} \big)
\varGamma \; ,  \hskip15pt  D (L + 1) \big(1 + L^{-2}\big) =
\big(1 + q^{-1}\big) \varGamma  \cr }  $$
for all  $ \; X \in {\big\{F_i, L^{\pm 1}, D, \varGamma, E_i
\big\}}_{i=1,\dots,n} \, $  and  $ \, i, j= 1, \dots, n \, $;
furthermore,  $ U_q^s(\gerh_n) $  is a Hopf subalgebra (over
$ R $),  with
  $$   \begin{matrix}
      \Delta(\varGamma) = \varGamma \otimes L^2 + L^{-2} \otimes \varGamma
\; , &  \hskip25pt  \epsilon(\varGamma) = 0  \; ,  &  \hskip25pt
S(\varGamma) = - \varGamma  \\
      \Delta(D) = D \otimes 1 + L \otimes D \; ,  &  \hskip25pt
\epsilon(D) = 0 \; ,  &  \hskip25pt  S(D) = - L^{-1} D \; .  \\
     \end{matrix}  $$
Moreover, from relations  $ \, L = 1 + (q-1) D \, $  and
$ \, L^{-1} = L^3 - \big( q - q^{-1} \big) L \varGamma \, $ 
it follows that
  $$  U_q^s(\gerh_n) \; = \; \hbox{$ R $--span  of}  \;\,
\bigg\{\, {\textstyle \prod\limits_{i=1}^n} F_i^{a_i} \! \cdot
\! D^b \varGamma^c \! \cdot \! {\textstyle \prod\limits_{i=1}^n}
E_i^{d_i} \,\bigg\vert\, a_i, b, c, d_i \in \N, \, \forall\,
i= 1, \dots, n \,\bigg\}   \eqno (7.1)  $$
                                                  \par
   The ``adjoint version'' of  $ \Bbb{U}_q^s(\gerh_n) $  is the
subalgebra  $ \Bbb{U}_q^a(\gerh_n) $  generated  by  $ \, F_i $,
$ K^{\pm 1} := L^{\pm 2} $,  $ E_i \; (i = 1, \dots, n) $,  which
is a {\sl Hopf\/}  subalgebra too.  It also has an  $ R $--integer 
form  $ U_q^a(\gerh_n) \, $,  the  $ R $--subalgebra  generated by 
$ \, F_1 \, $,  $ \dots $,  $ F_n \, $,  $ K^{\pm 1} \, $,  $ H :=
\displaystyle{\, K - 1 \, \over \, q - 1 \,} \, $,  $ \varGamma :=
\displaystyle{\, K - K^{-1} \, \over \, q - q^{-1} \,} \, $,  $ \,
E_1 \, $,  $ \dots $,  $ E_n  \, $:  this has relations  
  $$  \displaylines{
   \hskip7pt  H X = X H \; ,  \hskip25pt  K^{\pm 1} X = X K^{\pm 1} \; ,
\hskip25pt  \varGamma X = X \varGamma \; ,  \hskip25pt  E_i F_j - F_j
E_i = \delta_{i{}j} \varGamma  \cr  
  K = 1 + (q-1) H \; ,  \hskip15pt  K - K^{-1} = \big( q - q^{-1} \big)
\varGamma \; ,  \hskip15pt  H \big(1 + K^{-1}\big) =
\big( 1 + q^{-1} \big) \varGamma  \cr }  $$
for all  $ \; X \in {\big\{ F_i, K^{\pm 1}, H, \varGamma, E_i
\big\}}_{i=1,\dots,n} \, $  and  $ \, i, j= 1, \dots, n \, $, 
and Hopf operations given by
  $$   \begin{matrix}
   \Delta(E_i) = E_i \otimes 1 + K \otimes E_i \; ,  &  \hskip21pt
\epsilon(E_i) = 0 \; ,  &  \hskip21pt  S(E_i) = - K^{-1} E_i  \\
   \Delta\big(K^{\pm 1}\big) = K^{\pm 1} \otimes K^{\pm 1} \; ,  &
\hskip21pt  \epsilon\big(K^{\pm 1}\big) = 1 \; ,  &  \hskip21pt
S\big(K^{\pm 1}\big) = K^{\mp 1}  \\
   \Delta(H) = H \otimes 1 + K \otimes H \; ,  &  \hskip21pt
\epsilon(H) = 0 \; ,  &  \hskip21pt  S(H) = - K^{-1} H  \\
   \Delta(\varGamma) = \varGamma \otimes K^{-1} + K \otimes \varGamma
\; ,   &  \hskip21pt  \epsilon(\varGamma) = 0 \; ,  &  \hskip21pt
S(\varGamma) = - \varGamma  \\
   \Delta(F_i) = F_i \otimes K^{-1} + 1 \otimes F_i \; ,  &
\hskip21pt  \epsilon(F_i) = 0 \; ,  &  \hskip21pt  S(F_i) =
- F_i K^{+1}  \\
      \end{matrix}  $$
for all  $ \, i= 1, \dots, n \, $.  One can easily check that  $ \, U_q^s
(\gerh_n) \, $  is a QrUEA, with  $ \, U(\gerh_n) \, $  as semiclassical
limit: in fact, mapping the generators  $ \, F_i \mod (q-1) \, $,  $ \,
L^{\pm 1} \mod (q-1) \, $,  $ \, D \mod (q-1) \, $,  $ \, \varGamma \mod
(q-1) \, $,  $ \, E_i \mod (q-1) \, $  respectively to  $ \, f_i \, $,
$ 1 \, $,  $ h\big/2 \, $,  $ h \, $,  $ e_i \in U(\gerh_n) \, $  yields
a co-Poisson Hopf algebra isomorphism between  $ \, U_q^s(\gerh_n)
\Big/ (q-1) \, U_q^s(\gerh_n) \, $  and  $ U(\gerh_n) \, $.  Similarly,
$ \,U_q^a(\gerh_n) \, $  is a QrUEA too, again with limit  $ \, U(\gerh_n)
\, $,  for a co-Poisson Hopf algebra isomorphism between  $ \,
U_q^a(\gerh_n) \Big/ (q-1) \, U_q^a(\gerh_n) \, $  and  $ \, U(\gerh_n)
\, $  is given by mapping the generators  $ \, F_i \mod (q-1) \, $, 
$ \, K^{\pm 1} \mod (q-1) \, $,  $ \, H \mod (q-1) \, $,  $ \, \varGamma
\mod (q-1) \, $,  $ E_i \mod (q-1) \, $  respectively to  $ \, f_i \, $, 
$ 1 \, $,  $ h \, $,  $ h \, $,  $ e_i \in U(\gerh_n) \, $.   

\vskip7pt

  {\bf 7.3 Computation of  $ \, {U_q(\gerh_n)}' \, $  and
specialization  $ \, {U_q(\gerh_n)}' \,{\buildrel {q \rightarrow 1}
\over \llongrightarrow}\, F\big[{H_n}^{\!\star}\big] \, $.} \, Here
we compute  $ {U_q^s(\gerh_n)}' $  and  $ {U_q^a(\gerh_n)}' \, $, 
and their semiclassical limits, along the pattern of \S 5.3.
                                               \par
   Definitions give, for any  $ \, n \in \N \, $,  $ \, \Delta^n(E_i)
= \sum_{s=1}^n {(L^2)}^{\otimes (s-1)} \otimes E_i \otimes 1^{\otimes
(n-s)} $,  \, hence  $ \; \delta_n(E_i) = {(q-1)}^{n-1} \cdot D^{\otimes
(n-1)} \otimes E_i \; $  so  $ \; \delta_n\big( (q-1) E \big) \in
{(q-1)}^n \, U_q^s(\gerh_n) \setminus {(q-1)}^{n+1} \, U_q^s(\gerh_n)
\; $  whence  $ \, \dot{E}_i := (q-1) \, E_i \in {U_q^s(\gerh_n)}' $,
whereas  $ \, E_i \notin {U_q^s(\gerh_n)}' $;  similarly, we have  $ \;
\dot{F}_i := (q-1) \, F_i $,  $ L^{\pm 1} $,  $ \dot{D} := (q-1) \, D
= L - 1 $,  $ \dot{\varGamma} := (q-1) \, \varGamma \in {U_q^s(\gerh_n)}'
\setminus (q-1) \, {U_q^s(\gerh_n)}' $,  for all  $ \, i = 1, \dots, n
\, $.  Thus  $ {U_q^s(\gerh_n)}' $  contains the subalgebra  $ U' $
generated by  $ \, \dot{F}_i \, $,  $ L^{\pm 1} \, $,  $ \, \dot{D}
\, $,  $ \dot{\varGamma} \, $,  $ \, \dot{E}_i \, $;  we argue that 
$ \, {U_q^s(\gerh_n)}' = U' \, $:  this is easily seen   --- like for 
$ {SL}_2 $  and for  $ E_2 $  ---   using the formulas above along
with (7.1).  Therefore  $ {U_q^s(\gerh_n)}' $  is the unital 
$ R $--algebra  with generators  $ \, \dot{F}_1 \, $,  $ \dots $, 
$ \dot{F}_n \, $,  $ L^{\pm 1} \, $,  $ \dot{D} \, $,  $ \dot{\varGamma}
\, $,  $ \dot{E}_1 \, $,  $ \dots $,  $ \dot{E}_n \, $  and relations
  $$  \displaylines{
   \hskip7pt  \dot{D} \dot{X} = \dot{X} \dot{D} \; ,  \hskip25pt
L^{\pm 1} \dot{X} = \dot{X} L^{\pm 1} \; ,  \hskip25pt
\dot{\varGamma} \dot{X} = \dot{X} \dot{\varGamma} \; , 
\hskip25pt  \dot{E}_i \dot{F}_j - \dot{F}_j \dot{E}_i =
\delta_{i{}j} (q-1) \dot{\varGamma}  \cr
  L = 1 + \dot{D} \; ,  \hskip15pt  L^2 - L^{-2} = \big( 1 + q^{-1}
\big) \dot{\varGamma} \; ,  \hskip15pt  \dot{D} (L + 1) \big( 1 +
L^{-2}\big) = \big( 1 + q^{-1} \big) \dot{\varGamma}  \cr }  $$
for all  $ \; \dot{X} \in {\big\{ \dot{F}_i, L^{\pm 1}, \dot{D},
\dot{\varGamma}, \dot{E}_i \big\}}_{i=1,\dots,n} \, $  and  $ \,
i, j= 1, \dots, n \, $,  with Hopf structure given by
  $$   \begin{matrix}
   \Delta\big(\dot{E}_i\big) = \dot{E}_i \otimes 1 + L^2 \otimes \dot{E}_i
\; ,  &  \hskip7pt  \epsilon\big(\dot{E}_i\big) = 0 \; ,  &  \hskip7pt
S\big(\dot{E}_i\big) =  - L^{-2} \dot{E}_i  &  \hskip19pt  \forall
\; i = 1, \dots, n \; \phantom{.}  \\
   \Delta\big(L^{\pm 1}\big) = L^{\pm 1} \otimes L^{\pm 1} \; ,  &
\hskip7pt  \epsilon\big(L^{\pm 1}\big) = 1 \; ,  &  \hskip7pt
S\big(L^{\pm 1}\big) =  L^{\mp 1}  &  {}  \\
   \Delta\big(\dot{\varGamma}\big) = \dot{\varGamma} \otimes L^2 + L^{-2}
\otimes \dot{\varGamma} \; , &  \hskip7pt
\epsilon\big(\dot{\varGamma}\big) = 0  \; ,  &  \hskip7pt
S\big(\dot{\varGamma}\big) = - \varGamma  \\
   \end{matrix}  $$   
  $$   \begin{matrix}
   \Delta\big(\dot{D}\big) = \dot{D} \otimes 1 + L \otimes \dot{D} \; ,
&  \hskip7pt  \epsilon\big(\dot{D}\big) = 0 \; ,  &  \hskip7pt 
S\big(\dot{D}\big) = - L^{-1} \dot{D}  \\  
   \Delta\big(\dot{F}_i\big) = \dot{F}_i \otimes L^{-2} + 1 \otimes
\dot{F}_i \; ,  &  \hskip7pt  \epsilon\big(F_i\big) = 0 \; ,  & 
\hskip7pt  S\big(\dot{F}_i\big) = - \dot{F}_i L^2  &  \hskip19pt
\forall \; i = 1, \dots, n \; .  \\  
   \end{matrix}  $$   
   \indent   A similar analysis shows that  $ \, {U_q^a(\gerh_n)}' \, $
is the unital  $ R $--subalgebra  $ U'' $  of  $ U_q^a(\gerh_n) $ 
generated by  $ \, \dot{F}_i \, $,  $ \, K^{\pm 1} \, $,  $ \, \dot{H}
:= (q-1) \, H $,  $ \, \dot{\varGamma} \, $,  $ \dot{E}_i \, $  ($ i = 1,
\dots, n \, $);  in particular,  $ \, {U_q^a(\gerh_n)}' \subset {U_q^s
(\gerh_n)}' \, $.  Thus  $ {U_q^a(\gerh_n)}' $  is the unital associative 
$ R $--algebra  with generators  $ \, \dot{F}_1 \, $,  $ \dots $, 
$ \dot{F}_n \, $,  $ \dot{H} \, $,  $ K^{\pm 1} \, $,  $ \dot{\varGamma}
\, $,  $ \dot{E}_1 \, $,  $ \dots $,  $ \dot{E}_n \, $  and relations
  $$  \displaylines{  
   \hskip7pt  \dot{H} \dot{X} = \dot{X} \dot{H} \; ,  \hskip25pt
K^{\pm 1} \dot{X} = \dot{X} K^{\pm 1} \; ,  \hskip25pt  \dot{\varGamma}
\dot{X} = \dot{X} \dot{\varGamma} \; ,  \hskip25pt  \dot{E}_i \dot{F}_j
- \dot{F}_j \dot{E}_i = \delta_{i{}j} (q-1) \dot{\varGamma}  \cr
  K = 1 + \dot{H} \; ,  \hskip15pt  K - K^{-1} = \big( 1 + q^{-1} \big)
\dot{\varGamma} \; ,  \hskip15pt  \dot{H} \big( 1 + K^{-1} \big) =
\big( 1 + q^{-1} \big) \dot{\varGamma}  \cr }  $$
for all  $ \; \dot{X} \in {\big\{ \dot{F}_i, K^{\pm 1}, \dot{K},
\dot{\varGamma}, \dot{E}_i \big\}}_{i=1,\dots,n} \, $  and
$ \, i, j= 1, \dots, n \, $,  with Hopf structure given by
  $$   \begin{matrix}
   \Delta\big(\dot{E}_i\big) = \dot{E}_i \otimes 1 + K \otimes \dot{E}_i
\; ,  &  \hskip7pt  \epsilon\big(\dot{E}_i\big) = 0 \; ,  &  \hskip7pt
S\big(\dot{E}_i\big) =  - K^{-1} \dot{E}_i  &  \hskip19pt  \forall
\; i = 1, \dots, n \; \phantom{.}  \\
   \Delta\big(K^{\pm 1}\big) = K^{\pm 1} \otimes K^{\pm 1} \; ,  &
\hskip7pt  \epsilon\big(K^{\pm 1}\big) = 1 \; ,  &  \hskip7pt
S\big(K^{\pm 1}\big) = K^{\mp 1}  &  {}  \\
   \Delta\big(\dot{\varGamma}\big) = \dot{\varGamma} \otimes K +
K^{-1} \otimes \dot{\varGamma} \; , &  \hskip7pt
\epsilon\big(\dot{\varGamma}\big) = 0  \; ,  &  \hskip7pt
S\big(\dot{\varGamma}\big) = - \varGamma  \\
   \Delta\big(\dot{H}\big) = \dot{H} \otimes 1 + K \otimes \dot{H} \; ,
&  \hskip7pt  \epsilon\big(\dot{H}\big) = 0 \; ,  &  \hskip7pt
S\big(\dot{H}\big) = - K^{-1} \dot{H}  \\
   \Delta\big(\dot{F}_i\big) = \dot{F}_i \otimes K^{-1} + 1 \otimes
\dot{F}_i \; ,  &  \hskip7pt  \epsilon\big(F_i\big) = 0 \; ,  &
\hskip7pt S\big(\dot{F}_i\big) = - \dot{F}_i K  &  \hskip19pt
\forall \; i = 1, \dots, n \; .  \\
   \end{matrix}  $$
   \indent   As  $ \, q \rightarrow 1 \, $,  \, the presentation
above provides an isomorphism of Poisson Hopf algebras
 \vskip-11pt
  $$  {U_q^s(\gerh_n)}' \Big/ (q-1) \, {U_q^s(\gerh_n)}'
\;{\buildrel \cong \over \llongrightarrow}\;
F \big[ {}_s{H_n}^{\!*} \big]  $$
 \vskip5pt
\noindent   
given by  $ \; \dot{E}_i \mod (q-1) \mapsto \alpha_i \,
\gamma^{+1} \, $,  $ \; L^{\pm 1} \mod (q-1) \mapsto \gamma^{\pm 1} \, $, 
$ \, \dot{D} \mod (q-1) \mapsto \gamma - 1 \, $,  $ \, \dot{\varGamma}
\mod (q-1) \mapsto \big( \gamma^2 - \gamma^{-2} \big) \Big/ 2 \; $,
$ \; \dot{F}_i \mod (q-1) \mapsto \gamma^{-1} \beta_i \; $.  In other
words, the semiclassical limit of  $ \, {U_q^s(\gerh_n)}' \, $  is
$ \, F \big[ {}_s{H_n}^{\!*} \big] \, $,  as predicted by  Theorem
2.2{\it (c)\/}  for  $ \, p = 0 \, $.  Similarly, when considering
the ``adjoint case'', we find a Poisson Hopf algebra isomorphism
  $$  {U_q^a(\gerh_n)}' \Big/ (q-1) \, {U_q^a(\gerh_n)}'
\;{\buildrel \cong \over \llongrightarrow}\; F \big[ {}_a{H_n}^{\!*}
\big]  \quad  \Big( \subset F \big[ {}_s{H_n}^{\!*} \big] \Big)  $$
given by  $ \;\, \dot{E}_i \mod (q-1) \mapsto \alpha_i \, \gamma^{+1} \, $, 
$ \, K^{\pm 1} \mod (q-1) \mapsto \gamma^{\pm 2} \, $,  $ \, \dot{H} \mod
(q-1) \mapsto \gamma^2 - 1 \, $,  $ \, \dot{\varGamma} \mod (q-1) \mapsto
\big( \gamma^2 - \gamma^{-2} \big) \Big/ 2 \; $,  $ \; \dot{F}_i \mod (q-1)
\mapsto \gamma^{-1} \beta_i \; $.  That is to say,  $ \, {U_q^a(\gerh_n)}'
\, $  has semiclassical limit  $ \, F \big[ {}_a{H_n}^{\!*} \big] \, $,
\, as predicted by  Theorem 2.2{\it (c)\/}  for  $ \, p = 0 \, $.
                                         \par
   We stress the fact that  {\sl this analysis is characteristic-free}, 
so we get in fact that its outcome does hold for  $ \, p > 0 \, $  as
well, thus ``improving''  Theorem 2.2{\it (c)\/}  (like in \S\S 5--6).

\vskip7pt

  {\bf 7.4 The identity  $ \, {\big({U_q(\gerh_n)}'\big)}^{\!\vee} =
U_q(\gerh_n) \, $.} \,  In this section we verify the part of  Theorem
2.2{\it (b)}  claiming, for  $ \, p = 0 \, $,  that  $ \, H \in \QrUEA
\,\Longrightarrow\, {\big(H'\big)}^{\!\vee} = H \, $,  both for  $ \,
H = U_q^s(\gerh_n) \, $  and for  $ \, H = U_q^a(\gerh_n) \, $.  {\sl
In addition, the same arguments will prove such a result for  $ \,
p > 0 \, $  too.}
                                          \par
   To begin with, using (7.1) and the fact that  $ \, \dot{F}_i $,
$ \dot{D} $,  $ \dot{\varGamma} $,  $ \dot{E}_i \in \text{\sl
Ker} \Big( \epsilon \, \colon \, {U_q^s(\gerh_n)}' \,
\relbar\joinrel\twoheadrightarrow \, R \Big) \, $
we get that  $ \, J := \text{\sl Ker}\,(\epsilon) \, $
is the  $ R $--span  of  $ \, {\Bbb M} \setminus \{1\} \, $,
where  $ \, {\Bbb M} \, $  is the set in the right-hand-side
of (7.1).  
%
%
    Since  $ \, {\big({U_q^s(\gerh_n)}'\big)}^{\!\vee} :=
  \sum_{n \geq 0} \! {\Big( \! {(q\!-\!1)}^{-1} J \Big)}^n \, $,  we
have that  $ \, {\big( {U_q^s (\gerh_n)}' \big)}^{\!\vee}
\, $  is generated   --- as a unital  $ R $--subalgebra  of
$ \Bbb{U}_q^s(\gerh_n) $  ---   by  $ \, {(q-1)}^{-1} \dot{F}_i
= F_i \, $,  $ {(q-1)}^{-1} \dot{D} = D \, $,  $ {(q-1)}^{-1}
\dot{\varGamma} = \varGamma \, $,  $ {(q-1)}^{-1} \dot{E}_i
= E_i \, $  ($ i = 1, \dots, n $),  so it coincides with
$ U_q^s(\gerh_n) $,  q.e.d.  In the adjoint case the
procedure is similar: one changes  $ L^{\pm 1} $,
resp.~$ \dot{D} $,  with  $ K^{\pm 1} $,  resp.~$ \dot{H} $,
and everything works as before.

\vskip7pt

  {\bf 7.5 The quantum hyperalgebra  $ \hyp_q(\gerh_n) $.} \, Like in
\S\S 5.5 and 6.5, we can define ``quantum hyperalgebras'' associated
to  $ \gerh_n \, $.  Namely, first we define a Hopf subalgebra of
$ \Bbb{U}_q^s(\gerh_n) $  over  $ \Z \big[ q, q^{-1} \big] $  whose
specialization at  $ \, q = 1 \, $  is the natural Kostant-like
$ \Z $--integer  form  $ U_\Z(\gerh_n) $  of  $ U(\gerh_n) $
(generated by divided powers, and giving the hyperalgebra
$ \hyp(\gerh_n) $  over any field  $ \Bbbk $  by scalar
  extension),
%
%
   and then take its scalar extension
over  $ R \, $.
                                             \par
   To be precise, let  $ \hyp^{s,\Z}_q(\gerh_n) $  be the unital
$ \Z\big[q,q^{-1}\big] $--subalgebra  of  $ \Bbb{U}^s_q(\gerh_n) $
(defined like above  {\sl but over\/}  $ \Z\big[q,q^{-1}\big] $)
generated by the ``quantum divided powers''  
  $$  F_i^{(m)} \! := F_i^m \! \Big/ {[m]}_q! \; ,  \qquad 
\left( {{L \, ; \, c} \atop {m}} \right) \! := \prod_{r=1}^n {{\;
q^{c+1-r} L - 1 \,} \over {\; q^r - 1 \;}} \; ,  \qquad 
E_i^{(m)} \! := E_i^m \! \Big/ {[m]}_q!  $$  
(for all  $ \, m \in \N \, $,  $ \, c \in \Z \, $  and  $ \, i = 1, \dots,
n \, $,  with notation of \S 5.5)  and by  $ L^{-1} \, $.  Comparing with
the case of  $ \gersl_2 $   --- noting that for each  $ i $  the quadruple 
$ (F_i,L,L^{-1},E_i) $  generates a copy of  $ \Bbb{U}^s_q(\gersl_2) $ 
---   we see at once that this is a Hopf subalgebra of  $ \Bbb{U}_q^s
(\gerh_n) \, $,  and  $ \, \hyp^{s,\Z}_q(\gerh_n){\Big|}_{q=1} \! \cong
U_\Z(\gerh_n) \, $;  \, thus  $ \, \hyp^s_q(\gerh_n) := R \otimes_{\Z
[q,q^{-1}]} \hyp^{s,\Z}_q(\gerh_n) \, $  (for any  $ R $  like in \S 6.2,
with  $ \, \Bbbk := R \big/ \h \, R \, $  and  $ \, p := \Char(\Bbbk)
\, $)  specializes at  $ \, q = 1 \, $  to the  $ \Bbbk $--hyperalgebra
$ \hyp(\gerh_n) \, $.  Moreover, among all the  $ \left( {{L \, ; \, c}
\atop {n}} \right) $'s  it is enough to take only those with  $ \, c
= 0 \, $.  {\sl From now on we assume  $ \, p > 0 \, $.}
                                             \par
   Pushing forward the close comparison with the case of  $ \gersl_2 $
we also see that  $ \, {\hyp^s_q(\gerh_n)}' \, $  is the unital
$ R $--subalgebra  of  $ \hyp^s_q(\gerh_n) $  generated by  $ L^{-1} $
and the ``rescaled quantum divided powers''  $ \, {(q \! - \! 1)}^m
F_i^{(m)} \, $,  $ \, {(q\!-\!1)}^m \! \left( {{L \, ; \, 0} \atop {m}}
\right) \, $  and  $ \, {(q\!-\!1)}^m E_i^{(m)} \, $,  for all  $ \, m
\in \N \, $  and  $ \, i = 1, \dots, n \, $.  It follows that  $ \,
{\hyp^s_q(\gerh_n)}'{\Big|}_{q=1} \, $  is generated by the
specializations at  $ \, q = 1 \, $  of  $ \, {(q-1)}^{p^r}
F_i^{(p^r)} \, $,  $ \, {(q-1)}^{p^r} \! \left( {{L \, ; \, 0}
\atop {p^r}} \right) \, $  and  $ \, {(q-1)}^{p^r} E_i^{(p^r)} \, $,
for all  $ \, r \in \N \, $,  $ \, i = 1, \dots, n \, $:  \, this proves
directly that the spectrum of  $ \, {\hyp^s_q(\gerh_n)}'{\Big|}_{q=1}
\, $  has dimension 0 and height 1, and its cotangent Lie algebra
  $ \, J \Big/ J^{\,2} \, $
 (where  $ J $  is the augmentation ideal of  $ {\hyp^s_q(\gerh_n)}'
{\Big|}_{q=1} \, $)
 has
basis  $ \, \Big\{\, {(q\!-\!1)}^{p^r} F_i^{(p^r)}, \, {(q\!-\!1)}^{p^r}
\! \left( {{L \, ; \, 0} \atop {p^r}} \right) \! ,$  $\, {(q\!-\!1)}^{p^r}
E_i^{(p^r)} \, \mod (q\!-\!1) \, {\hyp^s_q(\gerg)}' \, \mod J^{\,2}
\;\Big|\; r \!\in\! \N \, , \, i = 1, \dots, n \,\Big\} \, $.  
%
%
Finally,  $ \, \big(
{\hyp^s_q(\gerh_n)}' \big)^{\!\vee} \, $  is generated by  $ \,
{(q-1)}^{p^r-1} F_i^{(p^r)} \, $,  $ \, {(q-1)}^{p^r-1} \left(
{{L \, ; \, 0} \atop {p^r}} \right) \, $,  $ L^{-1} $  and  $ \,
{(q-1)}^{p^r-1} E_i^{(p^r)} \, $  (for all  $ r $  and  $ i \, $): 
\, in particular  $ \, \big( {\hyp^s_q(\gerh_n)}' \big)^{\!\vee}
\! \subsetneqq \! \hyp^s_q(\gerh_n) \, $,  and  $ \, \big(
{\hyp^s_q(\gerh_n)}' \big)^{\!\vee}{\Big|}_{q=1} $  is generated
by the cosets modulo  $ (q-1) $  of these elements, which form a
basis of the restricted Lie bialgebra  $ \gerk $  such that  $ \,
\big( {\hyp^s_q(\gerh_n)}' \big)^{\!\vee}{\Big|}_{q=1} \! = \,
\u(\gerk) \; $.   
                                             \par
   The previous analysis stems from  $ \Bbb{U}_q^s(\gerh_n) $,  and
so gives ``simply connected quantum objects''.  Instead we can start
from  $ \Bbb{U}_q^a(\gerh_n) \, $,  thus getting ``adjoint quantum
objects'', moving along the same pattern but for replacing  $ L^{\pm 1} $ 
by  $ K^{\pm 1} $  throughout: apart from this, the analysis and its
outcome are exactly the same.  Like for  $ \, \gersl_2 $  (cf.~\S
5.5), all the adjoint quantum objects   ---  i.e.~$ \hyp^a_q(\gerh_n) $,
$ {\hyp^a_q(\gerh_n)}' $  and  $ \big( {\hyp^a_q(\gerh_n)}' \big)^{\!
\vee} $  ---   will be strictly contained in the corresponding simply
connected quantum objects.  However, the semiclassical limits will be
the same in the case of  $ \hyp_q(\gerg) $  (giving  $ \hyp(\gerh_n) $ 
in both cases) and in the case of  $ \big( {\hyp_q(\gerg)}' \big)^{\!
\vee} $  (always yielding  $ \u(\gerk) $), whereas the semiclassical
limit of  $ {\hyp_q(\gerg)}' $  in the simply connected case will
be a (countable) covering of the limit in the adjoint case.

\vskip7pt

  {\bf 7.6 The QFA  $ \, F_q[H_n] \, $.} \, Now we look at
Theorem 2.2 the other way round, i.e.~from QFAs to QrUEAs.
We begin by introducing a QFA for the Heisenberg group.
                                             \par
   Let  $ \, F_q[H_n] \, $  be the unital associative
\hbox{$R$--alge}bra  with generators  $ \, \text{a}_1 \, $, 
$ \dots $,  $ \text{a}_n \, $,  $ \text{c} \, $,  $ \text{b}_1 \, $,
$ \dots $,  $ \text{b}_n \, $,  and relations  (for all  $ \, i, j
= 1, \dots, n \, $)
  $$  \text{a}_i \text{a}_j = \text{a}_j \text{a}_i \; ,  \hskip9pt
\text{a}_i \text{b}_j = \text{b}_j \text{a}_i \; ,  \hskip9pt  \text{b}_i
\text{b}_j = \text{b}_j \text{b}_i \; ,  \hskip9pt  \text{c} \, \text{a}_i
= \text{a}_i \, \text{c} + (q-1) \, \text{a}_i \; ,  \hskip9pt  \text{c}
\, \text{b}_j = \text{b}_j \, \text{c} + (q-1) \, \text{b}_j  $$
with a Hopf algebra structure given by (for all  $ \, i,
j = 1, \dots, n \, $)
  $$  \displaylines{
   \Delta(\text{a}_i) = \text{a}_i \otimes 1 + 1 \otimes \text{a}_i
\; ,  \hskip13pt  \Delta(\text{c}) = \text{c} \otimes 1 + 1 \otimes
\text{c} + {\textstyle \sum\limits_{j=1}^{n}} \text{a}_\ell \otimes
\text{b}_\ell \; ,  \hskip13pt  \Delta(\text{b}_i) = \text{b}_i \otimes 1
+ 1 \otimes \text{b}_i  \cr
   \epsilon(\text{a}_i) = 0 \; ,  \hskip7pt  \epsilon(\text{c})
= 0 \; ,  \hskip7pt \epsilon(\text{b}_i) = 0 \; ,  \hskip21pt
S(\text{a}_i) = - \text{a}_i \; ,  \hskip7pt  S(\text{c}) = - \text{c}
+ {\textstyle \sum\limits_{j=1}^{n}} \text{a}_\ell \, \text{b}_\ell \; ,
\hskip7pt  S(\text{b}_i) = - \text{b}_i  \cr }  $$
and let also  $ \, \F_q[H_n] \, $  be the  $ F(R) $--algebra  obtained
from  $ F_q[H_n] $  by scalar extension.  Then  $ \, \Bbb{B} :=
\Big\{ \prod_{i=1}^n \text{a}_i^{a_i} \cdot \text{c}^c \cdot
\prod_{j=1}^n \text{b}_j^{b_j} \,\Big\vert\, a_i, c, b_j \in
\N \; \forall \, i, j \,\Big\} \, $  is an  $ R $--basis  of
$ F_q[H_n] \, $,  hence an  $ F(R) $--basis  of  $ \F_q[H_n] \, $. 
Moreover,  $ F_q[H_n] $  is a QFA  (at  $ \, \h = q \! - \! 1 $) 
with semiclassical limit  $ F[H_n] \, $.

\vskip7pt

  {\bf 7.7 Computation of  $ \, {F_q[H_n]}^\vee \, $  and
specialization  $ \, {F_q[H_n]}^\vee \,{\buildrel {q \rightarrow 1}
\over \llongrightarrow}\, U({\gerh_n}^{\!\times}) \, $.} \, This section
is devoted to compute  $ \, {F_q[H_n]}^\vee \, $  and its semiclassical
limit (at  $ \, q = 1 \, $).
                                             \par
   Definitions imply that  $ \, \Bbb{B} \setminus \{1\} \, $  is an
$ R $--basis  of  $ \, J := \hbox{\sl Ker} \, \Big( \epsilon \, \colon
\, F_q[H_n] \! \relbar\joinrel\relbar\joinrel\twoheadrightarrow \! R
\Big) \, $.  
%
%
  Therefore  $ \; {F_q[H_n]}^\vee = \sum_{n \geq 0} 
  {\Big( {(q-1)}^{-1} J \Big)}^n \; $  
is just the unital
$ R $--algebra  (subalgebra of  $ \F_q[H_n] $)  with generators
$ \, E_i := \displaystyle{\, \text{a}_i \, \over \, q - 1 \,} \, $,
$ \, H := \displaystyle{\, \text{c} \, \over \, q - 1\,} \, $,  and
$ \, F_i := \displaystyle{\, \text{b}_i \, \over \, q - 1\,} \, $  ($ \,
i = 1, \dots, n \, $)  and relations (for all  $ \, i, j = 1, \dots,
n \, $)
  $$  E_i E_j = E_j E_i \; ,  \hskip13pt  E_i F_j = F_j E_i \; ,
\hskip13pt  F_i F_j = F_j F_i \; ,  \hskip13pt  H E_i = E_i H + E_i \; ,
\hskip13pt  H F_j = F_j H + F_j  $$
with Hopf algebra structure given by (for all  $ \, i, j = 1,
\dots, n \, $)
 \vskip-13pt
  $$  \Delta(E_i) = E_i \otimes 1 + 1 \otimes E_i \, ,  \hskip5pt
\Delta(H) = H \otimes 1 + 1 \otimes H + (q-1) {\textstyle
\sum\limits_{j=1}^{n}} E_j \otimes F_j \, ,  \hskip5pt
\Delta(F_i) = F_i \otimes 1 + 1 \otimes F_i  $$
$$  \epsilon(E_i) = \epsilon(H) = \epsilon(F_i) = 0 \, ,
\hskip8pt  S(E_i) = \! - E_i \, ,  \hskip8pt  S(H) = \! - H
+ (q-1) {\textstyle \sum\limits_{j=1}^{n}} E_j F_j \, ,
\hskip8pt  S(F_i) = \! - F_i \, .  $$
 \vskip-5pt
   At  $ \, q = 1 \, $  this implies that  $ \, {F_q[H_n]}^\vee
\,{\buildrel \, q \rightarrow 1 \, \over \llongrightarrow}\;
U({\gerh_n}^{\!\star}) = U({\gerh_n}^{\! *}) \; $  as co-Poisson
Hopf algebras, for a co-Poisson Hopf algebra isomorphism
 \vskip-11pt
  $$  {F_q[H_n]}^\vee \Big/ (q-1) \, {F_q[H_n]}^\vee
\,{\buildrel \cong \over \llongrightarrow}\; U({\gerh_n}^{\! *})  $$
 \vskip1pt
\noindent   exists, given by  $ \;\, E_i \mod (q-1) \mapsto
\text{e}_i \, $,  $ \, H \mod (q-1) \mapsto \text{h} \, $,  $ \,
F_i \mod (q-1) \mapsto \text{f}_i \, $,  \; for all  $ \, i, j = 1,
\dots, n \, $.  Thus  $ \, {F_q[H_n]}^\vee \, $  specializes to  $ \,
U({\gerh_n}^{\! *}) \, $  {\sl as a co-Poisson Hopf algebra},  q.e.d.

\vskip7pt

  {\bf 7.8 The identity  $ \, {\big({F_q[H_n]}^\vee\big)}' = F_q[H_n]
\, $.} \, Finally, we check the validity of the part of  Theorem
2.2{\it (b)}  claiming, when  $ \, p = 0 \, $,  \, that  $ \; H
\in \QFA \,\Longrightarrow\, {\big(H^\vee\big)}' = H \; $  for the
QFA  $ \, H = F_q[H_n] \, $.  Once more the proof works for all
$ \, p \geq 0 \, $,  \, so we do improve  Theorem 2.2{\it (b)}.
                                         \par
 First of all, from definitions induction gives, for all
$ \, m \in \N \, $,
  $$  \Delta^m(E_i) = \hskip-4pt {\textstyle \sum\limits_{r+s=m-1}}
\hskip-4pt 1^{\otimes r} \otimes E_i \otimes 1^{\otimes s} \, ,
\hskip15pt  \Delta^m(F_i) = \hskip-4pt {\textstyle \sum\limits_{r+s
= m-1}}  \hskip-4pt 1^{\otimes r} \otimes F_i \otimes 1^{\otimes s}
\hskip25pt  \forall \; i = 1, \dots, n  $$
 \vskip-13pt
  $$  \Delta^m(H) = \hskip-4pt {\textstyle \sum\limits_{r+s=m-1}}
\hskip-4pt 1^{\otimes r} \otimes H \otimes 1^{\otimes s} +
{\textstyle \sum\limits_{i=1}^m} {\textstyle \sum\limits_{\substack{  j, k
= 1  \\   j < k   }}^m} 1^{\otimes (j-1)} \otimes E_i \otimes
1^{\otimes (k-j-1)} \otimes F_i \otimes 1^{\otimes (m-k)}  $$
 \vskip-1pt
\noindent   so that  $ \; \delta_m(E_i) = \delta_\ell(H) = \delta_m(F_i)
= 0 \; $  for all  $ \, m > 1 $,  $ \ell > 2 \, $  and  $ \, i = 1,
\dots, n \, $;  moreover, for  $ \; \dot{E}_i := (q-1) E_i = \text{a}_i
\, $,  $ \dot{H} := (q-1) H = \text{c} \, $,  $ \dot{F}_i := (q-1) F_i
= \text{b}_i \, $  ($ i = 1, \dots, n $)  one has
  $$  \displaylines{
   \delta_1\big( \dot{E}_i \big) \! = (q-1) E_i \, ,  \; \delta_1\big(
\dot{H} \big) \! = (q-1) H ,  \;  \delta_1\big( \dot{F}_i \big) \!
= (q-1) F_i \in (q-1) \, {F_q[H_n]}^\vee \setminus {(q-1)}^2
{F_q[H_n]}^\vee  \cr
   \delta_2\big( \dot{H} \big) \, = \, {(q-1)}^2 \, {\textstyle
\sum_{i=1}^n} \, E_i \otimes F_i \in {(q-1)}^2 {\big( {F_q[H_n]}^\vee
\big)}^{\otimes 2} \setminus {(q-1)}^3 {\big( {F_q[H_n]}^\vee
\big)}^{\otimes 2} \, .  \cr }  $$  
   \indent   The outcome is that  $ \, \dot{E}_i = \text{a}_i \, ,
\dot{H} = \text{c} \, , \dot{F}_i = \text{b}_i \in {\big( {F_q[H_n]}^\vee
\big)}' \, $,  \, so the latter algebra contains the one generated by
these elements, that is  $ F_q[H_n] \, $.  Even more,  $ {F_q[H_n]}^\vee $
is clearly the  $ R $--span  of the set  $ \, {\Bbb B}^\vee := \Big\{\,
\prod_{i=1}^n E_i^{a_i} \cdot H^c \cdot \prod_{j=1}^n F_j^{b_j}
\,\Big\vert\, a_i, c, b_j \in \N \; \forall \, i, j \,\Big\} \, $,
\, so from this and the previous formulas for  $ \Delta^n $  one
gets that  $ \, {\big( {F_q[H_n]}^\vee \big)}' = F_q[H_n] \, $,
\, q.e.d.    

\vskip2,1truecm

\end{document}